%% file: High_Order_Transition_Elements.tex
\documentclass[a4paper,3p,11pt,sort&compress]{elsarticle}
\usepackage{geometry}
\usepackage{lineno,hyperref}
\modulolinenumbers[5]
\usepackage{amsmath,amssymb}
\usepackage{subfig}
\usepackage{algorithm}
\usepackage[end]{algpseudocode}
\usepackage{booktabs}
\usepackage{longtable}
\usepackage{pdflscape}
\usepackage{cancel}
\usepackage{placeins}
\usepackage{paralist}
\usepackage{color}
\usepackage[bottom]{footmisc}
\usepackage{caption}
\usepackage{upgreek}
\usepackage{todonotes}
\usepackage{bbding}
\usepackage{relsize}
\newcommand*{\NewCirc}{\raisebox{1pt}{$\mathsmaller{\bigcirc}$}}
\newcommand*{\NewStar}{\raisebox{-2.5pt}{\FiveStarOpen}}
\definecolor{Matlab1}{rgb}{0.0,0.0,0.0}
\definecolor{Matlab2}{rgb}{0.878431372549020,0.015686274509804,0.015686274509804}
\definecolor{Matlab3}{rgb}{0.101960784313725,0.352941176470588,0.627450980392157}
\definecolor{Matlab4}{rgb}{0.6,0.6,0.6}
\definecolor{Matlab5}{rgb}{1.0,0.443137254901961,0.0}
\definecolor{Matlab6}{rgb}{0.274509803921569,0.560784313725490,1.0}
\definecolor{Matlab7}{rgb}{1.0,0.0,0.0}
\definecolor{Matlab8}{rgb}{0.0,0.0,1.0}
\definecolor{Matlab9}{rgb}{0.0,1.0,0.0}
\definecolor{xElem}{rgb}{0.8500,0.3250,0.0980}
\definecolor{xNyElem}{rgb}{0,0.4470,0.7410}
\definecolor{yElem}{rgb}{0.4940,0.1840,0.5560}
\definecolor{yNyElem}{rgb}{0.9290,0.6940,0.1250}
\journal{Comput Methods Appl Mech Eng}
\begin{document}

\begin{frontmatter}
\title{High order transition elements: The \emph{x\textbf{N}y}-element concept - Part I: Statics}

\author{S.~Duczek\corref{cor1}}
\ead{s.duczek@unsw.edu.au}

\author{A. A.~Saputra\corref{}}
\ead{a.saputra@unsw.edu.au}
\address{The University of New South Wales Sydney, School of Civil and Environmental Engineering, Sydney, NSW 2052, Australia.}
\cortext[cor1]{Corresponding author}

\author{H.~Gravenkamp\corref{}}
\ead{hauke.gravenkamp@uni-due.de}
\address{University of Duisburg-Essen, Department of Civil Engineering, Universit\"atsstra\ss{}e 15, 45141 Essen, Germany.}

\begin{abstract}
Advanced transition elements are of utmost importance in many applications of the finite element method (FEM) where a local mesh refinement is required. Considering problems that exhibit singularities in the solution, an adaptive \emph{hp}-refinement procedure must be applied. Even today, this is a very demanding task especially if only quadrilateral/hexahedral elements are deployed and consequently the hanging nodes problem is encountered. These element types, are, however, favored in computational mechanics due to the improved accuracy compared to triangular/tetrahedral elements. Therefore, we propose a compatible transition element -- \emph{x\textbf{N}y}-element -- which provides the capability of coupling different element types. The adjacent elements can exhibit different element sizes, shape function types, and polynomial orders. Thus, it is possible to combine independently refined \emph{h}- and \emph{p}-meshes. The approach is based on the \textit{transfinite mapping concept} and constitutes an extension/generalization of the \emph{pNh}-element concept. By means of several numerical examples, the convergence behavior is investigated in detail, and the asymptotic rates of convergence are determined numerically. Overall, it is found that the proposed approach provides very promising results for local mesh refinement procedures.
\end{abstract}

\begin{keyword}
Transition elements\sep \emph{p}-Version of the finite element method\sep Spectral element method\sep Transfinite mappings\sep Local mesh refinement.
\end{keyword}
\end{frontmatter}
\vspace*{-15pt}
\section*{Highlights}
\begin{compactitem}
\item Derivation of general transition elements coupling different
\begin{compactenum}
	\item element types,
	\item element sizes,
	\item polynomial orders.
\end{compactenum}
\item Verification that the elements are high order complete (high order patch test).
\item Proof that the asymptotic convergence rates depend on the minimum polynomial order $p_\mathrm{min}$.
\item Assessment of the performance of transition elements using examples with singularity.
\end{compactitem}
\rule[0pt]{\textwidth}{0.5pt}
\vspace*{-24pt}
\tableofcontents
%
%
%
%---------------------------------------------------------------------------%
\input{./tex/Intro.tex}
\input{./tex/HighOrderShapeFunctions.tex}
\input{./tex/Theory.tex}
\input{./tex/PatchTest.tex}
\input{./tex/ConvergenceRate.tex}
\input{./tex/SingularExamples.tex}
\input{./tex/Summary.tex}
%---------------------------------------------------------------------------%
%
%
%
\FloatBarrier
%
%\vspace*{24pt}
%
%\noindent\textbf{Acknowledgments}\\[3pt]
%SD would like to acknowledge the German Research Foundation (DFG) for its financial support under grant \mbox{DU 1613/1-1}.
%
\vspace*{24pt}
\addcontentsline{toc}{section}{References}
%\section*{References}
\bibliography{./bib/Literature_Duczek,./bib/paper_xNy_HG}
%----------------------------------Appendix---------------------------------%
%
\newpage
\input{./tex/Appendix.tex}
\end{document}

%% file: tex/Intro.tex
%
%
\section{Introduction}
\label{sec:Intro}
\begin{figure}[t!]
    \centering
    \includegraphics[clip,width=0.6\textwidth]{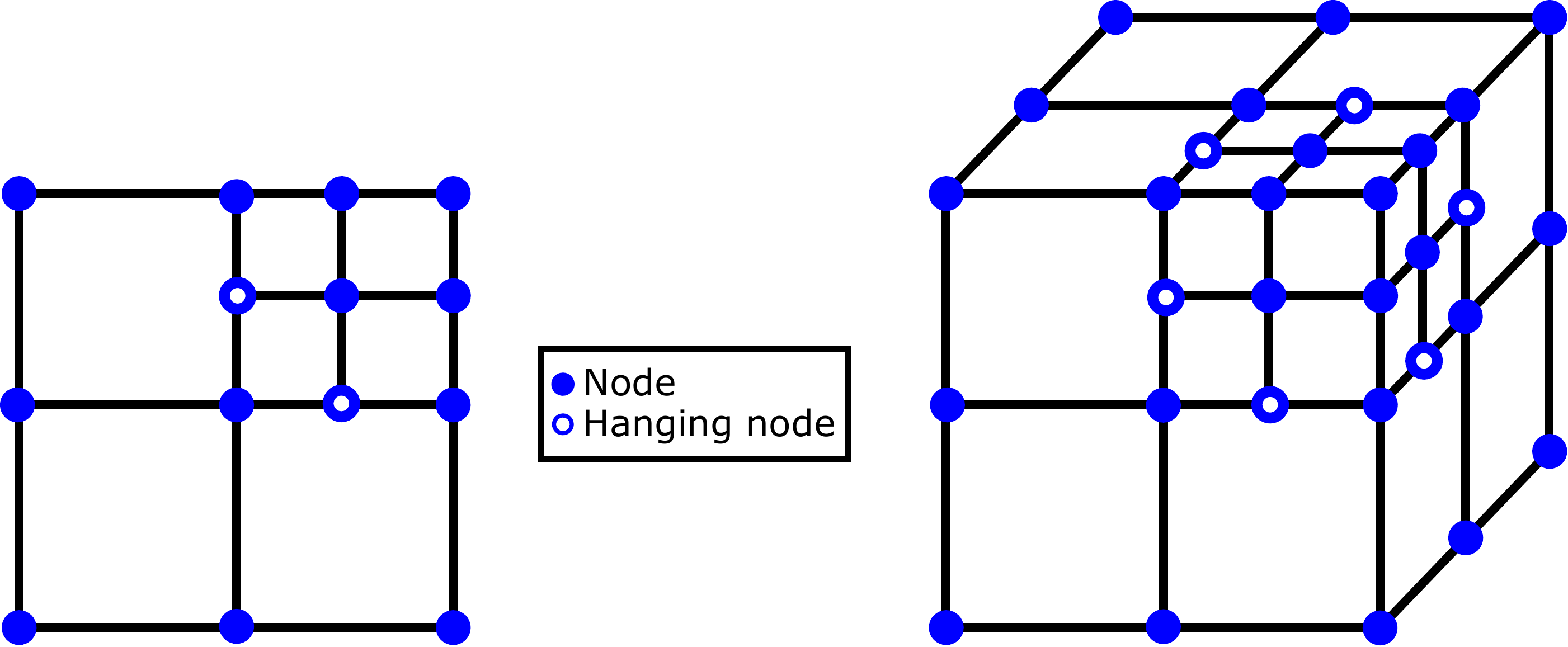}
	\caption{Hanging nodes problem in the standard FEM.}
	\label{fig:HangingNodes}
\end{figure}

Local mesh refinement procedures are of crucial importance if the solution to a partial differential equation (PDE) exhibits steep gradients or singular behavior. In these cases, it has been shown that high order finite element methods (FEMs) in conjunction with pure \emph{p}-refinement strategies (constant mesh size with increasing polynomial degree $p$) are not an optimal choice \cite{BookSzabo1991}, and only algebraic rates of convergence are obtained. In order to achieve the desired exponential convergence characteristics, a combination of \emph{h}- (constant polynomial degree $p$ with decreasing element size $h$) and \emph{p}-refinement is required \cite{BookSzabo1991, InbookDuester2018, BookDuester2002}. In computational mechanics, quadrilateral/hexahedral finite elements are generally favored over their triangular/tetrahedral counterparts. This preference is mainly related to the following two reasons: \textit{(i)}~quadrilateral/hexahedral finite elements provide more accurate results \cite{BookZienkiewicz2000a} and \textit{(ii)}~the formulation of high order shape functions for quadrilateral/hexahedral finite elements is significantly simplified due to the possibility of using simple tensor products of the one-dimensional shape functions \cite{BookPozrikidis2014, BookSzabo1991, PhDDuczek2014}. Considering triangular or tetrahedral (high order) finite elements, the formulation of suitable shape functions is much more elaborate \cite{ArticleDubiner1991, ArticleSherwin1995, ArticleBlyth2006, ArticleBlyth2006b, ArticleLuo2006, BookPozrikidis2014, BookKarniadakis2005}. The mentioned attributes explain why the analyst usually strives to utilize quadrilateral/hexahedral elements as much as possible. This, however, causes problems when a local \emph{h}-refinement is required to improve the solution. A conformal discretization can only be achieved if special element types are used or other measures are taken to prevent the so-called ``hanging node problem'' (see Fig.~\ref{fig:HangingNodes}). A matching interpolation between adjacent finite elements can be generated by employing the penalty method, Lagrange multipliers, Nitsche's method, or multi-point constraints \cite{BookBathe2002}. The problem with these methods is that they severely constrain the solution field (primary variables) and therefore, the quality of the numerical results deteriorates in a certain region around the coupling zone.
\begin{figure}[b!]
	\centering
	\includegraphics[clip,width=0.9\textwidth]{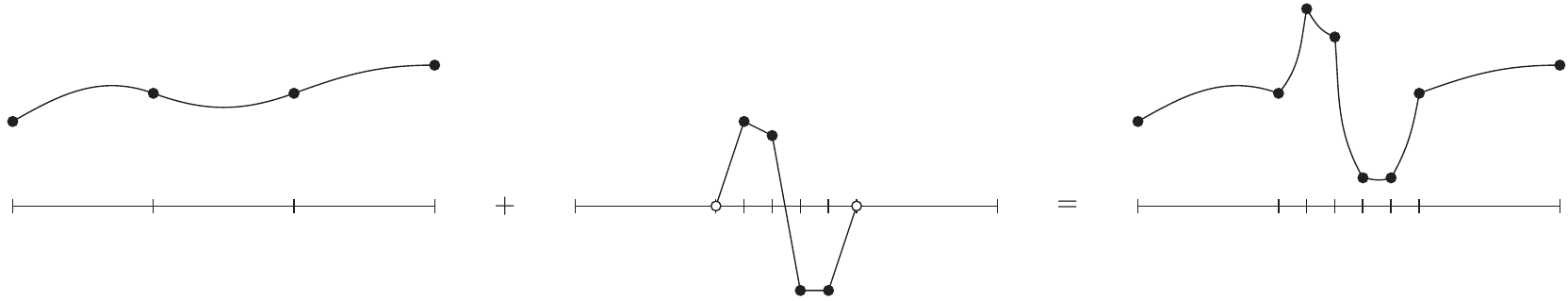}
	\caption{Mesh refinement by superposition (overlay meshes) \cite{PhDZander2017}. In the area of interest, the solution on the base mesh (left) is enhanced by results computed using an overlay mesh (middle) to obtain an accurate final solution (right).}
	\label{fig:OverlayMesh}
\end{figure}%

A different approach to ensure a continuous solution field is based on the partition of unity method (PUM) where the numerical results are improved by enriching the ansatz space \cite{ArticleFries2010}. Two prominent representatives of this class of methods are the extended finite element method (XFEM) and the generalized finite element method (GFEM)\footnote{Note that a PU is not a prerequisite for the application of the XFEM and GFEM, but PU shape functions are often used to approximate the enrichment function.}. The PUM is a widely used and very versatile approach, but depending on the enrichment functions that are utilized, problems with the numerical integration, conditioning of the system matrices, and blending elements\footnote{Blending elements are special elements in XFEM where not all nodes are enriched.} may arise \cite{ArticleBelytschko2009}. Closely related to this approach are schemes that use a coarse basis mesh, and the mesh refinement is executed by deploying an overlay mesh with a different element size or polynomial order in areas where the solution must be improved (see Fig.~\ref{fig:OverlayMesh}). The \emph{d}- \cite{ArticleRank1992}, \emph{s}- \cite{ArticleFish1992}, \emph{hpd}- \cite{ArticleRank1997} and multi-level \emph{hp}-versions of the FEM \cite{ArticleZander2015} fall into this category. Basically, the computational domain is decomposed into two regions with different solution characteristics: \textit{(i)}~smooth or \textit{(ii)}~non-smooth solution fields. In the first region, a coarse \emph{p}-mesh can be employed while a mesh refinement is needed in the second one. The mentioned approaches only differ in how the mesh refinement (overlay mesh) is set up. In the beginning of the 1990s, the \emph{d}- (d:~domain decomposition) and \emph{s}-versions (s:~superposition) of the FEM have been developed in parallel by Rank \cite{ArticleRank1992, InproceedingsRank1993} and Fish \cite{ArticleFish1992, ArticleFish1992b}, respectively. In the \emph{d}-version, a coarse \emph{p}-mesh is overlaid by a fine discretization using \emph{h}-elements \cite{InproceedingsRank1993}. Across the boundary of the refinement zone, a C\textsuperscript{0}-continuous approximation is achieved by imposing homogeneous Dirichlet boundary conditions (BCs). Since the coarse mesh is a subset of the fine discretization -- i.e., a certain number of small (refined) elements share the same elemental boundary with a large base element \cite{ArticleRank1992} -- the Dirichlet BCs can be applied in a strong sense. However, regarding the \emph{s}-FEM proposed by Fish \cite{ArticleFish1992}, the base and overlay meshes do not need to be conformal anymore. That is to say, C\textsuperscript{0}-continuity (compatibility) must be achieved by imposing Dirichlet BCs on the refinement domain in a weak sense. This can be done by means of standard approaches such as the penalty method, Lagrange multipliers, or Nitsche's method. The overlay mesh may consist of hierarchical high order finite elements in regions where the solution indicates high gradients. The \emph{hpd}-FEM is an extension of the \emph{d}-FEM \cite{ArticleRank1997, ArticleDuester2007} where the idea of several independent overlay meshes is introduced. The different overlay meshes may overlap, and linear shape functions are commonly employed. This idea can also be used when coupling dimensionally different elements (dimensional adaptivity) \cite{ArticleDuester2007}. In Refs. \cite{PhDSchillinger2012, ArticleSchillinger2012c}, several (hierarchical) levels of overlay meshes being conformal with the local spacetree decomposition of the integration domain have been proposed. With this approach, the discontinuities in the displacement field caused by material inclusions can be easily captured. Yet, the most general approach is the multi-level \emph{hp}-FEM, where several high order (hierarchical) overlay meshes are introduced \cite{ArticleZander2015, ArticleZander2016, ArticleZander2017, PhDZander2017}. In this method, an optimal combination of high order base elements and high or low order overlay elements can be set up. However, common to all introduced methods based on this \textit{refinement-by-superposition} technique is that extra care has to be taken to ensure the linear independence. A linearly dependent set of shape functions should be avoided at all costs, since it results in a rank-deficient stiffness matrix \cite{ArticleDuester2007}. In addition, the implementation of these methods in existing FE-codes is not straightforward since complex data structures have to be introduced when using high order (hierarchical) overlay meshes.

\begin{figure}[t!]
	\centering
	\subfloat[1-irregular (balanced) quadtree mesh]{\includegraphics[clip,width=0.47\textwidth]{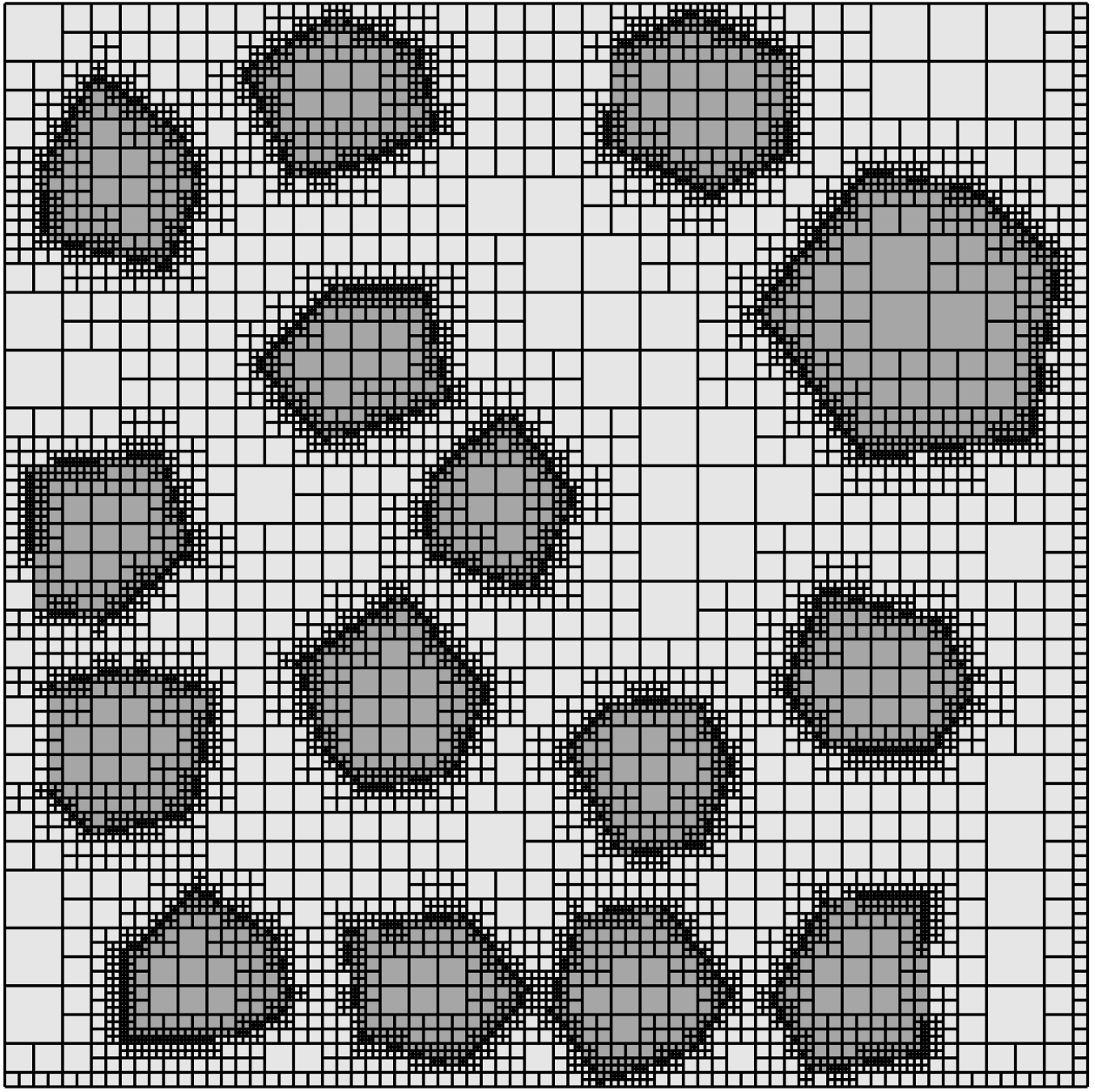}}\\
	\subfloat[Master elements]{\includegraphics[clip,width=0.75\textwidth]{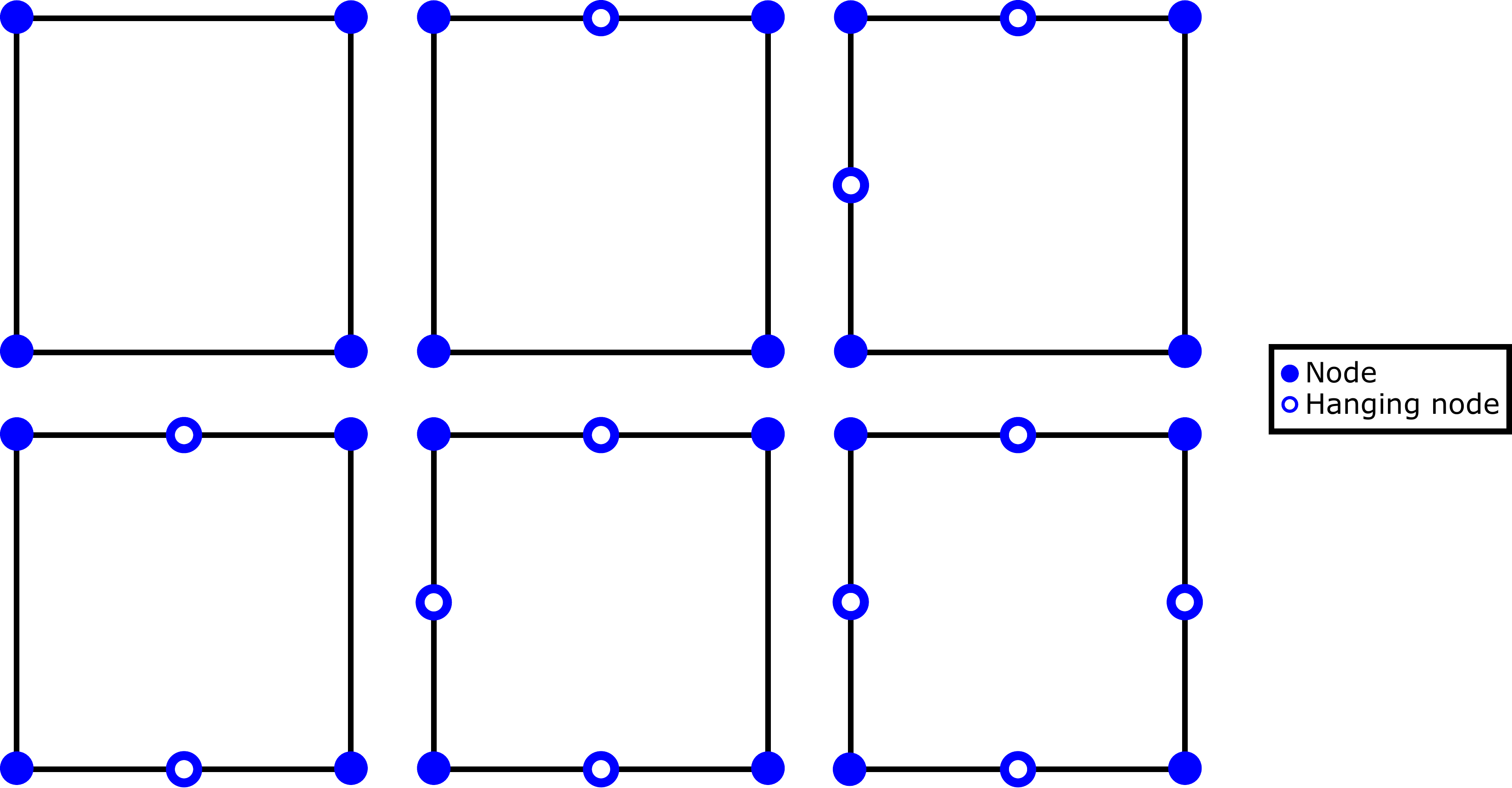}}
	\caption{Transition elements to achieve a balanced quadtree mesh (linear shape functions).}
	\label{fig:QuadtreeMesh}
\end{figure}

Another strategy to circumvent the hanging node problem, which will be followed in this article, is seen in the use of special transition elements. These elements are very helpful when local mesh refinement procedures are based on quadrilateral/hexahedral finite elements. Among those approaches, variable-node elements for 1-irregular/balanced meshes\footnote{A 1-irregular mesh is a discretization where on each edge only on hanging node may exist; i.e., in two-dimensional applications, one large element can be coupled to two smaller elements.} have been constructed by Gupta for two-dimensional applications \cite{ArticleGupta1978} (see Fig.~\ref{fig:QuadtreeMesh}), while Morton et al.\ published an extension for three dimensions \cite{ArticleMorton1995}. In Ref.~\cite{ArticleHuang2010}, an improved 5-node element is developed based on the work by Gupta \cite{ArticleGupta1978}, where the theoretically optimal rates of convergence are restored. These elements feature piecewise linear shape functions on their boundaries such that two smaller finite elements can be coupled conformally to a larger transition element. The derivation of these shape functions is described visually in Ref.~\cite{BookZienkiewicz2000a} (Chap.\ 8, p.\ 175). Note that this approach can also be used to derive Serendipity shape functions in a structured way\footnote{The derivation of Serendipity shape functions can also be achieved based on a tensor product formulation of the shape functions \cite{ArticleArnold2011, ArticleFloater2017}.} using some ``ingenuity'' \cite{ArticleScholz1985, BookZienkiewicz2000a}. To rely on ingenuity is, however, not a valid scientific approach and therefore, a generalization of this method is required. The sought-after generalization has already been proposed in the early 1970s and is closely related to the names Gordon and Hall \cite{ArticleGordon1971, ArtcileGordon1973a, ArtcileGordon1973b, BookGordon1982}. The works of Gordon and co-workers are probably best known in the context of geometry mapping in the \emph{p}-FEM \cite{BookSzabo1991, ArticleKiralyfalvi1997}. They developed a methodology which is known as transfinite interpolation or blending function method \cite{ArtcileGordon1973b}. In this framework, the term transfinite means that the interpolation of an arbitrary function matches the original function at a non-denumerable (transfinite\footnote{Note that the word transfinite can often be replaced by the word infinite without changing the meaning of the term. We, however, decided to stick to the terminology introduced by Gordon \cite{ArticleGordon1971}.}) set of points \cite{ArticleGordon1971}, i.e., the functions are identical in certain parts but not over the whole domain.
This idea is taken from Coons, who was the first to describe a method to construct an interpolatory surface which coincides with arbitrarily prescribed curves \cite{TechReportCoons1967}. It is obvious that this approach can also be used to generate complex shape functions. Transition elements which couple different element types are straightforward to derive by introducing piecewise polynomial shape functions on the element edges. As mentioned before, such elements are indispensable for local mesh refinement procedures based on quadrilateral/hexahedral discretizations. Scholz developed two- and three-dimensional transition elements with piecewise linear and quadratic shape functions for mesh refinement purposes \cite{ArticleScholz1985, PhDScholz1986, ArticleAltenbach1987}. Here, each edge can be divided into two parts and therefore, 1-irregular meshes can be generated without introducing hanging nodes. In addition to the mentioned element types, R\"ohr used the transfinite interpolation technique to derive  C\textsuperscript{1}-continuous transition elements based on cubic Hermite shape functions \cite{HabilRoehr1988}. In a series of papers, Cavendish used the blending function method to derive piecewise linear membrane elements, three-dimensional continuum elements, and C\textsuperscript{1}-continuous elements (cubic Hermite blending functions) \cite{ArticleBrikhoff1974, ArticleCavendish1975, ArticleCavendish1977, ArticleCavendish1984}. Each refined element is referred to as macro-element, and by  analogy with conventional substructuring techniques, the internal degrees of freedom (DOFs) are eliminated by static condensation \cite{ArticleCavendish1977}. In the literature, such elements are sometimes also called Coons-patch macro-elements, see Refs.~\cite{ArticleProvatidis2006, ArticleProvatidis2011, BookProvatidis2019}. In his works, Provatidis uses shape functions based on piecewise linear functions, cubic B-splines, or Lagrange polynomials to construct versatile transition elements \cite{ArticleProvatidis2006}. The role of interior nodes is investigated in Ref.~\cite{ArticleProvatidis2011}, where it is stated that ``a sufficient number of internal nodes is necessary to ensure convergence''. Note that these additional (internal) nodes ensure polynomial completeness of the shape functions. Provatidis also argues that this type of macro-element constitutes a precursor for the nowadays very successful isogeometric analysis (IGA) \cite{BookProvatidis2019}. In Refs.~\cite{PhDWeinberg1996, ArticleWeinberg1999, ArticleWeinberg2002}, Weinberg and Gabbert described the so-called \emph{pNh}-element concept (which uses an analogous ansatz as discussed above). The development of this element type was inspired by the difficulties that were encountered when solving contact problems with the \emph{p}-version of the FEM. To be able to employ standard contact algorithms that are widely implemented in commercial software tools, the idea to couple large \emph{p}-elements in areas with smooth solution characteristics to a fine mesh consisting of \emph{h}-elements where contact might occur (for contact problems the polynomial degree is $p_\mathrm{h}\,{=}\,1$, i.e. linear shape functions are generally preferred), was developed. The name derives from the fact that one \emph{p}-element can be coupled to $N$ \emph{h}-elements \cite{ArticleWeinberg2002}. In this formulation, it is, in principle, possible to divide each edge/face into a different number of subdomains \cite{PhDWeinberg1996} and also to assign different polynomial degrees to each element, although this possibility was not further explored. The approach being proposed in this article builds on the work of Weinberg and extends it in several aspects. The developed concept will be referred to as \emph{x\textbf{N}y}-element to emphasize that it is possible to couple different element families \emph{x} and \emph{y} (based on different shape functions).  Since it is basically possible to couple an arbitrary number of \emph{y}-elements to each edge/face of an \emph{x}-element, \emph{\textbf{N}} is denoted in boldface. That is to say, \emph{\textbf{N}} is generally vector-valued and contains four different components denoting the numbers of elements that are coupled to each of the four edges of a quadrilateral element. In this article, the focus is placed on a compatible coupling of spectral \cite{BookPozrikidis2014, BookKarniadakis2005} and \emph{p}-elements \cite{BookSzabo1991, InbookDuester2018}, although it is theoretically possible to couple arbitrary element types.

For the sake of completeness, other approaches that can be used for local mesh refinement in conjunction with quadrilateral/hexahedral finite elements are briefly introduced in the following.
Baitsch \& Hartmann also developed a transition element based on piecewise polynomial shape functions \cite{ArticleBaitsch2009}. This approach requires subdividing the refined element into a number of quadrilateral subdomains and consequently leads to a comparably large number of internal degrees of freedom within each macro-element. Therefore, a macro-element similar to the one proposed by Cavendish is obtained \cite{ArticleCavendish1977}. In Ref.~\cite{ArticleCho2005}, a variable-node transition element based on the moving least-square (MLS) approximation is proposed. This element is compatible with a piecewise quadratic interpolation on the element edges. Note that depending on the node configuration, a negative Jacobian (mapping from the local to the global domain) might be observed, which needs to be addressed by adding additional nodes \cite{ArticleCho2005}. This element formulation has been applied to contact problems using a node-to-node algorithm \cite{ArticleKim2008}. Based on the smoothed finite element method (SFEM), a polygonal variable-node element is devised \cite{ArticleKim2017}. This method can straightforwardly be used to couple elements of different sizes. Another approach that can be employed to create very versatile discretizations is the scaled boundary finite element method (SBFEM) \cite{BookSong2018}, in particular, the variant known as scaled boundary polygons \cite{Gravenkamp2018a, Chiong2014, Ooi2016}. This method allows the construction and coupling of star-convex polygonal elements consisting of an arbitrary number of edges, where the type and order of interpolation along each edge can be chosen independently. The SBFEM has been applied successfully to quadtree/octree decompositions \cite{Ooi2015c, Gravenkamp2017a, Saputra2017, Gravenkamp2017c}, as well as polytopal meshes \cite{Xing2018,Liu2017}. For linear static problems, the SBFEM achieves high order convergence while only requiring nodes on the element edges. However, for cases involving inertia terms, non-linearities, or arbitrary body loads, additional degrees of freedom have to be introduced to achieve high order convergence \cite{Song2009, Ooi2016, Bulling2018}. Extensive research in the area of transition elements, especially with respect to coupling solid (continuum) and shell elements has been conducted by Surana \cite{ArticleSurana1980a, ArticleSurana1980b, ArticleSurana1980c, ArticleSurana1982, ArticleSurana1983, ArticleSurana1986, ArticleSurana1987}. In Ref,~\cite{ArticleSurana1980a} Surana provides shape functions for transition elements which can be used to couple quadrilateral with beam-like or other quadrilateral elements for different polynomial orders ($p\,{=}\,1,2,3$). Extensions to axisymmetric and three-dimensional problems have been proposed in Refs.~\cite{ArticleSurana1980c, ArticleSurana1980b}. Thus, it is straightforwardly possible to couple hexahedral and (curved) shell-elements without the need for additional multi-point constraints, penalty, or mortar-coupling. Based on this work, geometrically non-linear transition elements are formulated using a total Lagrangian approach in Refs.~\cite{ArticleSurana1982, ArticleSurana1983}. Finally, it should be noted that similar elements are available for (Fourier) heat conduction problems where the nodal temperature and temperature gradients are included as primary variables \cite{ArticleSurana1986, ArticleSurana1987}.

%% file: tex/HighOrderShapeFunctions.tex
\section{High Order Shape Functions}
\label{sec:HO_ShapeFunctions}
One of the primary goals of this contribution is to study the performance of transition elements based on transfinite mappings when different element types are coupled. In this article, we focus on coupling the \emph{p}-version of the FEM with the spectral element method (SEM). Therefore, we will briefly recapitulate the derivation of the corresponding one-dimensional shape functions.
\subsection{Nodal shape functions: SEM}
\label{sec:SEM}
\begin{figure}[t!]
\centering
\subfloat[Shape functions]{\includegraphics[width=0.495\textwidth]{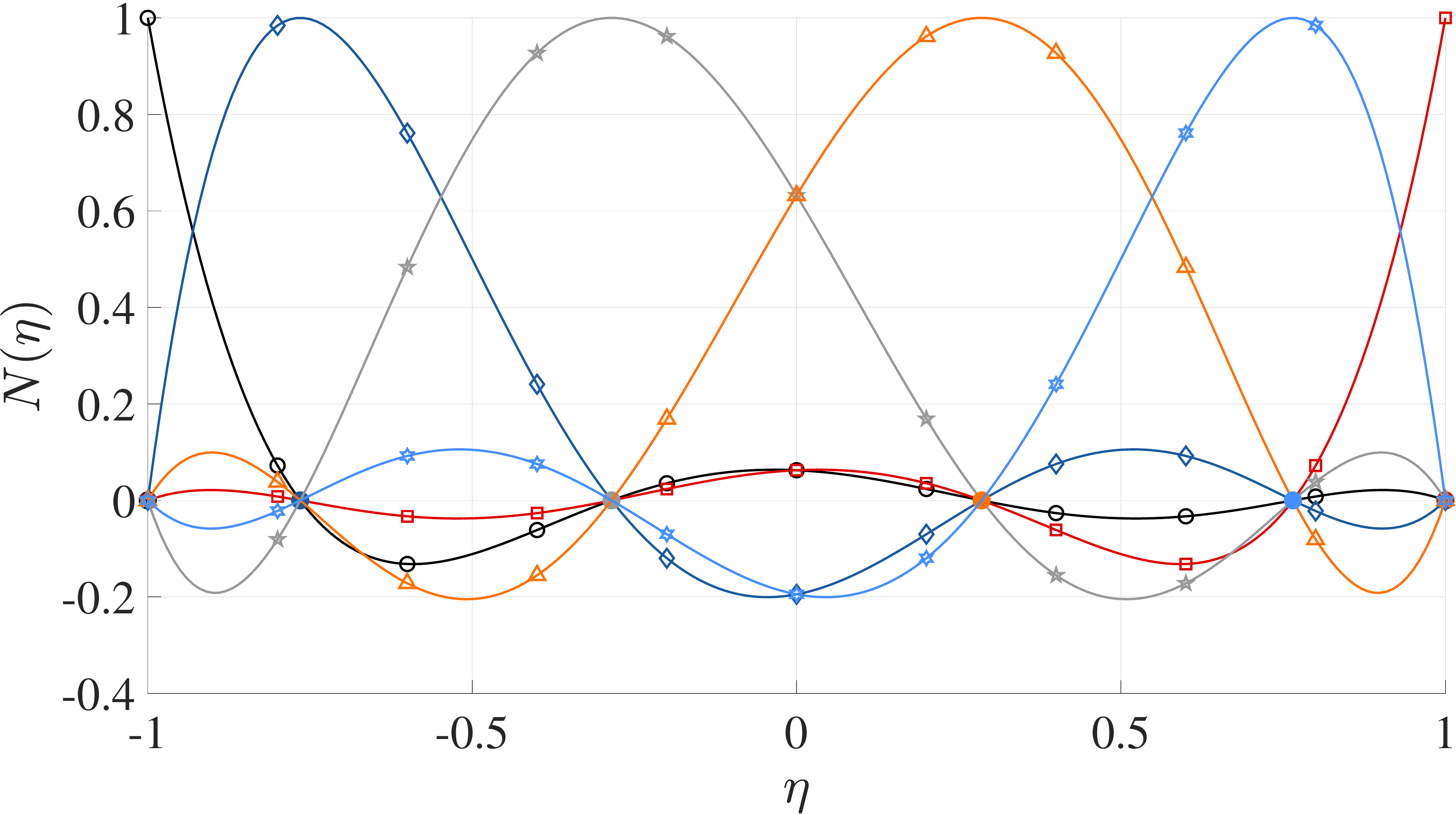}} 
\hfill
\subfloat[Derivatives]{\includegraphics[width=0.495\textwidth]{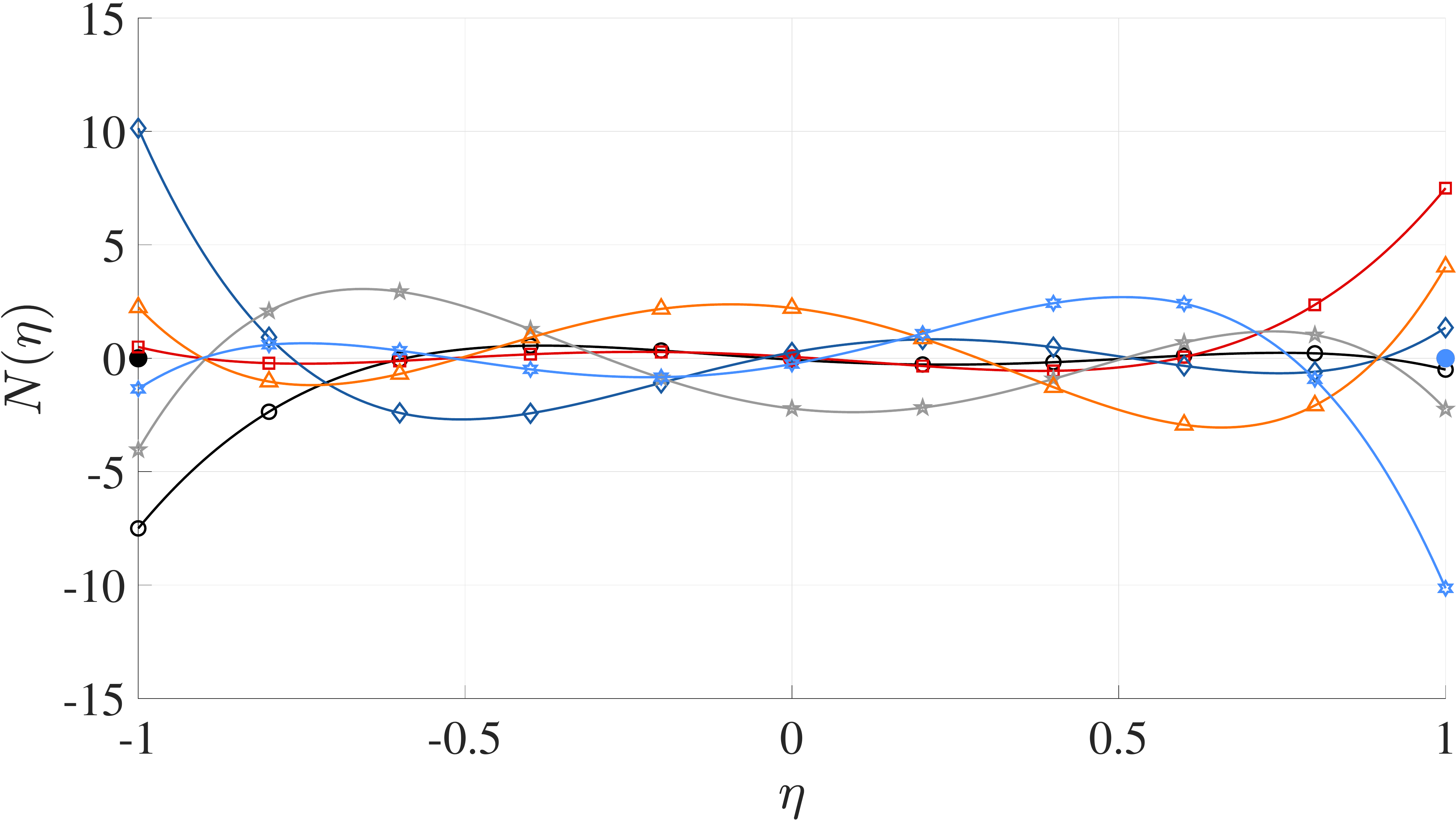}}
\caption{Lagrange shape functions defined on GLL points, $p\,{=}\,5$. Legend: \textcolor{Matlab1}{\rule[0.55ex]{5ex}{0.2ex}} $N_1(\eta)$ (\NewCirc), \textcolor{Matlab2}{\rule[0.55ex]{5ex}{0.2ex}} $N_2(\eta)$ ($\square$), \textcolor{Matlab3}{\rule[0.55ex]{5ex}{0.2ex}} $N_3(\eta)$ ($\Diamond$), \textcolor{Matlab4}{\rule[0.55ex]{5ex}{0.2ex}} $N_4(\eta)$ (\NewStar), \textcolor{Matlab5}{\rule[0.55ex]{5ex}{0.2ex}} $N_5(\eta)$ ($\bigtriangleup$), \textcolor{Matlab6}{\rule[0.55ex]{5ex}{0.2ex}} $N_6(\eta)$ ($\mathbf{\ast}$). \label{fig:shpFunLagrange}}
\end{figure}

The shape functions in the SEM are based on Lagrangian interpolation polynomials. Given a set of interpolation non-equidistant points/nodes, denoted by $\eta_i$ \cite{BookPozrikidis2014}, Lagrange polynomials of order $p$ are defined as \cite{BookAbramowitz1972}:
\begin{equation}
N_i(\eta) = \prod_{j=1,j\ne i}^{p+1}\cfrac{\eta-\eta_j}{\eta_i-\eta_j}\,, \qquad i=1,2,\ldots,p+1\,.
\label{eq:Lagrange}
\end{equation}
Such shape functions constitute a consistent extension of conventional low order (linear and quadratic) finite elements, inasmuch as they exhibit similar properties, in particular:
\begin{itemize}
\item \textit{Kronecker-delta-property:}
\begin{equation}\label{eq:kronecker}
N_i(\eta_j)=\delta_{ij}
\end{equation}
\item \textit{partition of unity:}
\begin{equation}\label{eq:partitionUnity}
\sum \limits_i N_i(\eta)=1 
\end{equation}
\end{itemize}
Equation~\eqref{eq:kronecker} also implies that these shape functions are `node-based' since the unknowns retain physical meaning. Therefore, the post-processing and the application of Dirichlet boundary conditions are straightforward. While Lagrange interpolation polynomials can, in principle, be defined for any set of distinct nodes, their properties depend critically on the nodal positions. In the wide body of literature, there are several possibilities to distribute these points within the reference interval $[-1,1]$, see \cite{BookPozrikidis2014}. In the framework of the SEM, the following distributions are commonly used: Gau\ss-Lobatto-Chebyshev (GLC) and Gau\ss-Lobatto-Legendre (GLL) points:
\begin{enumerate}
	\item GLC points \cite{ArtcilePatera1984}:
	\begin{itemize}
		\item $\eta_i^\mathrm{GLC} = -\cos\left(\cfrac{i-1}{p}\;\uppi\right)\,, \quad i = 1,2,\ldots,p\,{+}\,1\,.$
	\end{itemize}
	\item GLL points \cite{BookKarniadakis2005}:
	\begin{itemize}
		\item $\eta_i^\mathrm{GLL} = \begin{cases} -1 & \text{for } i\,{=}\,1\\ \hat{\eta}_i & \text{for } i\,{=}\,2,3,\ldots,p \\ +1 & \text{for } i\,{=}\,p\,{+}\,1 \end{cases}.$
	\end{itemize}
\end{enumerate}
In the above definition, $\hat{\eta}_i$ denotes the roots of the Lobatto polynomial $L_{\mathrm{o}_{p-1}}(\eta)$ of order $p-1$, which is defined as the first derivative of the Legendre polynomial $L_p(\eta)$ of order $p$~\cite{BookPozrikidis2014}
\begin{equation}
L_{\mathrm{o}_{p-1}}(\eta) = \cfrac{\mathrm{d}L_p(\eta)}{\mathrm{d}\eta}\,.
\end{equation}
Figure~\ref{fig:shpFunLagrange} depicts the spectral shape functions up to an order of $p=5$ (GLL nodal distribution). The most prominent drawback of this set of shape functions is that when changing the order $p$, all shape functions are modified and therefore, the coefficient matrices have to be computed anew. Lagrange interpolation polynomials are commonly used in the SEM for wave propagation phenomena~\cite{ArticleSeriani1994,ArticleDauksher2000}. In this context, the main advantage is the possibility to diagonalize the mass matrix (mass lumping) without loss of accuracy (retaining optimal rates of convergence) \cite{ArticleDuczek2019a,ArticleDuczek2019b}. 

\subsection{Hierarchic shape functions: p-FEM}
\label{sec:pFEM}
The concept of hierarchical shape functions is linked to the \emph{p}-version of the FEM \cite{BookSzabo1991,InbookDuester2018}. In this context, hierarchical means that all shape functions of order $p$ are contained in the set of shape functions of order $p\,{+}\,1$. The construction of these shape functions is based on the normalized integrals of the Legendre polynomials:
\begin{equation}
N_i(\eta) = \sqrt{\cfrac{2i-1}{2}} \int\limits_{-1}^{\eta} L_{i-1}(x)\mathrm{d}x\,, \quad i = 2,3,\ldots,p\,,
\label{eq:SF_pFEM_ho}
\end{equation}
with $L_{i-1}(x)$ denoting the Legendre polynomial of order $i\,{-}\,1$. The functions defined by Eq.~\eqref{eq:SF_pFEM_ho} contain polynomials of order 2 up to $p$ and vanish at the interval boundaries. To create a set of complete polynomials, these interpolants are augmented by the standard linear finite element shape functions:
\begin{subequations}
	\begin{align} \label{eq:N12}
	N_1(\eta) & = \cfrac{1}{2}(1-\eta)\,,\\
	N_{p+1}(\eta) & = \cfrac{1}{2}(1+\eta)\,.
	\end{align}
\end{subequations}
\begin{figure}[t!]
	\centering
	\subfloat[Shape functions]{\includegraphics[width=0.495\textwidth]{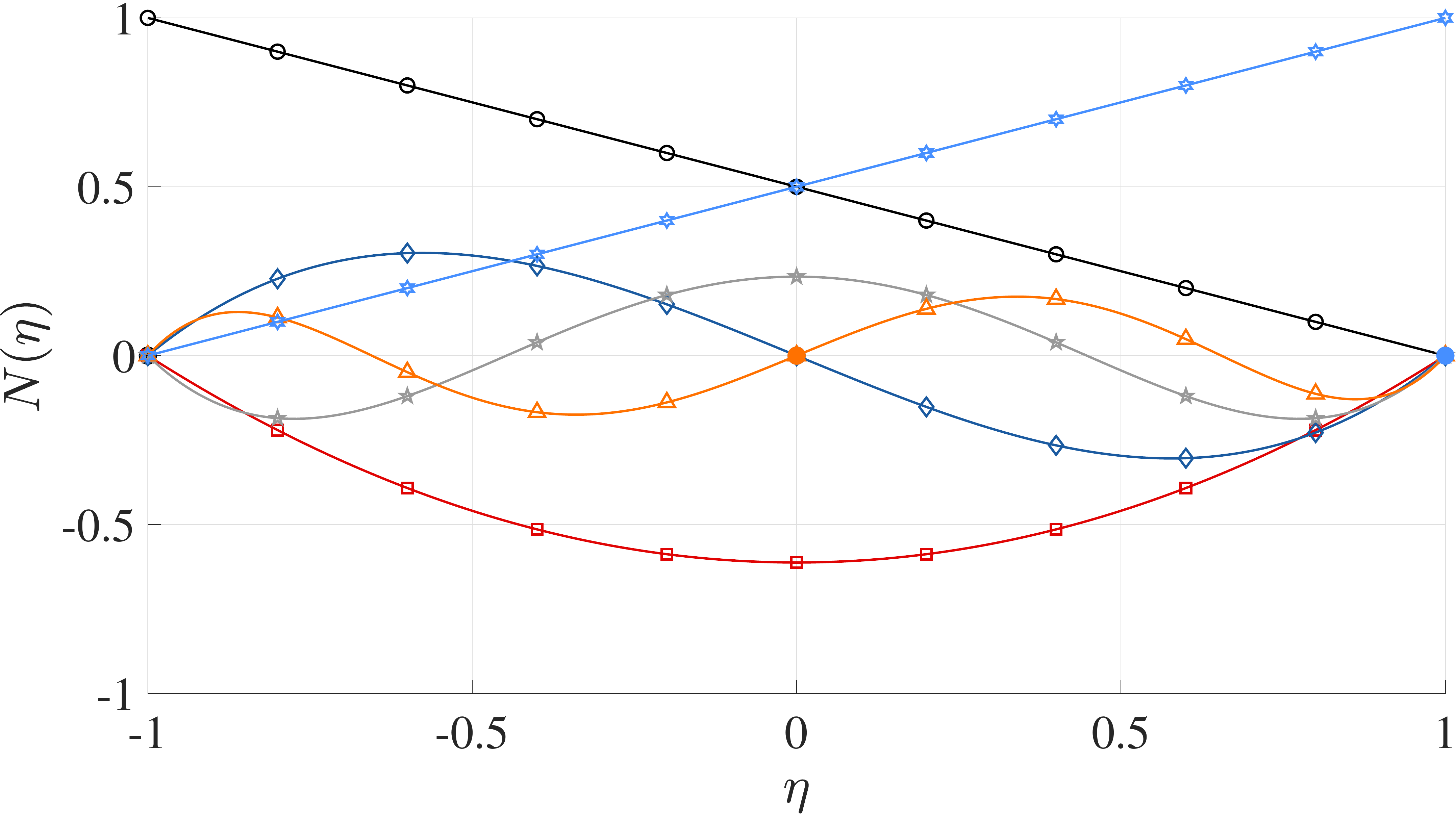}} 
	\hfill
	\subfloat[Derivatives]{\includegraphics[width=0.495\textwidth]{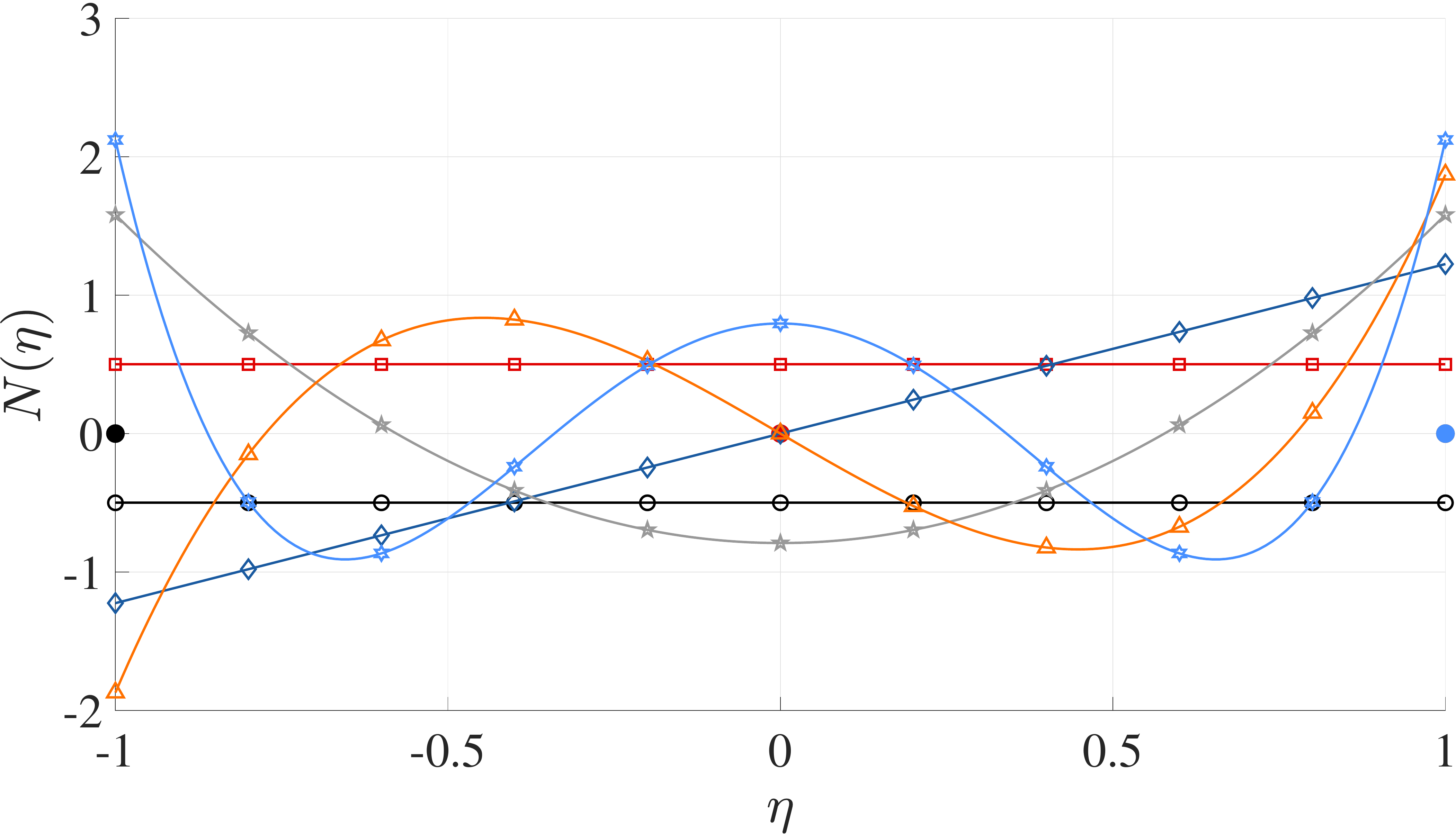}}
	\caption{\emph{p}-FEM shape functions, $p\,{=}\,5$. Legend: \textcolor{Matlab1}{\rule[0.55ex]{5ex}{0.2ex}} $N_1(\eta)$ (\NewCirc), \textcolor{Matlab2}{\rule[0.55ex]{5ex}{0.2ex}} $N_2(\eta)$ ($\square$), \textcolor{Matlab3}{\rule[0.55ex]{5ex}{0.2ex}} $N_3(\eta)$ ($\Diamond$), \textcolor{Matlab4}{\rule[0.55ex]{5ex}{0.2ex}} $N_4(\eta)$ (\NewStar), \textcolor{Matlab5}{\rule[0.55ex]{5ex}{0.2ex}} $N_5(\eta)$ ($\bigtriangleup$), \textcolor{Matlab6}{\rule[0.55ex]{5ex}{0.2ex}} $N_6(\eta)$ ($\mathbf{\ast}$). \label{fig:shpFunHier}}
\end{figure}%
Figure~\ref{fig:shpFunHier} shows the hierarchical shape functions up to an order of $p=5$.
The Legendre polynomials can be computed using Bonnet's recursion expression~\cite{BookDuester2002} or Rodriguez' formula
\begin{equation}
L_i(x) = \cfrac{1}{2^ii!} \cfrac{\mathrm{d}^i}{\mathrm{d}x^i}[(x^2-1)^i]\,, \quad i \in \mathbb{N}_0
\label{eq:Legendre}
\end{equation}
Note that only the points at the interval limits possess the Kronecker delta property and therefore, only DOFs associated with these functions (modes) retain a physical meaning, i.e., the linear shape functions alone possess both the Kronecker delta and the partition of unity properties. On the other hand, the higher order modes ensure the completeness of the ansatz for arbitrary polynomial orders. These shape functions are often referred to as modal functions~\cite{BookPozrikidis2014} and involve an additional step in the post-processing and the application of non-homogeneous Dirichlet boundary conditions (BCs). The hierarchical property implies that the order of interpolation can be increased by including additional shape functions without modifying the existing ones. The mass and stiffness matrices inherit this hierarchy \cite{BookDuester2002}.

%% file: tex/Theory.tex
\section{Derivation of Transition Elements}
\label{sec:TransitionElem}
The derivation of the shape functions for a general two-dimensional transition element, discussed in this section, is based on the concepts presented by Gordon and co-workers~\cite{ArticleGordon1971, ArtcileGordon1973a, ArtcileGordon1973b, BookGordon1982}. The basic idea is to construct bivariate interpolants that coincide with given univariate functions at the element's boundary. Gordon et al.\ referred  to this concept as \textit{transfinite interpolation}, indicating that the bivariate interpolant agrees exactly with the univariate functions at infinitely many points at the element edges. This approach is often used in the \emph{p}-FEM to ensure an accurate description of the geometry of the discretized structure \cite{BookSzabo1991, ArticleKiralyfalvi1997, PhDBroeker2001}. In the context of geometry mapping, the method is commonly referred to as \textit{blending function method}. The first use of this approach was reported by Coons \cite{TechReportCoons1967}. Based on this method, arbitrarily distorted element edges/faces can be realized. As an input, the parametric description of an element's boundary is required. If special transition elements are to be created, the functions on the individual edges need to be chosen such that a C\textsuperscript{0}-continuous approximation of the primary variables is achieved. Thus, standard isoparametric finite elements \cite{ArticleErgatoudis1968} can be combined with transition elements in a straightforward fashion. The following explanations are given for two-dimensional domains, while the extension to three-dimensional applications is straightforward. Although the general approach taken in this article is analogous to the one discussed by Gordon and co-workers \cite{ArticleGordon1971, ArtcileGordon1973a, ArtcileGordon1973b, BookGordon1982}, we extend this methodology to the conformal coupling of not only elements of different sizes but also elements with different types of shape functions. We will focus on a combination of nodal Lagrange basis functions \cite{BookKarniadakis2005, BookPozrikidis2014} and hierarchical (modal) ones based on the integrals of the Legendre polynomials \cite{BookSzabo1991, InbookDuester2018}, as already discussed in Sect.~\ref{sec:HO_ShapeFunctions}. A second important aspect of this contribution is a thorough assessment of the convergence properties of this class of transition elements. At least to the authors' knowledge, there is no concise investigation concerning the asymptotic rates of convergence related to these elements, which might also explain why they are rarely used despite their versatility.
\subsection{Transfinite interpolation}
\label{subsec:TransfiniteInterpolation}
Consider a given function $\Xi(\xi,\eta)$ that is to be interpolated over the domain $\Omega$:$\, [-1,1] \times [-1,1]$. An interpolation can generally be achieved by projecting the arbitrary bivariate function $\Xi(\xi,\eta)$ onto a chosen subspace of bivariate functions. Such a projection will be denoted as $\mathcal{P}[\Xi(\xi,\eta)]$. Let $\mathcal{P}_{\xi}[\Xi(\xi,\eta)]$ and $\mathcal{P}_{\eta}[\Xi(\xi,\eta)]$ be projectors that interpolate $\Xi(\xi,\eta)$ along the $\xi$- and $\eta$-directions, respectively. A suitable projection that interpolates $\Xi(\xi,\eta)$ over the domain $\Omega$ can then be defined as~\cite{ArtcileGordon1973b}
\begin{equation}
\mathcal{P}[\Xi(\xi,\eta)] = \mathcal{P}_{\xi}[\Xi(\xi,\eta)] \oplus \mathcal{P}_{\eta}[\Xi(\xi,\eta)] = \mathcal{P}_{\xi}[\Xi(\xi,\eta)] + \mathcal{P}_{\eta}[\Xi(\xi,\eta)] - \mathcal{P}_{\xi}[\mathcal{P}_{\eta}[\Xi(\xi,\eta)]]\,.
\label{eq:Projector}
\end{equation}
In Ref.~\cite{ArtcileGordon1973b}, the $\oplus$-operator introduced in Eq.~\eqref{eq:Projector} is referred to as \textit{Boolean sum}. In the following derivations, $\Xi(\xi,\eta)$ is assumed to be a scalar-valued function representing the primary variable (unknown). An analogous approach is taken for vector-valued functions, e.g., a displacement field. However, since it is not common practice to deploy different types of functions for the displacements in the different coordinate directions, it is sufficient to derive the necessary formulas based on a scalar function $\Xi(\xi,\eta)$. The last term in Eq.~\eqref{eq:Projector} -- $\mathcal{P}_{\xi}[\mathcal{P}_{\eta}[\square]]$ -- is often referred to as product/mixed projection operator. Note that the individual projection operators are both linear
\begin{equation}
\mathcal{P}_s[f(\xi,\eta) + g(\xi,\eta)] = \mathcal{P}_s[f(\xi,\eta)] + \mathcal{P}_s[g(\xi,\eta)]
\label{eq:ProjectorLinear}
\end{equation}
and idempotent
\begin{equation}
\mathcal{P}_s[\mathcal{P}_s[f(\xi,\eta)]] = \mathcal{P}_s[f(\xi,\eta)]\,,
\label{eq:ProjectorIdempotent}
\end{equation}
where $f(\xi,\eta)$ and $g(\xi,\eta)$ are (continuous) bivariate functions and the subscript $\square_s$ denotes the coordinate direction ($s\in \{\xi,\eta\}$).
\subsection{Transfinite bivariate Lagrange interpolation}
\label{subsec:TransfiniteLagrangeInterpolation}
\begin{figure}[t!]
    \centering
    \subfloat[Linear blending\label{fig:IntLinBlending}]{\includegraphics[clip,width=0.315\textwidth]{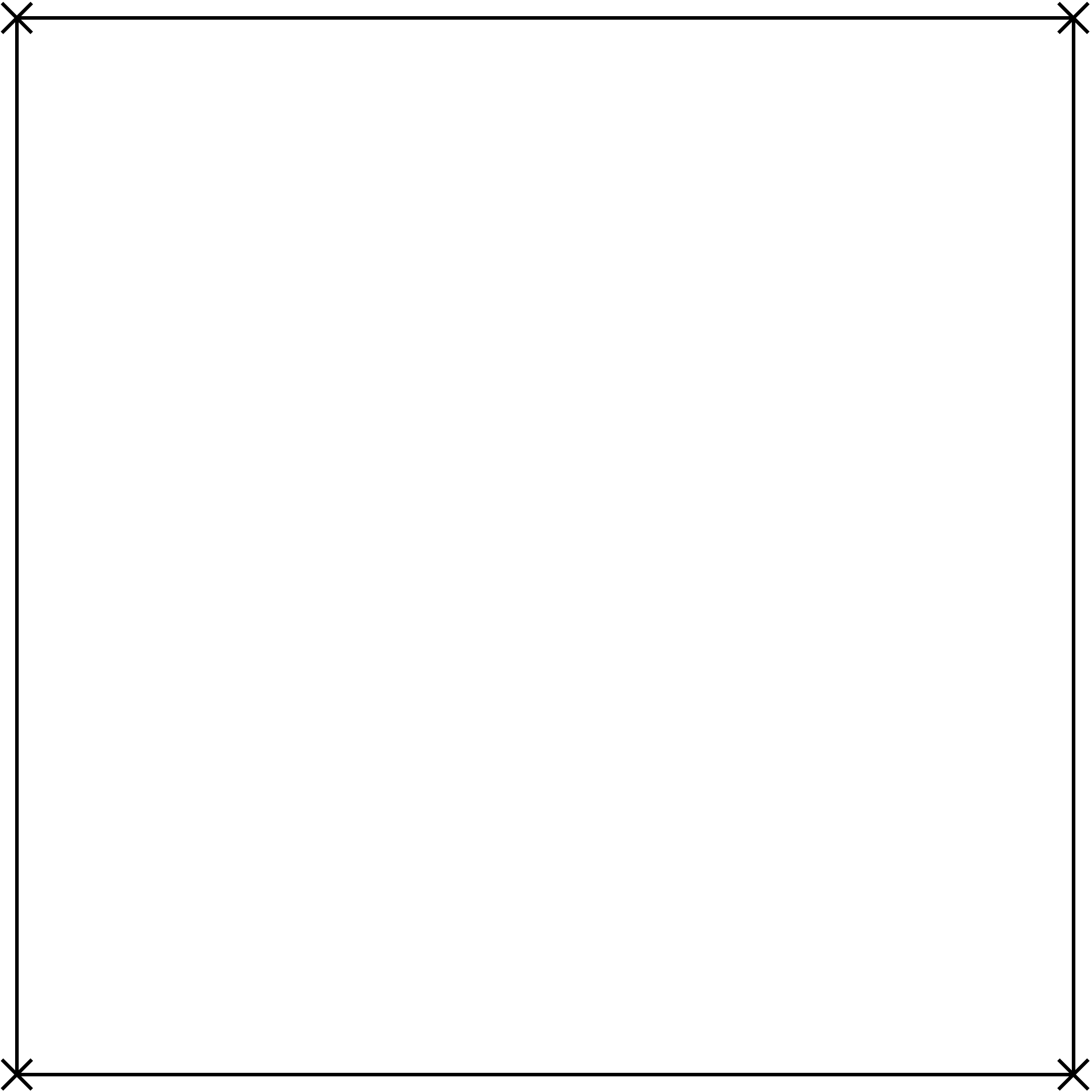}}
    \hfill
    \subfloat[Quadratic blending\label{fig:IntQuadBlending}]{\includegraphics[clip,width=0.315\textwidth]{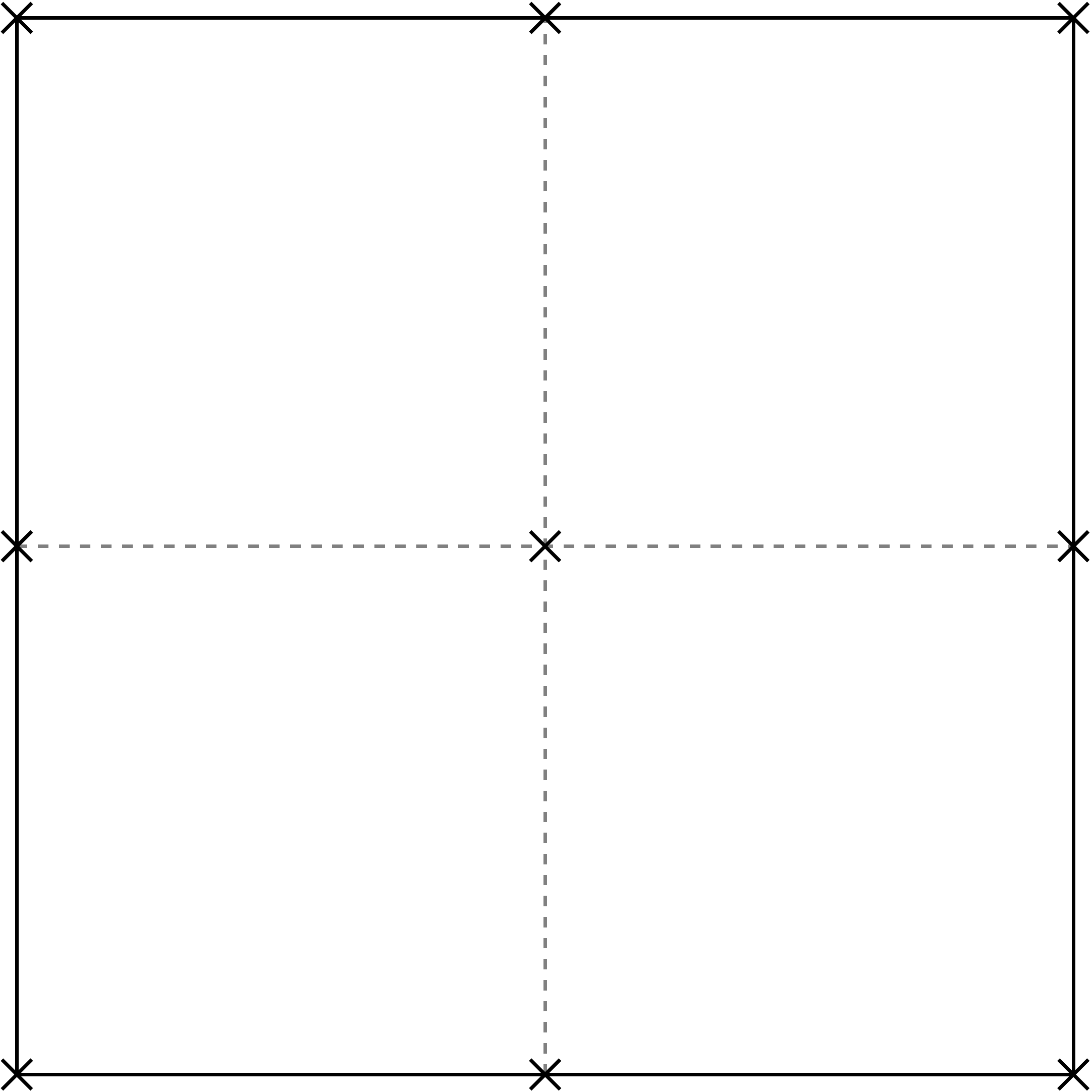}}
    \hfill
    \subfloat[Cubic blending\label{fig:IntCubBlending}]{\includegraphics[clip,width=0.315\textwidth]{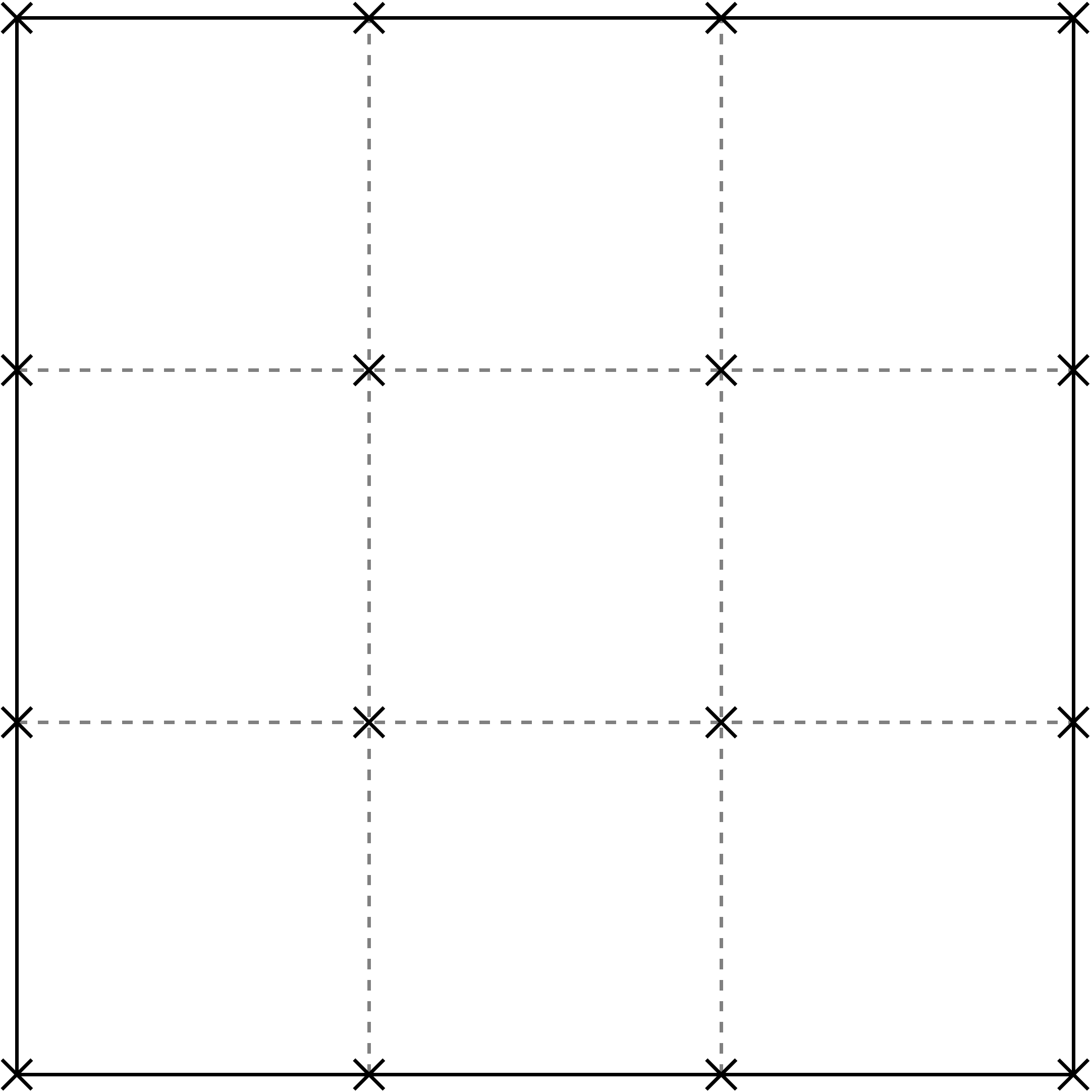}}
	\caption{Interpolation of the function $\bar{u}(\xi,\eta)$ for different blending functions (shown on the reference element in the $\xi\eta$-domain $\Omega$:$\, [-1,1] \times [-1,1]$). Depending on the polynomial degree of the blending functions, $\bar{u}(\xi,\eta)$ is not only interpolated at the boundary of the element but also at the dashed lines in the interior. This interpolation is achieved by the projection operators $\mathcal{P}_{\xi}[\bar{u}(\xi,\eta)]$ and $\mathcal{P}_{\eta}[\bar{u}(\xi,\eta)]$. The $\times$-markers indicate the locations where the functions are interpolated by the product projector $\mathcal{P}_{\xi}[\mathcal{P}_{\eta}[\bar{u}(\xi,\eta)]]$.}
	\label{fig:IntBlending}
\end{figure}

In this section, a special class of projection operators based on Lagrangian interpolation polynomials is introduced. The projectors in the individual coordinate directions are defined by the following expressions
\begin{align}
\mathcal{P}_{\xi}[\Xi(\xi,\eta)] & = \sum\limits_{i=1}^{p_{\xi}^\mathrm{p}+1} \varphi_i(\xi) \Xi(\xi_i,\eta)\,, \label{eq:ProjectorLagrangeXi}\\
\mathcal{P}_{\eta}[\Xi(\xi,\eta)] & = \sum\limits_{i=1}^{p_{\eta}^\mathrm{p}+1} \psi_i(\eta) \Xi(\xi,\eta_i)\,, \label{eq:ProjectorLagrangeEta}
\end{align}
where $\varphi_i(\xi)$ and $\psi_i(\eta)$ are referred to as blending functions which are given by univariate Lagrangian interpolation polynomials. We note that for the application of the projection operators, it is sufficient to know the function values at discrete curves (lines) in the interior of a finite element and at its boundary. Depending on the polynomial degree of the blending functions, the amount of information (data) that is required in the interior of an element differs (see also Fig.~\ref{fig:IntBlending}). In principle, the polynomial degree of the blending functions can be different in both directions ($p_{\xi}^\mathrm{p} \ne p_{\eta}^\mathrm{p}$). The definition of the product/mixed projection operator immediately follows from Eqs.~\eqref{eq:ProjectorLagrangeXi} and \eqref{eq:ProjectorLagrangeEta} as
\begin{equation}
\mathcal{P}_{\xi}[\mathcal{P}_{\eta}[\Xi(\xi,\eta)]] = \sum\limits_{i=1}^{p_{\xi}^\mathrm{p}+1} \sum\limits_{j=1}^{p_{\eta}^\mathrm{p}+1} \varphi_i(\xi)  \psi_j(\eta) \Xi(\xi_i,\eta_j)\,. \label{eq:ProjectorLagrangeProduct}
\end{equation}
From Eq.~\eqref{eq:ProjectorLagrangeProduct} we can infer that the product projector interpolates the function only at discrete nodal locations, whereas the individual projectors take certain lines into account (see Fig.~\ref{fig:IntBlending} -- interior: dashed lines, boundary: solid lines). In most applications, linear blending functions are deployed \cite{BookSzabo1991} and therefore, $\varphi_i(\xi)$ and $\psi_i(\eta)$ are given as
\begin{alignat}{2}\label{eq:linearBlending1}
\varphi_1(\xi) &= N_1^1(\xi) = \cfrac{1}{2}(1-\xi),\qquad \psi_1(\eta) &= N_1^1(\eta) = \cfrac{1}{2}(1-\eta)\,, \\ \label{eq:linearBlending2}
\varphi_2(\xi) &= N_2^1(\xi) = \cfrac{1}{2}(1+\xi),\qquad \psi_2(\eta) &= N_2^1(\eta) = \cfrac{1}{2}(1+\eta)\,.
\end{alignat}
The superscript $\square^1$ indicates that the polynomial degree is selected as $p\,{=}\,1$. Blending functions of arbitrary order $\square^p$ can easily be generated using standard Lagrangian interpolation techniques which have been discussed in Sect.~\ref{sec:SEM}. The nodal distribution that is used to derive the blending functions is principally arbitrary, but we assume an \textit{equidistant spacing} in the remainder of the article. 
\subsection{Derivation of shape functions for transition elements}
\label{subsec:ShapeFuncTransElem}
Based on the concept of transfinite bivariate Lagrange interpolation (see Sect.~\ref{subsec:TransfiniteLagrangeInterpolation}) we are now able to derive the shape functions for arbitrary transition elements. In the following, we will concentrate on two-dimensional displacement-type finite elements. A C\textsuperscript{0}-continuous displacement field between adjacent elements is achieved if the shape functions coincide at the common boundaries, i.e., on the elements' edges. Therefore, it suffices to employ linear blending functions in the derivation of the shape functions for arbitrary transition elements. As a result, the interpolated (projected) function, which is used to derive the shape functions, matches the original one (function that was projected) at the boundary of the element. This fact is illustrated in Fig.~\ref{fig:IntBlending}. Consequently, high order blending terms are not considered in the examples used in this contribution.

\begin{figure}[b!]
	\begin{minipage}[t]{0.45\textwidth}
		\centering
		\includegraphics[clip,width=0.733\textwidth]{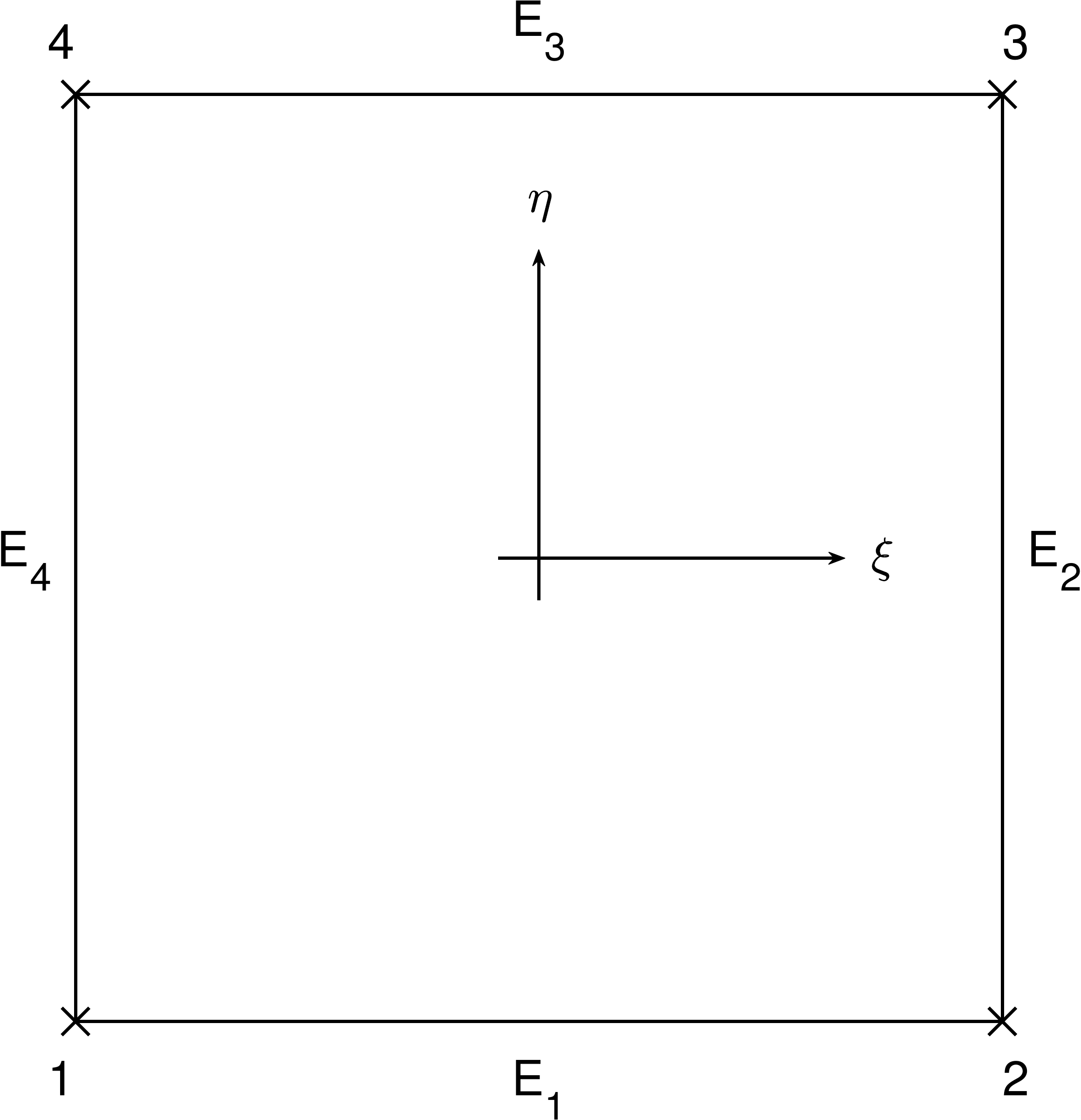}
		\caption{Reference element with edge numbering.}
		\label{fig:RefElem}
	\end{minipage}
	\hfill
	\begin{minipage}[t]{0.45\textwidth}
		\centering
		\includegraphics[clip,width=0.733\textwidth]{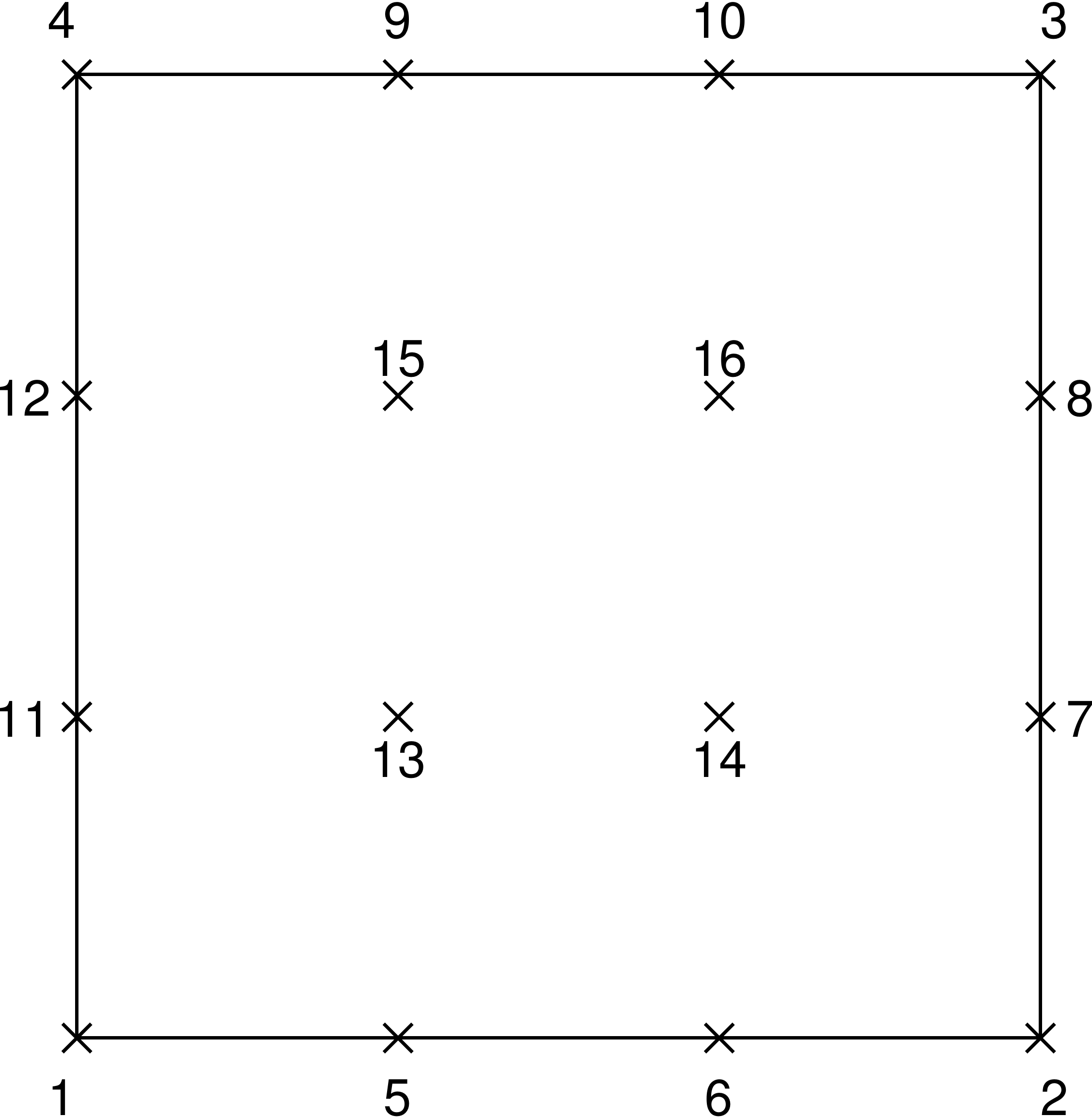}
		\caption{Numbering scheme of the DOFs (exemplarily given for a cubic finite element).}
		\label{fig:NumberingScheme}
	\end{minipage}
\end{figure}

The point of departure for the derivation of the shape functions of arbitrary transition elements are the shape functions of a standard high order finite element, which should be complete polynomials of order $p$. In our implementation, the basis element is either a \emph{p}-element \cite{BookSzabo1991} or a spectral element \cite{BookPozrikidis2014} based on a tensor product formulation. Keep in mind that only specific shape functions belonging to the boundary of an element are required for the coupling of either different types of elements or multiple elements to only one element. Therefore, it is sufficient to adjust the nodal (corner nodes) and edge shape functions that belong to the part of the boundary that is coupled while leaving the other edge and bubble/interior modes untouched. In the most general case, all four edges are involved and therefore, shape functions associated with all corner vertices and edges are subject to changes. As mentioned before, we need a description of the function that is to be interpolated on the boundary of the transition element. This function, in the remainder of the article referred to as $\bar{u}(\xi,\eta)$, is given in terms of one-dimensional shape functions on the four edges of a quadrilateral element (see Fig.~\ref{fig:RefElem})
\begin{equation}
\bar{u}(\xi,\eta) =
\begin{cases}
\bar{u}_\mathrm{E_1}(\xi)  = \sum\limits_{i=1}^{n_{\mathrm{E}_1}} \bar{N}_i(\xi) u_\mathrm{E_1}^{(i)},    & \text{for } \eta=-1 \,, \\
\bar{u}_\mathrm{E_2}(\eta) = \sum\limits_{i=1}^{n_{\mathrm{E}_2}} \hat{N}_i(\eta) u_\mathrm{E_2}^{(i)},   & \text{for } \xi=1   \,, \\
\bar{u}_\mathrm{E_3}(\xi)  = \sum\limits_{i=1}^{n_{\mathrm{E}_3}} \tilde{N}_i(\xi) u_\mathrm{E_3}^{(i)},  & \text{for } \eta=1  \,, \\
\bar{u}_\mathrm{E_4}(\eta) = \sum\limits_{i=1}^{n_{\mathrm{E}_4}} \breve{N}_i(\eta) u_\mathrm{E_4}^{(i)}, & \text{for } \xi=-1  \,,
\end{cases}
\label{eq:EdgeFunctions}
\end{equation}
where $\bar{N}_i(t)$, $\hat{N}_i(t)$, $\tilde{N}_i(t)$, and $\breve{N}_i(t)$ denote different standard or piecewise shape functions\footnote{Note that the edge shape functions do not need to be polynomial. In principle, it is possible to assume any type of shape functions such as trigonometric, splines, etc. on the boundary. In the remainder of the article, we will, however, concentrate on (complete) polynomial shape functions.}. This function describes only the displacements at the boundary of the element which is sufficient for the conformal/compatible coupling of different element types and numbers. The variables $n_{\mathrm{E}_i}$, where $i\in\{1,2,3,4\}$, denote the number of shape functions associated with the corresponding edge including the shape functions of the corner nodes. Observe that in the introduced numbering scheme, the shape functions associated with the corner nodes/vertices are always connected to the index values $i\,{=}\,1$ and $i\,{=}\,n_{\mathrm{E}_j}$. The parameters $u_\mathrm{E_i}^{(j)}$ denote the $j$\textsuperscript{th} degree of freedom belonging to the $i$\textsuperscript{th} edge of the finite element. In \ref{sec:12NodeElem}, we will come back to these one-dimensional shape functions and present a specific example to illustrate their use. The generic form of the displacement field including all shape functions of the transition element, i.e., the functions defined on the boundary and in the interior of the domain, is given as
\begin{equation}
u(\xi,\eta) = \underbrace{\mathcal{P}_{\xi}[\bar{u}(\xi,\eta)] \oplus \mathcal{P}_{\eta}[\bar{u}(\xi,\eta)]}_{u^\circ(\xi,\eta)} \;+\; \text{internal shape functions,}
\label{eq:DispTransitionElement}
\end{equation}
where $u^\circ(\xi,\eta)$ is the part of the displacement field related to the shape functions that are purely defined on the edges of the element and blended into the interior. Since we restrict ourselves to linear blending functions, $u^\circ(\xi,\eta)$ is given as
\begin{equation}
\begin{split}
u^\circ(\xi,\eta) = 
& \underbrace{N_2^1(\xi)\bar{u}_\mathrm{E_2}(\eta) + N_1^1(\xi)\bar{u}_\mathrm{E_4}(\eta)}_{\mathcal{P}_{\xi}[\bar{u}(\xi,\eta)]} +
  \underbrace{N_1^1(\eta)\bar{u}_\mathrm{E_1}(\xi) + N_2^1(\eta)\bar{u}_\mathrm{E_3}(\xi)}_{\mathcal{P}_{\eta}[\bar{u}(\xi,\eta)]} - \\
& \underbrace{\left[ N_1^1(\xi)N_1^1(\eta) u^{(1)} + N_2^1(\xi)N_1^1(\eta)u^{(2)} + N_2^1(\xi)N_2^1(\eta)u^{(3)} + N_1^1(\xi)N_2^1(\eta)u^{(4)} \right]}_{\mathcal{P}_{\xi}[\mathcal{P}_{\eta}[\bar{u}(\xi,\eta)]]}\,,
\end{split}
\label{eq:EdgeShapeFunctions}
\end{equation}
with the linear interpolation/blending functions defined in Eqs.~\eqref{eq:linearBlending1} and \eqref{eq:linearBlending2}. Here, the coefficients $u^{(i)}$ denote the values of the primary unknowns (structural analysis: displacements DOFs) at the four corner vertices. In the next step, we substitute the definition of the edge shape functions provided in Eq.~\eqref{eq:EdgeFunctions} into Eq.~\eqref{eq:EdgeShapeFunctions}
\begin{equation}
\begin{split}
u^\circ(\xi,\eta) = 
& \underbrace{N_2^1(\xi)\sum\limits_{i=1}^{n_{\mathrm{E}_2}} \hat{N}_i(\eta) u_\mathrm{E_2}^{(i)} + N_1^1(\xi)\sum\limits_{i=1}^{n_{\mathrm{E}_4}} \breve{N}_i(\eta) u_\mathrm{E_4}^{(i)}}_{\mathcal{P}_{\xi}[\bar{u}(\xi,\eta)]} + \underbrace{N_1^1(\eta)\sum\limits_{i=1}^{n_{\mathrm{E}_1}} \bar{N}_i(\xi) u_\mathrm{E_1}^{(i)} + N_2^1(\eta)\sum\limits_{i=1}^{n_{\mathrm{E}_3}} \tilde{N}_i(\xi) u_\mathrm{E_3}^{(i)}}_{\mathcal{P}_{\eta}[\bar{u}(\xi,\eta)]} - \\
& \underbrace{\left[ N_1^1(\xi)N_1^1(\eta) u^{(1)} + N_2^1(\xi)N_1^1(\eta)u^{(2)} + N_2^1(\xi)N_2^1(\eta)u^{(3)} + N_1^1(\xi)N_2^1(\eta)u^{(4)} \right]}_{\mathcal{P}_{\xi}[\mathcal{P}_{\eta}[\bar{u}(\xi,\eta)]]}\,.
\end{split}
\label{eq:EdgeShapeFunctions2}
\end{equation}
To extract the individual shape functions corresponding to each DOF we have to prescribe a unit value for that DOF and set all other DOFs to zero. Although the numbering scheme of the DOFs is arbitrary, we will introduce a specific order to facilitate the understanding. The chosen methodology is based on the common practice in the \emph{p}-FEM, i.e., first the DOFs of the corner vertices are numbered, then those belonging to the edges 1, 2, 3, and 4 (in this order), before the interior DOFs are taken care of. As an example, a cubic element with one DOF per node/mode is illustrated in Fig.~\ref{fig:NumberingScheme}. Special care has to be taken in the \emph{p}-FEM -- where hierarchic shape functions are employed which do not have the Kronecker delta property -- where the DOFs do not retain a physical meaning and are not connected to a certain nodal position. The illustration provided in Fig.~\ref{fig:NumberingScheme} only assigns positions to certain DOFs to facilitate the understanding. The DOFs associated with the corner nodes are always assigned the numbers  1, 2, 3, and 4. The DOFs connected to the different edges indicated by $u_\mathrm{E_i}^{(j)}$ are defined by the vectors $\mathbf{I}_i$ containing the indices of the local DOFs
\allowdisplaybreaks
\begin{alignat*}{2}
\mathbf{u}_\mathrm{E_1} &= \left[u^{(\mathbf{I}_1(1))}, u^{(\mathbf{I}_1(2))},\ldots, u^{(\mathbf{I}_1(n_{\mathrm{E}_1}))}\right]^\mathrm{T}  \quad\text{for } \mathbf{I}_1 = [&&1, 5, 6,\ldots, 4{+}(n_{\mathrm{E}_1}{-}2),2]\,, \\
\mathbf{u}_\mathrm{E_2} &= \left[u^{(\mathbf{I}_2(1))}, u^{(\mathbf{I}_2(2))},\ldots, u^{(\mathbf{I}_2(n_{\mathrm{E}_2}))}\right]^\mathrm{T}  \quad\text{for } \mathbf{I}_2 = [&&2, 4{+}(n_{\mathrm{E}_1}{-}1), 4{+}n_{\mathrm{E}_1},\ldots,\\
& && 4{+}(n_{\mathrm{E}_1}{-}2){+}(n_{\mathrm{E}_2}{-}2),3]\,, \\
\mathbf{u}_\mathrm{E_3} &= \left[(u^{\mathbf{I}_3(1))}, u^{(\mathbf{I}_3(2))},\ldots, u^{(\mathbf{I}_3(n_{\mathrm{E}_3}))}\right]^\mathrm{T}  \quad\text{for } \mathbf{I}_3 = [&&4, 4{+}(n_{\mathrm{E}_1}{-}2){+}(n_{\mathrm{E}_2}{-}1), 4{+}(n_{\mathrm{E}_1}{-}2){+}n_{\mathrm{E}_2},\ldots, \\
& && 4{+}(n_{\mathrm{E}_1}{-}2){+}(n_{\mathrm{E}_2}{-}2){+}(n_{\mathrm{E}_3}{-}2),3]\,, \\
\mathbf{u}_\mathrm{E_4} &= \left[u^{(\mathbf{I}_4(1))}, u^{(\mathbf{I}_4(2))},\ldots, u^{(\mathbf{I}_4(n_{\mathrm{E}_4}))}\right]^\mathrm{T}  \quad\text{for } \mathbf{I}_4 = [&&1, 4{+}(n_{\mathrm{E}_1}{-}2){+}(n_{\mathrm{E}_2}{-}2){+}(n_{\mathrm{E}_3}{-}1), \\
& && 4{+}(n_{\mathrm{E}_1}{-}2){+}(n_{\mathrm{E}_2}{-}2){+}n_{\mathrm{E}_3},\ldots,\\
& && 4{+}(n_{\mathrm{E}_1}{-}2){+}(n_{\mathrm{E}_2}{-}2){+}(n_{\mathrm{E}_3}{-}2){+}(n_{\mathrm{E}_4}{-}2),4]\,.
\end{alignat*}
Note that in this definition, the DOFs of the corner vertices are included in the description of the shape functions on the boundary of the transition element. Based on this definition, we can derive simple expressions for the shape functions corresponding to the corner nodes
\allowdisplaybreaks
\begin{alignat}{3}
N_1(\xi,\eta) &= N_1^1(\xi)\breve{N}_1(\eta) &&+ N_1^1(\eta)\bar{N}_1(\xi) &&- N_1^1(\xi)N_1^1(\eta)\,,\label{eq:SF_CornerNode1}\\
N_2(\xi,\eta) &= N_2^1(\xi)\hat{N}_1(\eta) &&+ N_1^1(\eta)\bar{N}_{n_{\mathrm{E}_1}}\!(\xi) &&- N_2^1(\xi)N_1^1(\eta)\,,\label{eq:SF_CornerNode2}\\
N_3(\xi,\eta) &= N_2^1(\xi)\hat{N}_{n_{\mathrm{E}_2}}\!(\eta) &&+ N_2^1(\eta)\tilde{N}_{n_{\mathrm{E}_3}}\!(\xi) &&- N_2^1(\xi)N_2^1(\eta) \,,\label{eq:SF_CornerNode3}\\
N_4(\xi,\eta) &= N_1^1(\xi)\breve{N}_{n_{\mathrm{E}_4}}\!(\eta) &&+ N_2^1(\eta)\tilde{N}_1(\xi) &&- N_1^1(\xi)N_2^1(\eta) \label{eq:SF_CornerNode4}\,,
\end{alignat}
and also for the shape functions that are associated with only one edge (edge modes)
\allowdisplaybreaks
\begin{align}
\text{E}_1:\; N_i(\xi,\eta) &= N_1^1(\eta) \bar{N}_j(\xi)\,, \label{eq:SF_Edge1}\\
\text{for } i &= 4 + (j-1)\,, \nonumber\\
j &= 2, 3,\ldots, n_{\mathrm{E}_1}{-}1\,, \nonumber\\
\text{E}_2:\; N_i(\xi,\eta) &= N_2^1(\xi) \hat{N}_j(\eta)\,, \label{eq:SF_Edge2}\\
\text{for } i &= 4 + (n_{\mathrm{E}_1}{-}2) + (j-1)\,, \nonumber\\
j &= 2, 3,\ldots, n_{\mathrm{E}_2}{-}1\,, \nonumber\\
\text{E}_3:\; N_i(\xi,\eta) &= N_2^1(\eta) \tilde{N}_j(\xi)\,, \label{eq:SF_Edge3}\\
\text{for } i &= 4 + (n_{\mathrm{E}_1}{-}2) + (n_{\mathrm{E}_2}{-}2) + (j-1)\,, \nonumber\\
j &= 2, 3,\ldots, n_{\mathrm{E}_3}{-}1\,, \nonumber\\
\text{E}_4:\; N_i(\xi,\eta) &= N_1^1(\xi) \breve{N}_j(\xi)\,, \label{eq:SF_Edge4}\\
\text{for } i &= 4 + (n_{\mathrm{E}_1}{-}2) + (n_{\mathrm{E}_2}{-}2) + (n_{\mathrm{E}_3}{-}2) + (j-1)\,, \nonumber\\
j &= 2, 3,\ldots, n_{\mathrm{E}_4}{-}1\,.\nonumber
\end{align}
We need to remark that the Kronecker delta property of nodal shape functions is lost for the interior nodes when this procedure is applied. This is because only the shape functions belonging to entities on the boundary (2D: edges, 3D: faces) are changed, while the interior ones remain unaltered. Consequently, similar approaches known from the post-processing procedures of the \emph{p}-FEM have to be used for the proposed transition elements. Furthermore, we notice once again that only the shape functions associated with the corner vertices are affected by the product projection operator $\mathcal{P}_{\xi}[\mathcal{P}_{\eta}[\square]]$. Shape functions that are connected to only one edge are  adjusted/changed by only one of the projectors $\mathcal{P}_{\xi}[\square]$ or $\mathcal{P}_{\eta}[\square]$.

\subsection{Guideline}
In the following, the methodology for the derivation of the shape functions of arbitrary transition elements is summarized in a few bullet points. As mentioned at the beginning of this section, we limit ourselves to linear blending functions as they are sufficient to guarantee a C\textsuperscript{0}-continuous approximation. Deploying high order blending functions is, however, very similar -- one difference being that the product projector is applied to a larger set of vertices/nodes. For the scheme based on linear blending functions, the general procedure is described as follows:
\begin{enumerate}
	\item Define the base finite element including its shape functions.
	\item Define the shape functions on the edges of the finite element:
	\begin{itemize}
		\item Provide the shape functions only for edges that are involved in the coupling (of different element types or numbers of elements).
		\item Edges that retain the standard basis functions are not included in the projection.
	\end{itemize}
	\item Apply the projection operators:
	\begin{itemize}
		\item The product projection operator $\mathcal{P}_{\xi}[\mathcal{P}_{\eta}[\square]]$ only affects the shape functions of the four corner nodes of a quadrilateral finite element. This operator is applied to all corner nodes that belong to a coupling edge.
		\item  The projection in the $\xi$-direction is achieved using $\mathcal{P}_{\xi}[\square]$. This operator is applied to all nodes/modes that belong to either edge $E_2$ or $E_4$, provided that the corresponding edge is a coupling edge.
		\item  The projection in the $\eta$-direction is achieved using $\mathcal{P}_{\eta}[\square]$. This operator is applied to all nodes/modes that belong to either edge $E_1$ or $E_3$, provided that the corresponding edge is a coupling edge.
	\end{itemize}
	\item Standard shape functions:
	\begin{itemize}
		\item All shape functions that belong to edges that are not used to couple different element types or numbers of elements are adopted from the base finite element.
		\item Interior (bubble) shape functions are not adjusted at all in the linear blending approach.
	\end{itemize}
\end{enumerate}
It is briefly mentioned at this point that, considering the geometry approximation, all standard concepts known from FEM such as iso- sub-, superparametric mapping \cite{BookZienkiewicz2000a} or blending functions \cite{BookSzabo1991} can be employed in a straightforward fashion.

The derivation of a 12-node/mode piecewise bi-quadratic transition element is discussed in detail in \ref{sec:12NodeElem}, where different shape function types are employed. The application of the devised methodology is presented for Lagrange-Lagrange, Lagrange-Legendre, and Legendre-Legendre transition elements.
\subsection{Numerical integration}
\label{subsec:NumInt}
\begin{figure}[b!]
    \centering
    \subfloat[\emph{x\textbf{N}y}-mesh]{\includegraphics[clip,width=0.433\textwidth]{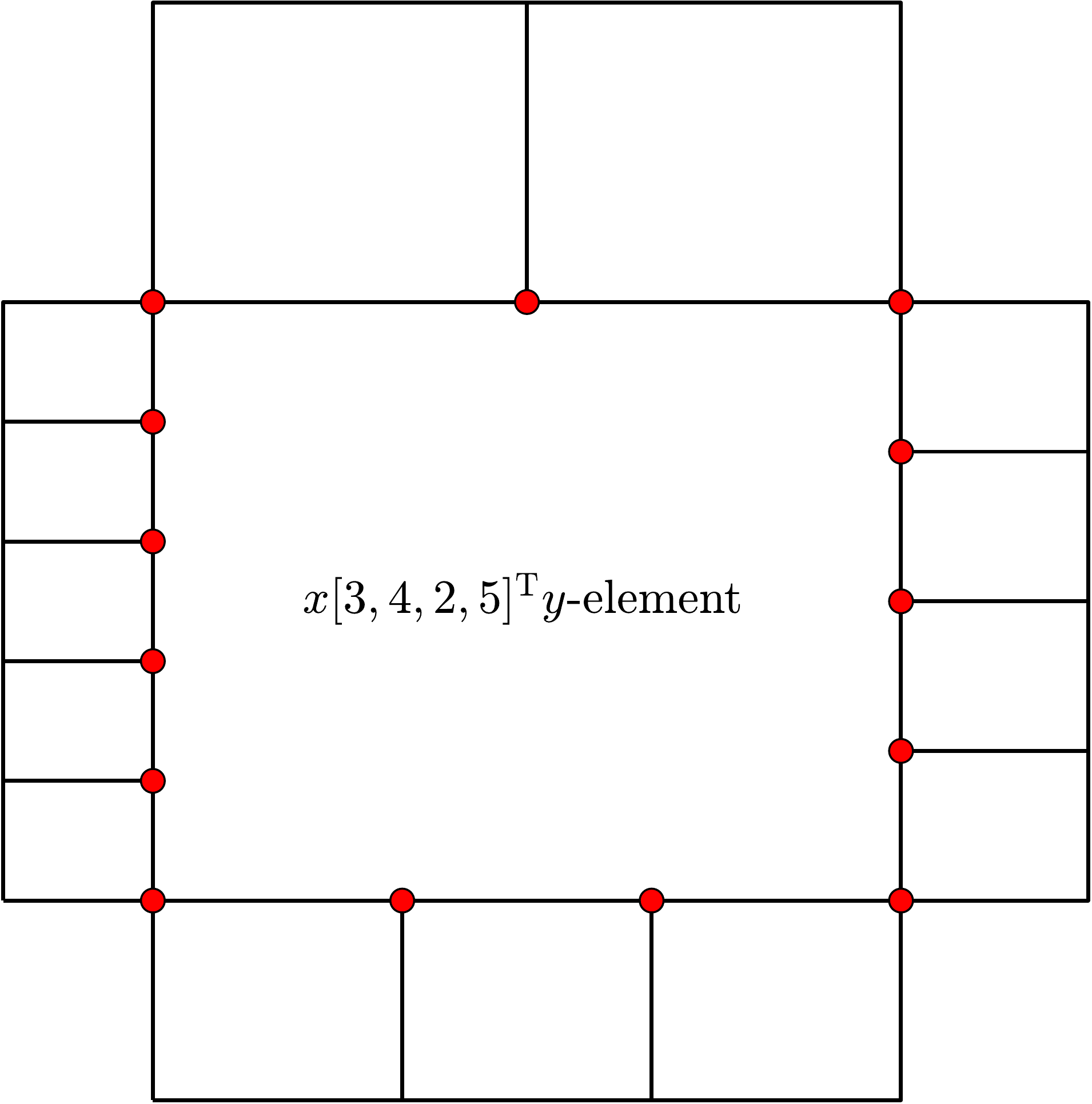}}
    \hfill
    \subfloat[Integration points and subdomains]{\includegraphics[clip,width=0.433\textwidth]{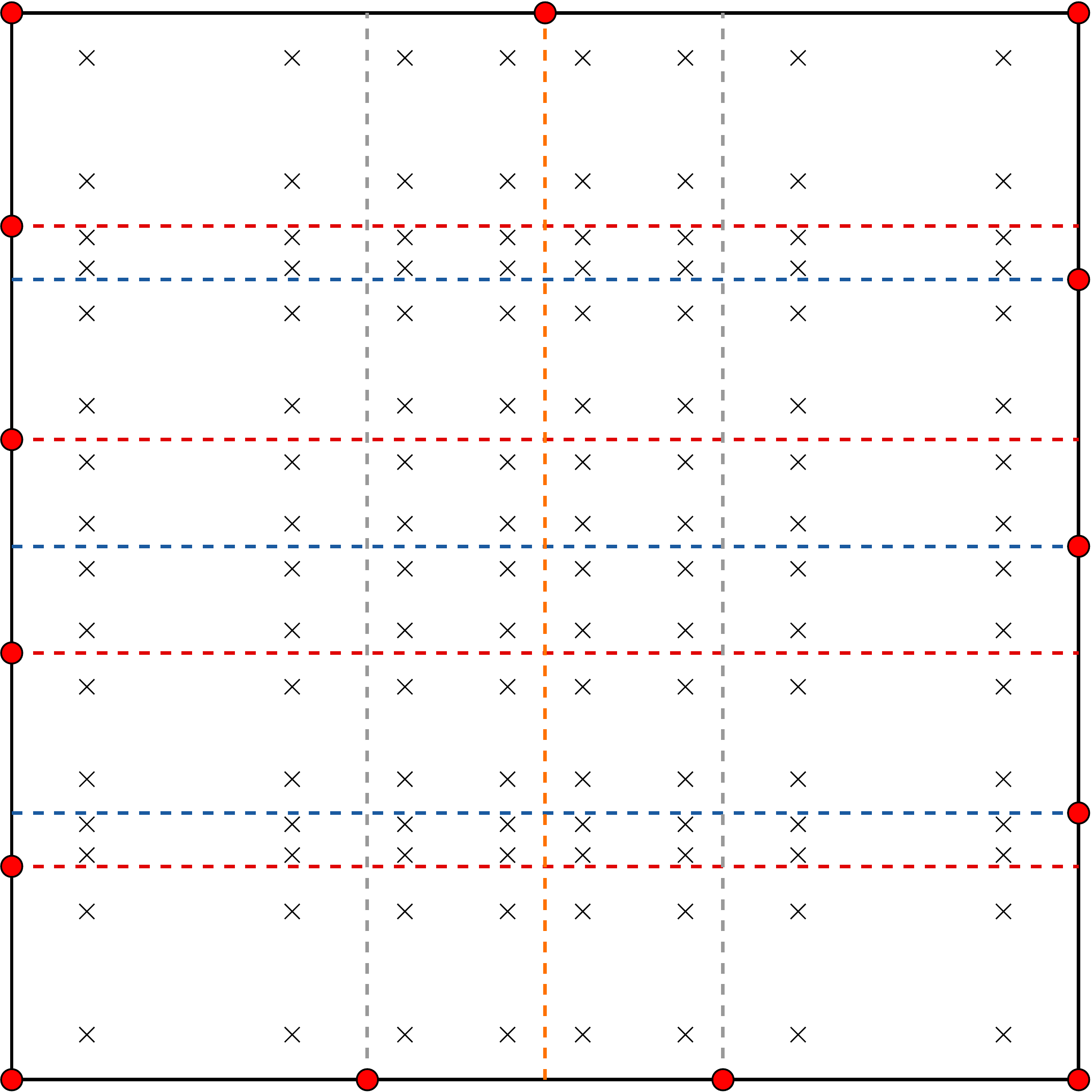}}
	\caption{Numerical integration approach for arbitrary \emph{x\textbf{N}y}-elements.}
	\label{fig:NumInt_xNy}
\end{figure}%
The numerical integration of transition elements proves to be a bit more involved compared to standard high order finite elements. Since the shape functions of \emph{x\textbf{N}y}-elements are only piecewise C\textsuperscript{0}-continuous within the element (weak discontinuity -- kink in the displacement field) their first derivatives are discontinuous. Thus, a standard Gau\ss{} quadrature over the reference element will not suffice and yield inaccurate results. The remedy for this problem is rather simple and can easily be implemented in any existing finite element code. It involves a subdivision of the integration domain according to the number of elements that are coupled to each edge of the \emph{x\textbf{N}y}-elements. Since the adjacent \emph{y}-elements are also only  C\textsuperscript{0}-continuous we know \textit{a priori} where the discontinuities arise and therefore, the element can be divided for integration purposes at those positions. The basic idea is illustrated in Fig.~\ref{fig:NumInt_xNy} for the most general case where each edge couples a different number of elements. The reference element is subdivided into several integration subdomains according to element sizes of the neighboring elements. In order to illustrate the subdivision process, the corner vertices of the \emph{y}-elements are denoted by red circles. Using this approach, it is ensured that in each integration subdomain, the shape functions do not exhibit any discontinuities and hence smooth functions are integrated. As indicated in Fig.~\ref{fig:NumInt_xNy}, a standard Gau\ss{} quadrature is applied for each subdomain and finally the results are added. For the sake of clarity, the numerical integration is indicated by a $2\times2$-point Gau\ss{} rule in each subdomain. Due to the subdivision of the finite element in the reference domain, the integration points need to be mapped from the canonical domain $[1\times1]^d$ to the integration subdomain. This mapping is achieved by a simple scaling procedure\footnote{Remark: The computation of the locations of the Gau\ss{} points requires a mapping from the canonical domain $[1\times1]^d$ to the integration domain. Both domains are rectangular and therefore, a linear scaling procedure is sufficient, and we do not need to map the points using the shape functions of the  \emph{x\textbf{N}y}-element.}. An analogous approach is taken in fictitious domain methods such as the FCM \cite{ArticleDuester2008} and in polygonal FEMs \cite{ArticleDuczek2016b}. The main advantage of this method is its robustness and suitability for all different cases of transition elements.

%% file: tex/PatchTest.tex
\section{Patch Test}
\label{sec:PatchTest}
In the finite element community, the patch test is a widely used and accepted approach to verify the implementation of displacement-type elements \cite{ArticleTaylor1986}. However, according to Zienkiewicz \& Taylor~\cite{BookZienkiewicz2000a} (Ch.~10), it can provide much more information about the element. The versatility of the patch test can briefly be summarized in the following points: It is
\begin{enumerate}
\item a necessary (sometimes sufficient) condition for assessing the convergence,
\item a means to investigate the (asymptotic) rate of convergence,
\item a check on the robustness of the element, and
\item a possibility to verify the correct programming.
\end{enumerate}
\subsection{Versions of the patch test}
\label{sec:PatchTestVersions}
The original idea of the patch test was introduced by Irons \cite{InproceedingsIrons1966} in order to check whether a finite element could exactly reproduce a constant stress state. Today, there are three different versions of the patch test, i.e., A, B, and C \cite{ArticleTaylor1986}:
\begin{itemize}
\item \textbf{Version A}: All nodes are constrained, i.e., the exact displacement field is prescribed at each node in the computational domain. To this end, we first need to compute an arbitrary solution to a partial differential equation. Then we insert these exact values into the system of equations and compute the difference between the external force vector and the internal force vector (due to the prescribed displacement field). If the method is correctly implemented, the residuum must be equal to zero. In a generic way, this can be expressed as
\begin{equation*}
K_{ij}u^\mathrm{pre}_j - f_i = 0\,,
\end{equation*}
where $\mathbf{K}$ is the stiffness matrix of the patch, $\mathbf{u}^\mathrm{pre}$ denotes the prescribed displacement field, and $\mathbf{f}$ represents the external force vector.
\item \textbf{Version B}: All nodes at the boundary are constrained, i.e., the displacement field is only prescribed at the boundary of the model. The known values are inserted in the system of equations and therefore, an additional force term exists. The system of equations is solved for the unknown displacements, and the numerical solution is compared to the exact/analytic one
\begin{equation*}
\tilde{\mathbf{u}} = \tilde{\mathbf{K}}^{-1} \left( \tilde{\mathbf{f}} - \bar{\mathbf{K}}\mathbf{u}^\mathrm{pre} \right)\,.
\end{equation*}
The tilde-operator $\tilde{\square}$ indicates that the prescribed boundary conditions have been incorporated in the system of equations, while the bar-operator $\bar{\square}$ denotes that only those columns of the stiffness matrix with prescribed boundary conditions are taken into account.
\item \textbf{Version C}: Minimum number of constraints to avoid rigid body motions are provided (global stiffness matrix is not singular; retains full rank). Neumann boundary conditions are applied such that the exact solution yields, if possible, a polynomial displacement field
\begin{equation*}
\mathbf{K}\mathbf{u} = \mathbf{f}\,.
\end{equation*}
The computed displacements are then compared against the exact solution. It can be stated that from versions A to C, the significance of the patch test increases, and version C is the most general test that provides an indication of the robustness of the finite element implementation.
\end{itemize}
\subsection{Model definition}
\label{sec:ModelDef}
\begin{figure}[t!]
	\centering
	\subfloat[$n_\mathrm{s}\,{=}\,1$, $n_\mathrm{y}\,{=}\,1$]{\includegraphics[clip,width=0.475\textwidth]{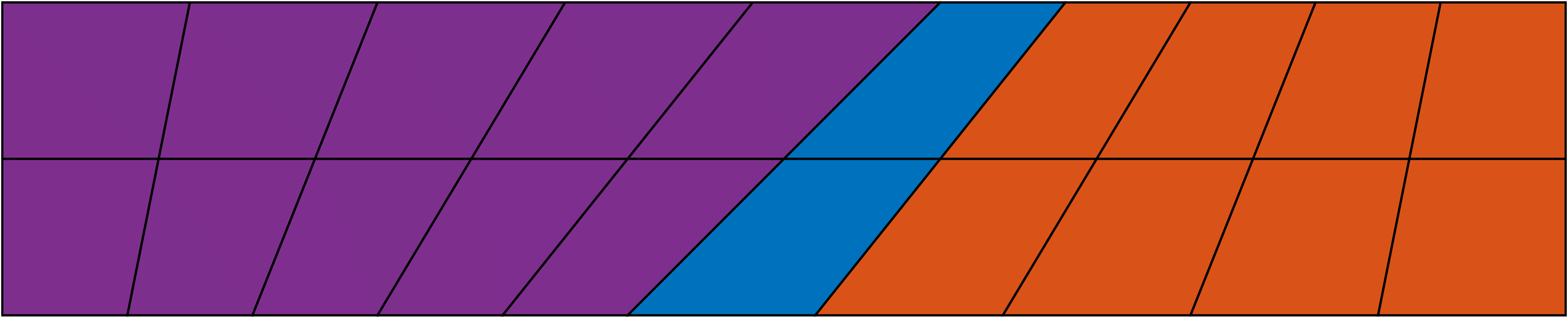}}
	\hfill
	\subfloat[$n_\mathrm{s}\,{=}\,1$, $n_\mathrm{y}\,{=}\,2$]{\includegraphics[clip,width=0.475\textwidth]{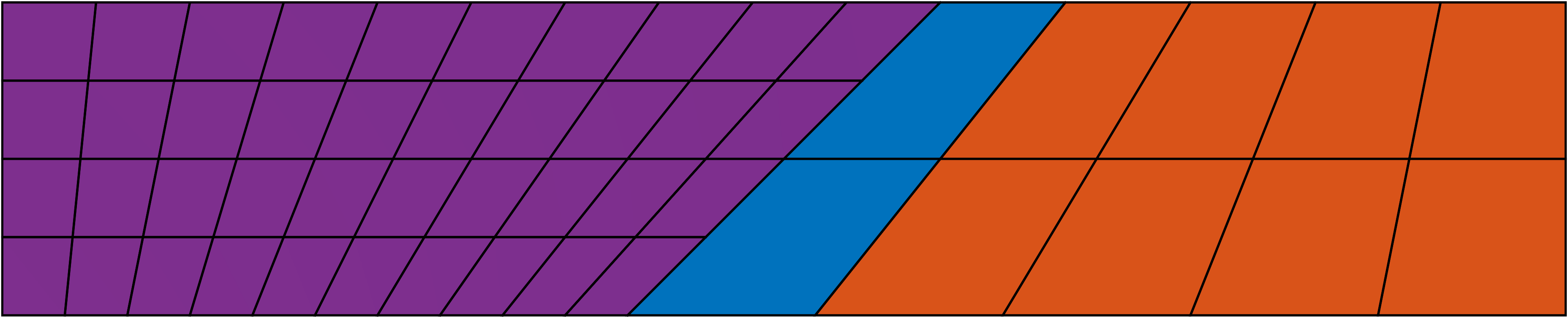}}\\
	\subfloat[$n_\mathrm{s}\,{=}\,1$, $n_\mathrm{y}\,{=}\,3$]{\includegraphics[clip,width=0.475\textwidth]{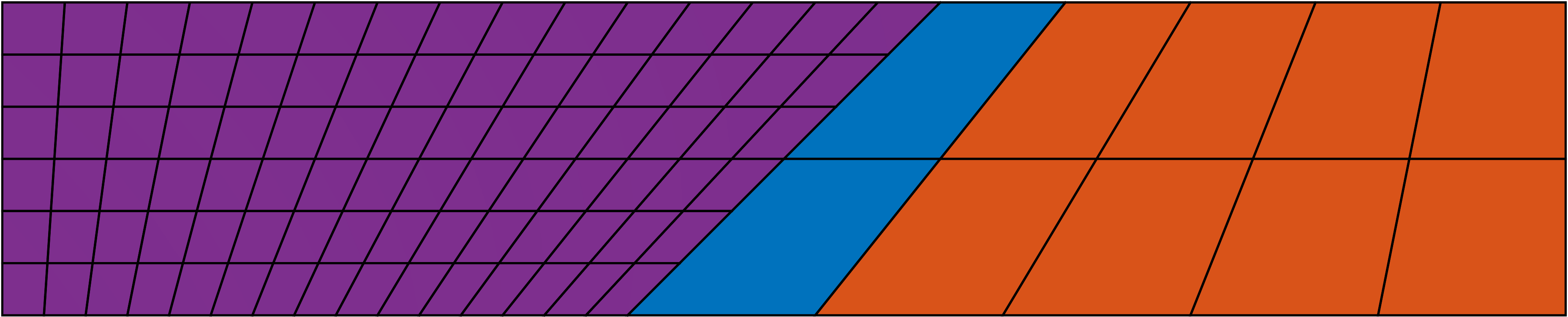}}
	\hfill
	\subfloat[$n_\mathrm{s}\,{=}\,1$, $n_\mathrm{y}\,{=}\,4$]{\includegraphics[clip,width=0.475\textwidth]{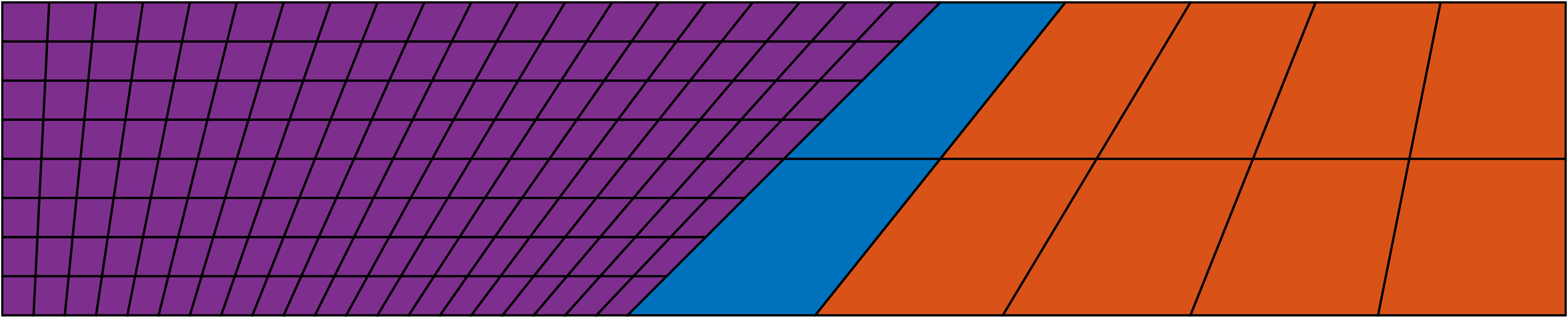}}\\
	\subfloat[$n_\mathrm{s}\,{=}\,2$, $n_\mathrm{y}\,{=}\,2$]{\includegraphics[clip,width=0.475\textwidth]{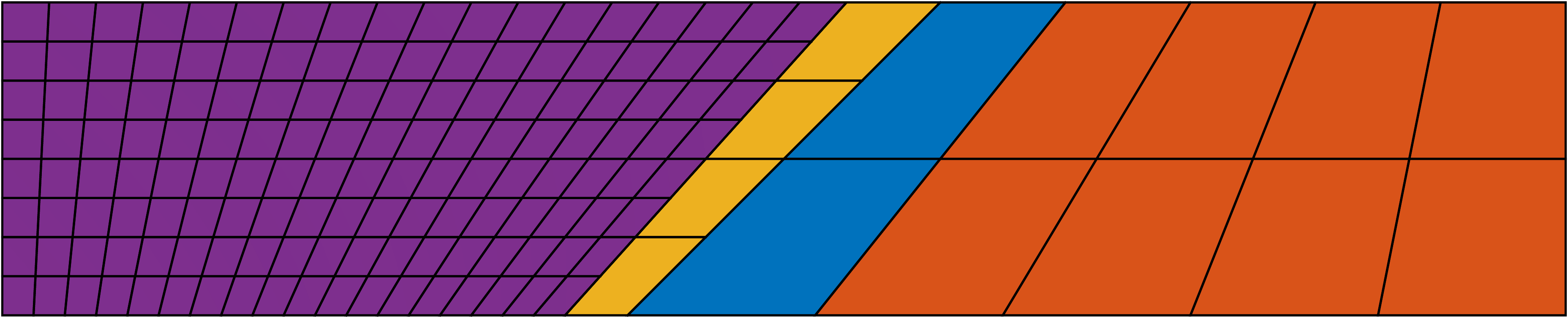}}
	\hfill
	\subfloat[$n_\mathrm{s}\,{=}\,2$, $n_\mathrm{y}\,{=}\,3$]{\includegraphics[clip,width=0.475\textwidth]{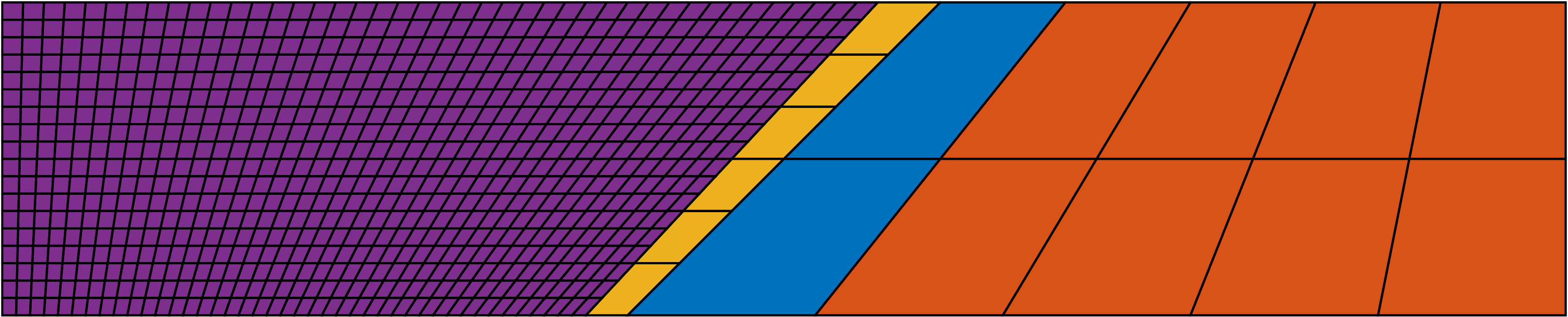}}\\
	\subfloat[$n_\mathrm{s}\,{=}\,2$, $n_\mathrm{y}\,{=}\,4$]{\includegraphics[clip,width=0.475\textwidth]{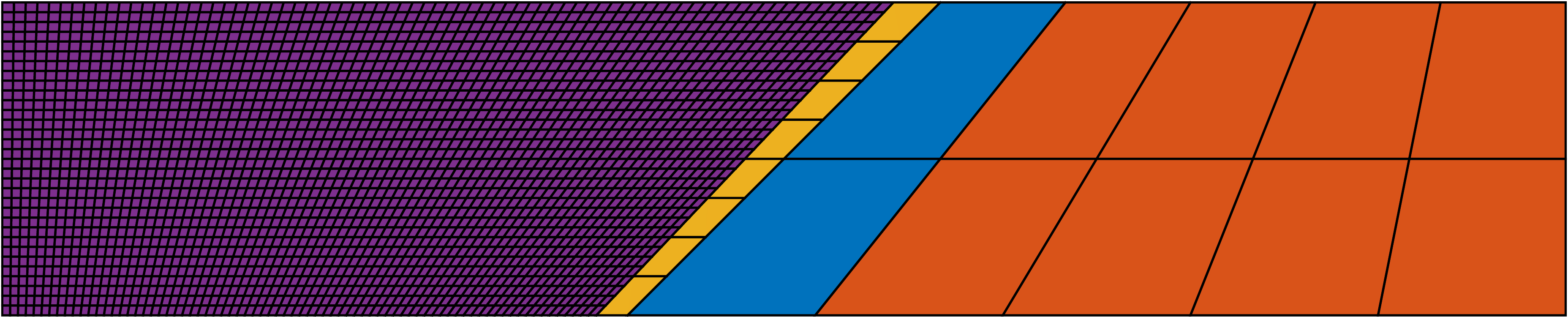}}
	\hfill
	\subfloat[$n_\mathrm{s}\,{=}\,3$, $n_\mathrm{y}\,{=}\,2$]{\includegraphics[clip,width=0.475\textwidth]{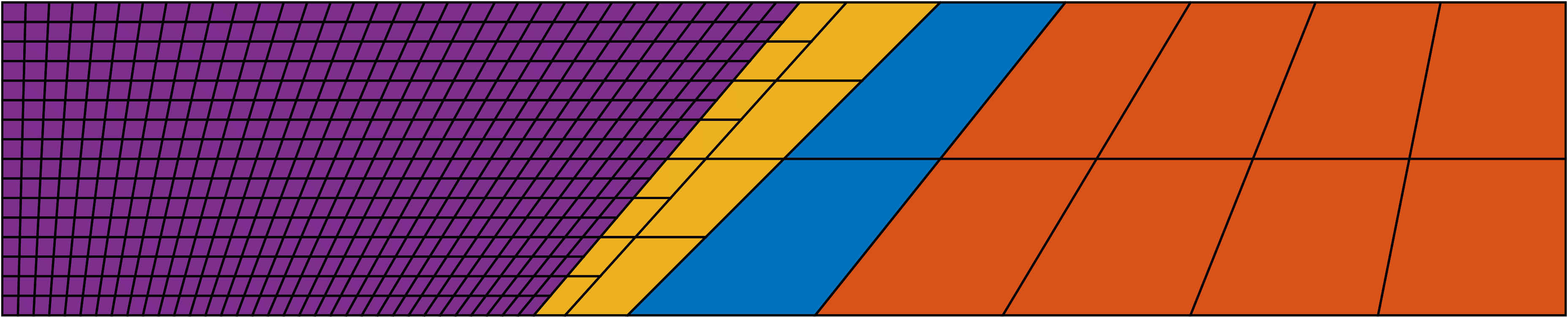}}\\
	\caption{Transition element meshes with different values of $n_\mathrm{s}$ and $n_\mathrm{y}$. The different element types are color-coded: \emph{x}-Elements \textcolor{xElem}{\rule[-0.25ex]{5ex}{2ex}}, \emph{y}-Elements \textcolor{yElem}{\rule[-0.25ex]{5ex}{2ex}}, \emph{x\textbf{N}y}-Elements \textcolor{xNyElem}{\rule[-0.25ex]{5ex}{2ex}}, \emph{y\textbf{N}y}-Elements \textcolor{yNyElem}{\rule[-0.25ex]{5ex}{2ex}}.}
	\label{fig:nsny}
\end{figure}

Versions B and C of the patch test are used in Sects.~\ref{sec:PatchTestLinear} to \ref{sec:PatchTestHigh} to investigate the performance of the proposed \emph{x\textbf{N}y}-element concept. The overarching goal is to show that the theoretically optimal rates of convergence can be retained independent of the set-up of these elements. However, before we start to investigate the properties of our approach, we have to introduce several parameters that will be varied in the parametric studies:
\begin{enumerate}
\item Shape functions
\begin{itemize}
    \item Lagrange (SEM)
    \item Legendre (\emph{p}-FEM)
\end{itemize}
\item $n_\mathrm{nd}$: Nodal distribution (Lagrange shape functions)
\begin{itemize}
    \item GLL
    \item GLC
\end{itemize}
\item $n_\mathrm{y}$: Number of \emph{y}-elements coupled to an edge of an \emph{x}-element
\item $n_\mathrm{s}$: Number of successive (hierarchic) mesh refinement steps
\item Polynomial degree of the shape functions\footnote{For the sake of simplicity, it is assumed that the polynomial degree of the shape functions is identical for all local coordinates.}
\begin{itemize}
    \item $p_\mathrm{x}$
    \item $p_\mathrm{y}$
\end{itemize}
\end{enumerate}
A few words regarding the parameters $n_\mathrm{y}$ and $n_\mathrm{s}$ are in order: $n_\mathrm{y}$ denotes the number of \emph{y}-elements that are coupled to an edge of an \emph{x}-element. In the current implementation, this means that a \emph{y}-element that is adjacent to an \emph{x}-element in the original mesh will be divided into $n_\mathrm{y} \times n_\mathrm{y}$ elements. This procedure constitutes one mesh refinement step. The number of such steps is denoted by $n_\mathrm{s}$. The process is referred to as hierarchical as the mesh refinement steps are executed in a sequential fashion. Once the first refinement step has been completed, all elements that are at the interface between the two different domains are subdivided in an identical way as described above. For a better understanding, the parameters $n_\mathrm{y}$ and $n_\mathrm{s}$ are illustrated in Fig.~\ref{fig:nsny}. Using different combinations of those parameters, we can either enforce a fast transition from a coarse to  a fine mesh or have a more gradual decrease in the element size. Although an isotropic subdivision of the \emph{y}-elements is shown, a graded mesh can also be implemented in a straightforward fashion.

There are different possibilities of how such a mesh refinement can be implemented. One idea is to start from a coarse discretization, in the following referred to as initial or base mesh, where the elements are divided into two sets, one corresponding to the \emph{x}-elements and the other one to the \emph{y}-elements. To characterize the created mesh, the number of elements in the base mesh along with the parameters $n_\mathrm{y}$ and $n_\mathrm{s}$ are specified. The \emph{y}-elements in the base mesh are subdivided according to the defined values of $n_\mathrm{y}$ and $n_\mathrm{s}$ (see~Fig.~\ref{fig:nsny}). A second idea would be to import two meshes with different element sizes and to merge both meshes by identifying the interface between \emph{x}- and \emph{y}-elements. Keep in mind that in both approaches to generate meshes with transition elements, the \emph{x}-elements at the common interface with the \emph{y}-elements must be converted into \emph{x\textbf{N}y}-elements in order to ensure a conformal/compatible coupling. If several successive subdivisions of the \emph{y}-elements are conducted, i.e., $n_\mathrm{s}\,{\ge}\,2$, we also have to introduce \emph{y\textbf{N}y}-elements because we need to couple the same element type.

The examples presented in the following sections are computed with all possible combinations of different shape function types for \emph{x}- and \emph{y}-elements (Lagrange-Lagrange, Lagrange-Legendre, Legendre-Lagrange, Legendre-Legendre), nodal distributions for the Lagrange shape functions (GLL, GLC), discretizations ($n_\mathrm{y} \in \{1,2,3,4\}$, $n_\mathrm{s} \in \{1,2,3,4\}$), and polynomial degrees ($p_\mathrm{x} \in \{2,3,\ldots,8\}$, $p_\mathrm{y} \in \{2,3,\ldots,8\}$). Thus, 6,272 different numerical models need to be analyzed for each benchmark problem. This procedure ensures that the implementation is correct and at the same time, we can investigate the performance of the transition elements in all its details.
\subsection{Linear patch test}
\label{sec:PatchTestLinear}
\begin{figure}[b!]
	\centering
	\includegraphics[clip,scale=0.72]{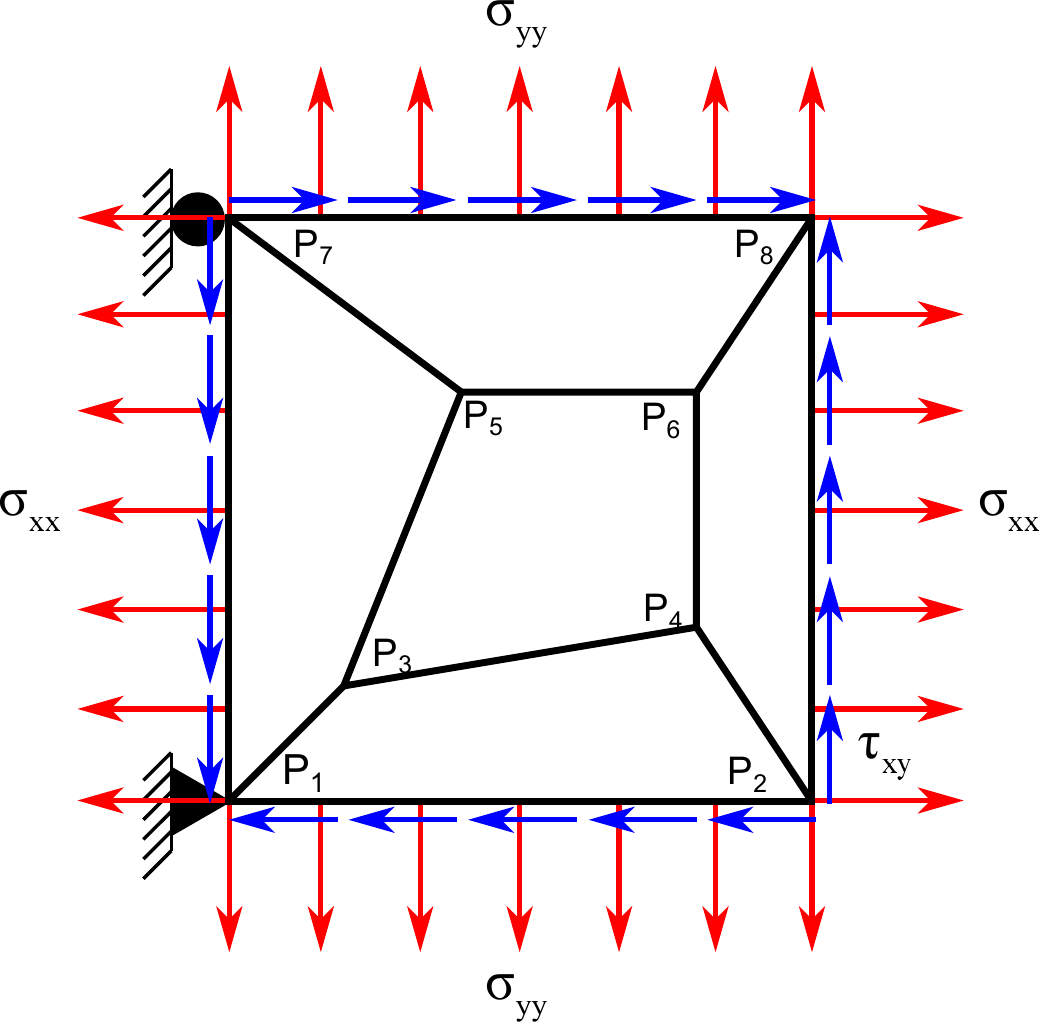}
	\caption{Geometry and boundary conditions for the linear patch test according to Ref.~\cite{BookBathe2002}.}
	\label{fig:PatchTestLinear}
\end{figure}%
In the linear patch test, it is checked whether the finite element is capable of exactly representing a constant strain state (linear displacement field), see Ref.~\cite{BookBathe2002} (Ch.~4.4, p.~264). This example belongs to the class of version C patch tests as discussed in the introduction to this section. Obviously, the number of possible constant strain states depends on the dimensionality of the model. In Fig.~\ref{fig:PatchTestLinear}, a typical patch of finite elements that is used in numerical investigations is depicted. The coordinates of the vertices for this model are listed in Table~\ref{tab:PatchTestLinear}. For this structure, the minimal number of essential (Dirichlet) boundary conditions has been enforced, and loads are applied such that three constant stress states arise. Note that each finite element has to pass the linear patch test in order to ensure convergence, i.e., the numerical solution has to be accurate up to machine precision.
\begin{table}[t!]
	\centering
	\caption{Coordinates of the vertices for the (linear) patch test according to Ref.~\cite{BookBathe2002}. \label{tab:PatchTestLinear}}
	\begin{tabular}{c|c|c||c|c|c}
		\toprule
		$\#$ Points & $x_i$ & $y_i$ & $\#$ Points & $x_i$ & $y_i$\\\hline
		1&	0	&	0	&	5	&	4	&	7 \\
		2&	10	&	0	&	6	&	8	&	7 \\
		3&	2	&	2	&	7	&	0	&	10\\
		4&	8	&	3	&	8	&	10	&	10\\
		\bottomrule
	\end{tabular}
\end{table}
\begin{figure}[t!]
	\centering
	\subfloat[Lagrange-Lagrange]{\includegraphics[clip,width=0.475\textwidth]{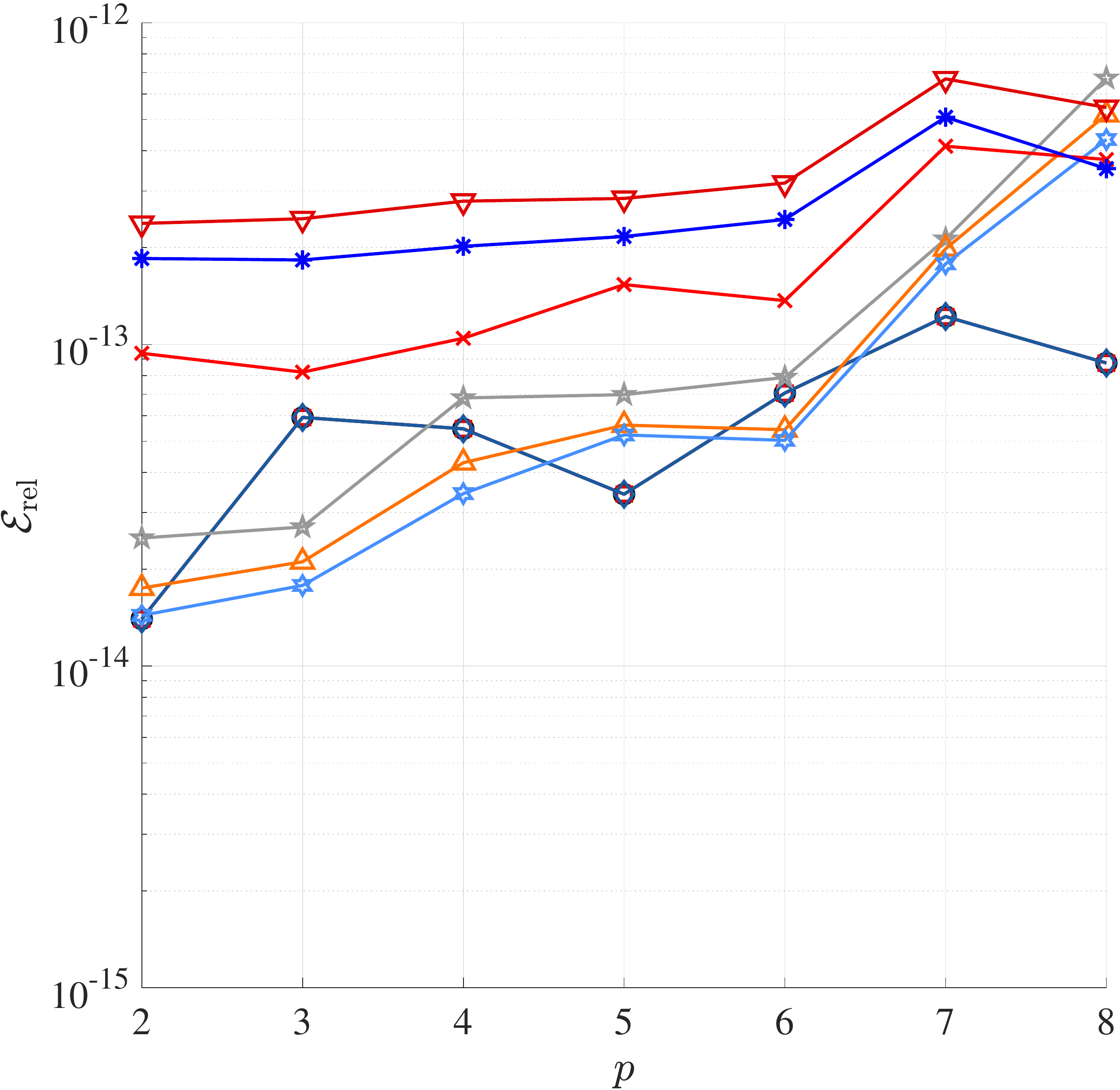}}
	\hfill
	\subfloat[Lagrange-Legendre]{\includegraphics[clip,width=0.475\textwidth]{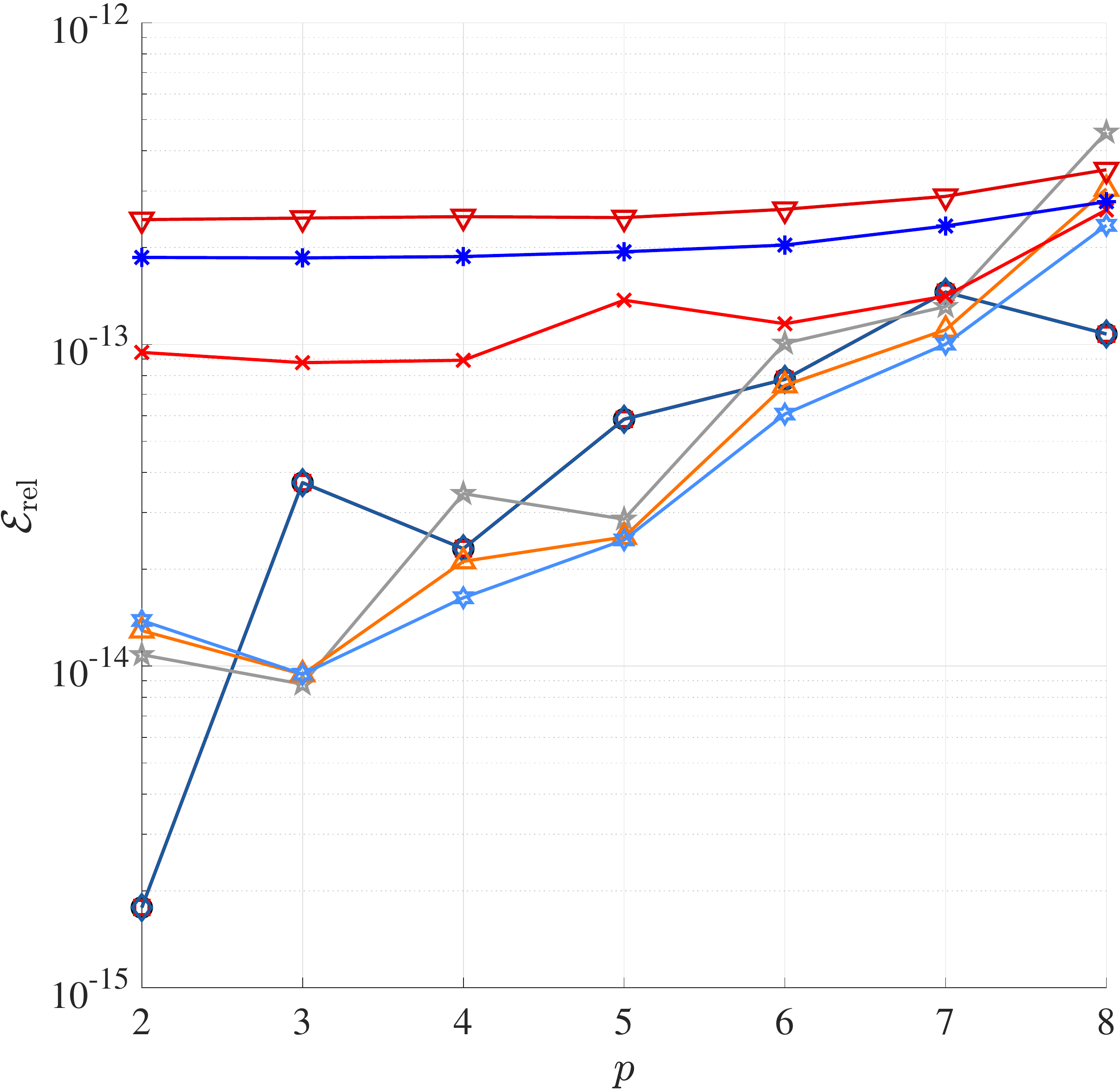}}\\
	\subfloat[Legendre-Legendre]{\includegraphics[clip,width=0.475\textwidth]{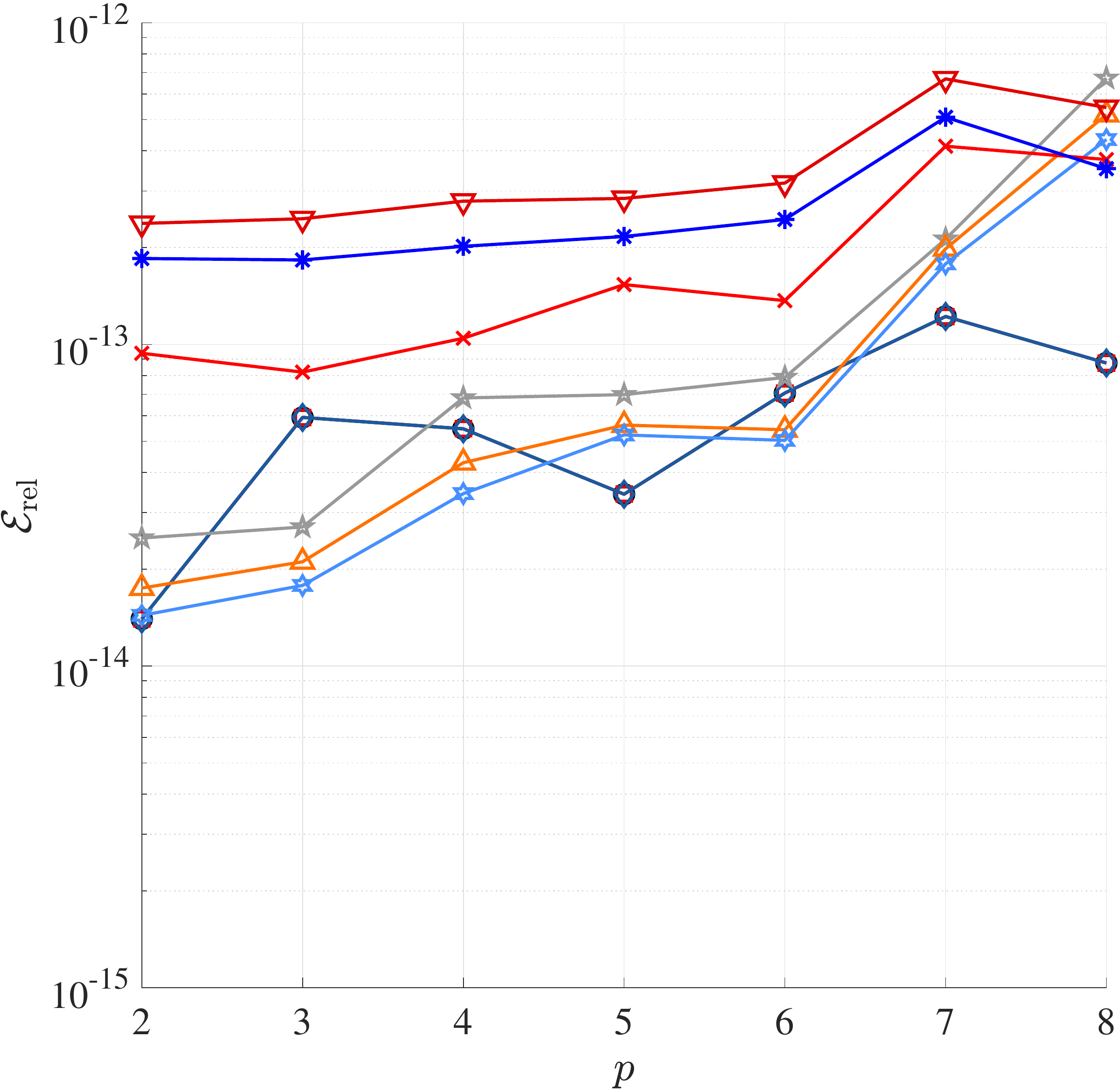}}
	\hfill
	\subfloat[Legendre-Lagrange]{\includegraphics[clip,width=0.475\textwidth]{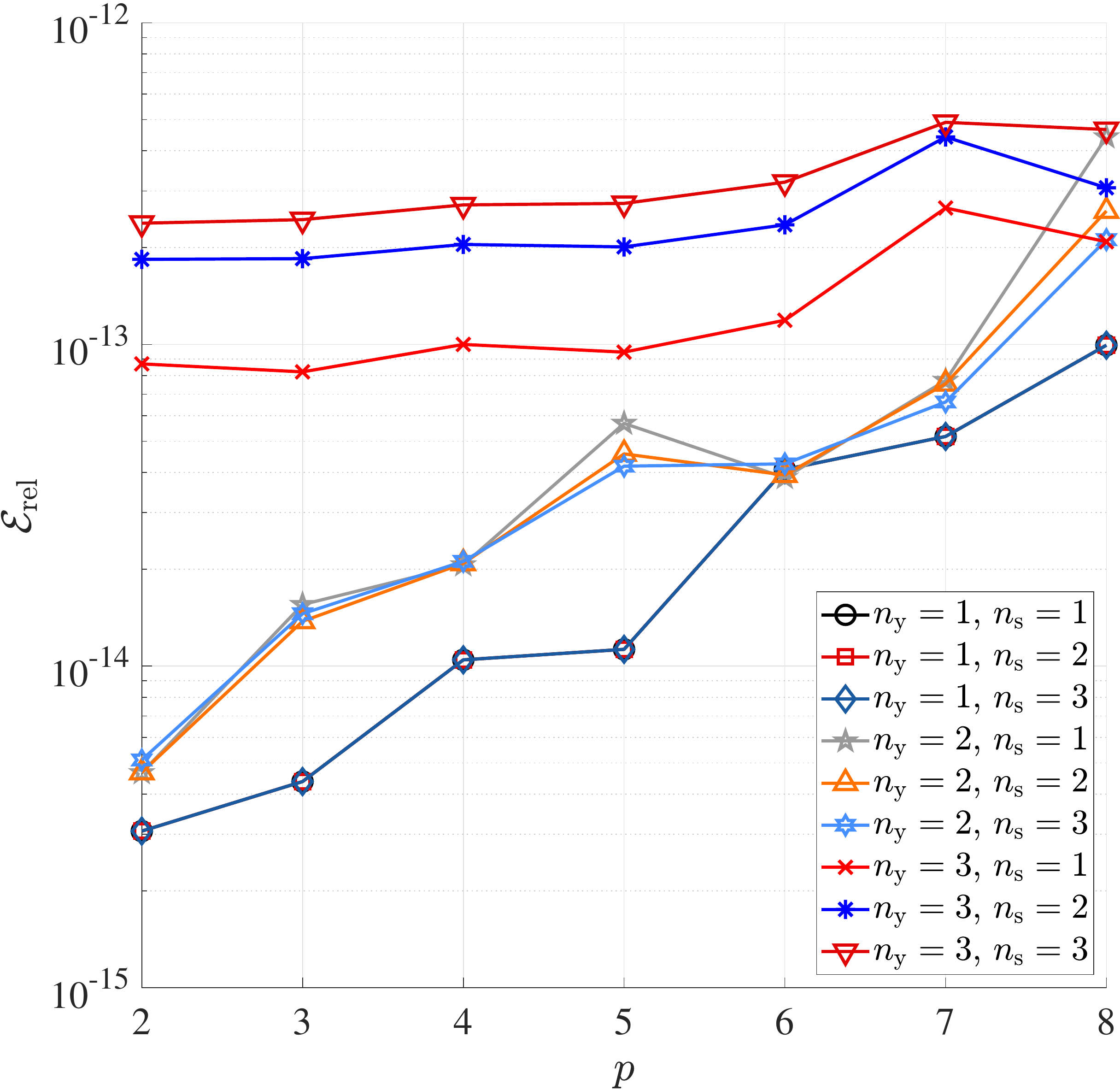}}\\
	\caption{Relative error for the linear patch test example ($p_\mathrm{x}\,{=}\,p_\mathrm{y}$, GLL nodal distribution).}
	\label{fig:ErrorLinPatchTest}
\end{figure}%

The base mesh for this example consists of five finite elements whereas the four outer elements are classified as \emph{x}-elements, while the fifth (interior) element is a \emph{y}-element. Hence, depending on the current values of  $n_\mathrm{s}$ and $n_\mathrm{y}$, the \emph{x}-elements are exchanged by suitable transition elements, and in the interior, a certain number of \emph{y}- and/or \emph{y\textbf{N}y}-elements is generated (see Fig.~\ref{fig:nsny} for the general procedure).

To assess the performance of the \emph{x\textbf{N}y}-elements, the three stress components are computed at specific locations in each element. Although it is generally accepted that the stress values are more accurate at the Gau\ss{} integration points of order $p\,{-}\,1$, -- with $p$ denoting the polynomial order of the corresponding shape functions of the finite elements -- we compute them at the nodes (Lagrange elements) or at points in an equidistant grid (hierarchic element). The reason being that each element must be capable of exactly reproducing a constant stress state independent on the location where the stresses are computed; otherwise, convergence is not ensured. As an error measure, we use the mean relative error of the computed quantities
\begin{equation}
\mathcal{E}_\mathrm{rel} = \cfrac{1}{n_\mathrm{P}} \sum\limits_{i}^{n_\mathrm{P}} \cfrac{|\square_\mathrm{ref} - \square_\mathrm{num}|}{|\square_\mathrm{ref}|}\,,
\end{equation}
where $\square_\mathrm{ref}$ and $\square_\mathrm{num}$ denote the reference (analytical) and numerical solution, respectively, and $n_\mathrm{P}$ is the number of evaluation points (nodes). For each of the computed models, the numerical results show relative errors of less than $10^{-12}$ in the stress components, see also Fig.~\ref{fig:ErrorLinPatchTest}. The largest error values are generally observed for the largest values of $n_\mathrm{y}$ and $n_\mathrm{s}$ and in cases where the polynomial degree of the shape functions $p_\mathrm{x}$ and $p_\mathrm{y}$ is comparably high (i.e., $p\,{\ge}\,7$). In Fig.~\ref{fig:ErrorLinPatchTest}, we exemplarily plot the error values for coupling different element types but with the same polynomial degree. This trend is rather universal and independent of whether $p_\mathrm{x}\,{=}\,p_\mathrm{y}$ or $p_\mathrm{x}\,{\ne}\,p_\mathrm{y}$. \textcolor{red}{The system of equations has been solved using Matlab's direct solver based on a Cholesky-factorization of the system matrix. The increasing error trend is still subject to investigations and it is so far not entirely clear what is causing this effect. Preliminary studies have revealed that the condition number of the stiffness matrix is growing slightly during the refinement process (about 3 orders of magnitude) and this might be one reason for the observed behavior.} Overall, we can still state that the tested transition elements pass the linear patch test and therefore, the convergence of the results is generally ensured, although it is not yet clear what convergence rates can be attained. In this example, the material properties of aluminum (Young's modulus $E\,{=}\,70\,$GPa, Poisson's ratio $\nu\,{=}\,0.3$) have been used. As we explicitly prescribe all stress components, the solution does not, however, depend on these properties. Therefore, the material data is only of importance if the displacement field needs to be computed and compared.
\subsection{Quadratic patch test}
\label{sec:PatchTestQuadratic}
\begin{figure}[t!]
    \centering
    \includegraphics[clip,scale=0.75]{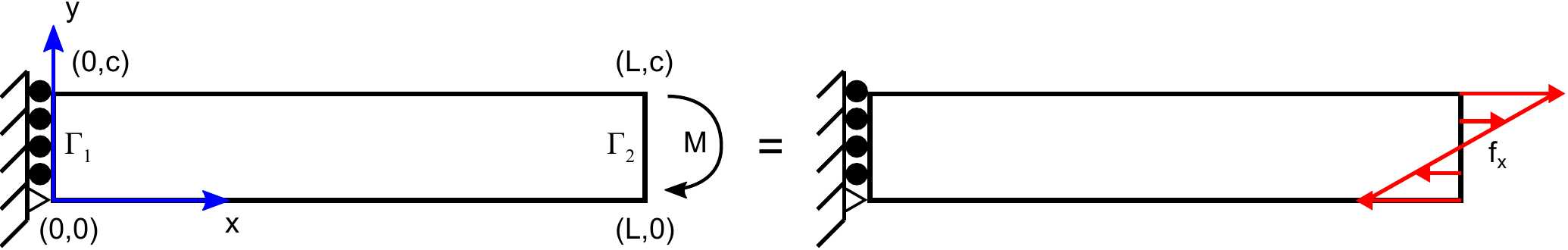}
	\caption{Quadratic patch test: Constant bending moment problem; Geometry ($L\,{=}\,10$, $c\,{=}\,2$) and boundary conditions \cite{ArticleLee1993}.}
	\label{fig:PatchTestQuadratic}
\end{figure}

A quadratic patch test has been proposed by Zienkiewicz \& Taylor \cite{BookZienkiewicz2000a} and Lee \& Bathe \cite{ArticleLee1993}. Here, a `constant bending-moment problem' is suggested which results in a quadratic displacement field. The model including Dirichlet and Neumann boundary conditions is depicted in Fig.~\ref{fig:PatchTestQuadratic}. The left edge ($\Gamma_1$) is constrained in $x$-direction, while the lower left vertex is fixed additionally in $y$-direction. These boundary conditions are sufficient to avoid rigid body motions of the structure. The load (Neumann boundary conditions) has been chosen as a distributed surface traction at the right edge ($\Gamma_2$), and the amplitude follows the analytical solution ($0 \le y \le c$, $0 \le x \le L$):
\begin{align}
\sigma_\mathrm{xx} &= \cfrac{240}{c}\,y - 120\,,\\
\sigma_\mathrm{yy} &= 0\,,\\
\tau_\mathrm{xy} &= 0\,.
\end{align}
Note that the stress distribution is independent of both the material's Young's modulus $E$ and Poisson's ratio $\nu$. These properties only influence the resulting displacement field. For this example, we can easily provide the exact solution for the displacement field, which can be derived in closed-from:
%
% Solution for nu = 0.3
%\begin{align}
%u_\mathrm{x} &= \cfrac{1}{E} \left( \cfrac{240}{c}\,xy - 120x \right)\,,\\
%u_\mathrm{y} &= \cfrac{1}{E} \left( -\cfrac{120}{c}\,x^2 - \cfrac{36}{c}\,y^2 + 36y \right)\,.
%\end{align}
%
\begin{align}
u_\mathrm{x} &= \cfrac{1}{E} \left( \cfrac{240}{c}\,xy - 120x \right)\,,\\
u_\mathrm{y} &= \cfrac{1}{E} \left( -\cfrac{120}{c}\,x^2 - \cfrac{120\nu}{c}\,y^2 + 120\nu y\, \right)\,.
\end{align}
\begin{figure}[b!]
    \centering
    \subfloat[$u_\mathrm{x}$]{\includegraphics[clip,width=0.315\textwidth]{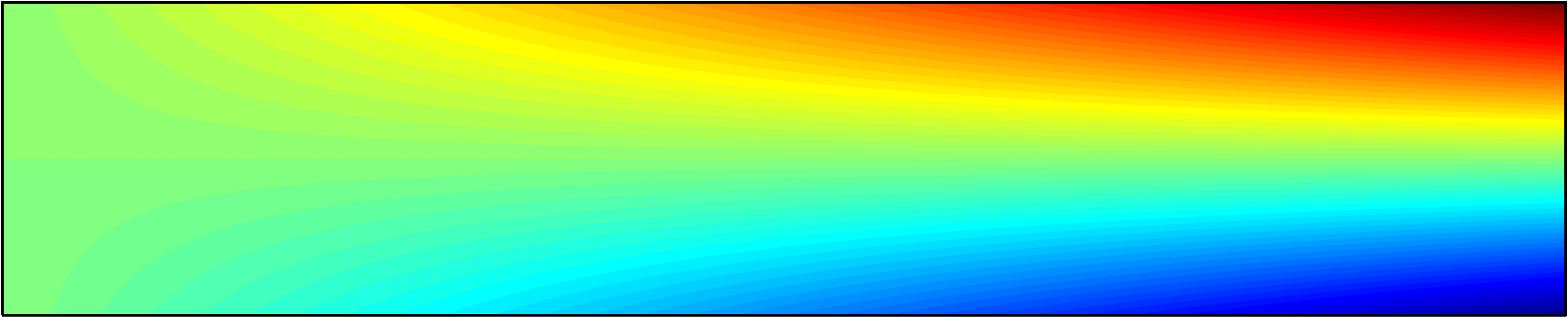}}
    \hfill
    \subfloat[$u_\mathrm{y}$]{\includegraphics[clip,width=0.315\textwidth]{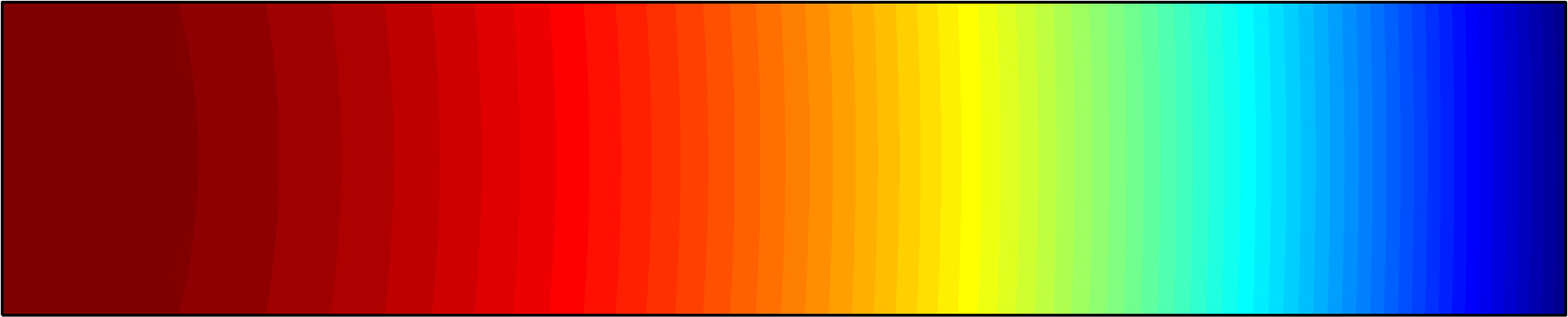}}
    \hfill
    \subfloat[$u_\mathrm{mag}$]{\includegraphics[clip,width=0.315\textwidth]{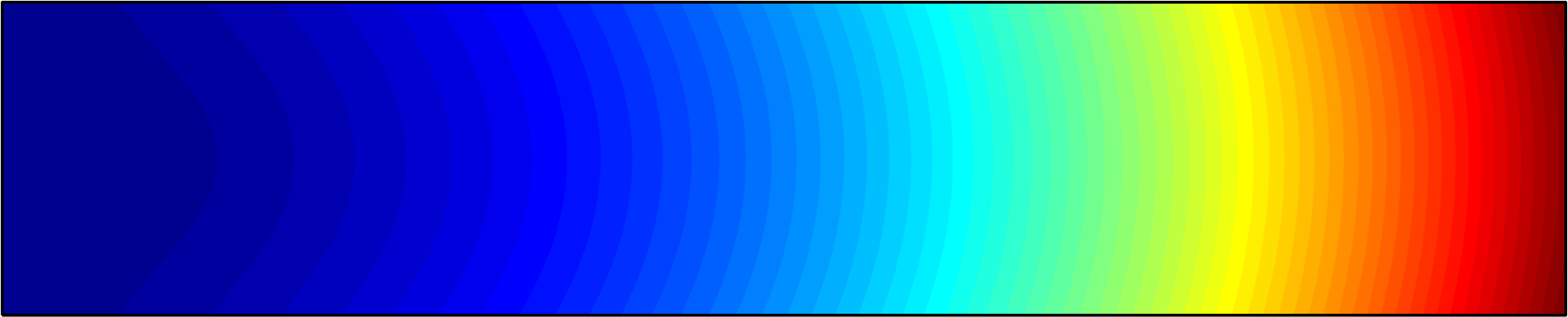}}
    \\
    \subfloat[$\sigma_\mathrm{xx}$]{\includegraphics[clip,width=0.315\textwidth]{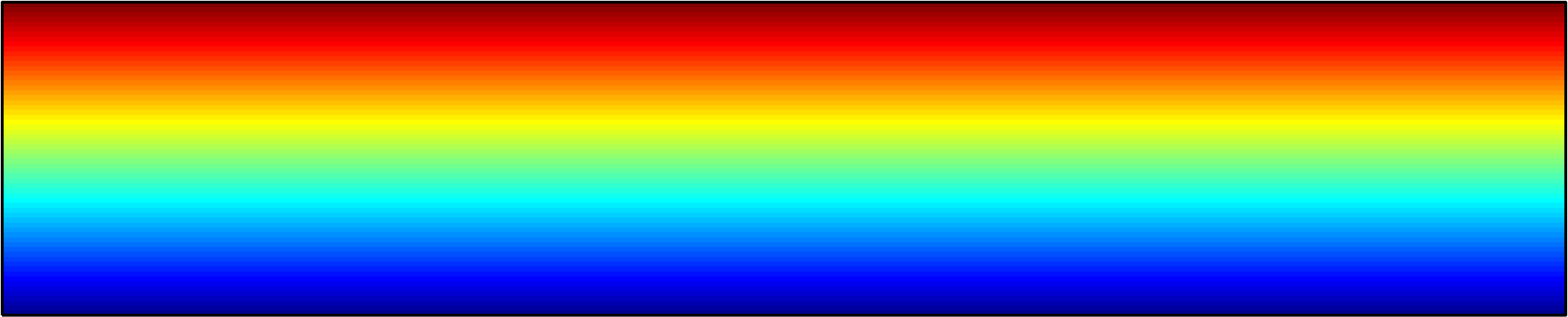}}
    \hfill
    \subfloat[$\sigma_\mathrm{yy}$]{\includegraphics[clip,width=0.315\textwidth]{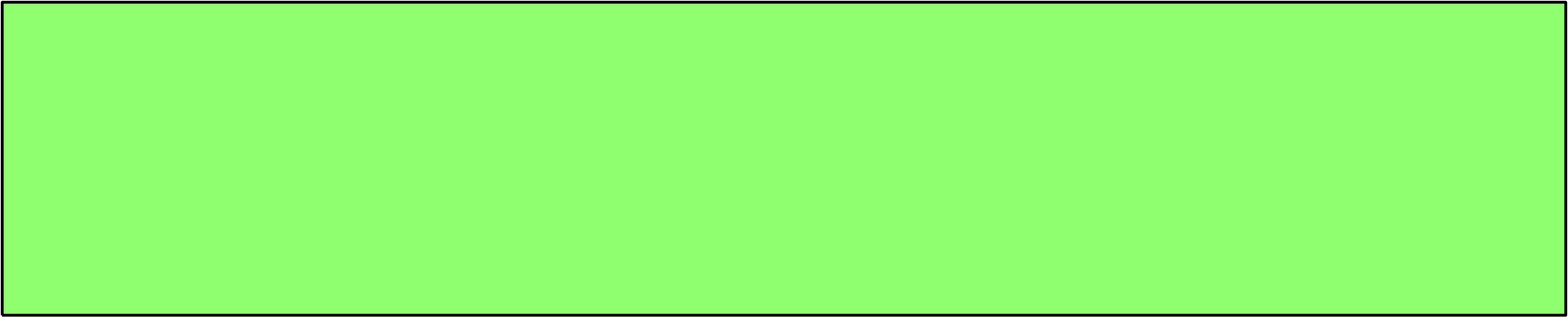}}
    \hfill
    \subfloat[$\tau_\mathrm{xy}$]{\includegraphics[clip,width=0.315\textwidth]{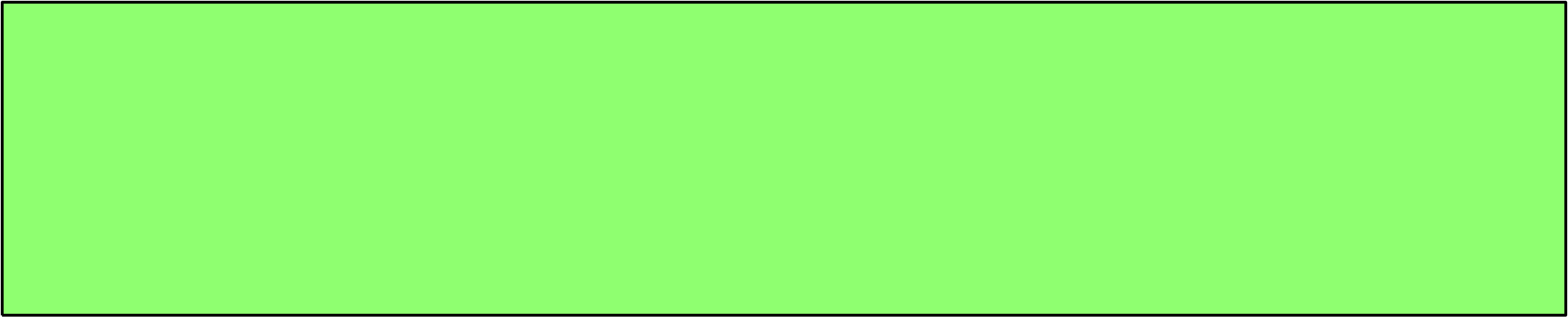}}
\caption{Displacement and stress solution for the quadratic patch test \cite{ArticleLee1993}.}
\label{fig:QuadPatchTestUS}
\end{figure}%
\begin{figure}[t!]
    \centering
    \subfloat[Lagrange-Lagrange]{\includegraphics[clip,width=0.475\textwidth]{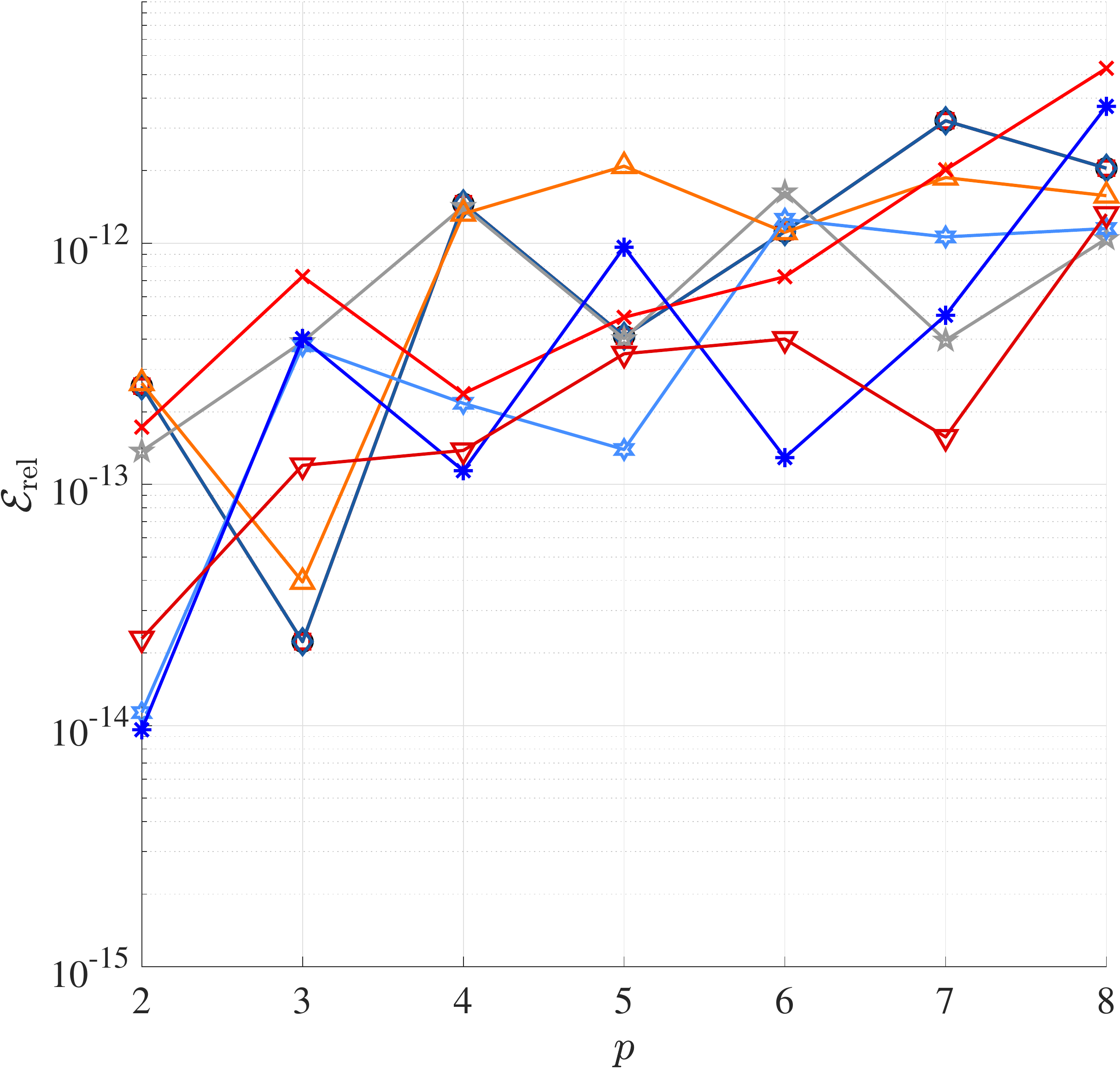}}
    \hfill
    \subfloat[Lagrange-Legendre]{\includegraphics[clip,width=0.475\textwidth]{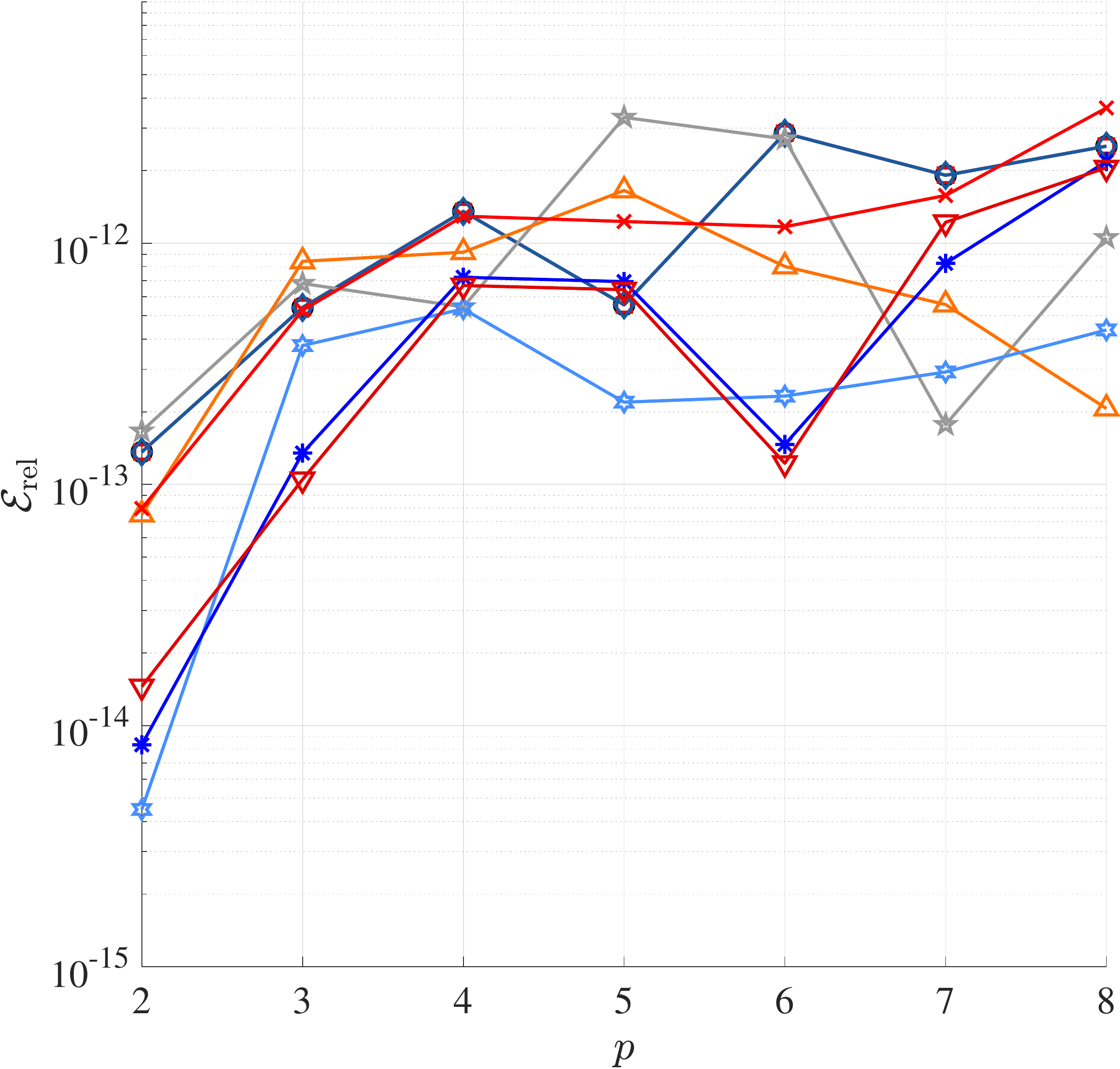}}\\
    \subfloat[Legendre-Legendre]{\includegraphics[clip,width=0.475\textwidth]{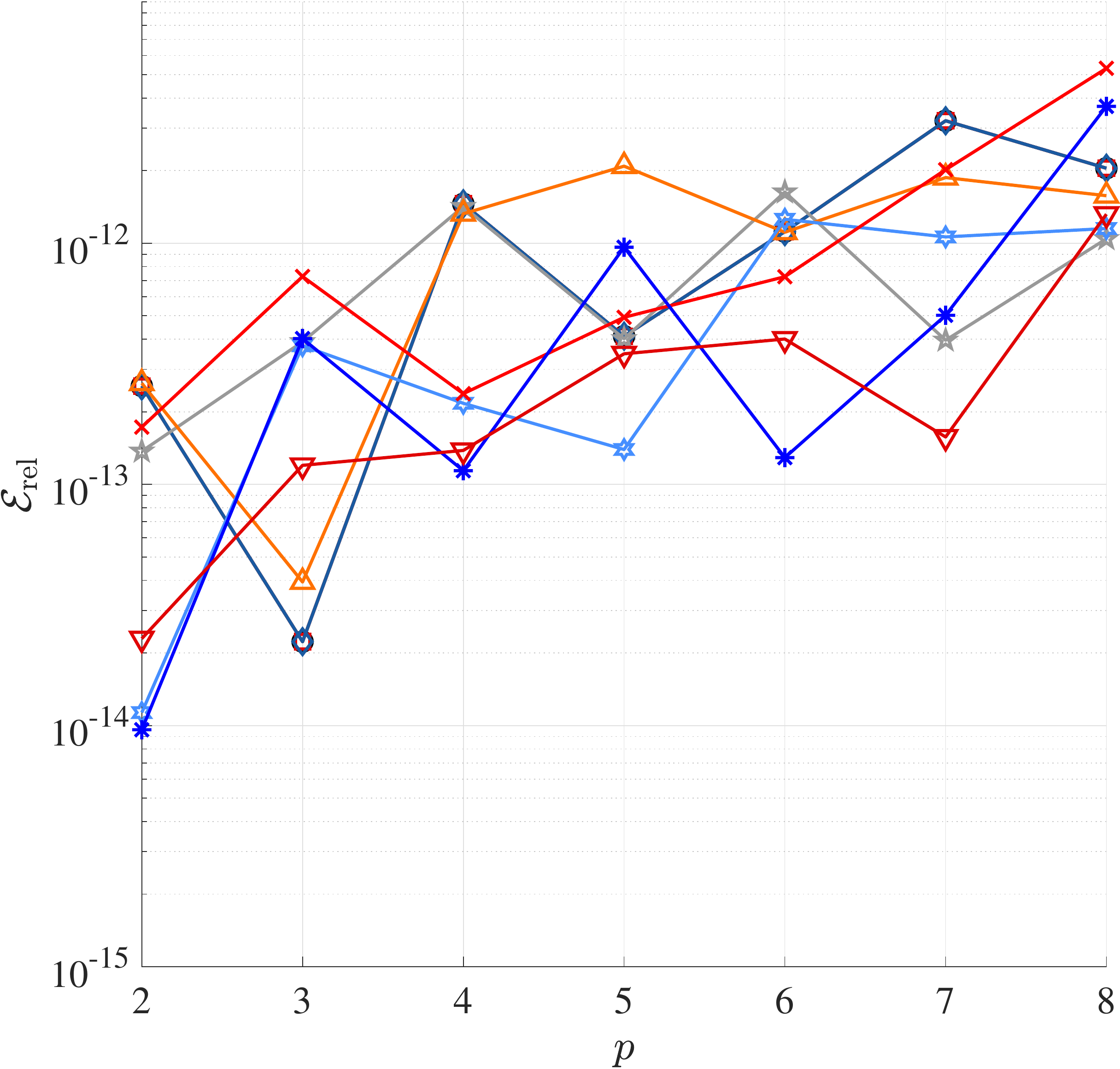}}
    \hfill
    \subfloat[Legendre-Lagrange]{\includegraphics[clip,width=0.475\textwidth]{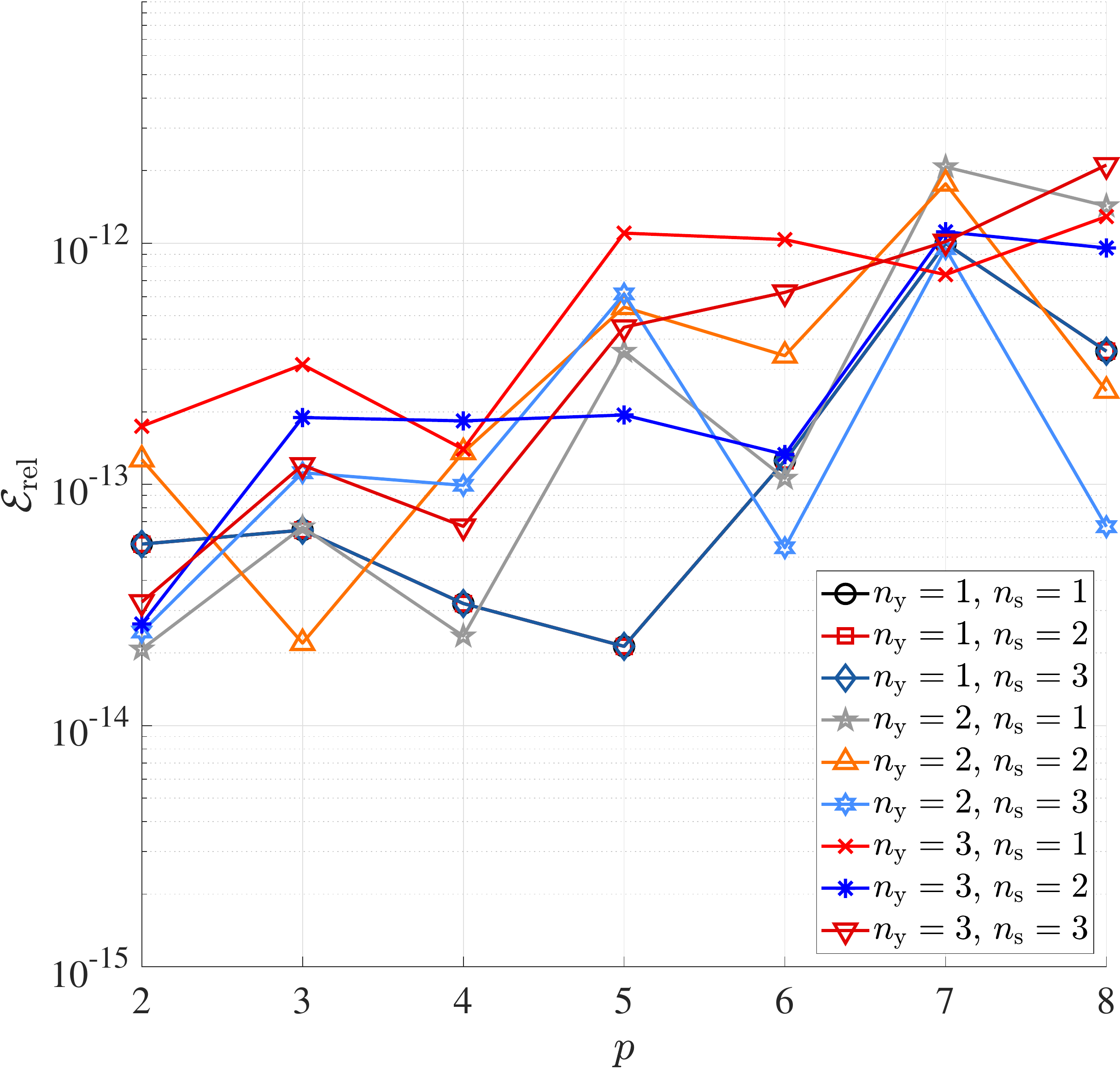}}
	\caption{Relative error for the quadratic patch test example ($p_\mathrm{x}\,{=}\,p_\mathrm{y}$, GLL nodal distribution).}
	\label{fig:ErrorQuadPatchTest}
\end{figure}%
The theoretical displacement field and stress distributions are depicted in Fig.~\ref{fig:QuadPatchTestUS} for a Poisson's ratio of $\nu\,{=}\,0.3$. To achieve these results in the numerical analysis, we apply the exact tractions at the boundary $\Gamma_2$ as external loads. Therefore, a line load acting in the global $x$-direction is applied to introduce a constant bending moment. According to the analytical solution, the Neumann boundary condition at the right edge takes the following form
\begin{equation*}
f_{\mathrm{x},\Gamma_2} = \cfrac{240}{c}\,y - 120\,.
\end{equation*}
As discussed in Sect.~\ref{sec:TransitionElem}, the shape functions of the transition elements are based on tensor product formulations of \emph{p}- or spectral elements. Therefore, the minimum polynomial degree
\begin{equation}
    p_\mathrm{min} = \min\{[p_\mathrm{x},\, p_\mathrm{y}]\}
\end{equation}
determines the order of the highest complete polynomial in the function space. Consequently, it is expected that the asymptotic rate of convergence depends on $p_\mathrm{min}$ only. Since all models that are investigated in this article feature a minimum polynomial degree of 2, the numerical solution should be exact. The obtained results for a two-element base mesh are depicted in Fig.~\ref{fig:ErrorQuadPatchTest}. The beam is divided by a skew line with the starting point (4,\,0) and the end point (6,\,2). Thus, two elements with angular distortion are created. It has to be mentioned that the used mesh is actually a coarser version of the one depicted in Fig.~\ref{fig:nsny}. By agglomerating all \emph{x}-/\emph{x\textbf{N}y}-elements and all \emph{y}-/\emph{y\textbf{N}y}-elements, we obtain only two elements which constitutes the base mesh for this example. The observed behavior of the transition elements for this example is slightly different from the one discussed in Sect.~\ref{sec:PatchTestLinear}. For all models, relative errors of less than $10^{-11}$ in the stresses have been obtained for the quadratic (constant bending) patch test. We can observe that higher polynomial orders of the shape functions tend to also result in increased error levels. Still, a sufficient accuracy with 11 significant digits can be achieved.
\subsection{Cubic patch test}
\label{sec:PatchTestCubic}
\begin{figure}[b!]
    \centering
    \includegraphics[clip,scale=0.75]{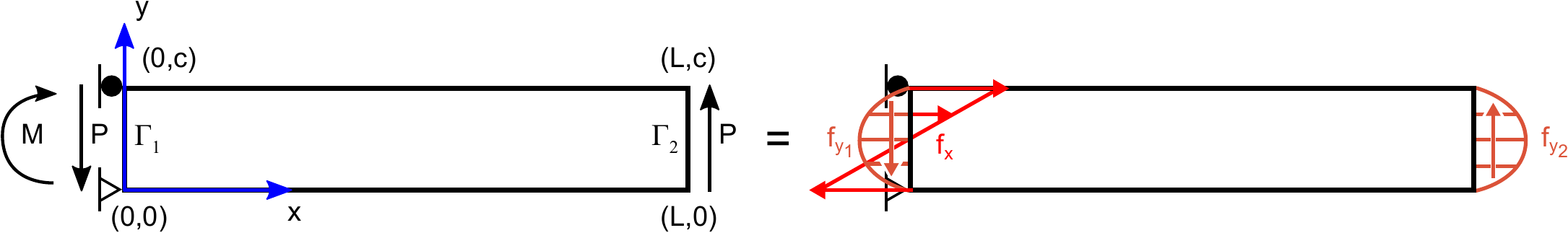}
	\caption{Cubic patch test: Linear bending moment problem; Geometry ($L\,{=}\,10$, $c\,{=}\,2$) and boundary conditions \cite{ArticleLee1993}.}
	\label{fig:PatchTestCubic}
\end{figure}

A cubic patch test has been suggested by Lee \& Bathe in their seminal work on the influence of distortion on the accuracy of finite element simulations \cite{ArticleLee1993}. They propose the `linear bending-moment problem' as illustrated in Fig.~\ref{fig:PatchTestCubic}. The minimum Dirichlet boundary conditions to avoid rigid body motion have been applied at the vertices of the left edge, while at the left ($\Gamma_1$) and right end ($\Gamma_2$) of the beam distributed loads are acting on the structure. The Neumann boundary conditions have been imposed according to the analytical solution of the problem ($0 \le x \le c$, $0 \le y \le L$):
\begin{align}
\sigma_\mathrm{xx} &= \cfrac{240}{cL}\,xy - \cfrac{120}{L}\,x - \cfrac{240}{c}\,y + 120\,,\\
\sigma_\mathrm{yy} &= 0\,,\\
\tau_\mathrm{xy} &= \cfrac{1}{1+\nu}\left(-\cfrac{138+60\nu}{cL}\,y^2 + \cfrac{138+60\nu}{L}\,y \right)\,.
\end{align}
Again, it is possible to derive a closed-form solution of the displacement field for the current example:
%
% Solution for nu = 0.3
%\begin{align}
%u_\mathrm{x} &= \cfrac{1}{E} \left( \cfrac{120}{cL}\,x^2y - \cfrac{92}{cL}\,y^3 - \cfrac{60}{L}\,x^2 - %\cfrac{240}{c}\,xy + \cfrac{138}{L}\,y^2 + 120x - \cfrac{46c}{L}\,y \right)\,,\\
%u_\mathrm{y} &= \cfrac{1}{E} \left( -\cfrac{40}{cL}\,x^3 - \cfrac{36}{cL}\,xy^2 + \cfrac{120}{c}\,x^2 %+ \cfrac{36}{L}\,xy + \cfrac{36}{c}\,y^2 + \cfrac{46c}{L}\,x - 36y \right)\,.
%\end{align}
%
\begin{align}
u_\mathrm{x} &= \cfrac{1}{E} \left( \cfrac{120}{cL}\,x^2y - \cfrac{92}{cL}\,y^3 - \cfrac{60}{L}\,x^2 - \cfrac{240}{c}\,xy + \cfrac{138}{L}\,y^2 + 120x - \cfrac{46c}{L}\,y \right)\,,\\
u_\mathrm{y} &= \cfrac{1}{E} \left( -\cfrac{40}{cL}\,x^3 - \cfrac{120\nu}{cL}\,xy^2 + \cfrac{120}{c}\,x^2 + \cfrac{120\nu}{L}\,xy + \cfrac{120\nu}{c}\,y^2 + \cfrac{46c}{L}\,x - 120\nu y \right)\,.
\end{align}
\begin{figure}[b!]
    \centering
    \subfloat[$u_\mathrm{x}$]{\includegraphics[clip,width=0.315\textwidth]{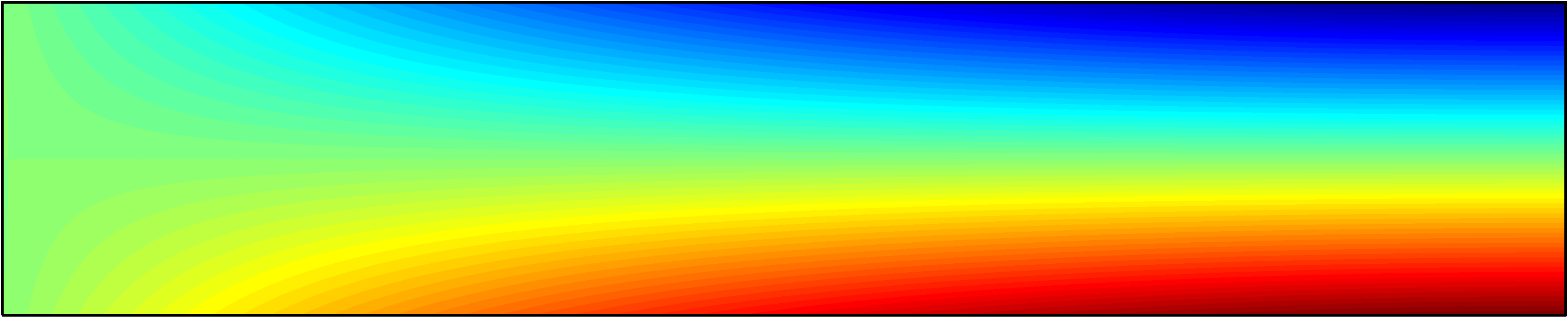}}
    \hfill
    \subfloat[$u_\mathrm{y}$]{\includegraphics[clip,width=0.315\textwidth]{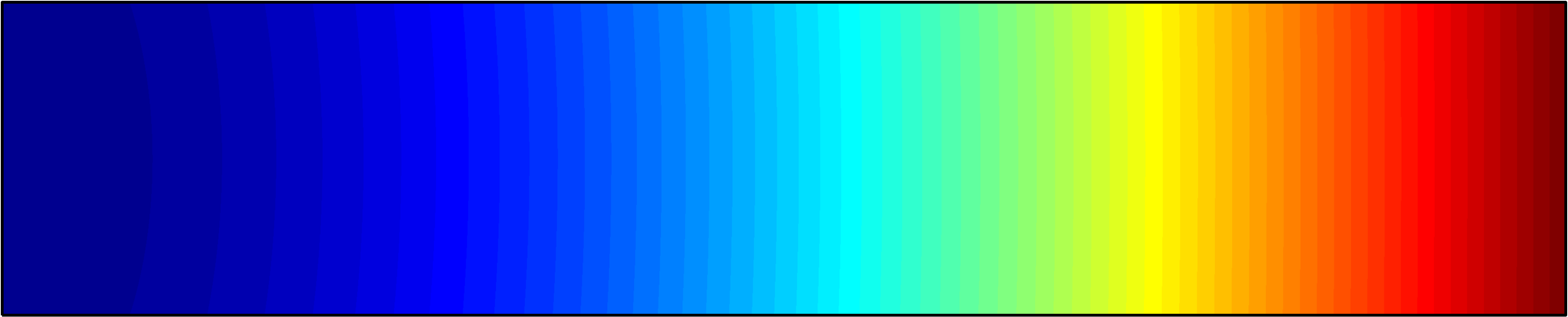}}
    \hfill
    \subfloat[$u_\mathrm{mag}$]{\includegraphics[clip,width=0.315\textwidth]{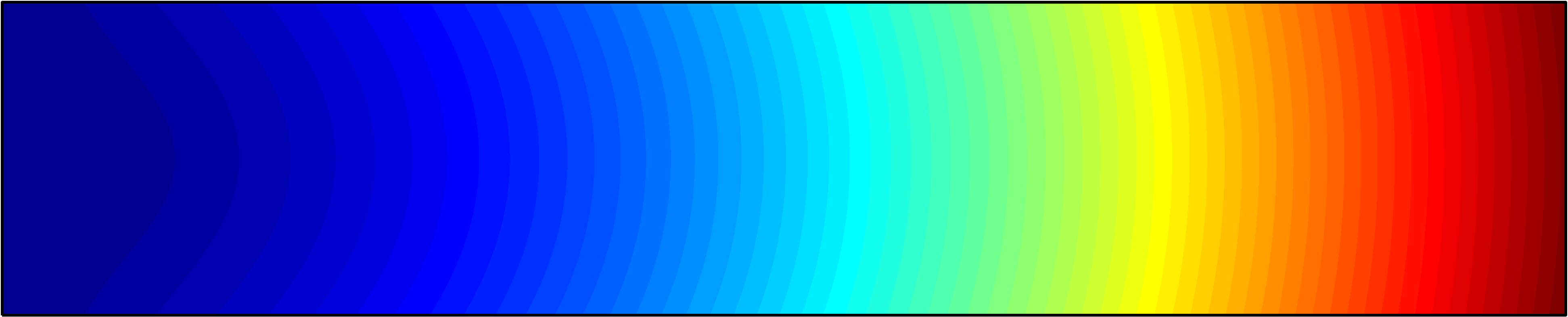}}
    \\
    \subfloat[$\sigma_\mathrm{xx}$]{\includegraphics[clip,width=0.315\textwidth]{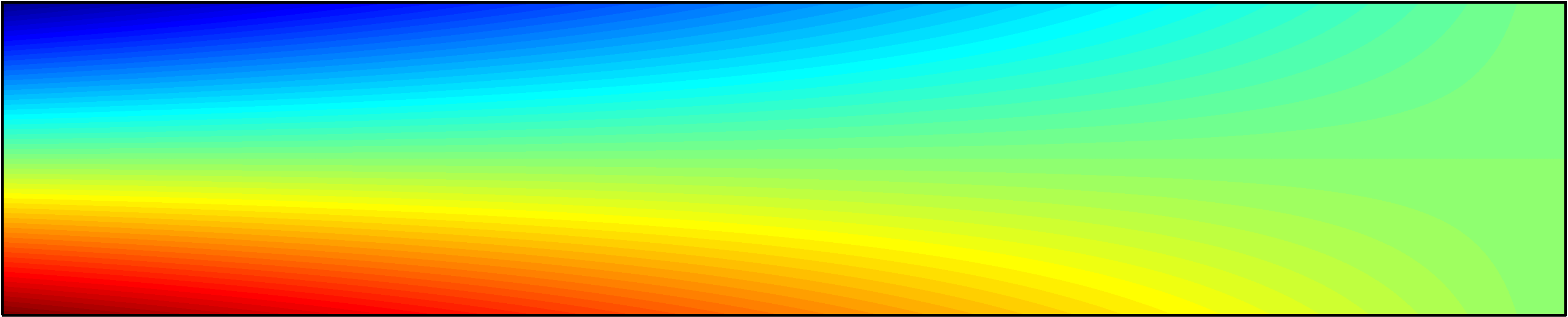}}
    \hfill
    \subfloat[$\sigma_\mathrm{yy}$]{\includegraphics[clip,width=0.315\textwidth]{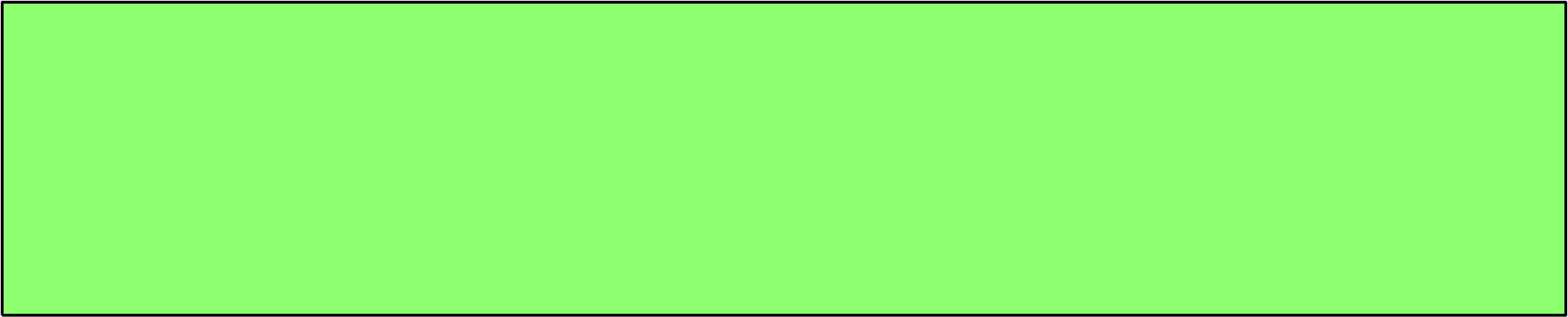}}
    \hfill
    \subfloat[$\tau_\mathrm{xy}$]{\includegraphics[clip,width=0.315\textwidth]{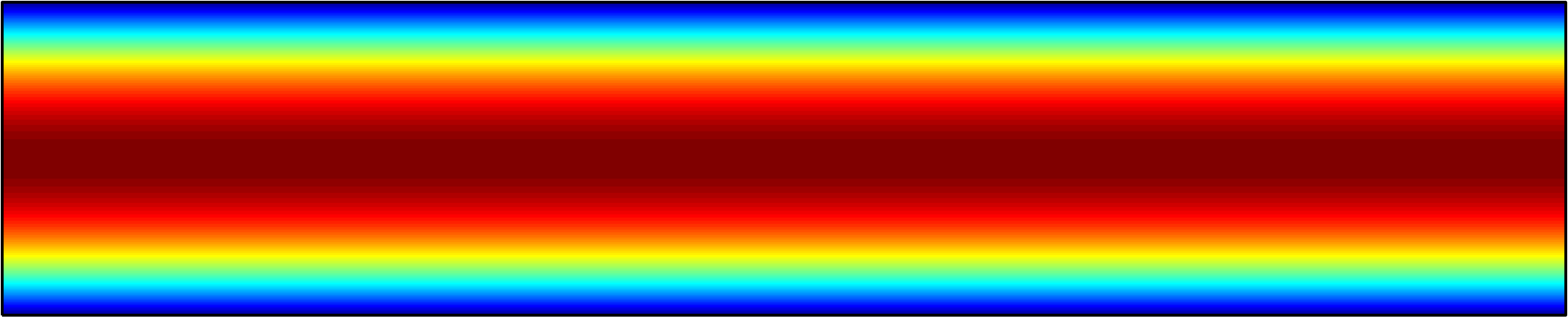}}
\caption{Displacement and stress solution for the cubic patch test \cite{ArticleLee1993}. Solution for $\nu\,{=}\,0.3$.}
\label{fig:CubPatchTestUS}
\end{figure}%
\begin{figure}[t!]
    \centering
    \subfloat[Lagrange-Lagrange]{\includegraphics[clip,width=0.475\textwidth]{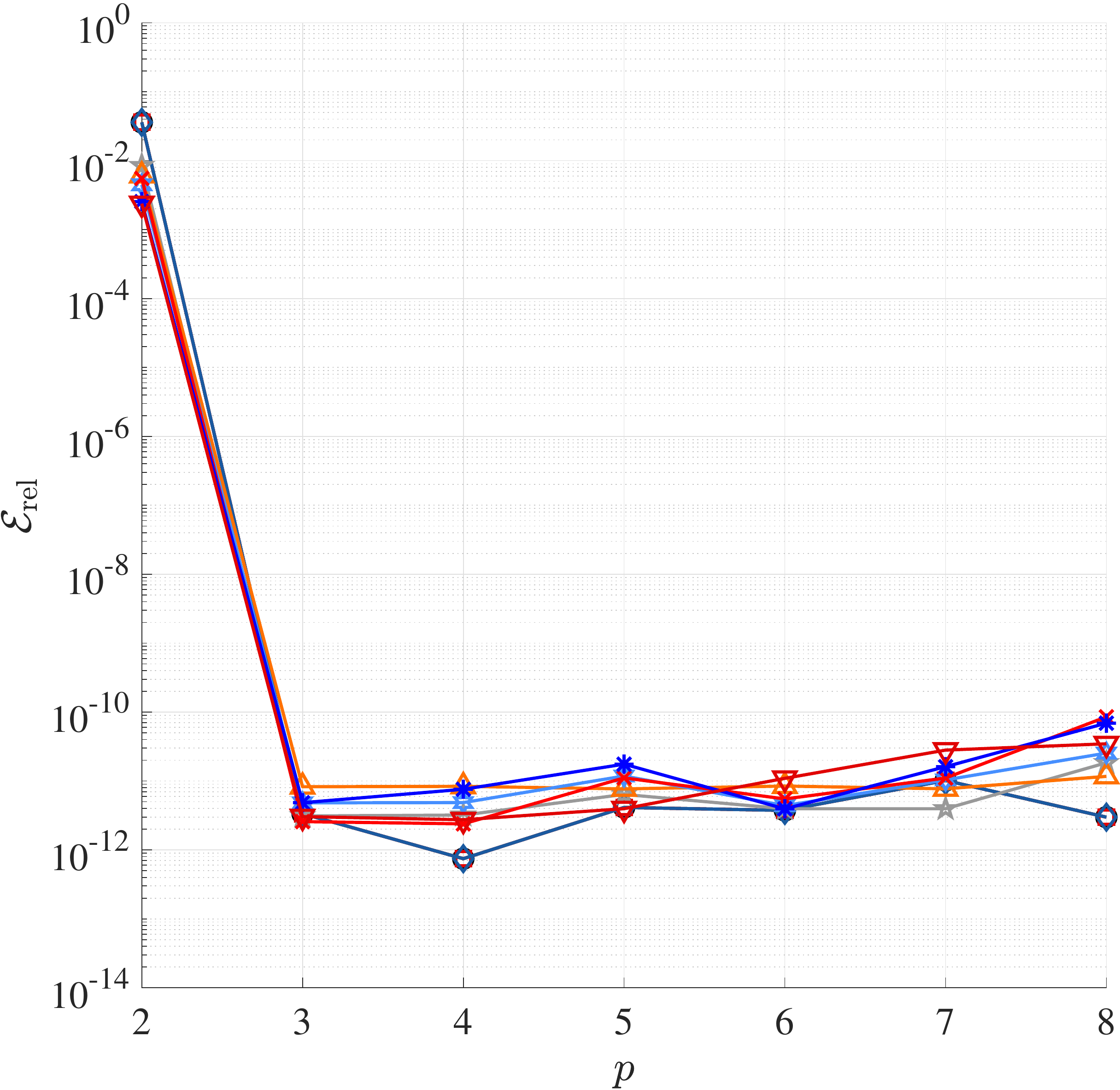}}
    \hfill
    \subfloat[Lagrange-Legendre]{\includegraphics[clip,width=0.475\textwidth]{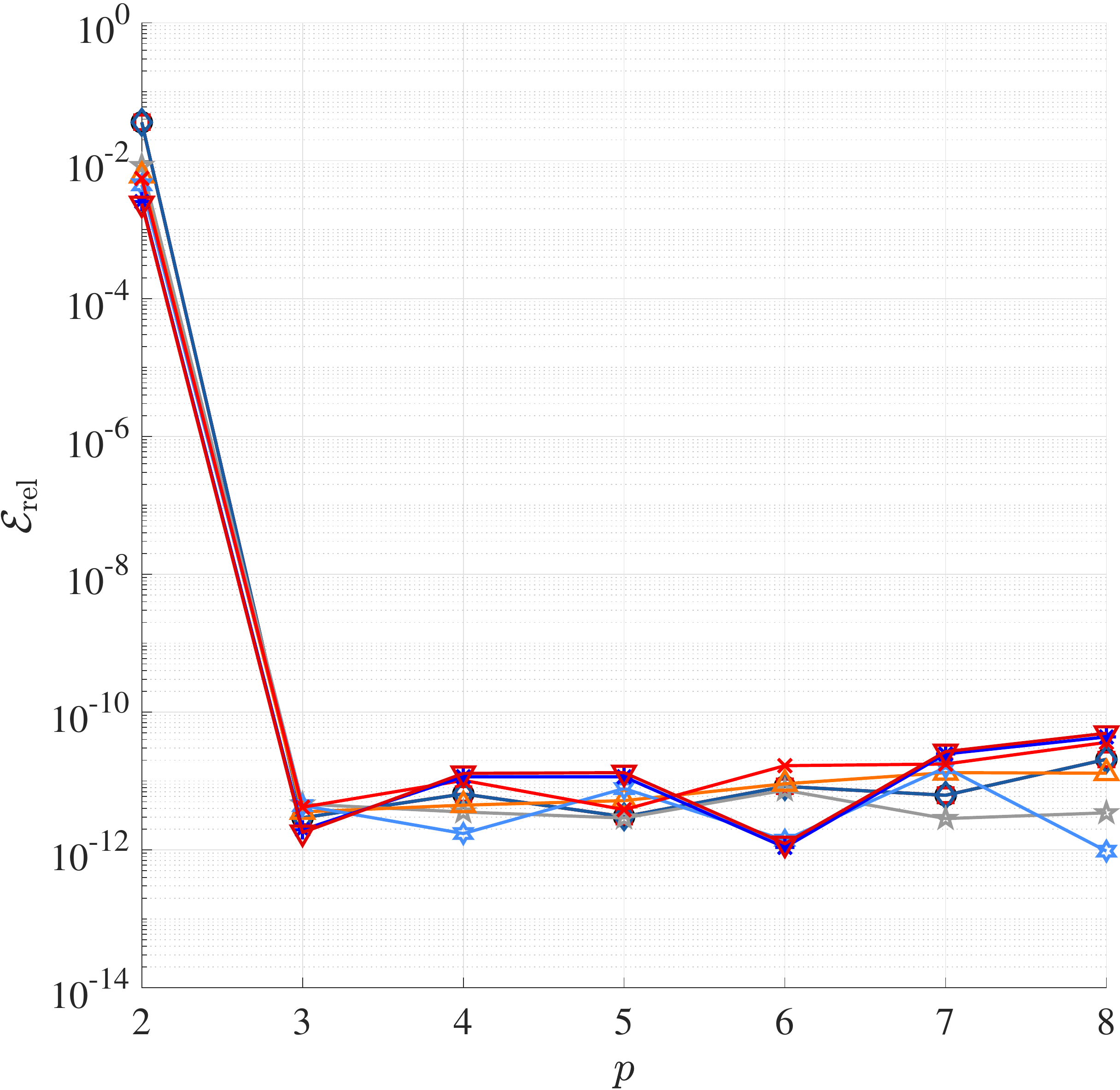}}\\
    \subfloat[Legendre-Legendre]{\includegraphics[clip,width=0.475\textwidth]{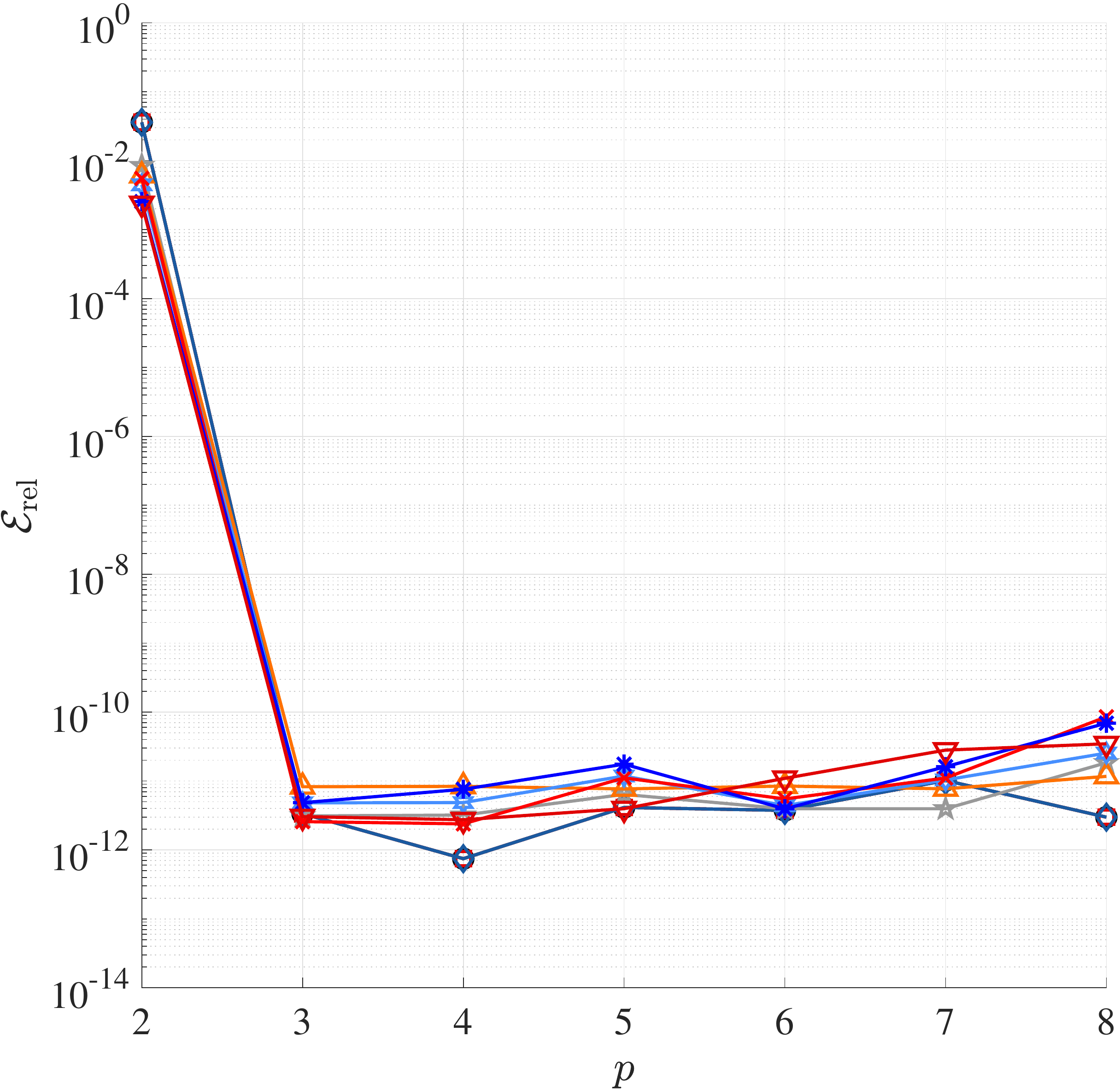}}
    \hfill
    \subfloat[Legendre-Lagrange]{\includegraphics[clip,width=0.475\textwidth]{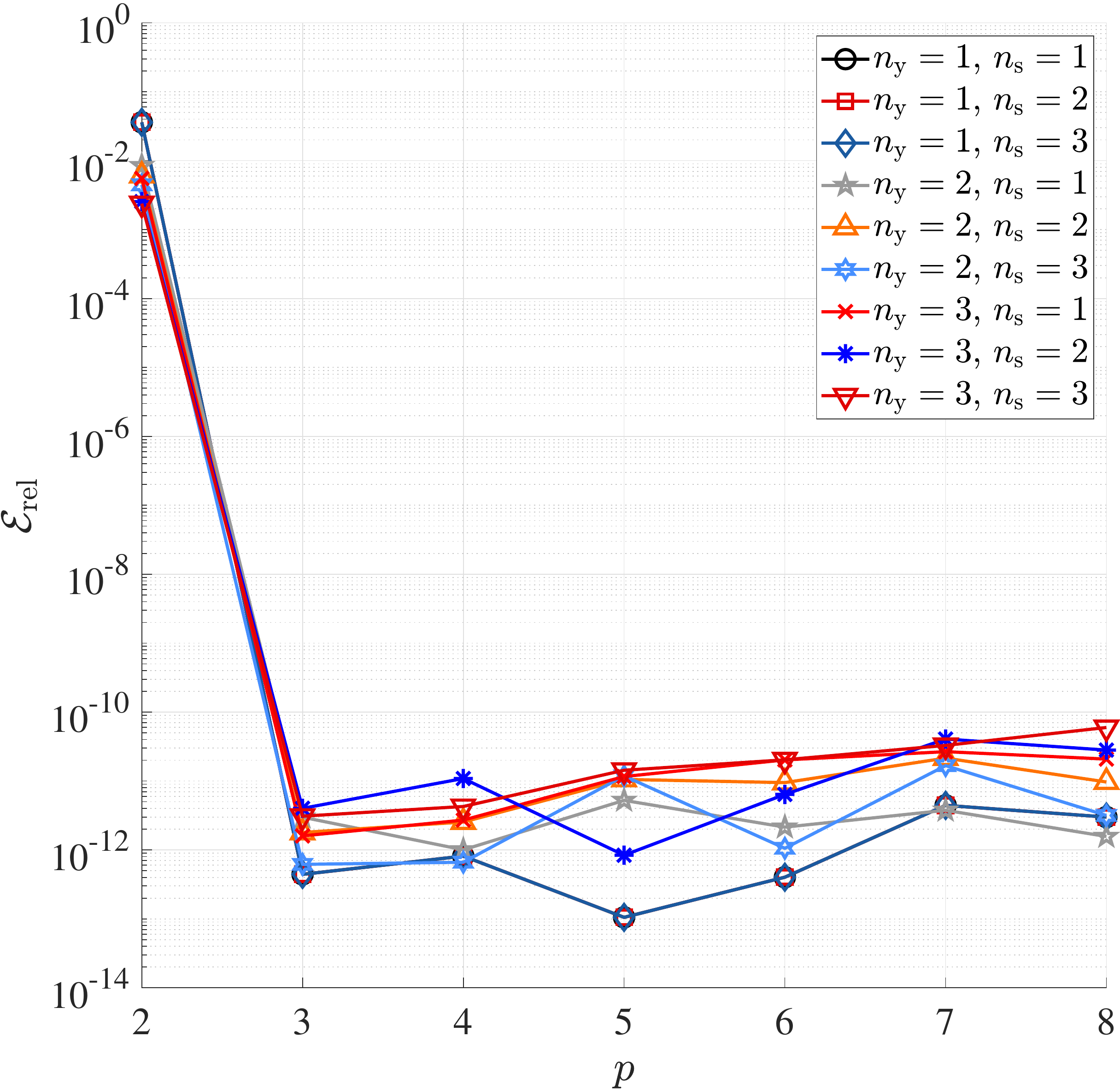}}\\
	\caption{Relative error for the cubic patch test example ($p_\mathrm{x}\,{=}\,p_\mathrm{y}$, GLL nodal distribution).}
	\label{fig:ErrorCubPatchTest}
\end{figure}%
The theoretical displacement field and stress distributions are depicted in Fig.~\ref{fig:CubPatchTestUS}. The external load, shear force, and bending moment are modeled as a pressure load and a surface traction. For the numerical test, we assume a Poisson's ratio of $\nu\,{=}\,0.3$. Hence, at the left edge of the beam, we apply
\begin{align}\nonumber
f_{\mathrm{x},\Gamma_1} &= \cfrac{240}{c}\,y - 120\,,\\\nonumber
f_{\mathrm{y},\Gamma_1}  &= -\cfrac{120}{L}\,y + \cfrac{120}{cL}\,y^2\,,
\end{align}
while
\begin{equation*}
f_{\mathrm{y},\Gamma_2} = \cfrac{120}{L}\,y - \cfrac{120}{cL}\,y^2\,.
\end{equation*}

The cubic patch test as introduced above can only be passed by transition elements that have a minimum polynomial degree of 3. Therefore, all combinations that include either $p_\mathrm{x}\,{=}\,2$ or $p_\mathrm{y}\,{=}\,2$ will fail this test and exhibit significantly increased error values. The discretization for this example is identical to the one discussed in Sect.~\ref{sec:PatchTestQuadratic}. For all models, relative errors of less than $10^{-10}$ in the stresses have been obtained for the cubic (linear bending) patch test, see Fig.~\ref{fig:ErrorCubPatchTest}. Again, we observe increasing error values for higher orders of the shape functions. Nevertheless, the error levels for transition elements that only feature a quadratically complete ansatz space are several orders of magnitude higher than for those that contain all terms of a complete polynomial of order 3. In the case of a cubically complete shape function, the error is reduced by eight orders of magnitude (error: $10^{-2}$ vs. $10^{-10}$; see Fig.~\ref{fig:ErrorCubPatchTest}).
\subsection{Higher order patch test}
\label{sec:PatchTestHigh}
It is well-known that when using higher order polynomial basis functions, a higher rate of convergence can be achieved \cite{BookSzabo1991}. To assess the (asymptotic) rate of convergence, a higher order patch test can be employed \cite{ArticleZienkiewicz1997}. Here, we prescribe a displacement field of the form
\begin{equation}
\mathbf{u} = \sum\limits_{i}^{} \mathbf{a}_iP^{p}_i(x,y)\,, \quad \forall\; i \in\{1,2,\ldots,n_\mathrm{p}\}
\label{eq:HighOrderDisplacementField}
\end{equation}
\begin{figure}[t!]
	\centering
	\includegraphics[clip,scale=0.8]{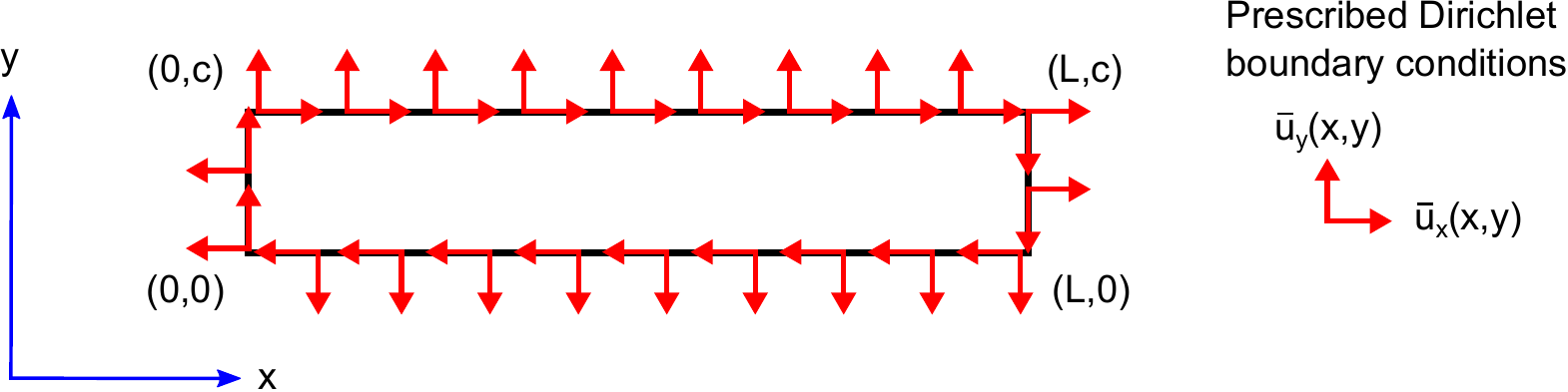}
	\caption{Model for the higher order patch test; Geometry and boundary conditions.}
	\label{fig:PatchTestHigherOrder}
\end{figure}%
\begin{figure}[b!]
	\centering
	\subfloat[Lagrange-Lagrange]{\includegraphics[clip,width=0.475\textwidth]{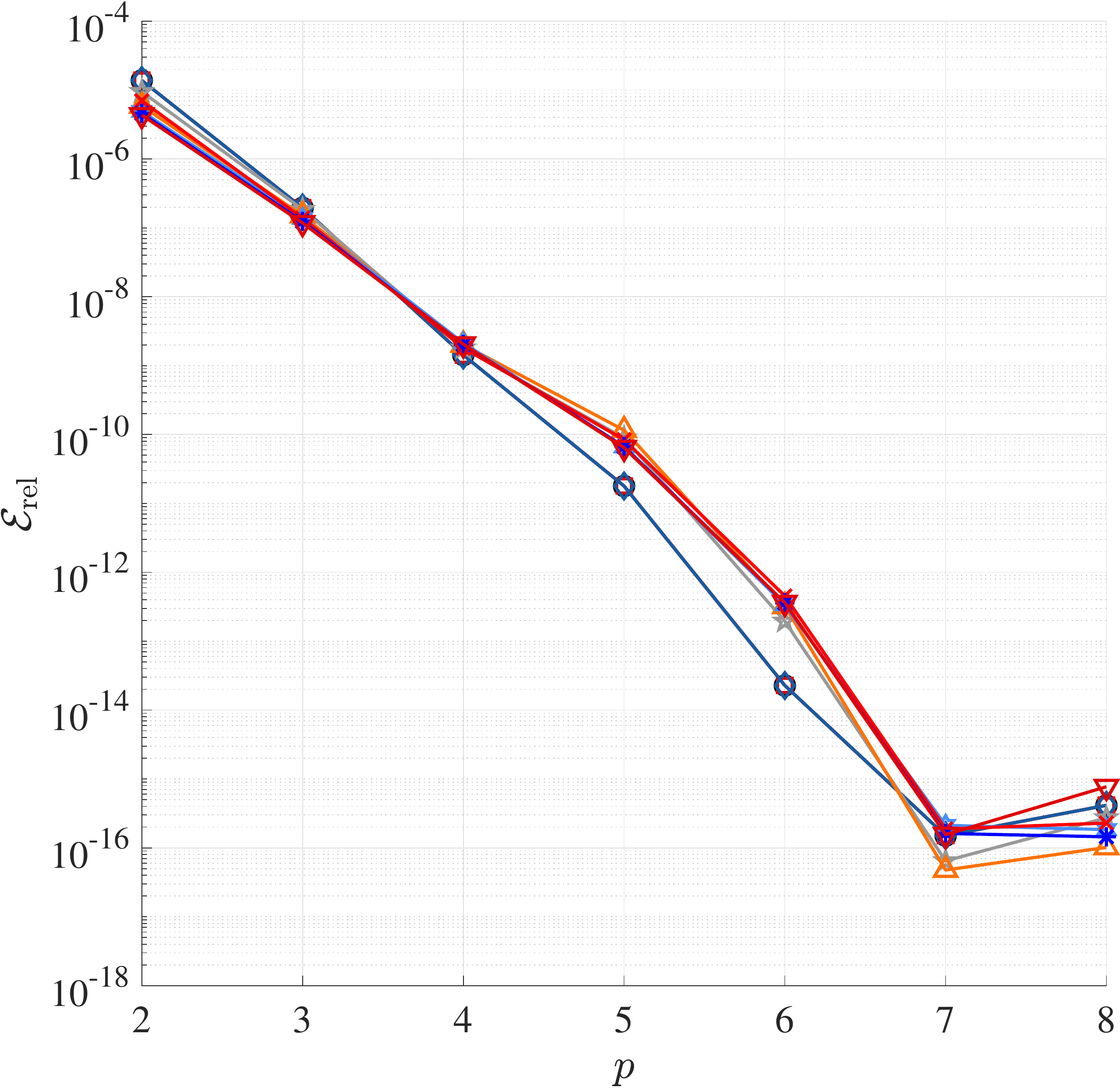}}
	\hfill
	\subfloat[Lagrange-Legendre]{\includegraphics[clip,width=0.475\textwidth]{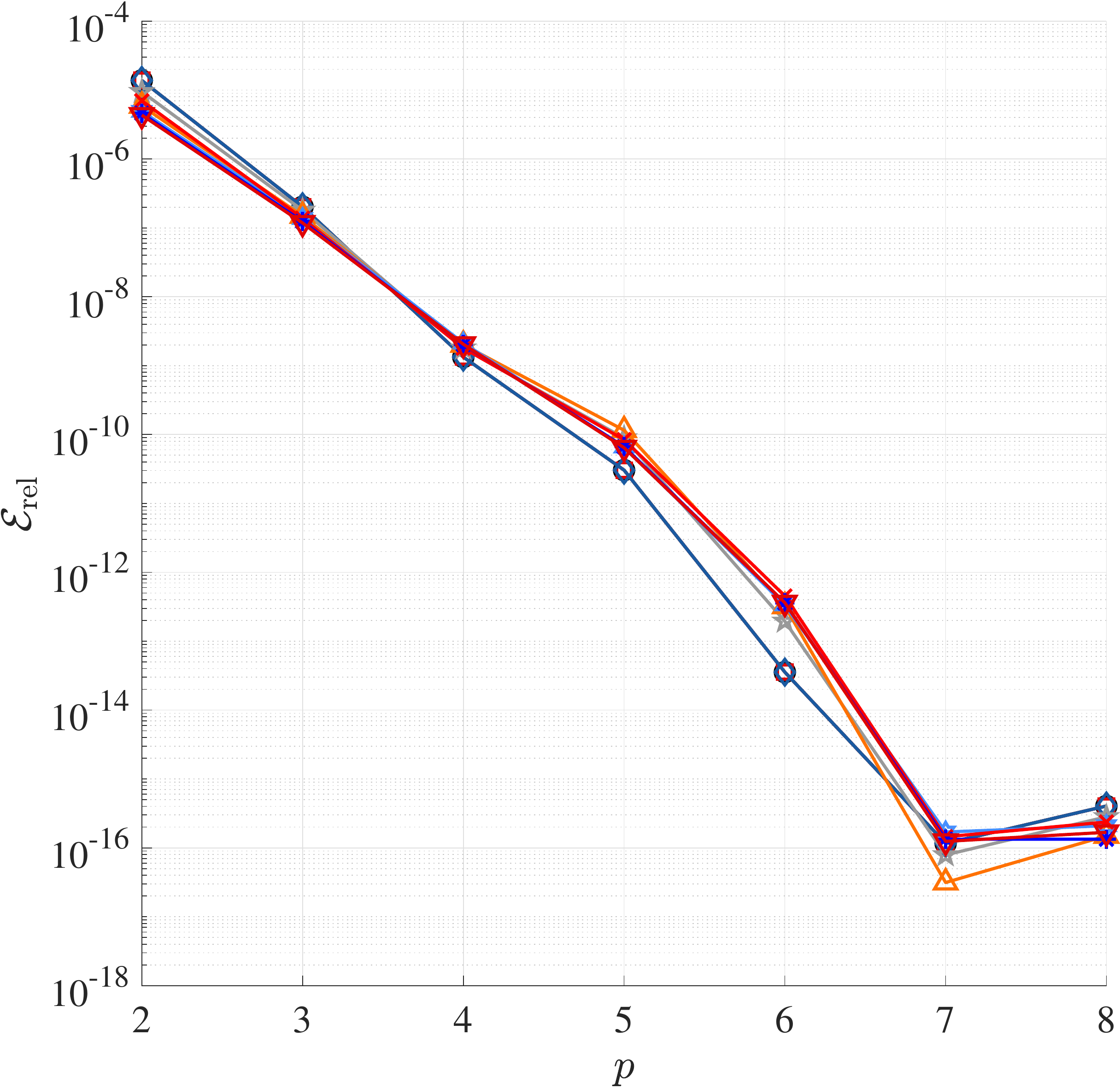}}\\
	\subfloat[Legendre-Legendre]{\includegraphics[clip,width=0.475\textwidth]{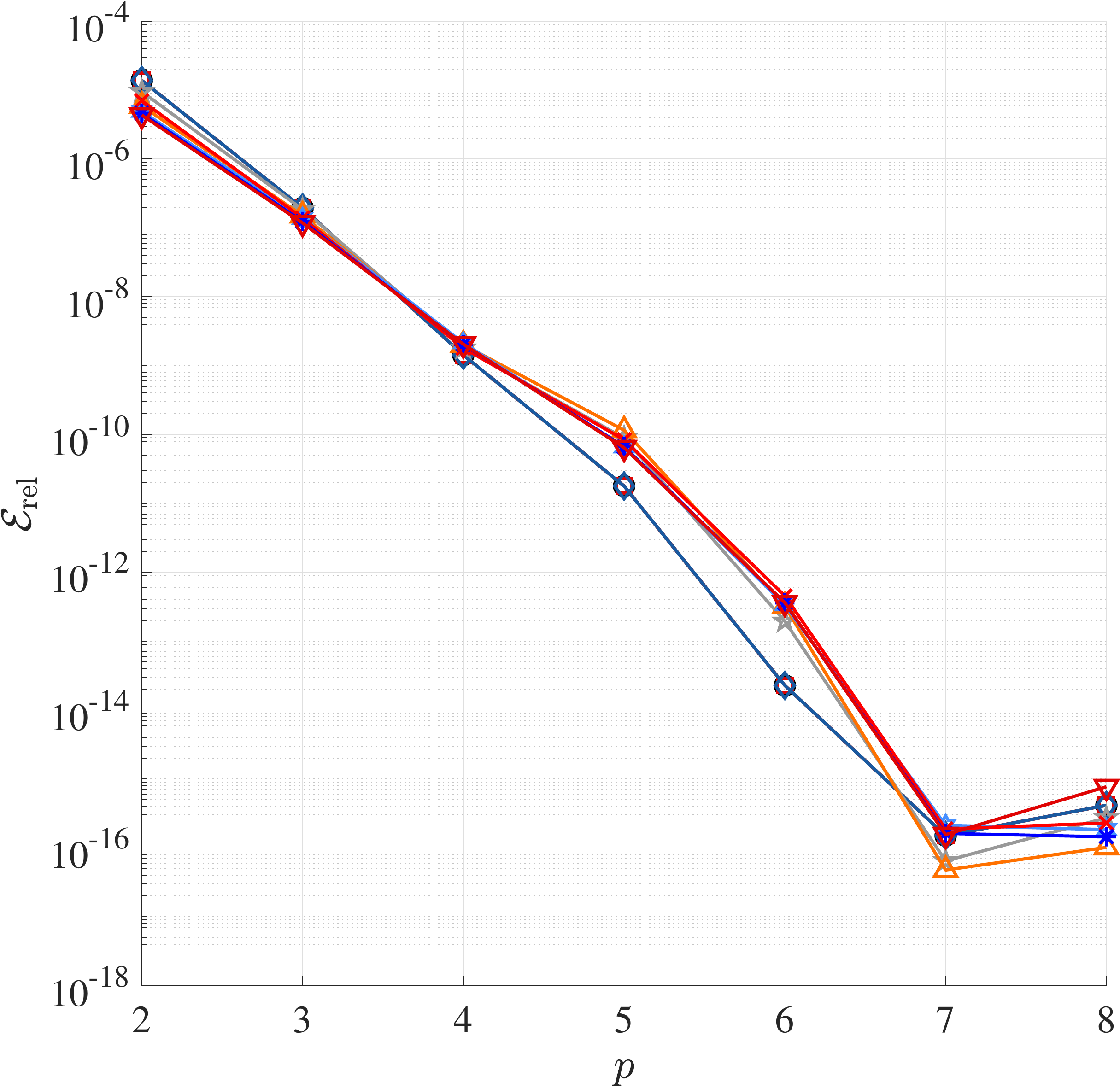}}
	\hfill
	\subfloat[Legendre-Lagrange]{\includegraphics[clip,width=0.475\textwidth]{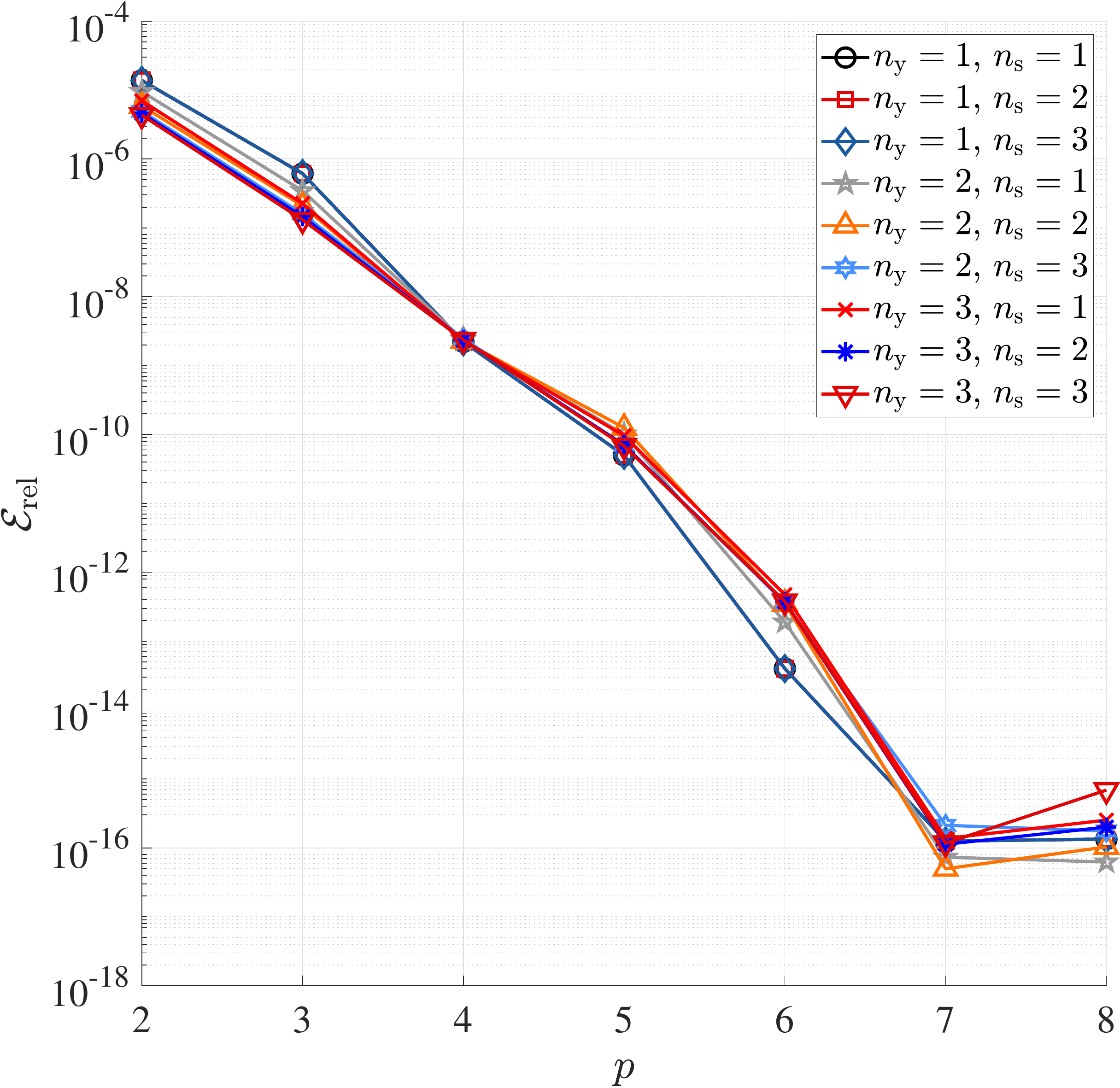}}\\
	\caption{Relative error for the high order patch test example ($p_\mathrm{x}\,{=}\,p_\mathrm{y}$, GLL nodal distribution, seventh-order displacement field).}
	\label{fig:ErrorHighPatchTest}
\end{figure}%
at the boundary of the computational domain (see Fig.~\ref{fig:PatchTestHigherOrder}), where $\mathbf{a}_i$ denotes unknown parameters of the polynomial functions, and $n_\mathrm{p}$ represents the number of monomials in the extension
$$n_\mathrm{p} = \cfrac{(p+1)(p+2)}{2}\,,$$
while $\mathbf{P}^{p}$ denotes the vector containing all monomials that construct a complete polynomial of order $p$ for the higher order patch test. For a linear patch test it is sufficient to employ $\mathbf{P}^1\,{=}\,[1,x,y]^\mathrm{T}$ and a cubic patch test requires $\mathbf{P}^3\,{=}\,[1,x,y,x^2,xy,y^2,x^3,x^2y,xy^2,y^3]^\mathrm{T}$. Once a finite element cannot pass a patch test with a higher order polynomial displacement field as the solution, we have obtained an estimate for the (asymptotic) rate of convergence. This approach can be employed for the version B patch test, where we prescribe the displacement at the whole boundary of the structure. Note that the high order displacement field cannot be chosen arbitrarily but needs to fulfill the equilibrium conditions in order to be admissible for the intended purpose, i.e., the coefficients $\mathbf{a}_i$ in Eq.~\eqref{eq:HighOrderDisplacementField} have to be determined from
\begin{equation}
\mathrm{div}(\boldsymbol{\sigma}) = \mathbf{0}\,,
\end{equation}
with $\boldsymbol{\sigma}$ denoting the stress tensor. A detailed derivation of the prescribed displacement fields is provided in \ref{sec:PatchTestHigh_Disp}.

The numerical results for the different admissible displacement fields clearly demonstrate that the asymptotic convergence rate  of the proposed transition elements depends only on $p_\mathrm{min}$, i.e., the error in the (strain) energy norm should decrease with
\begin{equation}
\mathcal{E}_\mathrm{SE} \propto Ch^{p_\mathrm{min}}\,,
\end{equation}
where is $C$ a positive constant and $h$ denotes the element size. As an example, we plotted the mean error in the nodal displacements for a seventh-order displacement field in Fig.~\ref{fig:ErrorHighPatchTest}. In this example, we observe the convergence to the exact solution for a \emph{p}-refinement. The error decreases monotonically until machine precision is reached for $p_\mathrm{x}\,{=}\,p_\mathrm{y}\,{=}\,7$. Since the prescribed displacement field is smooth (without the presence of singularities), optimal rates of convergence are expected. It is noted that the convergence behavior is almost identical for all different transition elements which leads to the conclusion that the methodology has been successfully implemented. With this particular example, it is shown for the first time that the optimal (expected) rates of convergence can be recovered using the transfinite mapping approach. To date, several applications of related approaches have been discussed, but according to the authors' knowledge, the convergence properties have not been addressed before.

In Sect.~\ref{sec:PatchTest}, we have shown that the proposed transition elements are capable of passing various patch tests. It is generally accepted in the computational mechanics community that this means convergence will be ensured. By means of the high order patch test we were successful in demonstrating that the theoretically predicted convergence rates can be obtained for problems with smooth solutions. In Sect.~\ref{sec:Singularity}, we apply the proposed transition elements to benchmark examples that feature singularities (re-entrant corners).

%% file: tex/ConvergenceRate.tex
\section{Convergence Rate}
\label{sec:ConvRate}
Passing the high order patch test is a valuable means to determine the asymptotic rate of convergence for standard finite elements (as discussed in Sect.~\ref{sec:PatchTest}) where the associated shape functions are continuously differentiable within the element domain \cite{BookZienkiewicz2000a}. This assumption is made use of when deriving the error estimates. However, considering transition elements, we must be aware of the fact that the shape functions exhibit kinks at the coupling edges which impact the whole element. Therefore, two examples are selected in this section to demonstrate that the convergence rates, predicted in Sect.~\ref{sec:PatchTestHigh}, are indeed correct. The error in the displacement field is measured in the $L_2$-norm
\begin{equation}
\mathcal{E}_{\mathrm{rel}}^{\mathrm{L}_2} = \cfrac{||\mathbf{u}_\mathrm{ex} - \mathbf{u}_\mathrm{num}||_{\mathrm{L}_2(\Omega)}}{||\mathbf{u}_\mathrm{ex}||_{\mathrm{L}_2(\Omega)}}\,.
\end{equation}
\subsection{Prescribed eighth order polynomial displacement field}
In this subsection, the example of the high order displacement field is revisited, and an \emph{h}- refinement is executed. The model is the same as reported in Sect.~\ref{sec:PatchTestHigh}. In order to compute the convergence rates, we partition the model at $x\,{=}\,L/2$ and discretize it initially with one \emph{x}- and one \emph{y}-element\footnote{Remark: The initial or base mesh for this method is generated using the commercial FE-software ABAQUS. In the pre-processing stage the elements are divided into two sets, one corresponding to the class of \emph{x}- and the other to the class of \emph{y}-elements. This mesh is the basis for all simulations conducted in this paper. The mesh data is in the next step imported to MATLAB, and \emph{x\textbf{N}y}-elements are set up at the interface between the two domains associated with the different element types. As a last step in the mesh generation process, the \emph{y}-elements are divided according to the parameters $n_\mathrm{y}$ and $n_\mathrm{s}$ which also entails setting up \emph{y\textbf{N}y}-elements if $n_\mathrm{s}\,{\ge}\,2.$}. This two-element mesh is refined six times by subdividing each quadrilateral element into four smaller ones. Hence, the finest base mesh contains 8,192 finite elements, where one half are \emph{x}- and the other \emph{y}-elements. For the purpose of this study, we choose the following parameters: $n_\mathrm{y}\,{=}\,2$, $n_\mathrm{s}\,{=}\,1$, $p_\mathrm{x},p_\mathrm{y}\,{=}\,1,2,\ldots,8$. The \emph{x}-elements are of the spectral type (Lagrange shape functions), while the \emph{y}-elements are \emph{p}-elements (Legendre shape functions).
\begin{figure}[b!]
	\centering
	\subfloat[$p_\mathrm{x}\,{=}\,8$]{\includegraphics[clip,width=0.485\textwidth]{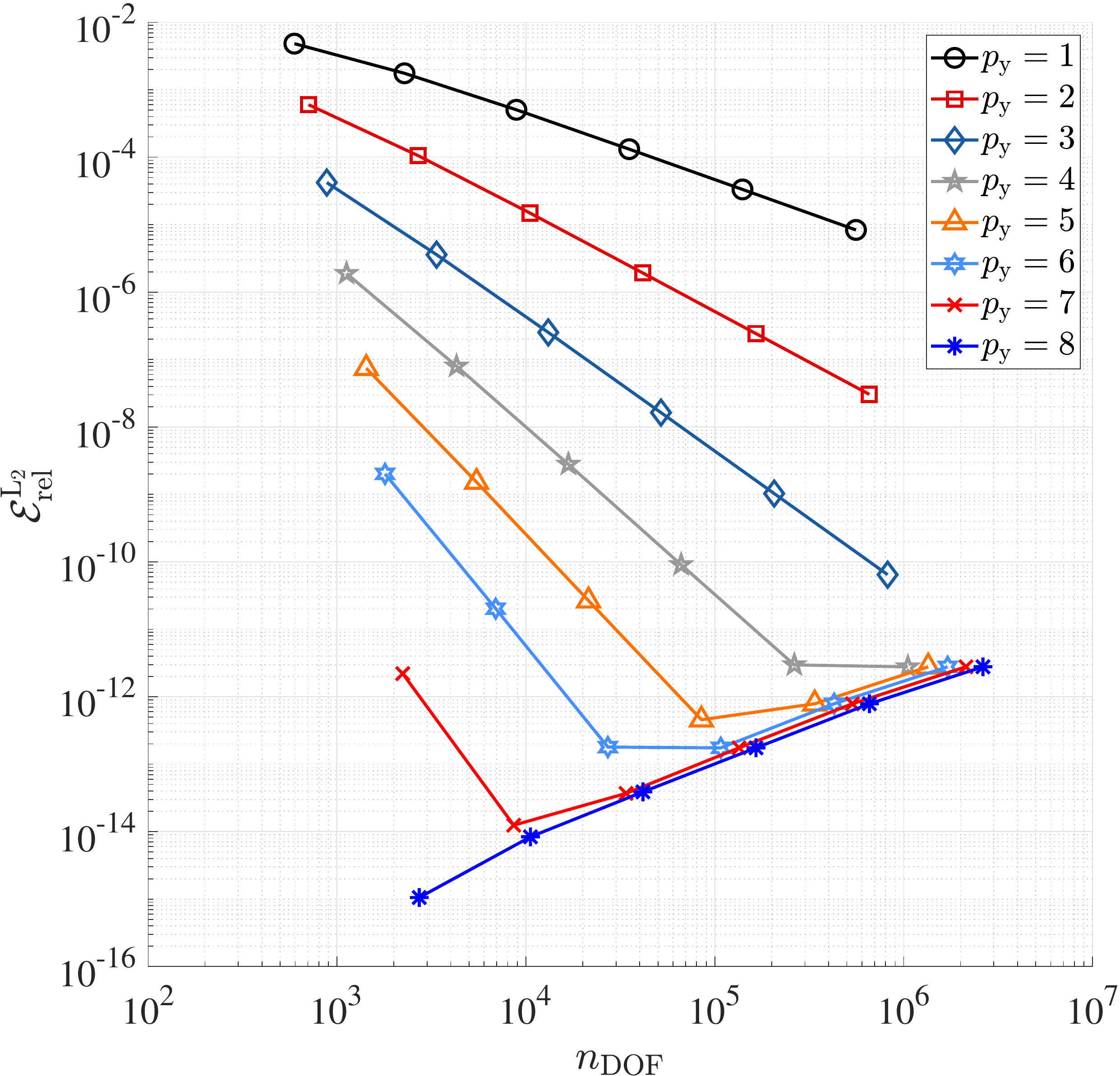}}
	\hfill
	\subfloat[$p_\mathrm{x}\,{=}\,p_\mathrm{y}$]{\includegraphics[clip,width=0.475\textwidth]{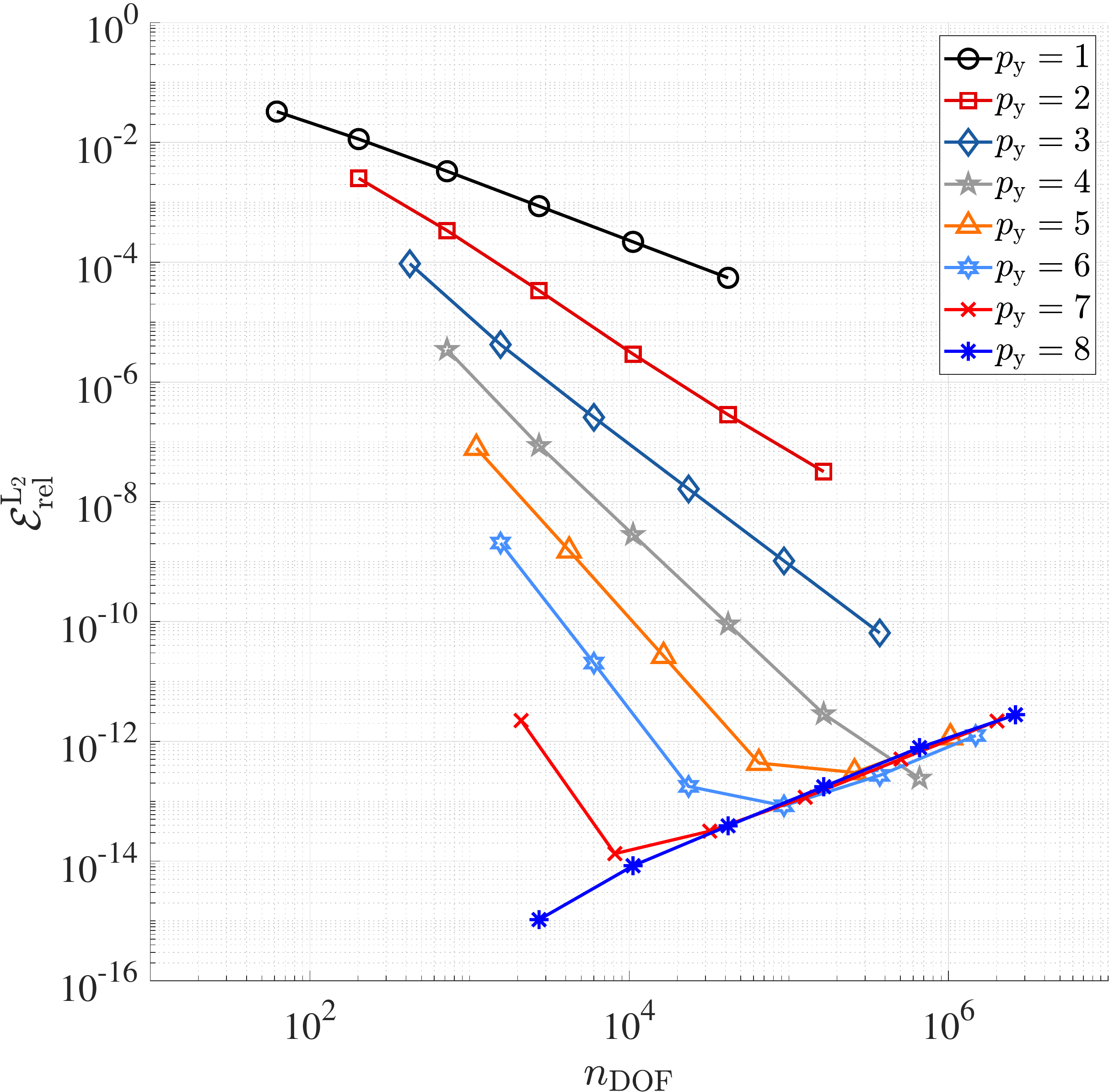}}
	\caption{Relative displacement error in the $L_2$-norm under \emph{h}-refinement -- prescribed high order displacement field.}
	\label{fig:ConvRate1}
\end{figure}%
\begin{table}[b!]
	\centering
	\caption{Convergence rates -- high order displacement field. \label{tab:ConvRate1}}
	\begin{tabular}{c|c|c|c}
		\toprule
		$p_\mathrm{y}$ & $p_\mathrm{x}\,{=}\,8$ & $p_\mathrm{x}\,{=}\,p_\mathrm{y}$ & theoretical \\\hline
		1 &	-0.9716 &	-0.9818	& 	-1.0\\
		2 &	-1.4819 &	-1.6817	& 	-1.5\\
		3 &	-1.9580 &	-2.0946	& 	-2.0\\
		4 &	-2.4477 &	-2.5789	& 	-2.5\\
		5 &	-2.9369 &	-2.9720	& 	-3.0\\
		6 &	-3.4339 &	-3.4503	& 	-3.5\\
		7 &	-3.8258 &	-3.7860	& 	-4.0\\
		8 &	--     &	--   	& 	-4.5\\
		\bottomrule
	\end{tabular}
\end{table}

The results of the convergence analysis are depicted in Fig.~\ref{fig:ConvRate1} for two different set-ups. In one example, only the polynomial degree of the \emph{y}-elements is varied, while the polynomial degree of the \emph{x}-elements is fixed at $p_\mathrm{x}\,{=}\,8$. The other idea is to increase both  $p_\mathrm{x}$ and  $p_\mathrm{y}$ simultaneously. Overall, a similar behavior as in the patch test examples is seen. For very fine meshes, the implementation seems to be suboptimal in that round-off errors accumulate and thus dominate the overall error. Bear in mind that the numerical rate of convergence can only be determined for models where at least one of the polynomial degrees is less than eight, as otherwise the exact solution is already contained in the ansatz space. The theoretically optimal convergence rate (with respect to the number of degrees of freedom $n_\mathrm{DOF}$) is predicted as
\begin{equation}
\chi = \cfrac{\log_{10}\left(\mathcal{E}_{\mathrm{rel,2}}^{\mathrm{L}_2}/\mathcal{E}_{\mathrm{rel,1}}^{\mathrm{L}_2}\right)}{\log_{10}\left(n_\mathrm{DOF,2}/n_\mathrm{DOF,1}\right)} = \cfrac{p_\mathrm{min}+1}{2}\,.
\end{equation}
In Fig.~\ref{fig:ConvRate1}, we clearly observe the expected fan-shape of all curves. The rates of convergence are listed in Table~\ref{tab:ConvRate1} and confirm our statements provided in Sect.~\ref{sec:PatchTestHigh}, i.e., optimal rates of convergence are attained.
\subsection{Infinite two-dimensional domain with circular hole under uni-axial loading}
\begin{figure}[b!]
	\centering
	\includegraphics[clip,scale=0.75]{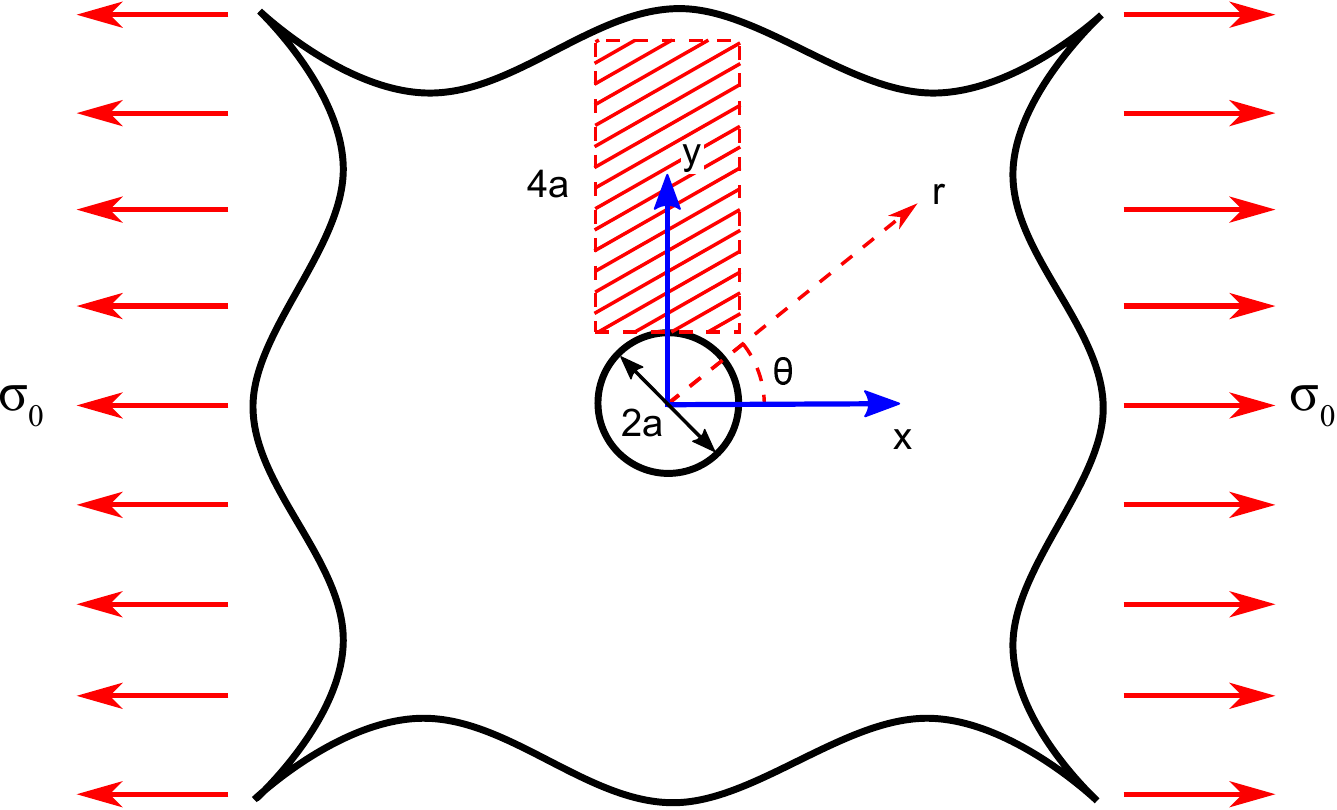}
	\caption{Circular hole in an infinite domain under uni-axial tension.}
	\label{fig:ModelInfinite}
\end{figure}

The second example to study the convergence rates of the proposed transition elements is the well-known infinite domain with a circular hole under uni-axial loading as depicted in Fig.~\ref{fig:ModelInfinite}. The reason behind this choice is twofold: \textit{(i)} existence of an analytical solution and \textit{(ii)} smooth trigonometric displacement field, i.e., the exact solution is not contained in the ansatz of the \emph{x\textbf{N}y}-elements. The analytical solution is given in polar coordinates as
\begin{align}
u_\mathrm{x}(r,\theta) & = \cfrac{\sigma_0 a}{8G} \left[ \cfrac{r}{a}\,(\kappa+1)\cos(\theta) + \cfrac{2a}{r}\,[(1+\kappa)\cos(\theta)+\cos(3\theta)] - \cfrac{2a^3}{r^3}\,\cos(3\theta) \right]\,, \\
u_\mathrm{y}(r,\theta) & = \cfrac{\sigma_0 a}{8G} \left[ \cfrac{r}{a}\,(\kappa-3)\sin(\theta) + \cfrac{2a}{r}\,[(1-\kappa)\sin(\theta)+\sin(3\theta)] - \cfrac{2a^3}{r^3}\,\sin(3\theta) \right]\,.
\end{align}
Here, $\sigma_0$ is a uniform uni-axial loading, $G$ denotes the shear modulus, and $\kappa$ is Kolosov's constant\footnote{Remark: The parameter $\kappa$ takes different values depending on whether a plane stress or strain problem is investigated. In the first case, we obtain $\kappa = \cfrac{3-\nu}{1+\nu}$, while in the second one the following definition is valid $\kappa = 3-4\nu$ with Poisson's ratio $\nu$.}. The material parameters for this example correspond to that of aluminum ($E\,{=}\,70\,$GPa, $\nu\,{=}\,0.3$). The radius of the circular hole is given by $a$, whereas $r$ and $\theta$ represent the polar coordinates. The transformation between the Cartesian coordinate system used for the simulation and the polar one is given as 
\begin{alignat}{3}
x &= r\cos{\theta}\,, \quad&&\quad y&&=r\sin(\theta)\,,\\
r&=\sqrt{x^2+y^2}\,, \quad&&\quad \theta&&=\mathrm{atan2}\left( \cfrac{y}{x} \right)\,.
\end{alignat}
\begin{figure}[t!]
	\centering
	\subfloat[$p_\mathrm{x}\,{=}\,8$]{\includegraphics[clip,width=0.485\textwidth]{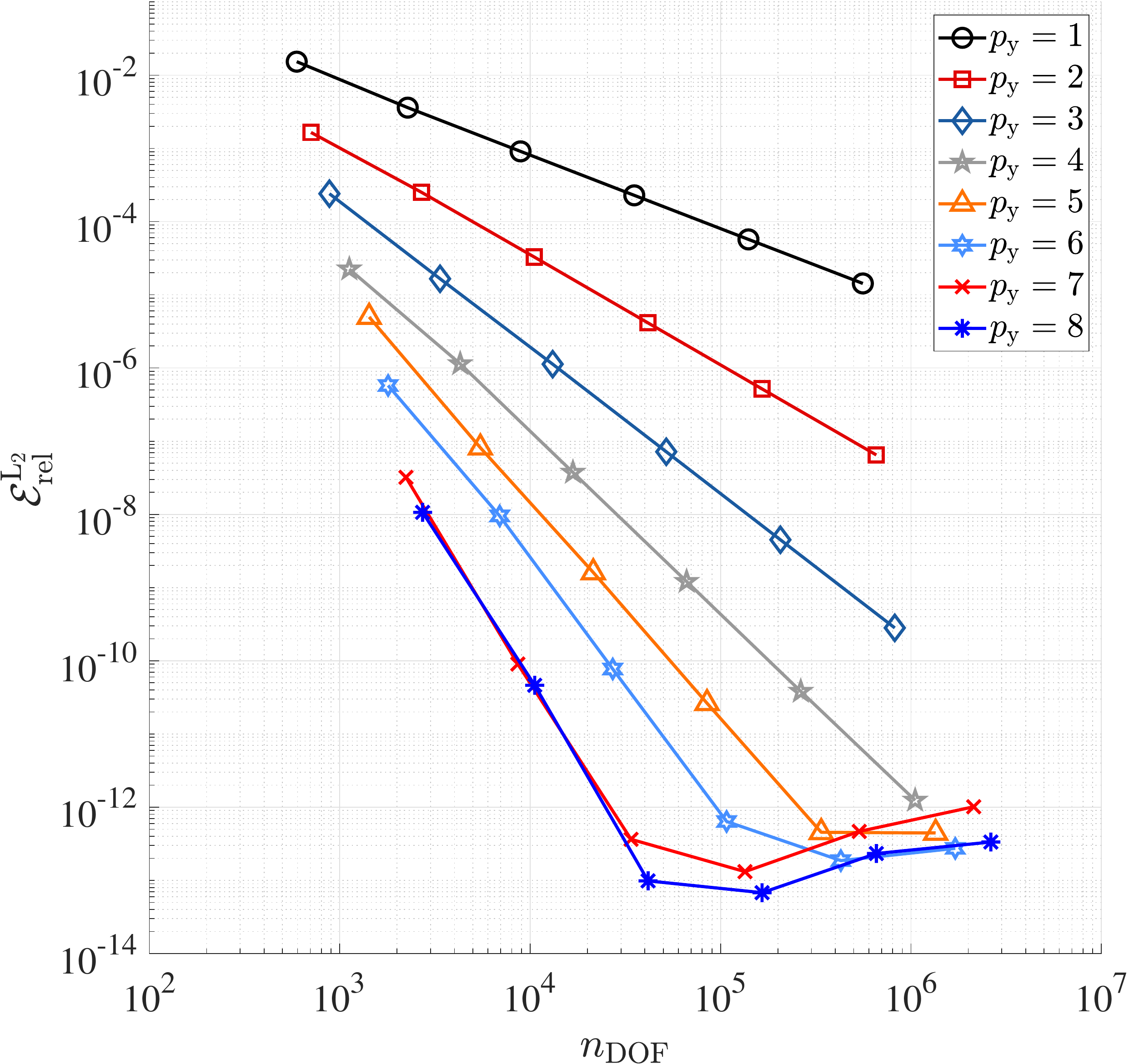}}
	\hfill
	\subfloat[$p_\mathrm{x}\,{=}\,p_\mathrm{y}$]{\includegraphics[clip,width=0.475\textwidth]{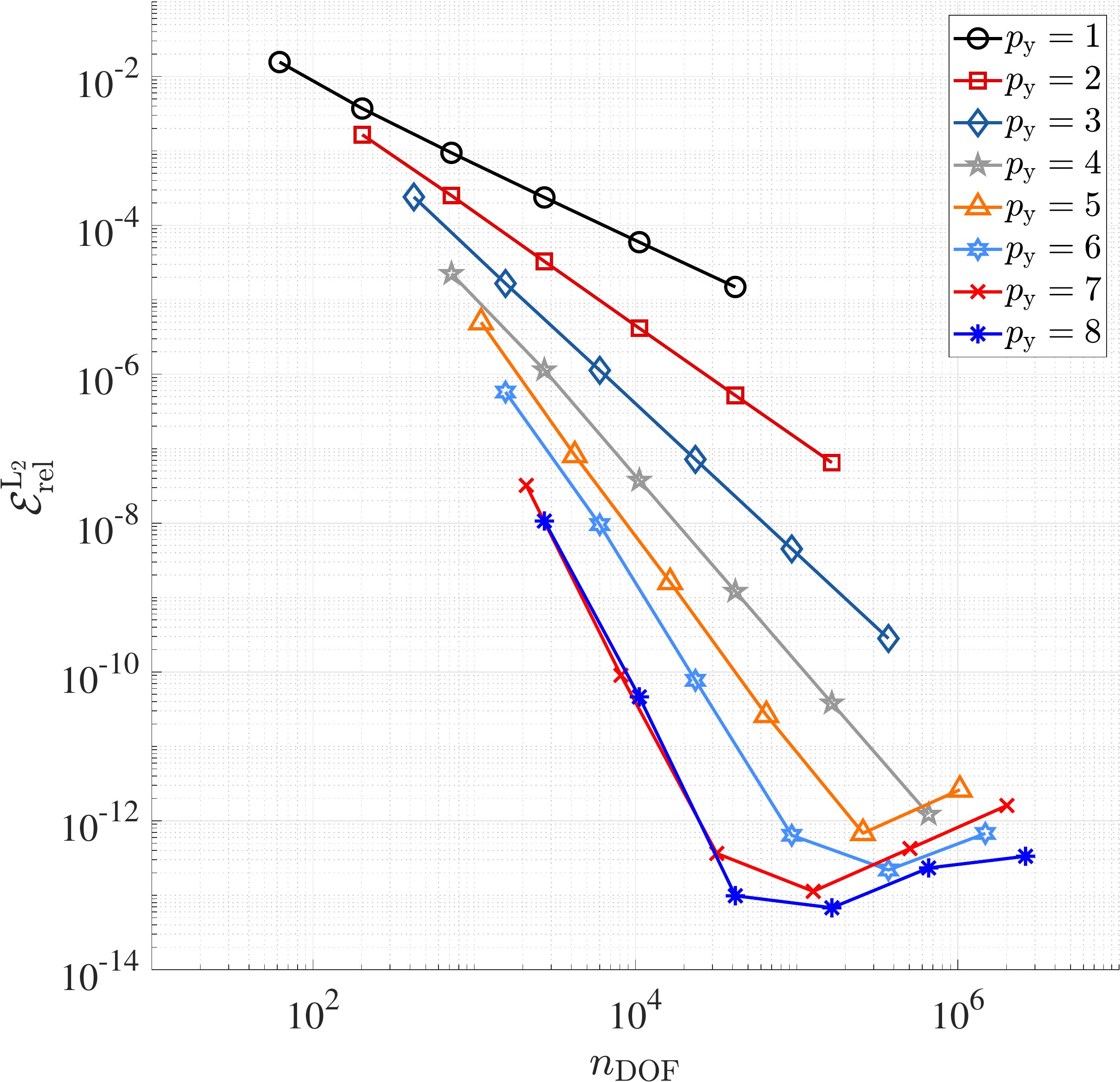}}
	\caption{Relative displacement error in the $L_2$-norm under \emph{h}-refinement - infinite domain.}
	\label{fig:ConvRate2}
\end{figure}%
\begin{table}[t!]
	\vspace*{-6pt}
	\centering
	\caption{Convergence rates - infinite domain. \label{tab:ConvRate2}}
	\begin{tabular}{c|c|c|c}
		\toprule
		$p_\mathrm{y}$ & $p_\mathrm{x}\,{=}\,8$ & $p_\mathrm{x}\,{=}\,p_\mathrm{y}$ & theoretical \\\hline
		1 &	-1.0189 &	-1.0688	& 	-1.0\\
		2 &	-1.4839 &	-1.5126	& 	-1.5\\
		3 &	-1.9976 &	-2.0150	& 	-2.0\\
		4 &	-2.4417 &	-2.4543	& 	-2.5\\
		5 &	-2.9648 &	-2.8970	& 	-3.0\\
		6 &	-3.3504 &	-3.2884	& 	-3.5\\
		7 &	-4.1831 &	-4.1849	& 	-4.0\\
		8 &	-4.2511 &	-4.2511 & 	-4.5\\
		\bottomrule
	\end{tabular}
\end{table}%
In order to avoid any effect related to the geometry approximation on the convergence rates, only the red hatched part is modeled (see Fig.~\ref{fig:ModelInfinite}) and the exact displacements are prescribed at the boundary of the domain analogously to the previous examples. The base mesh for the \emph{h}-refinement consists of two elements, i.e., one \emph{x}-element and one \emph{y}-element. \textcolor{red}{Note that when using a linear or quadratic geometry approximation of the circular hole, the geometry error will dominate the convergence behavior \cite{BookDuester2002}. Thus, an exact or at least very accurate description of the circle by means of the blending function method or quasi-regional mapping would be necessary \cite{ArticleKiralyfalvi1997}. This, however, is out of the scope of the present article and therefore, we decided for a simple domain with a structured mesh.} To this end, elements are square-shaped, the computational domain is rectangular and has the dimensions $2a\times 4a$. The obtained numerical results are illustrated in Fig.~\ref{fig:ConvRate2} and compiled in Table~\ref{tab:ConvRate2}. Also in this example, the theoretically optimal rates of convergence are recovered.

%% file: tex/SingularExamples.tex
\section{Application to Problems with Singularities}
\label{sec:Singularity}
It is well-known that exponential rates of convergence can only be attained by a pure \emph{p}-refinement approach if the solution is smooth within the computational domain \cite{BookSzabo1991, BookDuester2002}. However, if singularities are present, only the combination of \emph{h}- and \emph{p}-refinement strategies can achieve high (theoretically optimal) convergence rates. Unfortunately, such an \emph{hp}-adaptive procedure is rather difficult to implement \cite{BookDemkowicz2006, BookDemkowicz2008}, and standard procedures result in so-called hanging nodes that need to be constrained. This drawback can be circumvented either by employing the \textit{refine-by-superposition} approach presented by Zander \cite{PhDZander2017} instead of the classical \textit{refine-by-replacement} methodology or by using the proposed \emph{x\textbf{N}y}-transition elements. As discussed in the previous sections, conformal coupling along the element boundary is easily achieved without the need for imposing additional constraints. Thus, a local mesh refinement as required for problems with non-smooth solution characteristics can be implemented in a straightforward fashion. In the remainder of this section, we use the error in the energy norm to assess the performance of the transition elements in comparison to spectral elements based on a GLL nodal distribution. Thus, the error measure is defined as
\begin{equation}
\mathcal{E}_\mathrm{SE} = \sqrt{\cfrac{\mathbf{u}_\mathrm{ov}^\mathrm{T} \tilde{\mathbf{K}}\mathbf{u}_\mathrm{ov} - \mathbf{u}_\mathrm{num}^\mathrm{T} \mathbf{K}\mathbf{u}_\mathrm{num}}{\mathbf{u}_\mathrm{ov}^\mathrm{T} \tilde{\mathbf{K}}\mathbf{u}_\mathrm{num}}}\,,
\end{equation}
where $\mathbf{K}$ is the stiffness matrix and $\mathbf{u}$ the displacement vector. The subscripts $\square_\mathrm{ov}$ and $\square_\mathrm{num}$ denote the reference (overkill solution: fine mesh with high polynomial degree of the shape functions) and numerical solutions, respectively.
\subsection{Cantilever beam}
\label{sec:Cantilever}
\begin{figure}[t!]
	\centering
	\includegraphics[clip,scale=0.75]{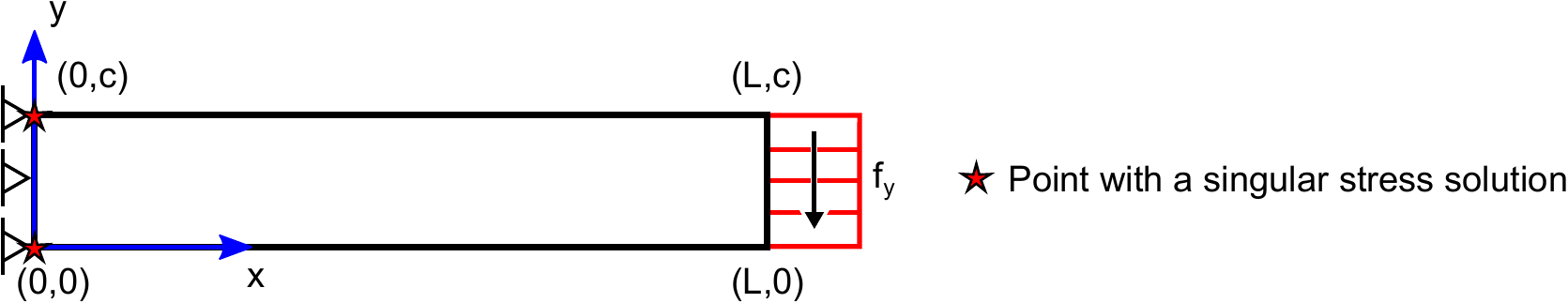}
	\caption{Model of the cantilever beam; Geometry and boundary conditions.}
	\label{fig:CantileverBeam}
\end{figure}

The first example that exhibits singular points in the solution is the cantilever beam shown in Fig.~\ref{fig:CantileverBeam}. In order to illustrate the performance of the proposed transition element, its results are compared to an established numerical method. To this end, the cantilever beam is also analyzed using the SEM. Considering the proposed transition element, a coupling of spectral elements has been chosen due to the fact the different types of coupling have shown a very similar performance in the patch test discussed in Sect.~\ref{sec:PatchTest}.

\begin{figure}[t!]
	\centering
	\subfloat[\emph{x\textbf{N}y}-mesh, $n_\mathrm{s}\,{=}\,2$, $n_\mathrm{y}\,{=}\,2$, 32 elements \label{fig:CantileverBeamMesh1}]{\includegraphics[clip,width=0.45\textwidth]{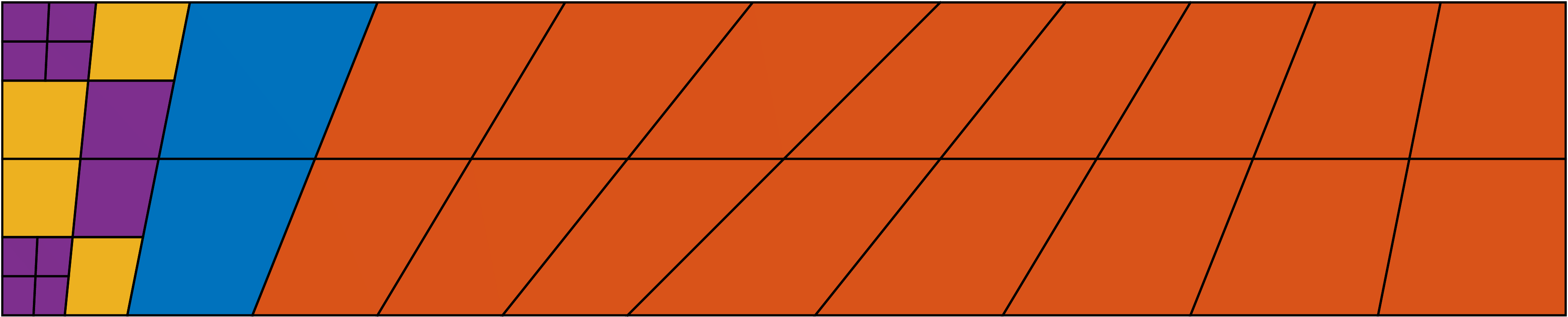}}
	\hfill
	\subfloat[\emph{x\textbf{N}y}-mesh, $n_\mathrm{s}\,{=}\,8$, $n_\mathrm{y}\,{=}\,2$, 68 elements \label{fig:CantileverBeamMesh2}]{\includegraphics[clip,width=0.45\textwidth]{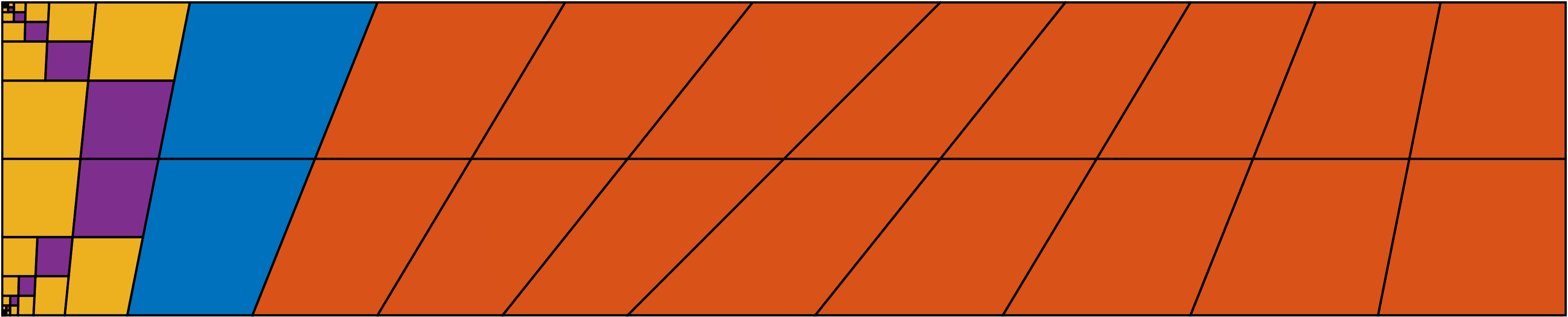}}\\
	\subfloat[Spectral element mesh, 320 elements \label{fig:CantileverBeamMesh3}]{\includegraphics[clip,width=0.45\textwidth]{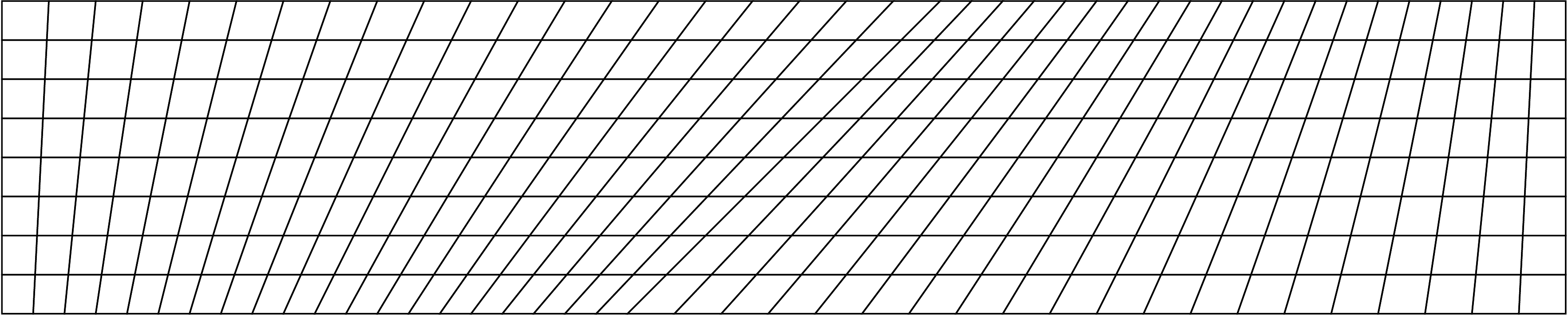}}
	\caption{Discretization of the cantilever beam model. The different element types are color-coded: \emph{x}-Elements \textcolor{xElem}{\rule[-0.25ex]{5ex}{2ex}}, \emph{y}-Elements \textcolor{yElem}{\rule[-0.25ex]{5ex}{2ex}}, \emph{x\textbf{N}y}-Elements \textcolor{xNyElem}{\rule[-0.25ex]{5ex}{2ex}}, \emph{y\textbf{N}y}-Elements \textcolor{yNyElem}{\rule[-0.25ex]{5ex}{2ex}}.}
	\label{fig:CantileverBeamMesh}
\end{figure}

\begin{figure}[t!]
	\centering
	\subfloat[{Detail view of the mesh at a singular point ($x\in[0,2] $, $y\in[0,1]$)} \label{fig:CantileverBeamMeshRef1}]{\includegraphics[clip,width=0.485\textwidth]{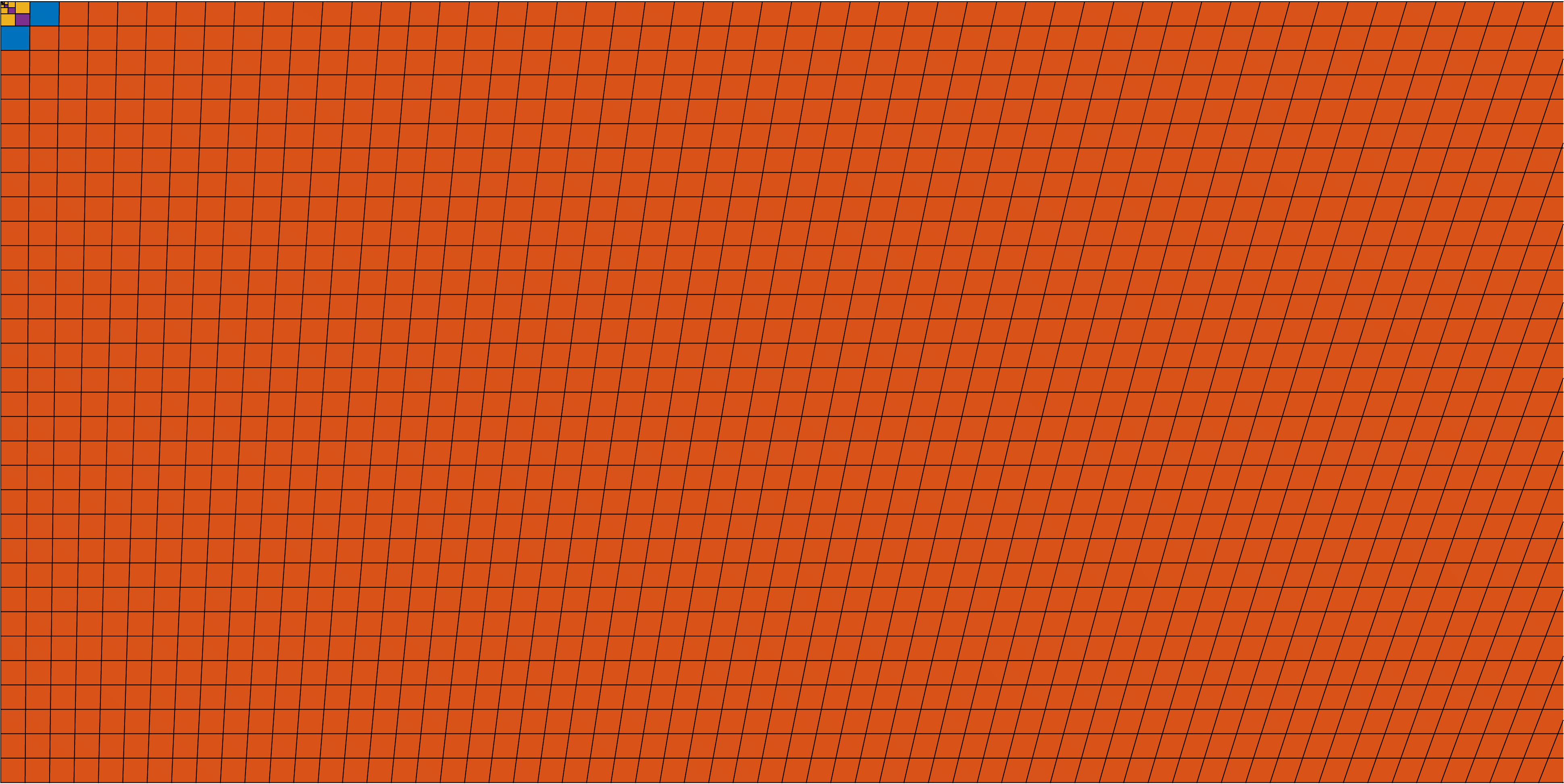}}
	\hfill
	\subfloat[{Detail view of the mesh at a singular point ($x\in[0,0.1] $, $y\in[0.9,1]$)} \label{fig:CantileverBeamMeshRef2}]{\includegraphics[clip,width=0.485\textwidth]{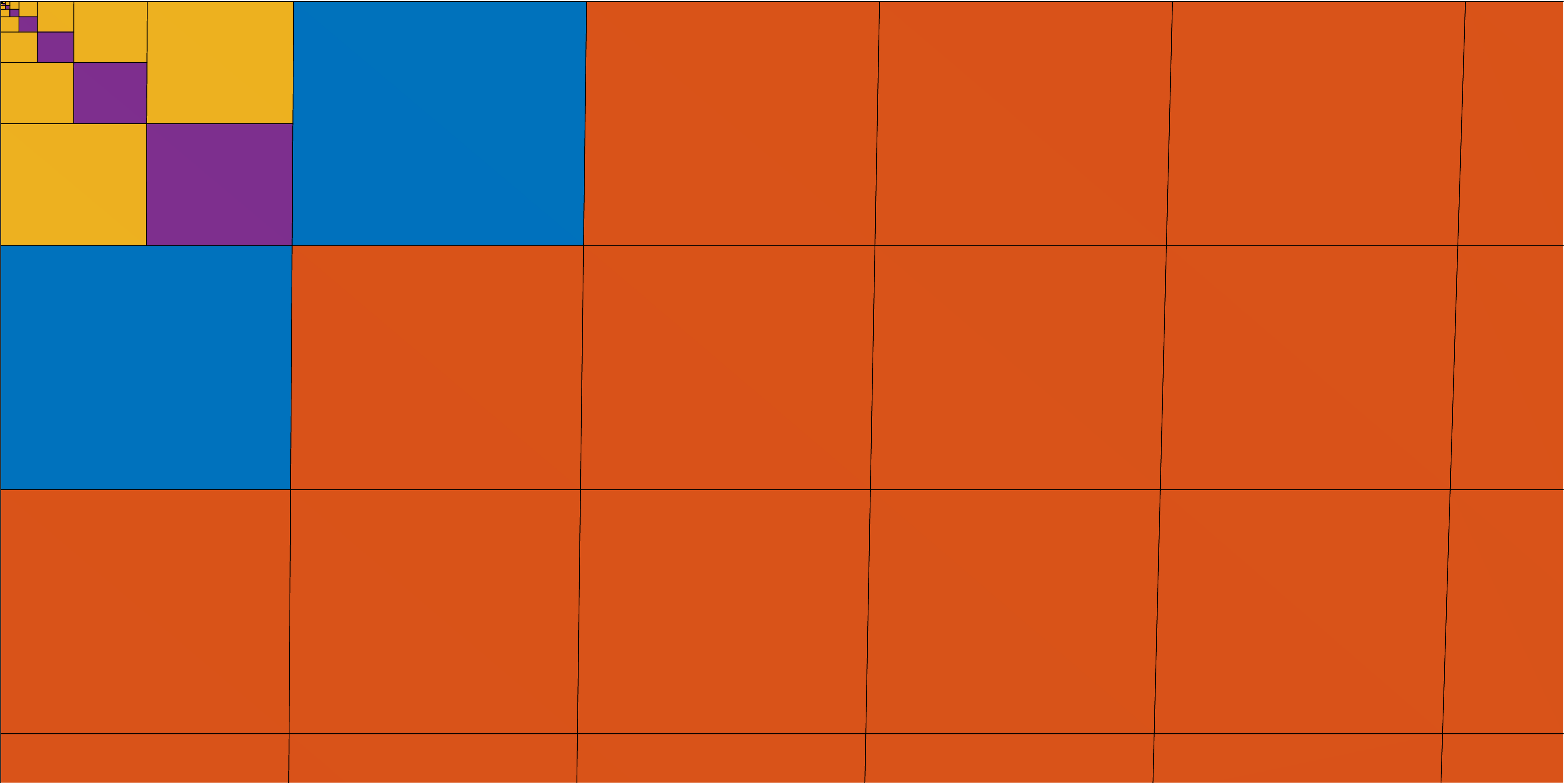}}
	\caption{Discretization of the cantilever beam model (overall mesh: 20,528 elements) - reference model (\emph{x\textbf{N}y}-mesh, $n_\mathrm{s}\,{=}\,8$, $n_\mathrm{y}\,{=}\,2$). The different element types are color-coded: \emph{x}-Elements \textcolor{xElem}{\rule[-0.25ex]{5ex}{2ex}}, \emph{y}-Elements \textcolor{yElem}{\rule[-0.25ex]{5ex}{2ex}}, \emph{x\textbf{N}y}-Elements \textcolor{xNyElem}{\rule[-0.25ex]{5ex}{2ex}}, \emph{y\textbf{N}y}-Elements \textcolor{yNyElem}{\rule[-0.25ex]{5ex}{2ex}}.}
	\label{fig:CantileverBeamMeshRef}
\end{figure}

\textcolor{red}{The discretization of the numerical model is based on the coarse mesh depicted in Fig.~\ref{fig:nsny} consisting of 20 (angularly) distorted quadrilateral finite elements. Finer meshes are constructed from this base discretization by dividing each element into four smaller ones. Thus, the spectral element mesh is constructed from the base mesh by refining it twice (see Fig.~\ref{fig:CantileverBeamMesh3}). The SEM results are obtained by executing a simple \emph{p}-refinement procedure. Two different discretizations that are utilized in conjunction with the spectral--spectral transition elements are depicted in Figs.~\ref{fig:CantileverBeamMesh1} and \ref{fig:CantileverBeamMesh2}. Here, the initial discretization is chosen, and only the two elements adjacent to the corner singularities are refined. The first mesh is constructed using the parameter combination $n_\mathrm{s}\,{=}\,2$ and $n_\mathrm{y}\,{=}\,2$, while the second mesh corresponds to $n_\mathrm{s}\,{=}\,8$ and $n_\mathrm{y}\,{=}\,2$. It is noted that the reference solution for this analysis was obtained using a fine \emph{x\textbf{N}x}-discretization -- the base mesh has been refined five times resulting in 20,480 elements -- where spectral elements were coupled to spectral elements. Additionally, the two elements near the singular points have been successively refined eight times ($n_\mathrm{s}\,{=}\,8$, $n_\mathrm{y}\,{=}\,2$). The polynomial degree of the shape functions was set to $p\,{=}\,10$. Thus, the final analysis mesh for the reference solution contains 20,528 elements (see Fig.~\ref{fig:CantileverBeamMeshRef}) with 4,113,992 DOFs.}


\begin{figure}[t!]
	\centering
	\includegraphics[clip,width=0.475\textwidth]{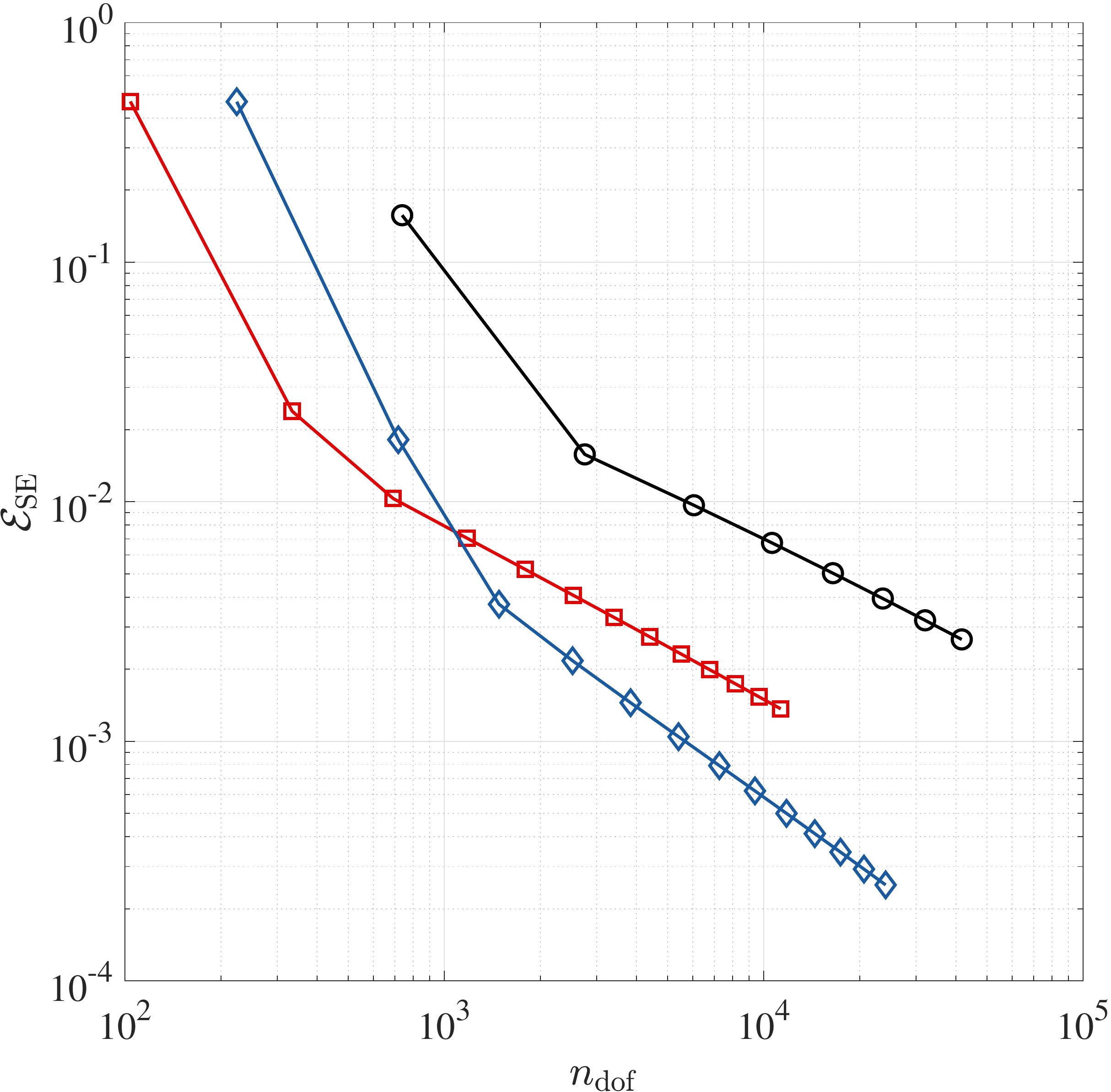}
	\caption{Relative error for the cantilever beam example -- \emph{p}-refinement ($p \in \{1,2,3,\ldots\}$). Legend: \textcolor{Matlab1}{\rule[0.55ex]{5ex}{0.2ex}} Spectral elements (\NewCirc), \textcolor{Matlab2}{\rule[0.55ex]{5ex}{0.2ex}} \emph{x\textbf{N}y}-elements, $n_\mathrm{s}\,{=}\,2$ ($\square$), \textcolor{Matlab3}{\rule[0.55ex]{5ex}{0.2ex}} \emph{x\textbf{N}y}-elements, $n_\mathrm{s}\,{=}\,8$ ($\Diamond$).}
	\label{fig:ErrorCantileverBeam}
\end{figure}

In Fig.~\ref{fig:ErrorCantileverBeam}, the advantages of a local mesh refinement are clearly observed. A uniform \emph{p}-refinement on the spectral element mesh leads to an algebraic convergence which is related to the two singular points. By refining the mesh near the singularities, we achieve an exponential convergence rate in the pre-asymptotic range as the error in the smooth part of the solution is significantly decreased. In the asymptotic range, an algebraic convergence that is limited by the singularity is observed analogous to the SEM curve. Another point to notice is that more layers of refinement near the singular points result in a slightly increased number of DOFs but also in a lower error. By using more hierarchic refinement steps, the pre-asymptotic range of convergence is extended. In Ref. \cite{BookSzabo1991}, it is shown that if a proper \emph{hp}-refinement strategy is employed, we are able to recover exponential rates of convergence during the whole refinement process. To this end, the mesh has to be graded with a ratio of 0.15:0.85 and the polynomial degree must be decreased near the singularity.
\subsection{L-shaped domain}
\label{sec:L}
\begin{figure}[b!]
	\centering
	\includegraphics[clip,scale=0.75]{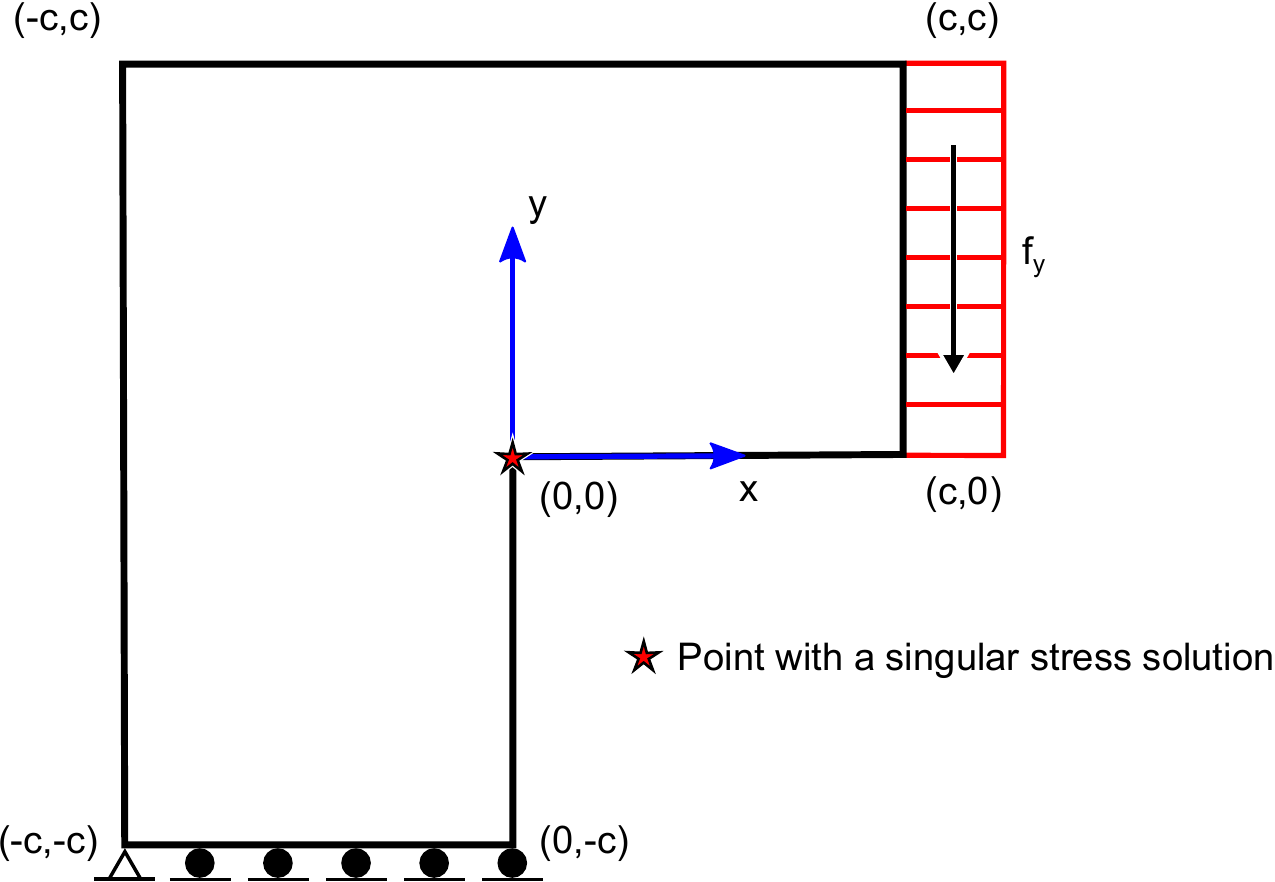}
	\caption{Model of the L-shaped domain; Geometry and boundary conditions.}
	\label{fig:L}
\end{figure}

The L-shaped domain is a classical example to demonstrate the adverse effects of a singular point due to the existence of a re-entrant corner. The geometry including all boundary conditions is depicted in Fig.~\ref{fig:L}. For the numerical example presented in this section the dimension is chosen as: $c\,{=}\,20$. 
\begin{figure}[b!]
	\centering
	\subfloat[\emph{x\textbf{N}y}-mesh, $n_\mathrm{s}\,{=}\,2$, $n_\mathrm{y}\,{=}\,2$, 30 elements \label{fig:LDomainMesh1}]{\includegraphics[clip,width=0.345\textwidth]{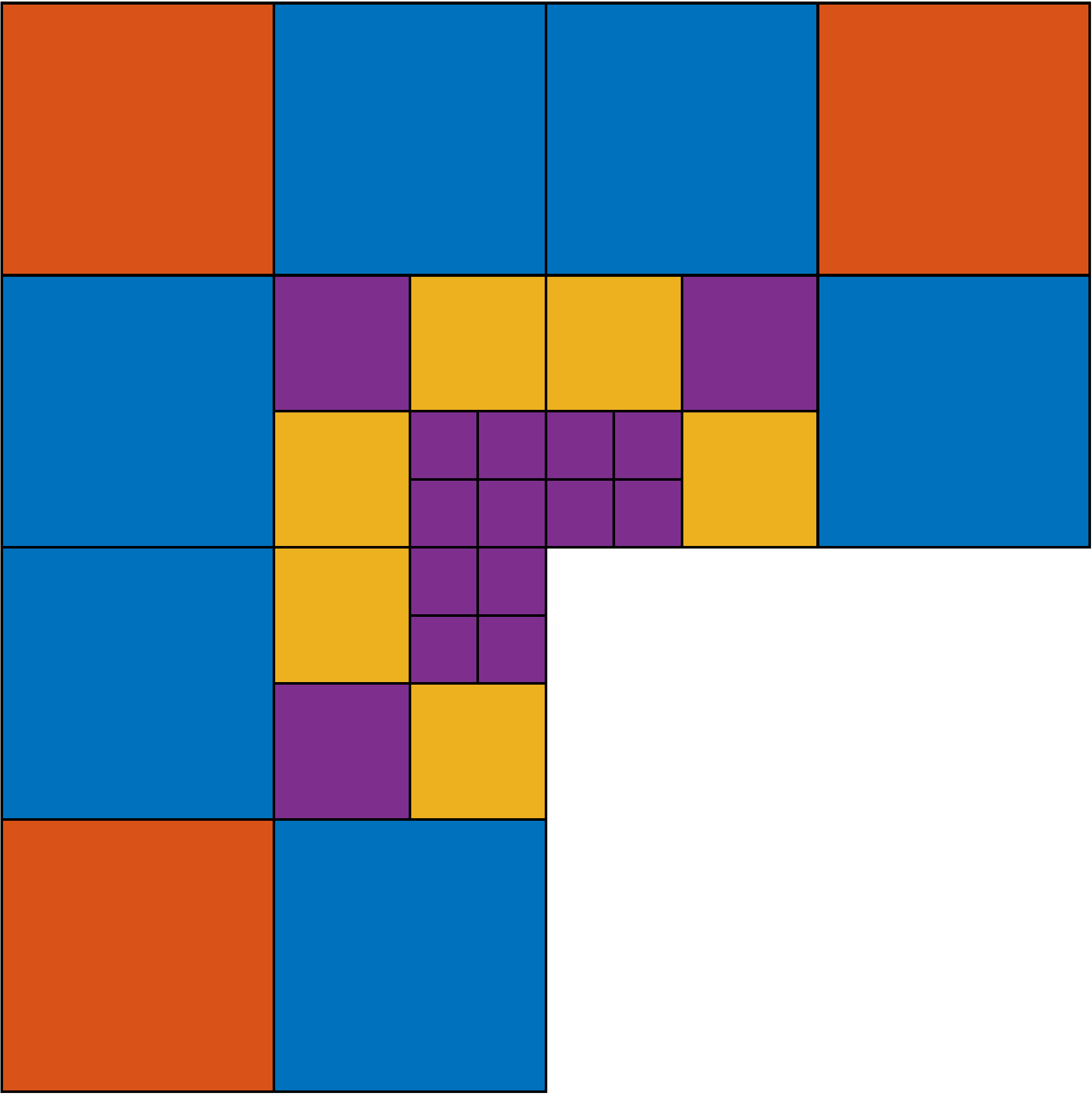}}
	\hspace*{64pt}
	\subfloat[\emph{x\textbf{N}y}-mesh, $n_\mathrm{s}\,{=}\,8$, $n_\mathrm{y}\,{=}\,2$, 84 elements \label{fig:LDomainMesh2}]{\includegraphics[clip,width=0.345\textwidth]{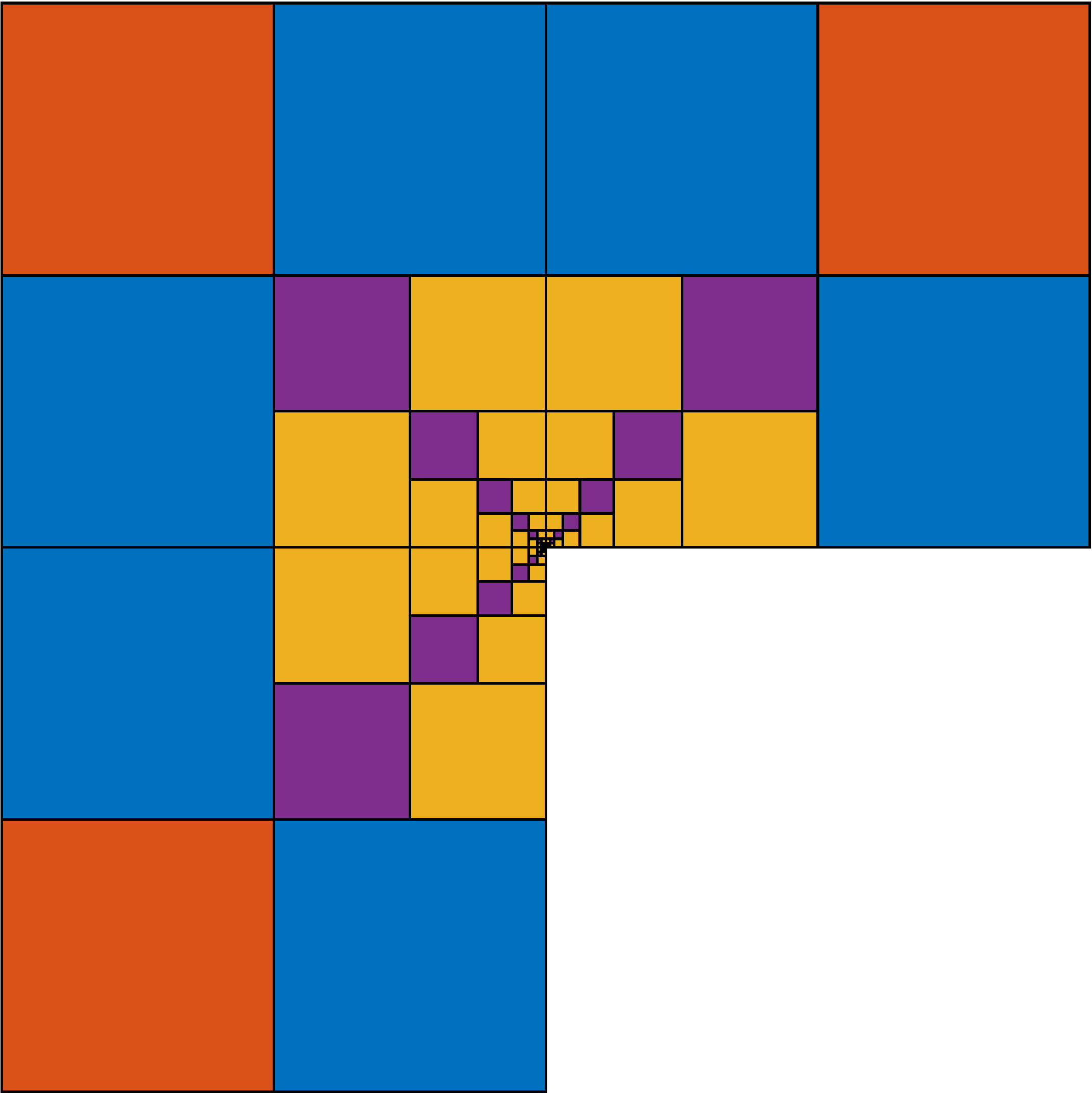}}\\
	\subfloat[Spectral element mesh, 192 elements \label{fig:LDomainMesh3}]{\includegraphics[clip,width=0.345\textwidth]{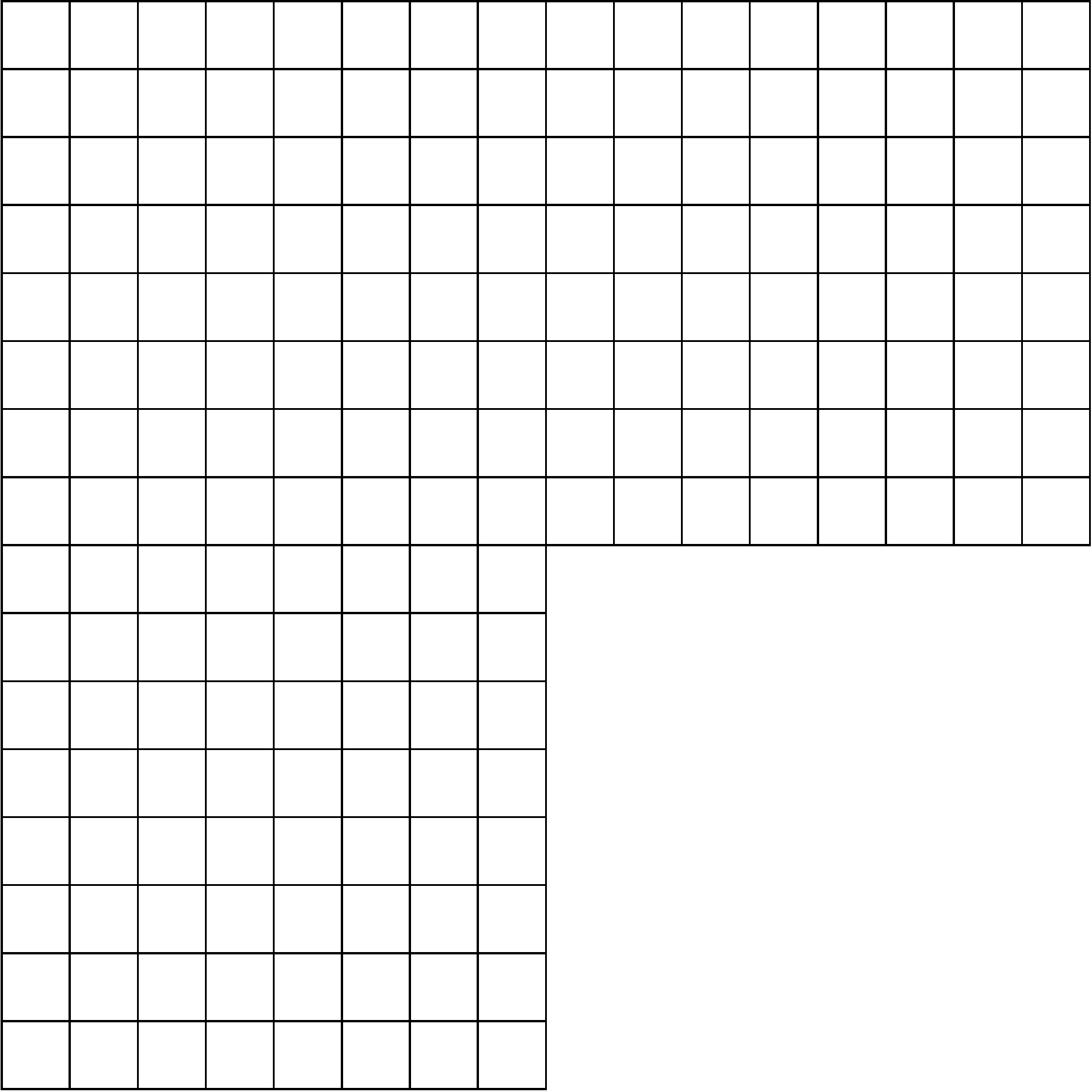}}
	\caption{Discretization of the L-shaped domain model. The different element types are color-coded: \emph{x}-Elements \textcolor{xElem}{\rule[-0.25ex]{5ex}{2ex}}, \emph{y}-Elements \textcolor{yElem}{\rule[-0.25ex]{5ex}{2ex}}, \emph{x\textbf{N}y}-Elements \textcolor{xNyElem}{\rule[-0.25ex]{5ex}{2ex}}, \emph{y\textbf{N}y}-Elements \textcolor{yNyElem}{\rule[-0.25ex]{5ex}{2ex}}.}
	\label{fig:LDomainMesh}
\end{figure}%

\textcolor{red}{The discretization of this structure is set up in a similar fashion as the one reported in the previous example. The initial mesh consists of 12 square-shaped quadrilateral finite elements and finer discretizations are generated by dividing each element into four. Also in this example, the results obtained based on a simple \emph{p}-refinement in combination with the SEM are compared with spectral--spectral transition elements. The different discretizations that are the basis for the \emph{p}-refinement are depicted in Fig.~\ref{fig:LDomainMesh}, where 2 and 8 mesh refinement levels ($n_\mathrm{s}$) have been employed for the transition elements and the SEM grid has been generated by refining the initial mesh twice.} 
\begin{figure}[t!]
	\centering
	\includegraphics[clip,width=0.475\textwidth]{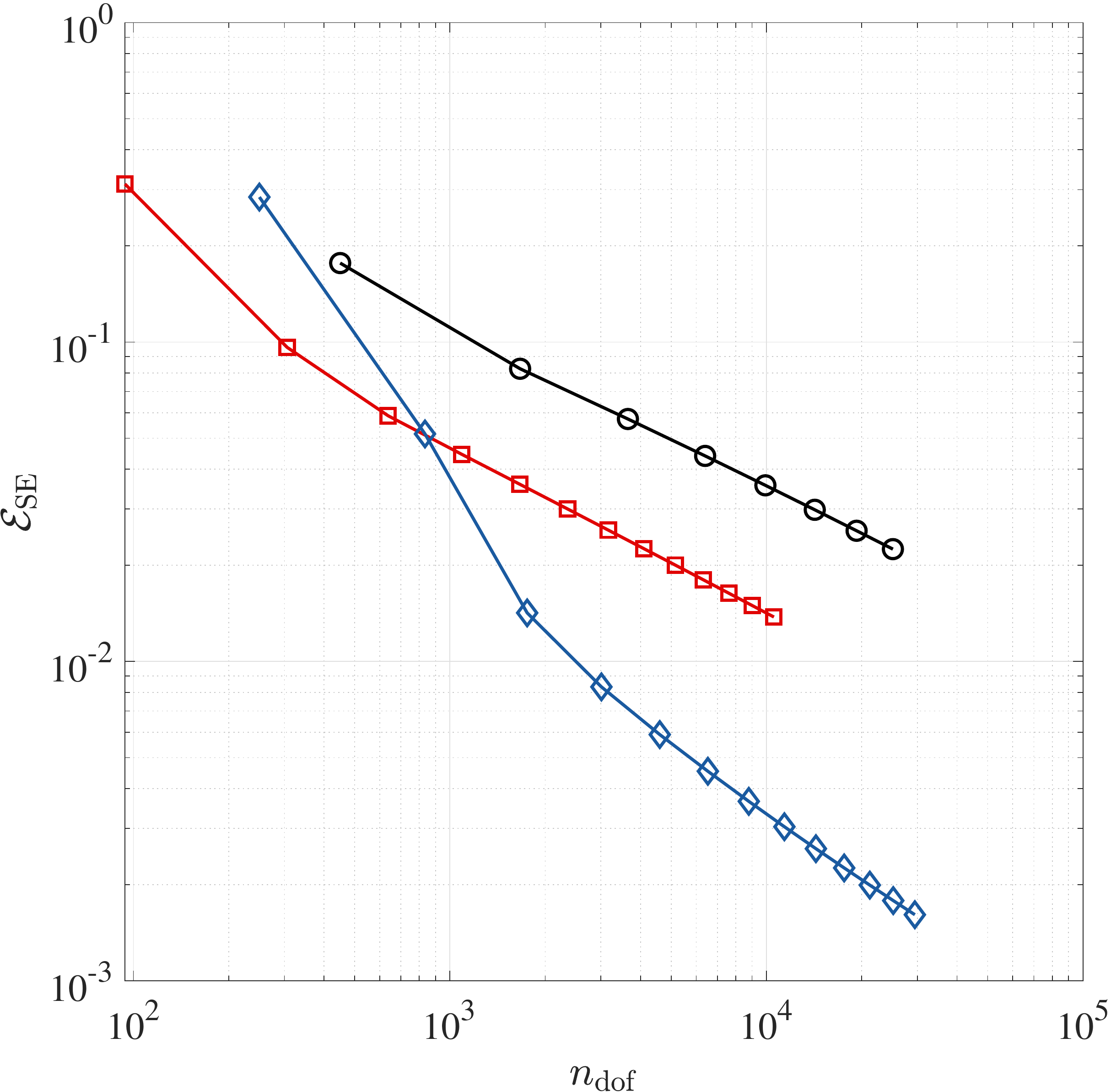}
	\caption{Relative error for the L-shaped domain example -- \emph{p}-refinement ($p \in \{1,2,3,\ldots\}$). Legend: \textcolor{Matlab1}{\rule[0.55ex]{5ex}{0.2ex}} Spectral elements (\NewCirc), \textcolor{Matlab2}{\rule[0.55ex]{5ex}{0.2ex}} \emph{x\textbf{N}y}-elements, $n_\mathrm{s}\,{=}\,2$ ($\square$), \textcolor{Matlab3}{\rule[0.55ex]{5ex}{0.2ex}} \emph{x\textbf{N}y}-elements, $n_\mathrm{s}\,{=}\,8$ ($\Diamond$).}
	\label{fig:ErrorLDomain}
\end{figure}%
\begin{figure}[b!]
	\centering
	\subfloat[{Detail view of the mesh at the re-entrant corner ($x\in[5,15] $, $y\in[5,15]$)} \label{fig:LDomainMeshRef1}]{\includegraphics[clip,width=0.46\textwidth]{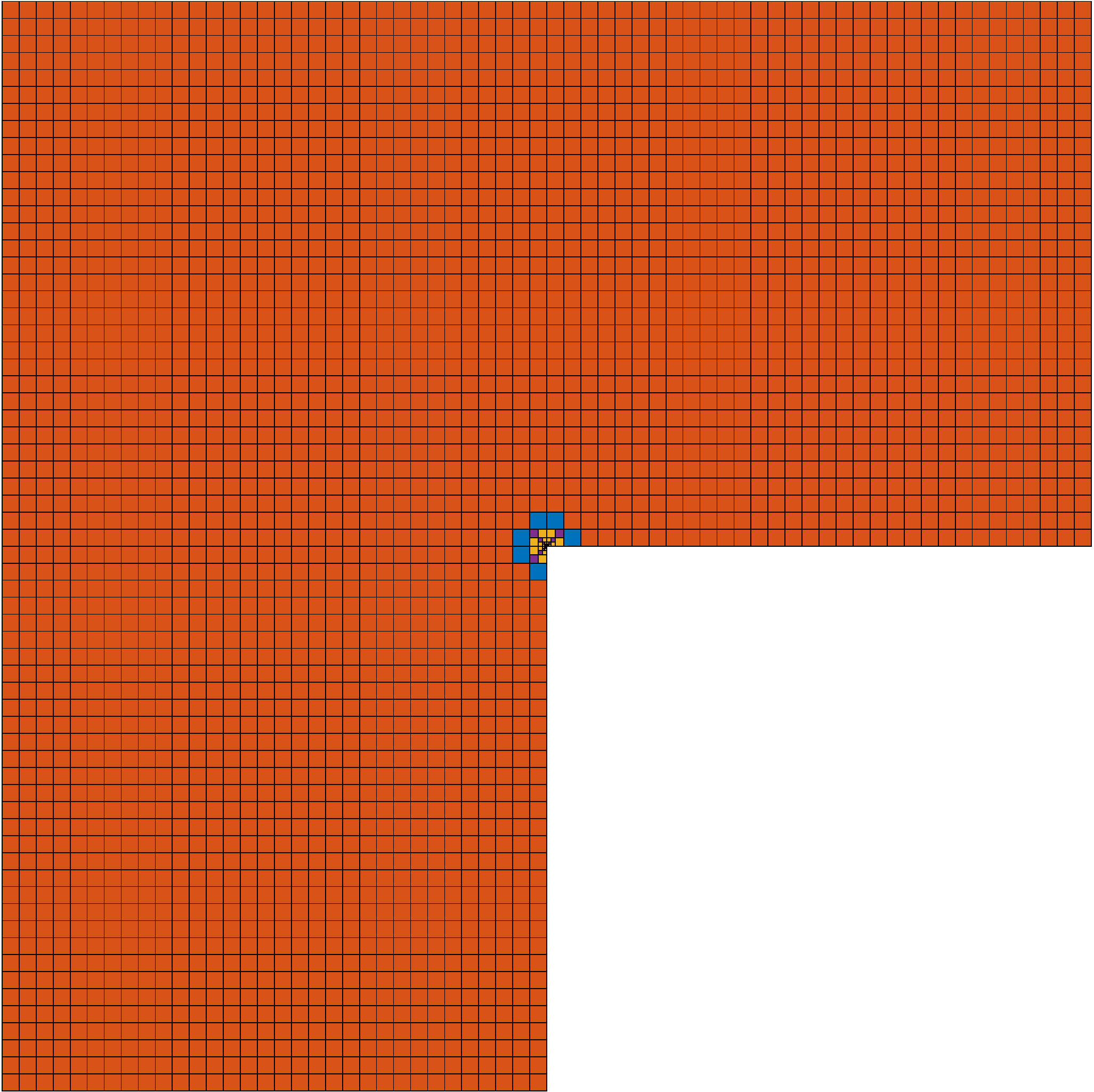}}
	\hfill
	\subfloat[{Detail view of the mesh at the re-entrant corner ($x\in[9.5,10.5] $, $y\in[9.5,10.5]$)} \label{fig:LDomainMeshRef2}]{\includegraphics[clip,width=0.46\textwidth]{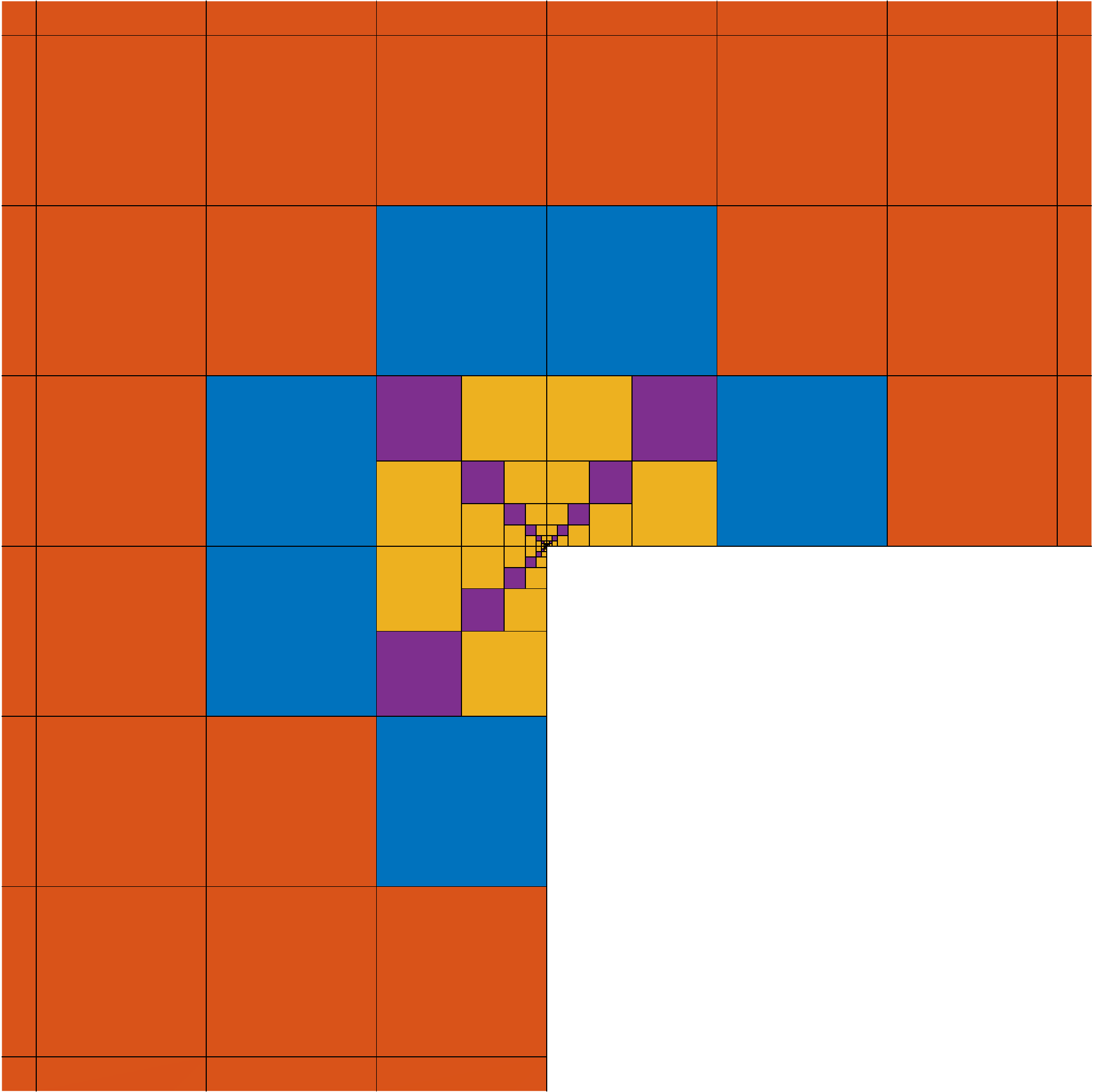}}
	\caption{Discretization of the L-shaped domain model (overall mesh: 12,300 elements) - reference model (\emph{x\textbf{N}y}-mesh, $n_\mathrm{s}\,{=}\,8$, $n_\mathrm{y}\,{=}\,2$). The different element types are color-coded: \emph{x}-Elements \textcolor{xElem}{\rule[-0.25ex]{5ex}{2ex}}, \emph{y}-Elements \textcolor{yElem}{\rule[-0.25ex]{5ex}{2ex}}, \emph{x\textbf{N}y}-Elements \textcolor{xNyElem}{\rule[-0.25ex]{5ex}{2ex}}, \emph{y\textbf{N}y}-Elements \textcolor{yNyElem}{\rule[-0.25ex]{5ex}{2ex}}.}
	\label{fig:LDomainMeshRef}
\end{figure}%

The numerical results show a similar behavior to that discussed in Sect.~\ref{sec:Cantilever}. Again, the advantages of a local mesh refinement near the re-entrant corner are observed (see Fig.~\ref{fig:ErrorLDomain}). Using the same number of DOFs, the error in the energy norm can be decreased by up to two orders of magnitude. Note that the reference solution for this analysis was obtained using a fine \emph{x\textbf{N}x}-discretization where spectral elements were coupled with spectral elements. \textcolor{red}{The fine mesh was generated from the initial one by refining it five times resulting in 12,288 finite elements.} Additionally, the two elements near the singular points have been refined eight times successively ($n_\mathrm{s}\,{=}\,8$, $n_\mathrm{y}\,{=}\,2$). The polynomial degree of the shape functions was set to $p\,{=}\,10$. Thus, the final analysis mesh for the reference solution contains 12,360 elements (see Fig.~\ref{fig:LDomainMeshRef}) with 2,477,762 DOFs.

Although the obtained results are very promising, they can be significantly improved by implementing an automatic mesh generation capable of providing a geometrically graded mesh. A further improvement can also be obtained by means of an adaptive refinement scheme based on a suitable error estimator \cite{Song_2018}. Additionally, the polynomial degree of the shape functions near the singularity should be decreased with each hierarchical sub-division of the elements. Thus, fewer DOFs are ``wasted'' near the re-entrant corner where the solution is inaccurate due to the assumption of the numerical model.

%% file: tex/Summary.tex
\section{Summary}
\label{sec:Summary}
In this contribution, we extended the concept of transfinite mapping and transition elements to achieve a compatible coupling strategy between various element types and elements of different sizes and polynomial orders. This procedure is especially useful if a local mesh refinement is required in order to adapt the discretization to the requirements of the solution. The obtained numerical results clearly demonstrate that the theoretical optimal rates of convergence can be achieved, i.e., the rate of convergence is proportional to the minimum polynomial degree $p_\mathrm{min}$ of the \emph{x}- or \emph{y}-shape functions. Therefore, the \emph{x\textbf{N}y}-element concept can be employed to locally refine the discretization where needed without loss of accuracy, which is observed when using multi-point constraints or other coupling methods.
\\\\
At the current stage, there is still the problem of an efficient mesh generation for these transition elements. The discretizations in this article have been generated by using Abaqus as a pre-processor. To this end, a coarse base mesh is generated and imported to Matlab, where specified elements are successively refined according to the selected values of $n_\mathrm{s}$ and $n_\mathrm{y}$. Another idea would be to create two meshes with different element sizes and to merge them in regions where a mesh refinement is needed. In this case, it is also straightforward to identify the \emph{x}- and \emph{y}-elements. In a further step, \emph{x}-elements that are adjacent to \emph{y}-elements are changed to \emph{x\textbf{N}y}-elements. Especially for complex geometries and three-dimensional models, this procedure is very cumbersome and not feasible. Consequently, powerful mesh generators that are capable of accounting for transition elements are one prerequisite for the successful use of this concept.
\\\\
Future research activities are directed at investigating the performance of this concept in dynamic analyses. Here, it is of great interest to see if these elements can also be used in conjunction with explicit time integration approaches. To this end, it has to be possible to diagonalize the mass matrix \cite{ArticleDuczek2019a, ArticleDuczek2019b}. A natural area of application of the proposed method would be the analysis of wave propagation phenomena in layered media where the local mesh refinement capability could be exploited.
\\\\
It is, furthermore, intended to extend the \emph{x\textbf{N}y}-approach to fictitious domain methods such as the finite cell method (FCM) \cite{ArticleParvizian2007, ArticleDuester2008}. Although the discretization is not geometry-conforming, a certain number of elements are still required to capture the complex displacement/stress field near geometrically complex details of the microstructure \cite{ChapterWuerkner2018, ArticleJomo2019}.

%% file: tex/Appendix.tex
% Appendix
\appendix
\setcounter{figure}{0}
\setcounter{table}{0}
\section{Shape functions of a 12-node/mode transition element}
\label{sec:12NodeElem}
In Sect.~\ref{sec:HO_ShapeFunctions} a detailed description of the derivation of high order transition elements is given. These elements are either based on nodal or modal shape functions and have the ability to conformally couple different finite elements at specific edges. In this section, it is shown how transition element shape functions are derived for particular cases. As an illustrative example, we derive the shape functions for a piecewise bi-quadratic transition element with 12 nodes/modes (see Fig.~\ref{fig:12NodeElem}). Since this element does not feature interior nodes/modes, it can be seen as a Serendipity-type/trunk space transition element.
\subsection{One-dimensional standard shape functions}
The basis for the derivation of the shape functions for the 12-node/mode transition element are the one-dimensional quadratic Lagrange shape functions \cite{BookPozrikidis2014}:
\allowdisplaybreaks
\begin{subequations}
	\label{eq:1DQuadShapeFunc_La}
	\begin{align}
	^\mathrm{La}N^\mathrm{2}_{1}(\xi)  &= \cfrac{1}{2}(\xi^2-\xi)\,, \label{eq:1DQuadShapeFuncNode1_La}\\
	^\mathrm{La}N^\mathrm{2}_{2}(\xi)  &= (1-\xi^2)\,, \label{eq:1DQuadShapeFuncNode2_La}\\
	^\mathrm{La}N^\mathrm{2}_{3}(\xi)  &= \cfrac{1}{2}(\xi^2+\xi)\,, \label{eq:1DQuadShapeFuncNode3_La}
	\end{align}
\end{subequations}
and the one-dimensional hierarchic shape functions (based on the normalized integrals of the Legendre polynomials) of order 2 (\emph{p}-FEM) \cite{BookSzabo1991}:
\allowdisplaybreaks
\begin{subequations}
	\label{eq:1DQuadShapeFunc_Le}
	\begin{align}
	^\mathrm{Le}N^\mathrm{2}_{1}(\xi)  &= \cfrac{1}{2}(1-\xi)\,, \label{eq:1DQuadShapeFuncNode1_Le}\\
	^\mathrm{Le}N^\mathrm{2}_{2}(\xi)  &= \cfrac{1}{4}\sqrt{6}(\xi^2-1)\,, \label{eq:1DQuadShapeFuncNode2_Le}\\
	^\mathrm{Le}N^\mathrm{2}_{3}(\xi)  &= \cfrac{1}{2}(1+\xi)\,. \label{eq:1DQuadShapeFuncNode3_Le}
	\end{align}
\end{subequations}
These are the one-dimensional standard shape functions that describe the variation of the multi-dimensional ones at the edges of an element and are, therefore, of utmost importance in the derivation of transition elements.
\subsection{Two-dimensional standard shape functions}
\begin{figure}[t!]
	\begin{minipage}[t]{0.45\textwidth}
		\centering
		\includegraphics[clip,width=0.75\textwidth]{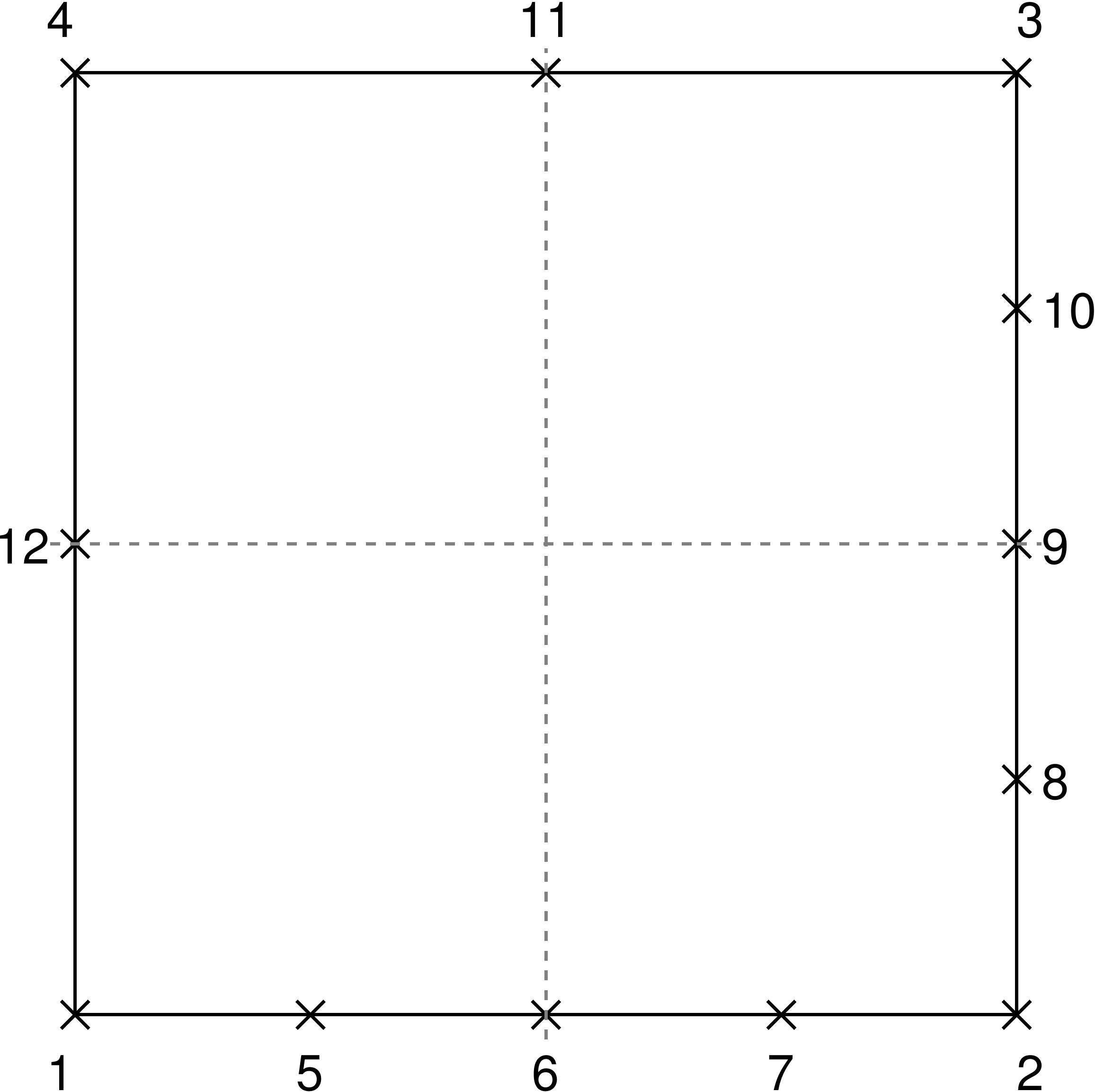}
		\caption{12-node bi-quadratic transition finite element (nodal numbering).}
		\label{fig:12NodeElem}
	\end{minipage}
	\hfill
	\begin{minipage}[t]{0.45\textwidth}
		\centering
		\includegraphics[clip,width=0.733\textwidth]{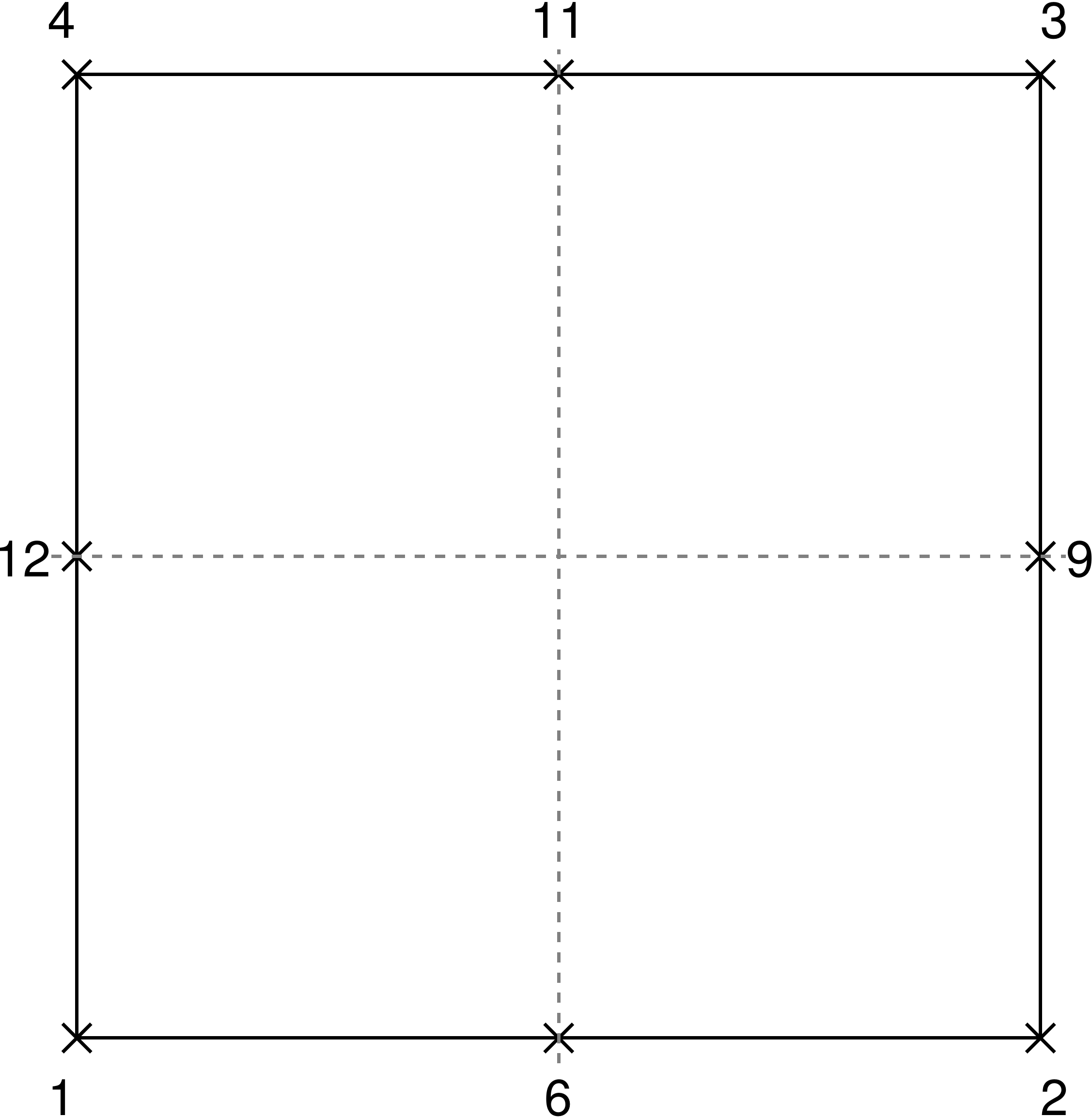}
		\caption{8-node bi-quadratic finite element (nodal numbering).}
		\label{fig:8NodeElem}
	\end{minipage}
\end{figure}

In addition to the one-dimensional shape functions defined above, we also need the shape functions of the well-known 8-node quadrilateral element (see Fig.~\ref{fig:8NodeElem}) \cite{BookZienkiewicz2000a}:
\allowdisplaybreaks
\begin{subequations}
	\label{eq:QuadShapeFunc_Serendipity}
	\begin{align}
	^\mathrm{s}N^{*}_{1}(\xi,\eta)  &= -\cfrac{1}{4}(1-\xi)(1-\eta)(1+\xi+\eta)\,, \label{eq:QuadShapeFuncNode1}\\
	^\mathrm{s}N^{*}_{2}(\xi,\eta)  &= -\cfrac{1}{4}(1+\xi)(1-\eta)(1-\xi+\eta)\,, \label{eq:QuadShapeFuncNode2}\\
	^\mathrm{s}N^{*}_{3}(\xi,\eta)  &= -\cfrac{1}{4}(1+\xi)(1+\eta)(1-\xi-\eta)\,, \label{eq:QuadShapeFuncNode3}\\
	^\mathrm{s}N^{*}_{4}(\xi,\eta)  &= -\cfrac{1}{4}(1-\xi)(1+\eta)(1+\xi-\eta)\,, \label{eq:QuadShapeFuncNode4}\\
	^\mathrm{s}N^{*}_{6}(\xi,\eta)  &=  \cfrac{1}{2}(1-\xi^2)(1-\eta)\,, \label{eq:QuadShapeFuncNode5}\\
	^\mathrm{s}N^{*}_{9}(\xi,\eta)  &=  \cfrac{1}{2}(1+\xi)(1-\eta^2)\,, \label{eq:QuadShapeFuncNode6}\\
	^\mathrm{s}N^{*}_{11}(\xi,\eta) &=  \cfrac{1}{2}(1-\xi^2)(1+\eta)\,, \label{eq:QuadShapeFuncNode7}\\
	^\mathrm{s}N^{*}_{12}(\xi,\eta) &=  \cfrac{1}{2}(1-\xi)(1-\eta^2)\,, \label{eq:QuadShapeFuncNode8}
	\end{align}
\end{subequations}
and the shape functions of the quadrilateral \emph{p}-element which is based on the trunk space formulation (similar to the Serendipity ansatz space) \cite{BookSzabo1991}:
\allowdisplaybreaks
\begin{subequations}
	\label{eq:QuadShapeFunc_Trunk}
	\begin{align}
	^\mathrm{t}N^{*}_{1}(\xi,\eta)  &= \,N_{1}^{1}(\xi) N_{1}^{1}(\eta)\,, \label{eq:QuadShapeFuncNode1_trunk}\\
	^\mathrm{t}N^{*}_{2}(\xi,\eta)  &= \,N_{2}^{1}(\xi) N_{1}^{1}(\eta)\,, \label{eq:QuadShapeFuncNode2_trunk}\\
	^\mathrm{t}N^{*}_{3}(\xi,\eta)  &= \,N_{2}^{1}(\xi) N_{2}^{1}(\eta)\,, \label{eq:QuadShapeFuncNode3_trunk}\\
	^\mathrm{t}N^{*}_{4}(\xi,\eta)  &= \,N_{1}^{1}(\xi) N_{2}^{1}(\eta)\,, \label{eq:QuadShapeFuncNode4_trunk}\\
	^\mathrm{t}N^{*}_{6}(\xi,\eta)  &= \,^\mathrm{Le}N^\mathrm{2}_{2}(\xi)  N_{1}^{1}(\eta)\,, \label{eq:QuadShapeFuncNode5_trunk}\\
	^\mathrm{t}N^{*}_{9}(\xi,\eta)  &= \,^\mathrm{Le}N^\mathrm{2}_{2}(\eta) N_{2}^{1}(\xi)\,, \label{eq:QuadShapeFuncNode6_trunk}\\
	^\mathrm{t}N^{*}_{11}(\xi,\eta) &= \,^\mathrm{Le}N^\mathrm{2}_{2}(\xi)  N_{2}^{1}(\eta)\,, \label{eq:QuadShapeFuncNode7_trunk}\\
	^\mathrm{t}N^{*}_{12}(\xi,\eta) &= \,^\mathrm{Le}N^\mathrm{2}_{2}(\eta) N_{1}^{1}(\xi)\,. \label{eq:QuadShapeFuncNode8_trunk}
	\end{align}
\end{subequations}
Studying the composition of the shape functions of an 8-node serendipity element corresponding to the edge mid-side nodes (6, 9, 11, 12) given by Eq.~\eqref{eq:QuadShapeFunc_Serendipity}, we observe that these are derived using the linear blending discussed in Sect.~\ref{sec:TransitionElem}. The one-dimensional shape function for a mid-side node is just multiplied by a linear blending term
\begin{subequations}
	\begin{align}
	N_1^1(\xi) & = \cfrac{1}{2}(1-\xi)\,,\\
	N_2^1(\xi) & = \cfrac{1}{2}(1+\xi)\,.
	\end{align}
\end{subequations}
Note that, in the following, the numbering scheme indicated in Figs.~\ref{fig:12NodeElem} and \ref{fig:8NodeElem} is assumed.
\subsection{Lagrange-Lagrange: Nodal shape functions}
\label{app:LaLa}
The first example to demonstrate the derivation of transition element shape functions is a bi-quadratic 12-node transition element. Here, only Lagrange-type nodal shape functions are employed. The basis functions for the derivation are given by Eqs.~\eqref{eq:1DQuadShapeFunc_La} and \eqref{eq:QuadShapeFunc_Serendipity}, accordingly.
\subsubsection{Shape functions on the boundary}
\begin{figure}[t!]
	\begin{minipage}[t]{0.49\textwidth}
		\centering
		\includegraphics[clip,width=1\textwidth]{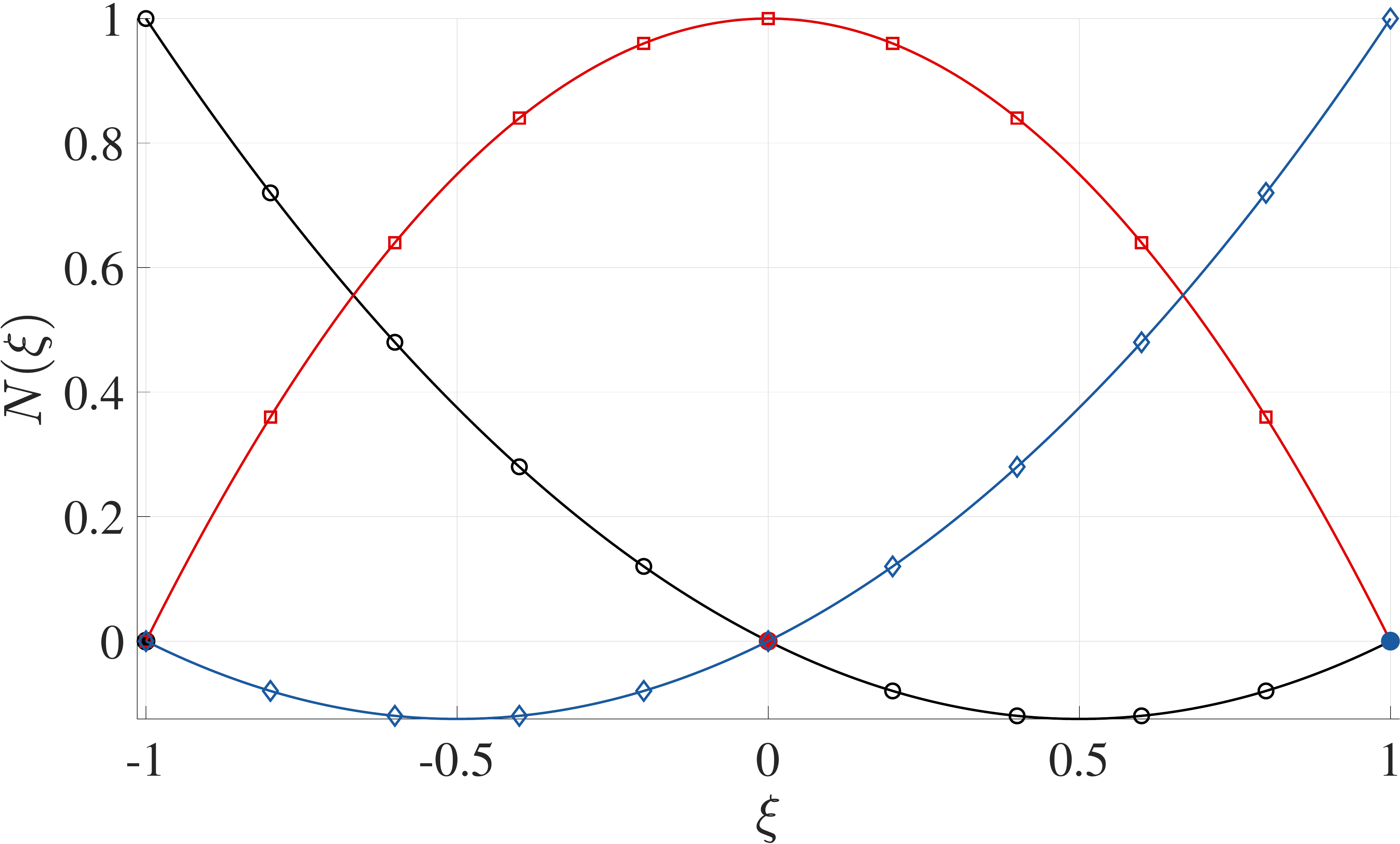}
		\caption{One-dimensional quadratic shape functions (3-node Lagrange element). Legend: \textcolor{Matlab1}{\rule[0.55ex]{5ex}{0.2ex}} $N_1(\xi)$ (\NewCirc), \textcolor{Matlab2}{\rule[0.55ex]{5ex}{0.2ex}} $N_2(\xi)$ ($\square$), \textcolor{Matlab3}{\rule[0.55ex]{5ex}{0.2ex}} $N_3(\xi)$ ($\Diamond$).}
		\label{fig:ShapeFuncQuad1d}
	\end{minipage}
	\hfill
	\begin{minipage}[t]{0.49\textwidth}
		\centering
		\includegraphics[clip,width=1\textwidth]{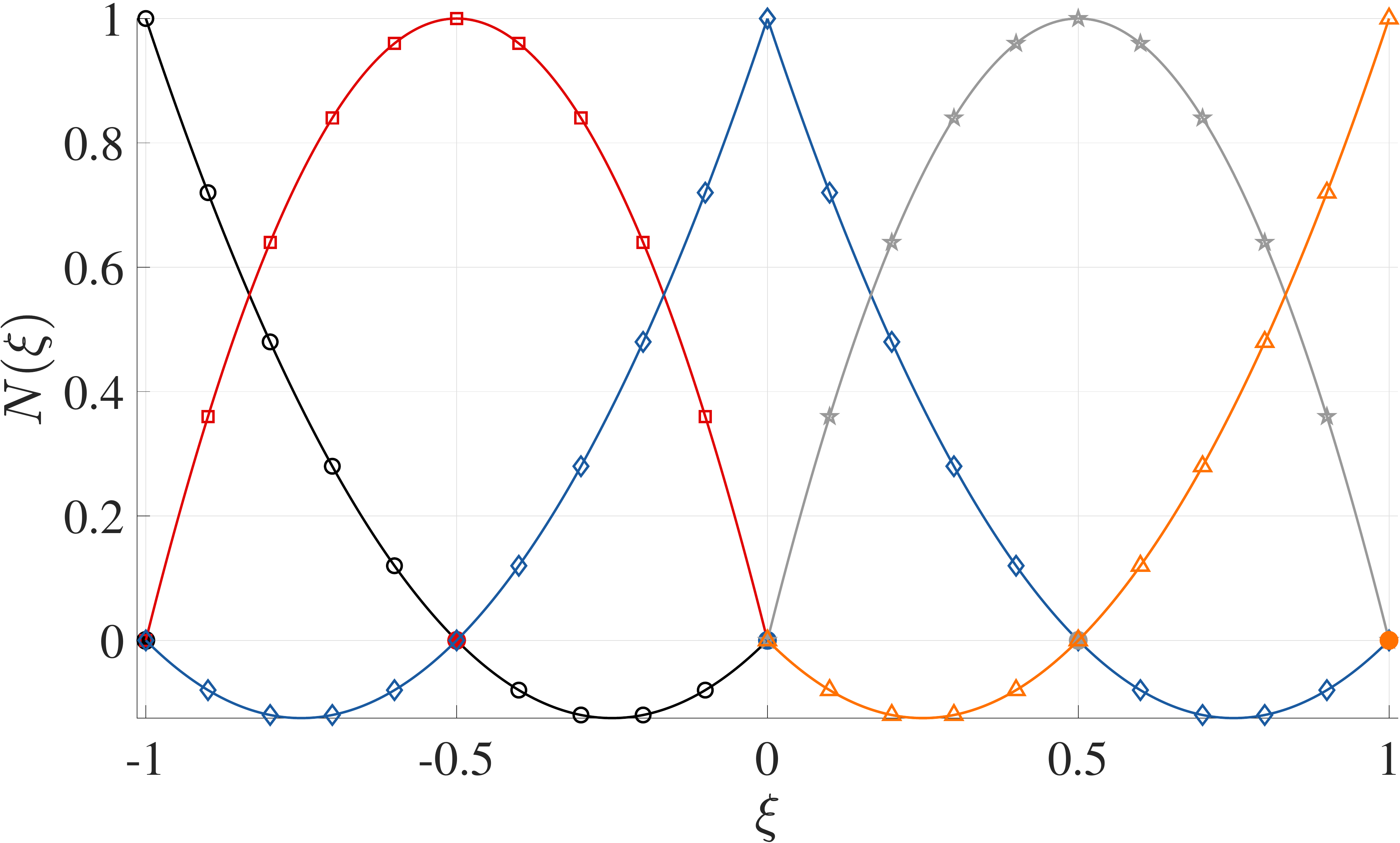}
		\caption{One-dimensional piecewise quadratic shape functions (based on 3-node Lagrange elements). Legend: \textcolor{Matlab1}{\rule[0.55ex]{5ex}{0.2ex}} $N_1(\xi)$ (\NewCirc), \textcolor{Matlab2}{\rule[0.55ex]{5ex}{0.2ex}} $N_2(\xi)$ ($\square$), \textcolor{Matlab3}{\rule[0.55ex]{5ex}{0.2ex}} $N_3(\xi)$ ($\Diamond$), \textcolor{Matlab4}{\rule[0.55ex]{5ex}{0.2ex}} $N_4(\xi)$ (\NewStar), \textcolor{Matlab5}{\rule[0.55ex]{5ex}{0.2ex}} $N_5(\xi)$ ($\bigtriangleup$).}
		\label{fig:ShapeFuncPiecewiseQuad1d}
	\end{minipage}
\end{figure}

In the first step, the function $\bar{u}(\xi,\eta)$, describing the displacement field on the boundary of the element, needs to be constructed -- see Eq.~\eqref{eq:EdgeFunctions}. The shape functions of edges $E_1$ and $E_2$ (see the definition of the reference element provided in Fig.~\ref{fig:RefElem}) are generated as piecewise quadratic polynomials, while $E_3$ and $E_4$ retain the standard quadratic shape functions. By means of this example, we will show that merely shape functions on the ``coupling'' edges need to be adjusted while the other shape functions can be adopted from the basis element. Therefore, only a limited number of shape functions are modified. Accordingly, the interpolations along the four element edges are given as
\allowdisplaybreaks
\begin{align}
\text{E}_1:\; \bar{N}_i(\xi) &=
\begin{cases}
^\mathrm{La}N^\mathrm{2}_{1}(\check{\xi}_1) & \text{for } i = 1, \text{ node \#1} \,,\\
^\mathrm{La}N^\mathrm{2}_{2}(\check{\xi}_1) & \text{for } i = 2, \text{ node \#5} \,,\\
\begin{cases}
^\mathrm{La}N^\mathrm{2}_{3}(\check{\xi}_1) \\
^\mathrm{La}N^\mathrm{2}_{1}(\check{\xi}_2) \\
\end{cases}
& \text{for } i = 3, \text{ node \#6} \,, \\
^\mathrm{La}N^\mathrm{2}_{2}(\check{\xi}_2) & \text{for } i = 4, \text{ node \#7} \,, \\
^\mathrm{La}N^\mathrm{2}_{3}(\check{\xi}_2) & \text{for } i = 5, \text{ node \#2} \,,
\end{cases}
\\
\text{E}_2:\; \hat{N}_i(\eta) &=
\begin{cases}
^\mathrm{La}N^\mathrm{2}_{1}(\check{\eta}_1) & \text{for } i = 1, \text{ node \#2} \,,\\
^\mathrm{La}N^\mathrm{2}_{2}(\check{\eta}_1) & \text{for } i = 2, \text{ node \#8} \,,\\
\begin{cases}
^\mathrm{La}N^\mathrm{2}_{3}(\check{\eta}_1) \\
^\mathrm{La}N^\mathrm{2}_{1}(\check{\eta}_2) \\
\end{cases}
& \text{for } i = 3, \text{ node \#9} \,,\\
^\mathrm{La}N^\mathrm{2}_{2}(\check{\eta}_2) & \text{for } i = 4, \text{ node \#10} \,, \\
^\mathrm{La}N^\mathrm{2}_{3}(\check{\eta}_2) & \text{for } i = 5, \text{ node \#3} \,,
\end{cases}
\\
\text{E}_3:\; \tilde{N}_i(\xi) &=
\begin{cases}
^\mathrm{La}N^\mathrm{2}_{1}(\xi) & \text{for } i = 1, \text{ node \#4} \,,\\
^\mathrm{La}N^\mathrm{2}_{2}(\xi) & \text{for } i = 2, \text{ node \#11} \,,\\
^\mathrm{La}N^\mathrm{2}_{3}(\xi) & \text{for } i = 3, \text{ node \#3} \,,
\end{cases}
\\
\text{E}_4:\; \breve{N}_i(\eta) &=
\begin{cases}
^\mathrm{La}N^\mathrm{2}_{1}(\eta) & \text{for } i = 1, \text{ node \#1} \,,\\
^\mathrm{La}N^\mathrm{2}_{2}(\eta) & \text{for } i = 2, \text{ node \#12} \,,\\
^\mathrm{La}N^\mathrm{2}_{3}(\eta) & \text{for } i = 3, \text{ node \#4} \,.
\end{cases}
\end{align}
Considering the expressions to compute the shape functions on the divided edges $E_1$ and $E_2$, we introduced new coordinates $\check{\xi}_1$, $\check{\xi}_2$, $\check{\eta}_1$, and $\check{\eta}_2$. These coordinates are defined by a simple linear transformation from $[-1,0]$ or $[0,1]$ to the reference domain $[-1,1]$. Consequently, $\check{\xi}_1$ and $\check{\xi}_2$ are given as
\begin{align}
\check{\xi}_1 &= 2\xi + 1\,, \quad \xi \in [-1,0]\,,\\
\check{\xi}_2 &= 2\xi - 1\,, \quad \xi \in [0,1]\,.
\end{align}
The same formulae also apply to $\check{\eta}_1$ and $\check{\eta}_2$. In cases where the edge is divided into more elements, the following general mapping functions can be used
\begin{equation}
\check{\xi}_k = \cfrac{2\xi - \xi_k - \xi_{k+1}}{\xi_{k+1} - \xi_k}\,, \quad \xi \in [\xi_k,\xi_{k+1}]\,.
\end{equation}
This expression can be straightforwardly derived by inverting the one-dimensional linear mapping using standard finite element shape functions
\begin{equation}
\xi =  N_1^1(\check{\xi}_k) \xi_k + N_2^1(\check{\xi}_k) \xi_{k+1}\,.
\end{equation}
Note that the shape functions corresponding to the edge mid-side nodes \#6 and \#9 are piecewise quadratic polynomials. Consequently, a weak discontinuity (kink) in the displacement field is observed. This, however, results in a jump in the first derivative which is directly related to the computed stress field. The shape functions for the edges are exemplarily sketched in Figs.~\ref{fig:ShapeFuncQuad1d} and \ref{fig:ShapeFuncPiecewiseQuad1d}.
\subsubsection{Application of the projection operators}
After having discussed the functions on the edges of the transition element, we can project these functions into the interior of the element using the projection operators $\mathcal{P}_{\xi}[\square]$, $\mathcal{P}_{\eta}[\square]$, and $\mathcal{P}_{\xi}[\mathcal{P}_{\eta}[\square]]$. The general operators are simplified significantly when linear blending functions are deployed and therefore, Eqs.~\eqref{eq:SF_CornerNode1} to \eqref{eq:SF_Edge4} are used in this section.

\begin{figure}[t!]
	\centering
	\subfloat[Corner vertex \#1\label{fig:12NodeElem_SF1}]{\includegraphics[clip,width=0.24\textwidth]{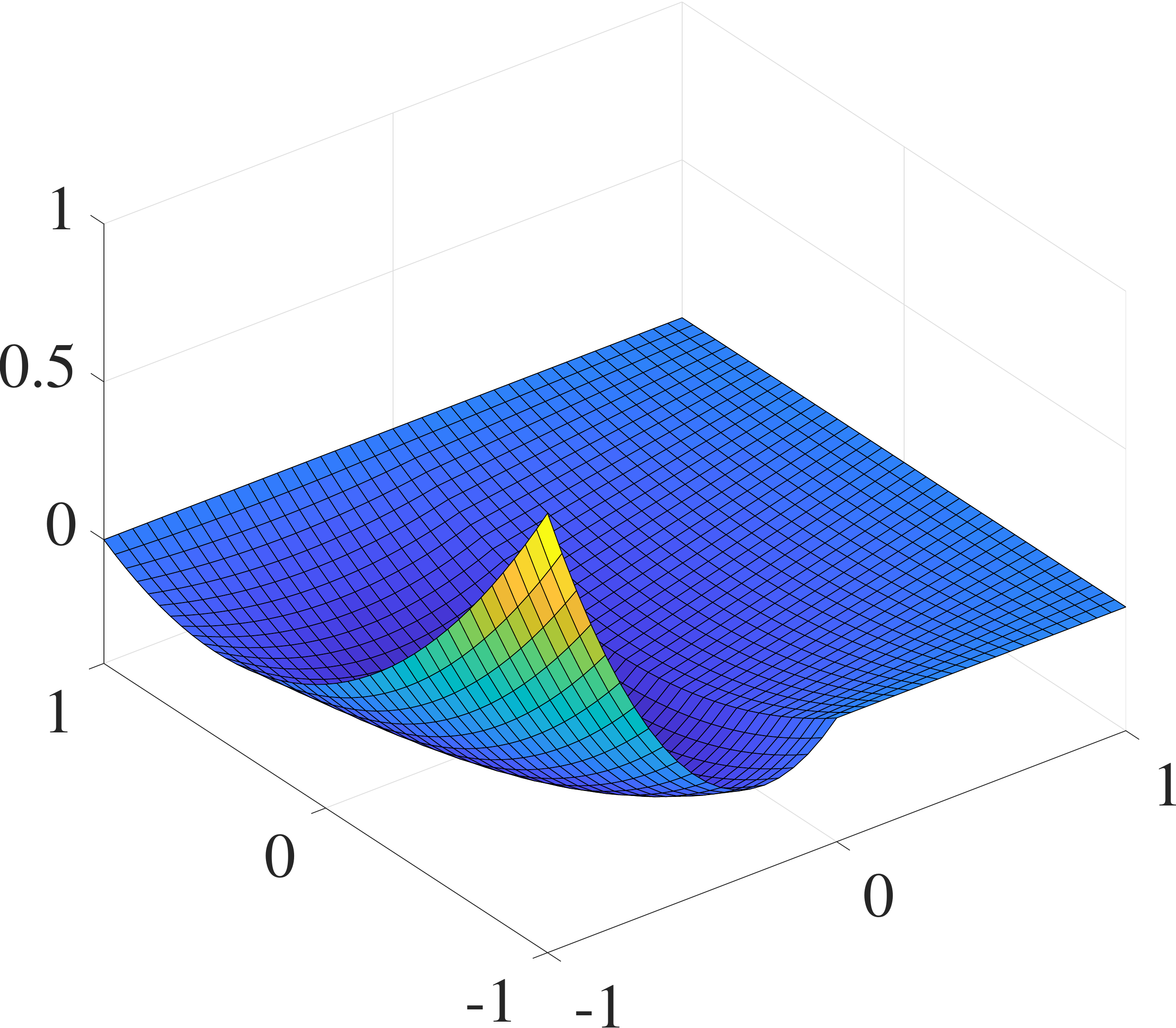}}
	\hfill
	\subfloat[Corner vertex \#2\label{fig:12NodeElem_SF2}]{\includegraphics[clip,width=0.24\textwidth]{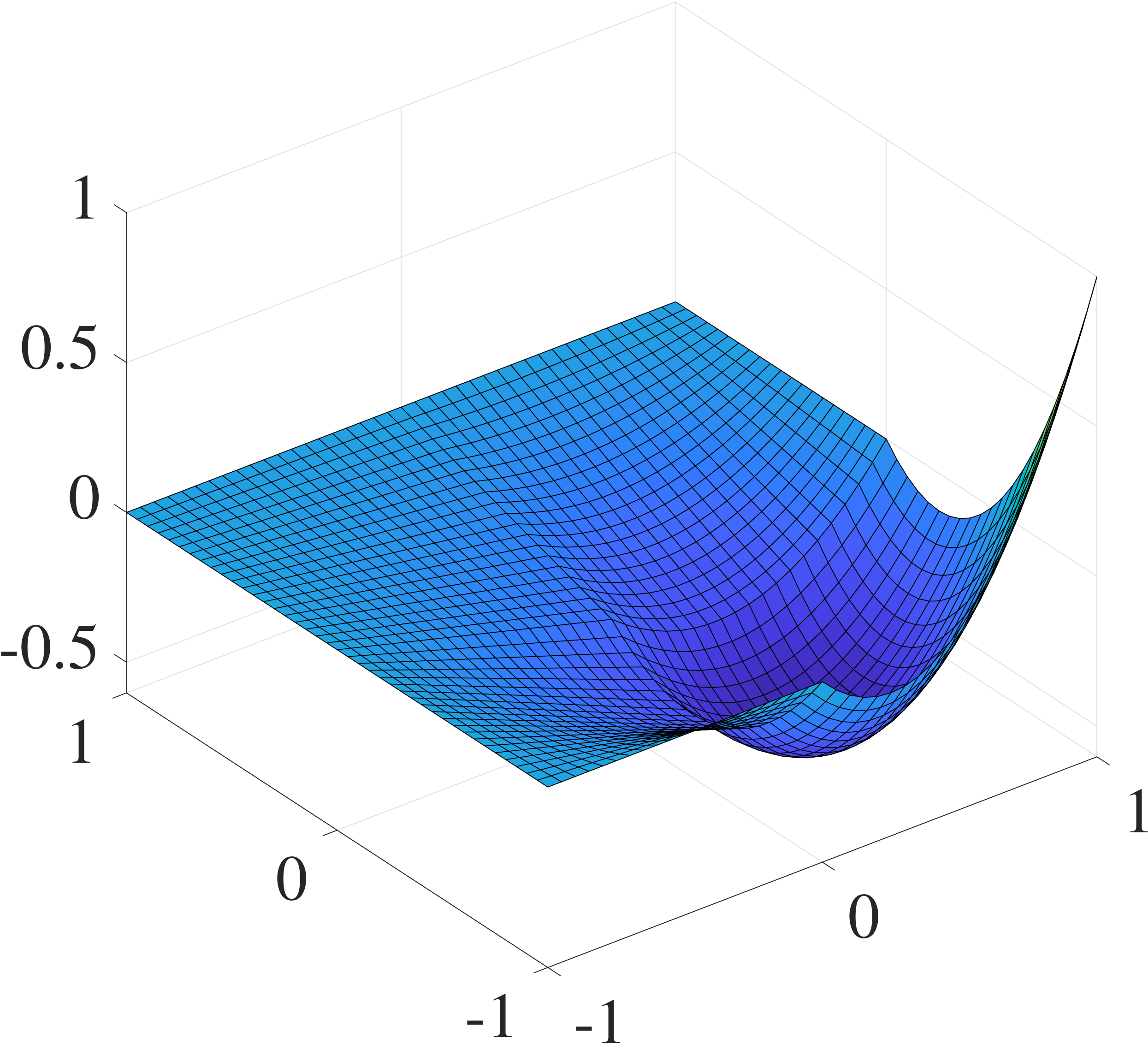}}
	\hfill
	\subfloat[Corner vertex \#3\label{fig:12NodeElem_SF3}]{\includegraphics[clip,width=0.24\textwidth]{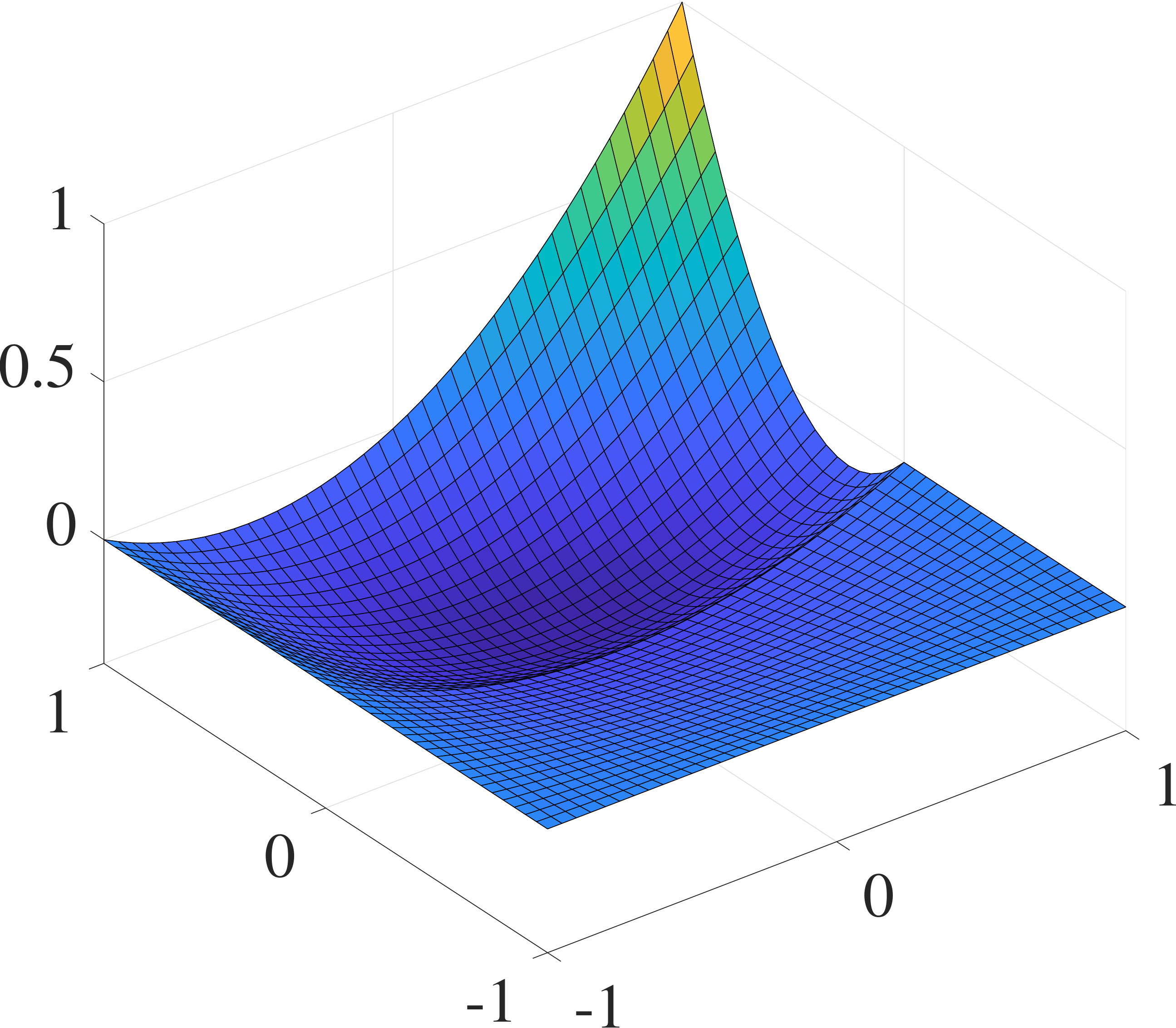}}
	\hfill
	\subfloat[Corner vertex \#4\label{fig:12NodeElem_SF4}]{\includegraphics[clip,width=0.24\textwidth]{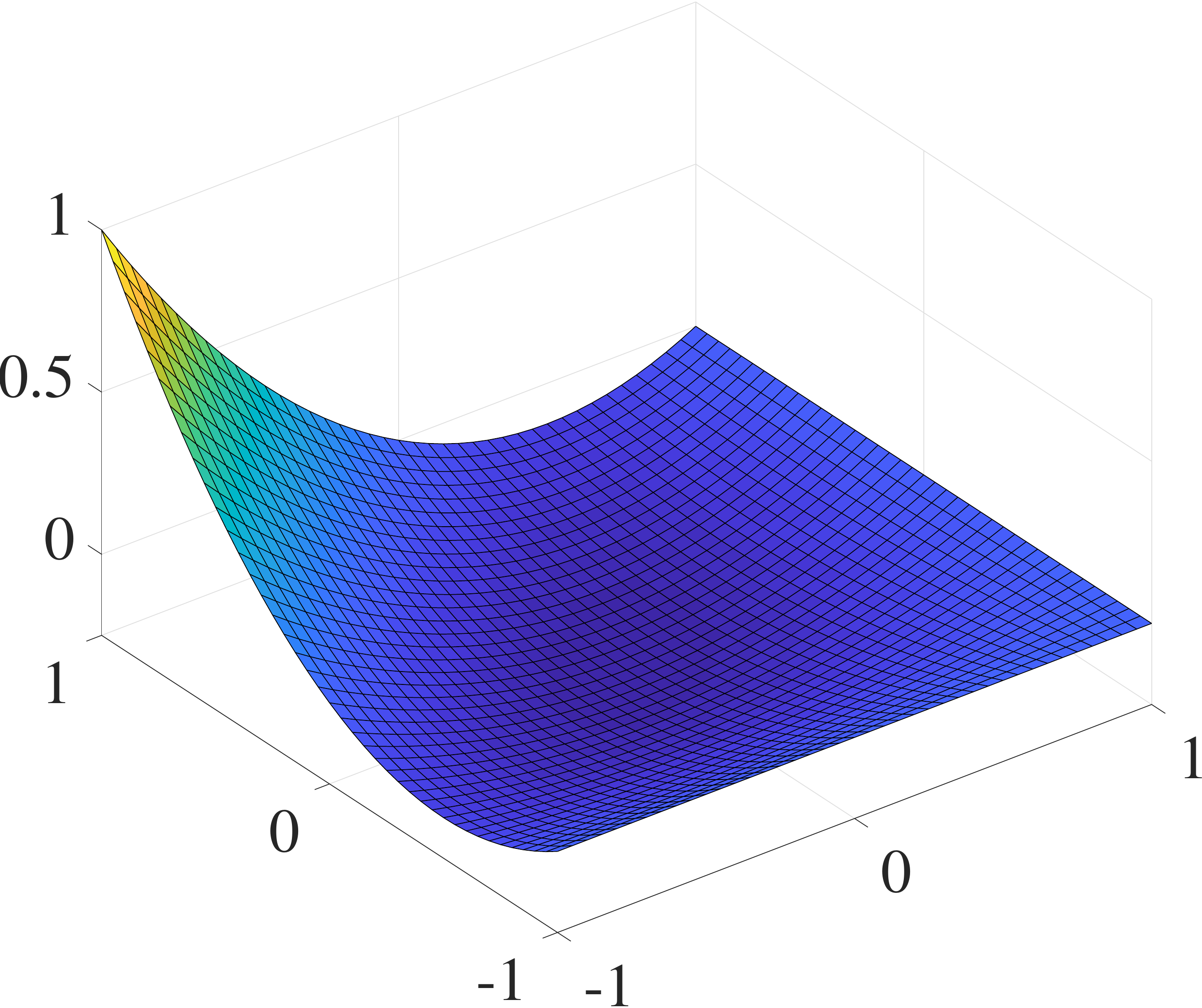}}\\
	\subfloat[Edge vertex \#5\label{fig:12NodeElem_SF5}]{\includegraphics[clip,width=0.24\textwidth]{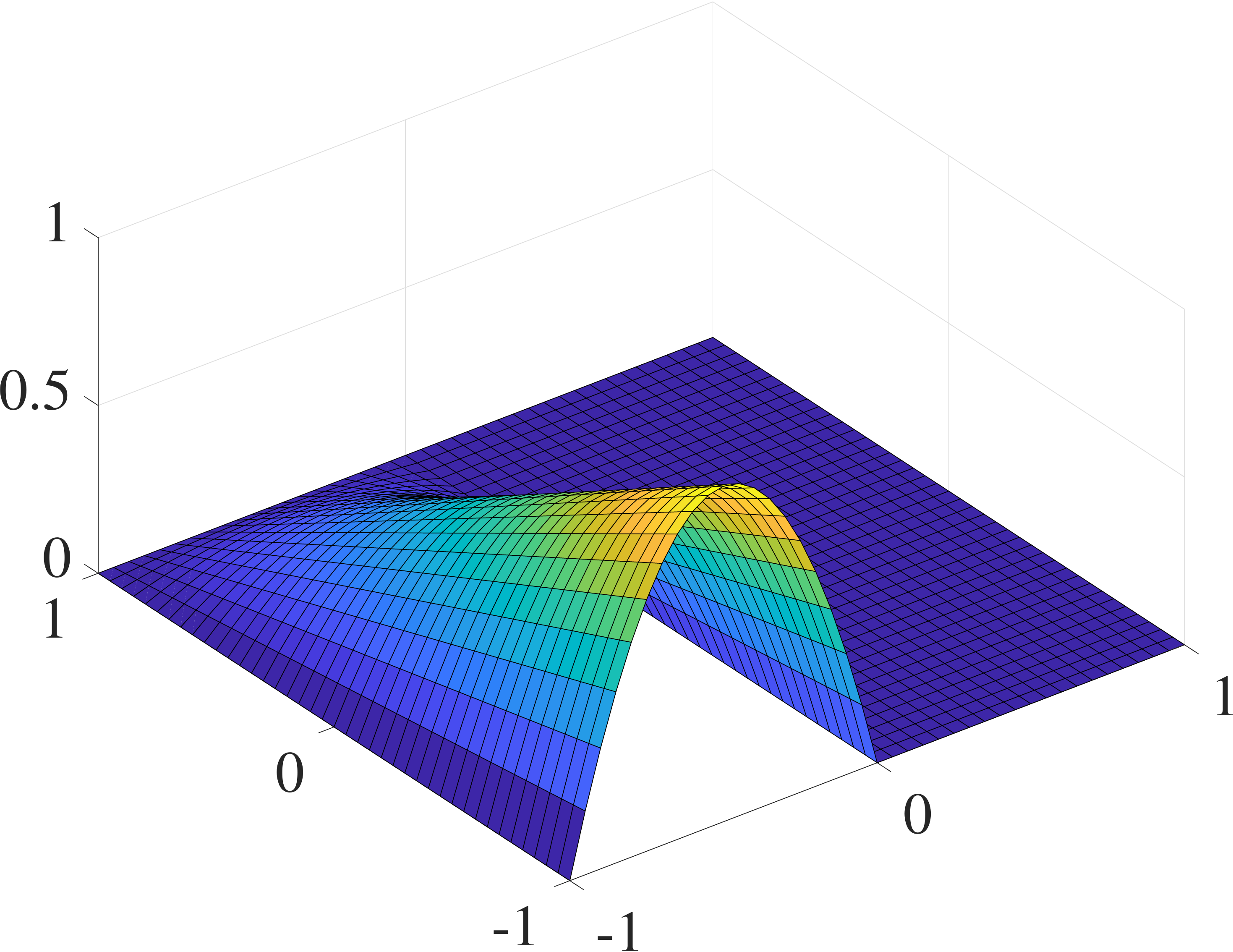}}
	\hfill
	\subfloat[Edge vertex \#6\label{fig:12NodeElem_SF6}]{\includegraphics[clip,width=0.24\textwidth]{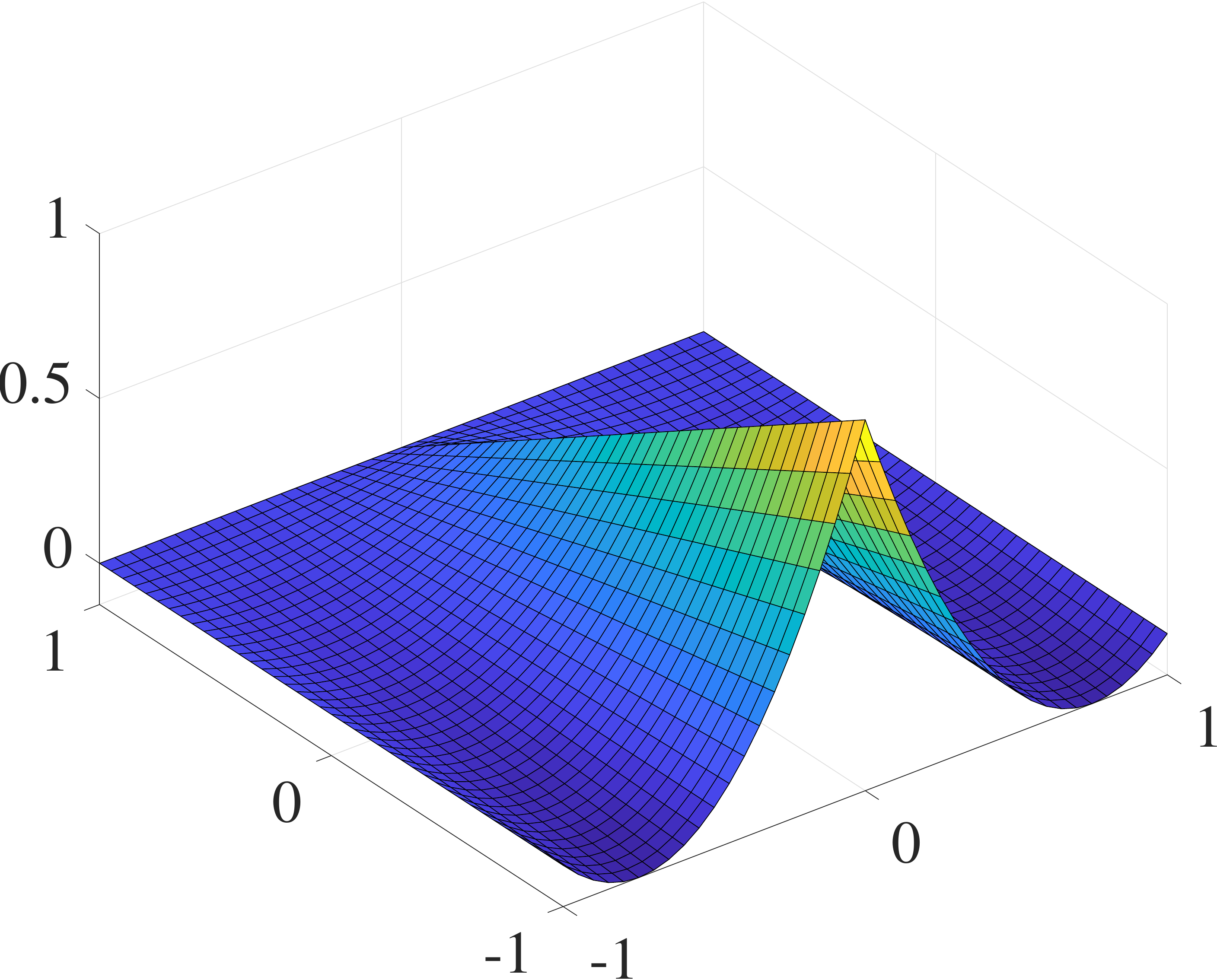}}
	\hfill
	\subfloat[Edge vertex \#7\label{fig:12NodeElem_SF7}]{\includegraphics[clip,width=0.24\textwidth]{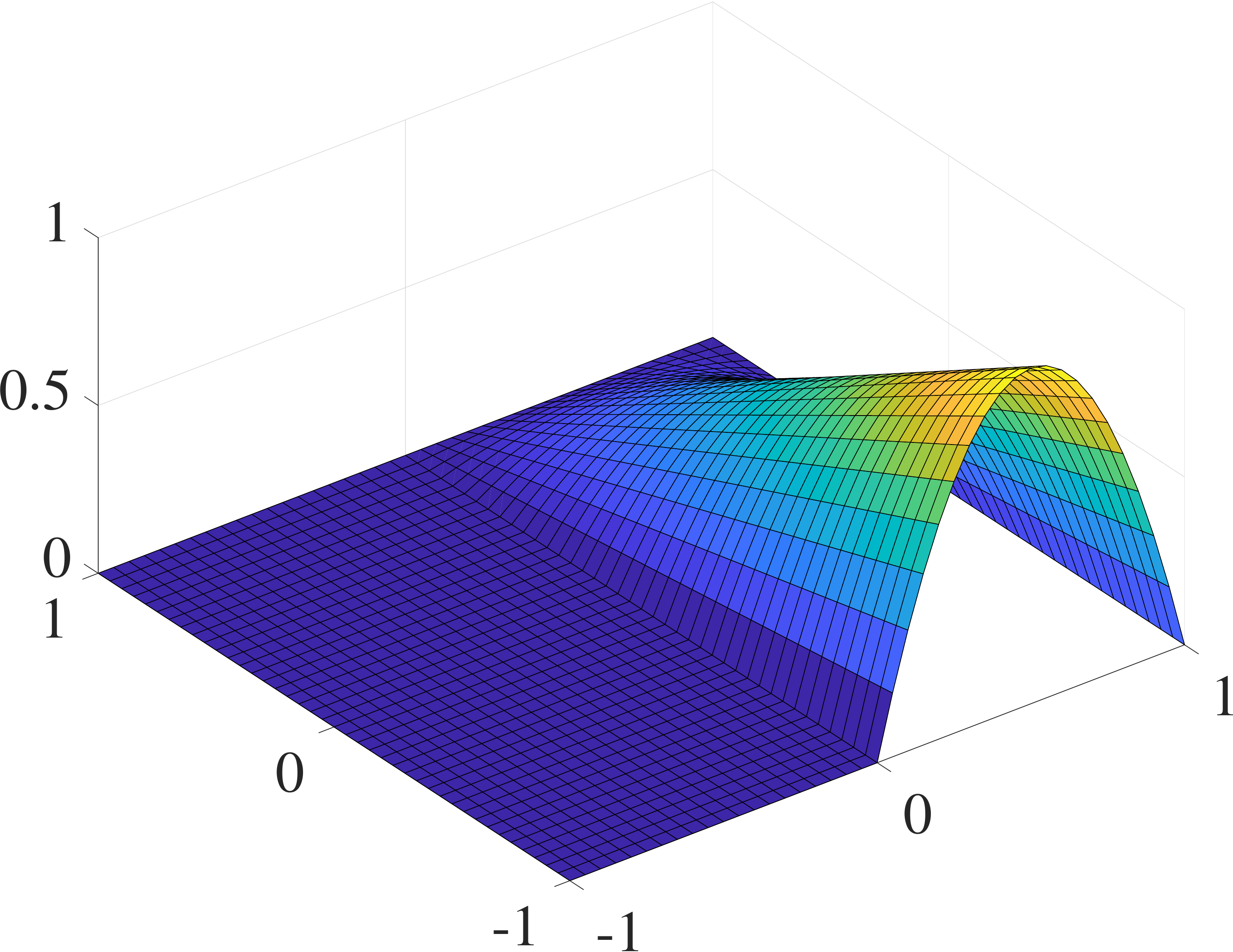}}
	\hfill
	\subfloat[Edge vertex \#8\label{fig:12NodeElem_SF8}]{\includegraphics[clip,width=0.24\textwidth]{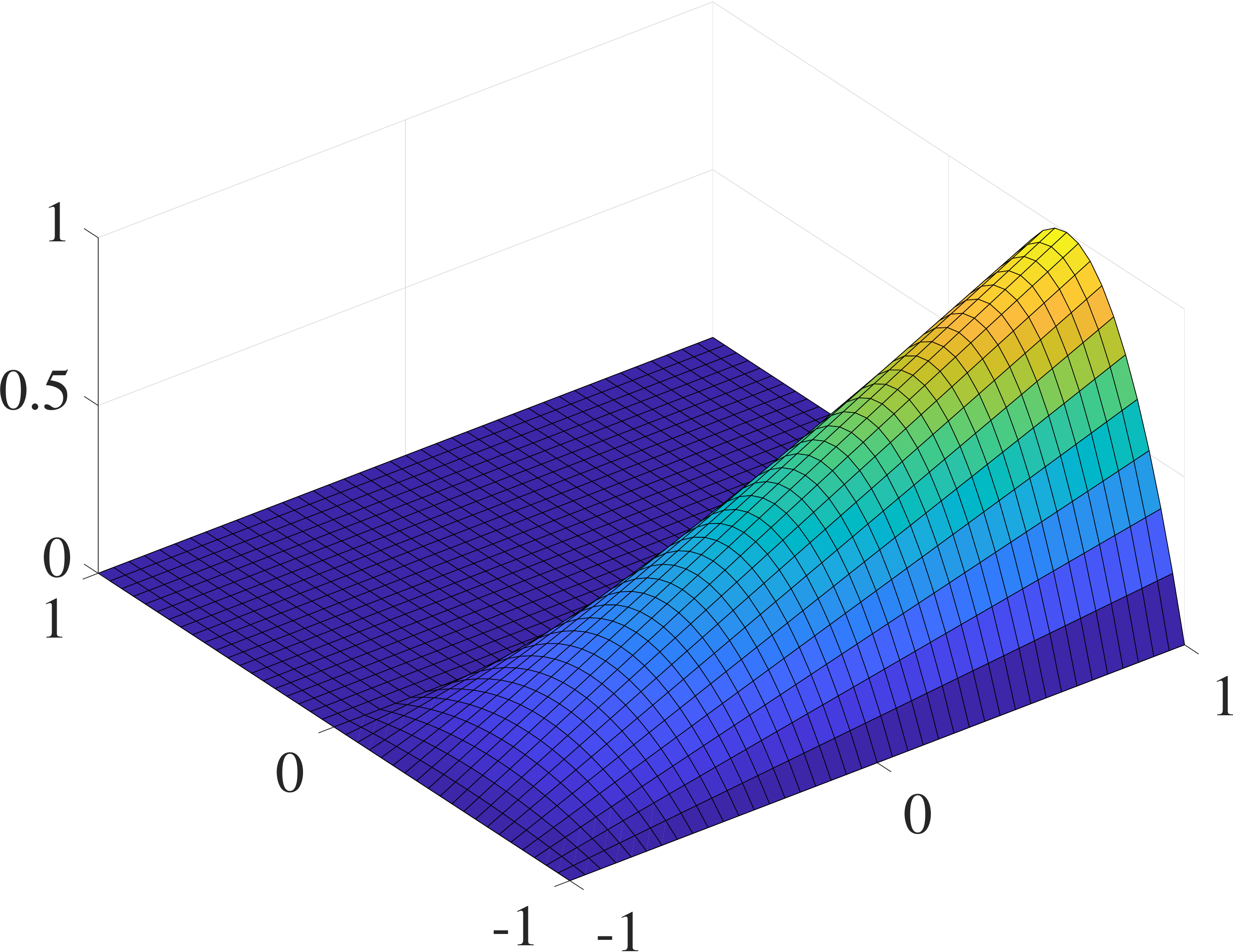}}\\
	\subfloat[Edge vertex \#9\label{fig:12NodeElem_SF9}]{\includegraphics[clip,width=0.24\textwidth]{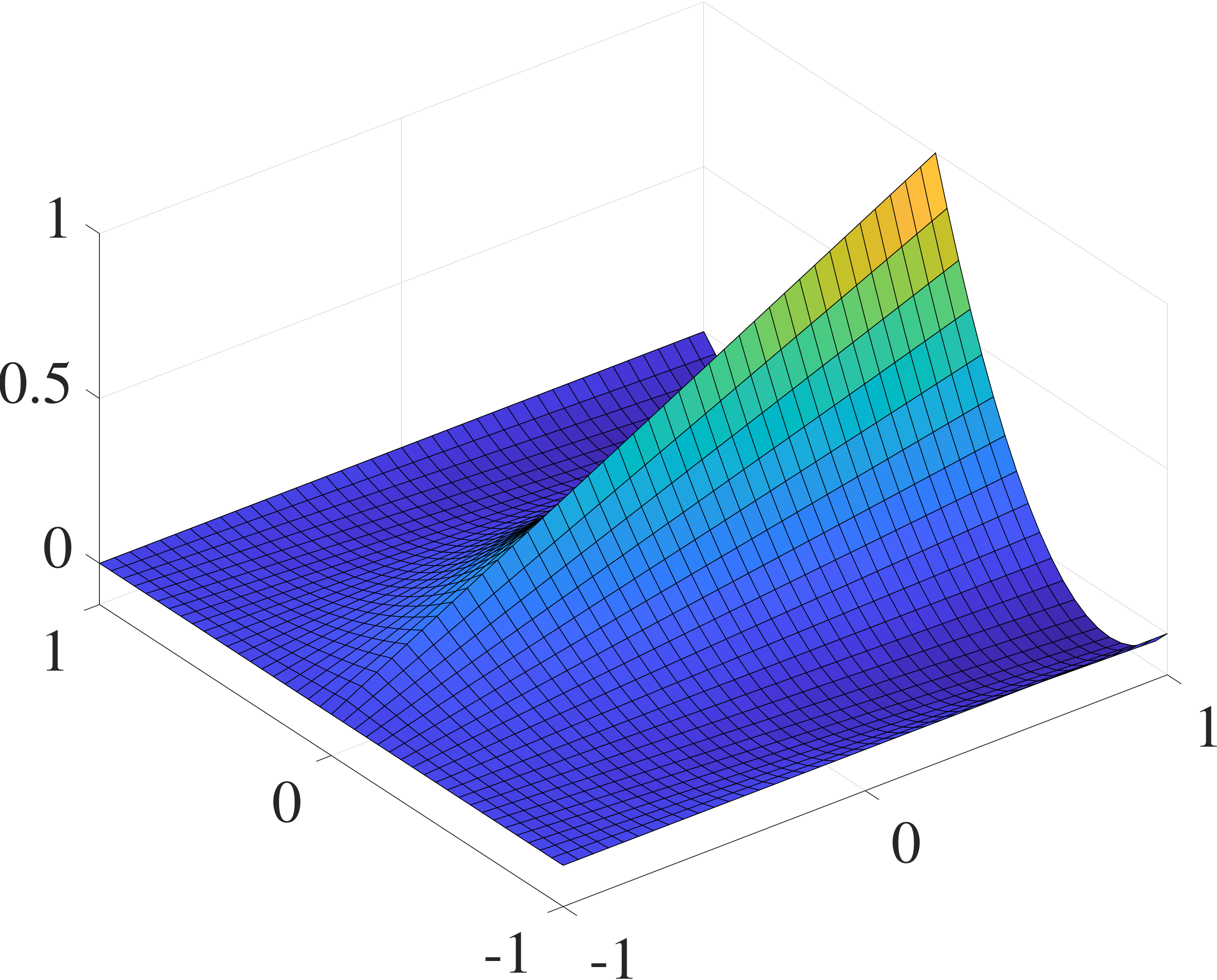}}
	\hfill
	\subfloat[Edge vertex \#10\label{fig:12NodeElem_SF10}]{\includegraphics[clip,width=0.24\textwidth]{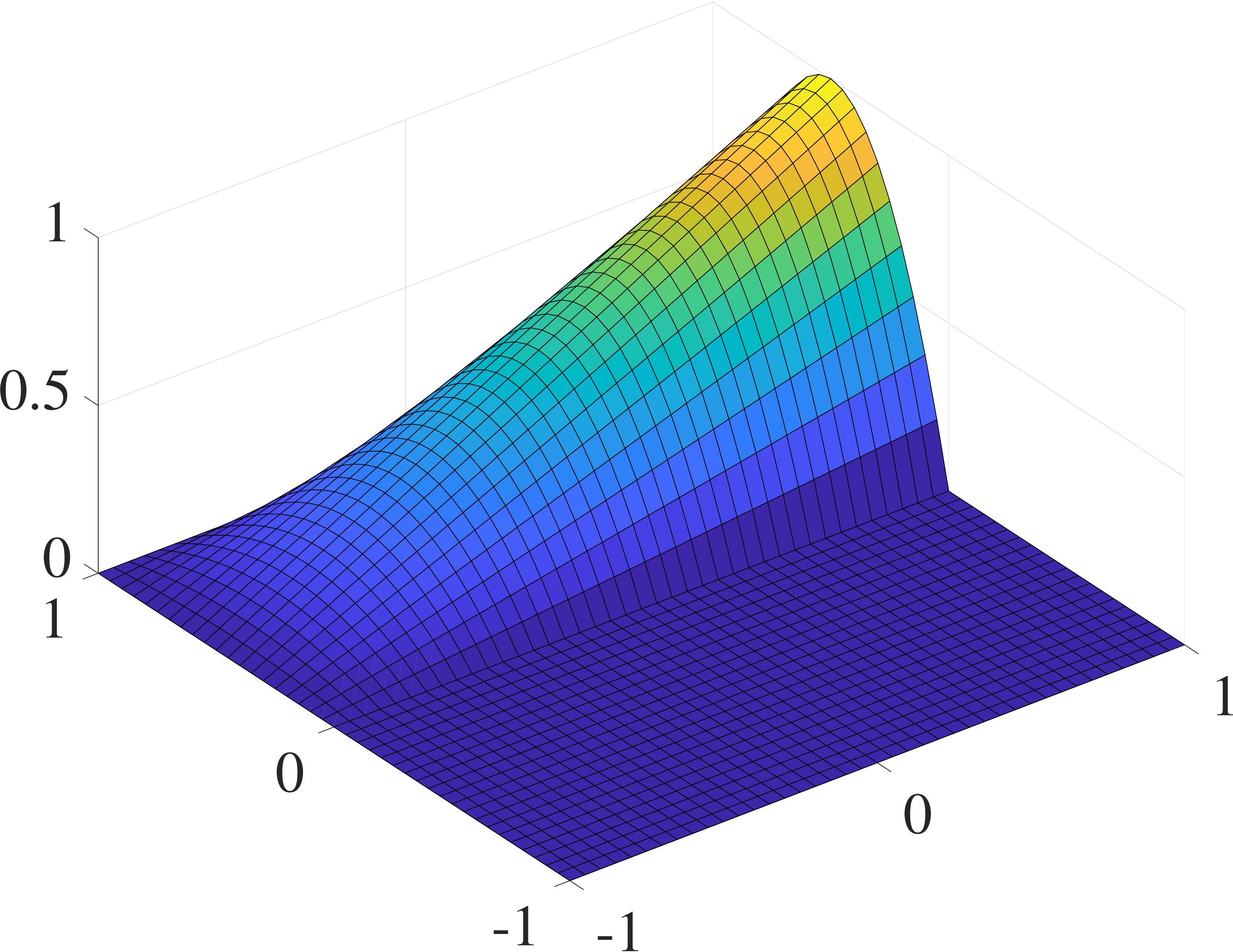}}
	\hfill
	\subfloat[Edge vertex \#11\label{fig:12NodeElem_SF11}]{\includegraphics[clip,width=0.24\textwidth]{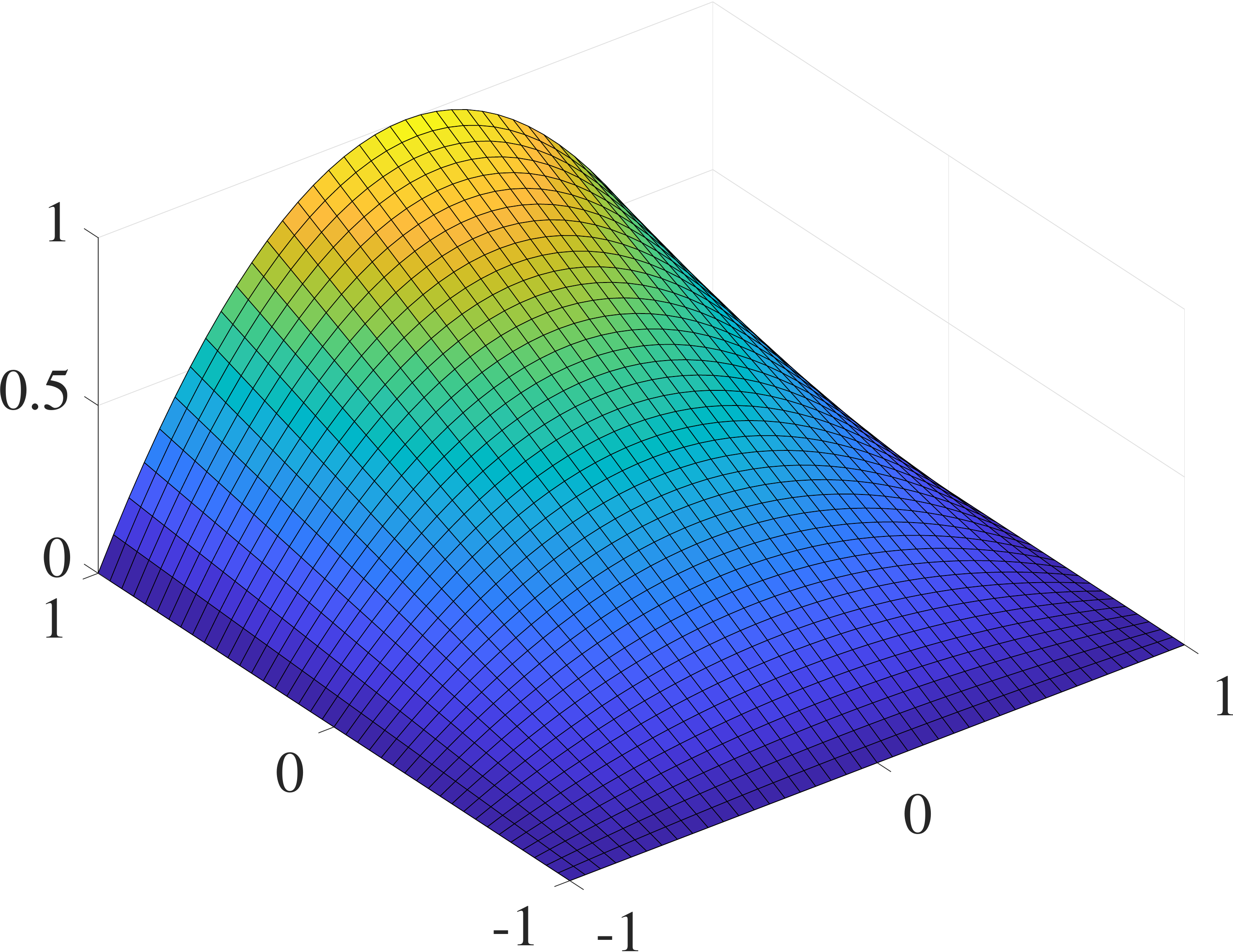}}
	\hfill
	\subfloat[Edge vertex \#12\label{fig:12NodeElem_SF12}]{\includegraphics[clip,width=0.24\textwidth]{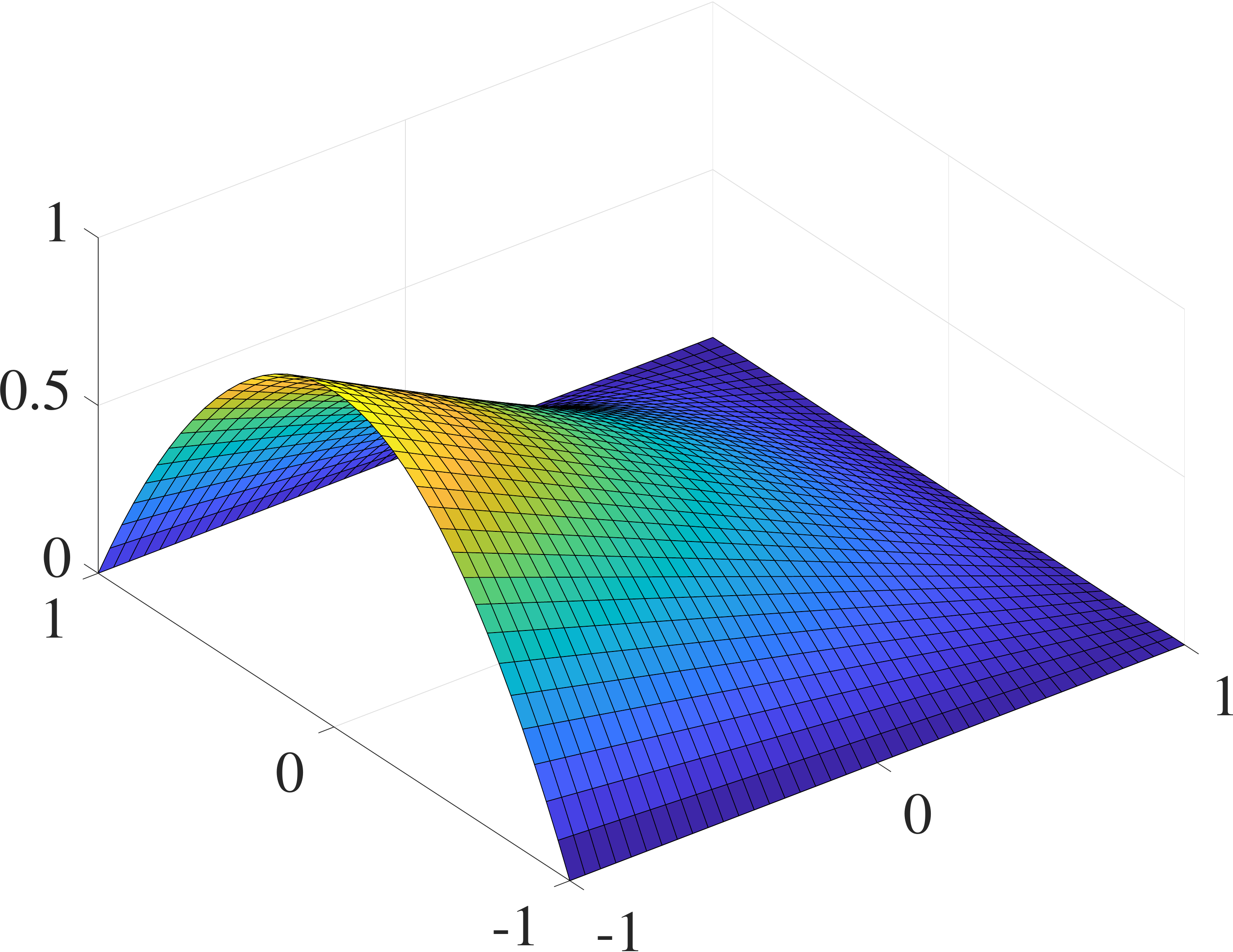}}\\
	\caption{Shape functions of the 12-node bi-quadratic transition element -- Lagrange- to Lagrange-based shape functions.}
	\label{fig:12NodeElem_SF}
\end{figure}

\paragraph{Shape functions of the corner vertices}
As mentioned before, when using linear blending functions, only the shape functions of the corner vertices are influenced by the product projector $\mathcal{P}_{\xi}[\mathcal{P}_{\eta}[\square]]$. According to Eqs.~\eqref{eq:SF_CornerNode1} to \eqref{eq:SF_CornerNode4} these shape functions are
\begin{alignat}{3}
N_1(\xi,\eta) &= N_{1}^{1}(\xi) \,^\mathrm{La}N^\mathrm{2}_{1}(\eta) &&+ N_{1}^{1}(\eta) \,^\mathrm{La}N^\mathrm{2}_{1}(\check{\xi}_1) &&- N_{1}^{1}(\xi)N_{1}^{1}(\eta)\,, \\
N_2(\xi,\eta) &= N_{2}^{1}(\xi) \,^\mathrm{La}N^\mathrm{2}_{1}(\check{\eta}_1) &&+ N_{1}^{1}(\eta) \,^\mathrm{La}N^\mathrm{2}_{3}(\check{\xi}_2) &&- N_{2}^{1}(\xi)N_{1}^{1}(\eta)\,, \\
N_3(\xi,\eta) &= N_{2}^{1}(\xi) \,^\mathrm{La}N^\mathrm{2}_{3}(\check{\eta}_2) &&+ N_{2}^{1}(\eta) \,^\mathrm{La}N^\mathrm{2}_{3}(\xi) &&- N_{2}^{1}(\xi)N_{2}^{1}(\eta)\,, \\
N_4(\xi,\eta) &= N_{1}^{1}(\xi) \,^\mathrm{La}N^\mathrm{2}_{3}(\eta) &&+ N_{2}^{1}(\eta) \,^\mathrm{La}N^\mathrm{2}_{1}(\xi) &&- N_{1}^{1}(\xi)N_{2}^{1}(\eta) = \!^\mathrm{s}N^{*}_{4}(\xi,\eta)  \,.
\end{alignat}
It is left to the reader to ascertain that the shape function belonging to node \#4 is indeed the function given in Eq.~\eqref{eq:QuadShapeFuncNode4}. The shape functions of the other nodes are altered/adjusted due to the presence of the piecewise quadratic polynomials featured on edges $E_1$ and $E_2$. The definition of the individual shape functions already indicates a difference in the intervals in which they are defined by the use of the mapped coordinates $\check{\square}$. The shape functions for the corner vertices are visualized in Figs.~\ref{fig:12NodeElem_SF1} to \ref{fig:12NodeElem_SF4}. The shape functions corresponding to the nodes that are associated with only one edge are obtained by applying either the projection operator $\mathcal{P}_{\xi}[\square]$ to shape functions living on edges $E_2$ and $E_4$ or the projection operator $\mathcal{P}_{\eta}[\square]$ to shape functions living on edges $E_1$ and $E_3$, respectively. Therefore, the shape functions are given as
\begin{alignat}{2}
&N_5(\xi,\eta)    &&= N_{1}^{1}(\eta) \,^\mathrm{La}N^\mathrm{2}_{2}(\check{\xi}_1) = \!^\mathrm{s}N^{*}_{6}(\check{\xi}_1,\eta) \,, \\
&N_6(\xi,\eta)    &&= N_{1}^{1}(\eta) 
\begin{cases}
^\mathrm{La}N^\mathrm{2}_{3}(\check{\xi}_1)\,, \\
^\mathrm{La}N^\mathrm{2}_{1}(\check{\xi}_2)\,, \\
\end{cases} \\
&N_7(\xi,\eta)    &&= N_{1}^{1}(\eta) \,^\mathrm{La}N^\mathrm{2}_{2}(\check{\xi}_2) = \!^\mathrm{s}N^{*}_{6}(\check{\xi}_2,\eta) \,, \\
&N_8(\xi,\eta)    &&= N_{2}^{1}(\xi) \,^\mathrm{La}N^\mathrm{2}_{2}(\check{\eta}_1) = \!^\mathrm{s}N^{*}_{9}(\xi,\check{\eta}_1) \,, \\
&N_9(\xi,\eta)    &&= N_{2}^{1}(\xi)
\begin{cases}
^\mathrm{La}N^\mathrm{2}_{3}(\check{\eta}_1)\,, \\
^\mathrm{La}N^\mathrm{2}_{1}(\check{\eta}_2)\,, \\
\end{cases} \\
&N_{10}(\xi,\eta) &&= N_{2}^{1}(\xi) \,^\mathrm{La}N^\mathrm{2}_{2}(\check{\eta}_2) = \!^\mathrm{s}N^{*}_{9}(\xi,\check{\eta}_2) \,, \\
&N_{11}(\xi,\eta) &&= N_{2}^{1}(\eta) \,^\mathrm{La}N^\mathrm{2}_{2}(\xi) = \!^\mathrm{s}N^{*}_{11}(\xi,\eta) \,, \\
&N_{12}(\xi,\eta) &&= N_{1}^{1}(\xi) \,^\mathrm{La}N^\mathrm{2}_{2}(\eta) = \!^\mathrm{s}N^{*}_{12}(\xi,\eta) \,.
\end{alignat}
All the shape functions of the derived transition element are illustrated in Fig.~\ref{fig:12NodeElem_SF}. Note that if we were to derive the shape functions on the basis of a bi-quadratic tensor product finite element featuring 9 nodes,  we would have to utilize quadratic blending functions as indicated in Fig.~\ref{fig:IntQuadBlending}. Otherwise, the shape functions for the edge mid-side nodes and the interior/bubble node would be incorrect\footnote{Remark: The presented approach is very general, and the standard Serendipity and spectral elements can be easily obtained as special cases. For spectral elements, we simply choose the same interpolation along all four edges. To derive transition elements, it is also possible to only set up the required shape functions at the edges to couple different element types without introducing interior nodes/modes. Those are, however, needed to ensure completeness of the polynomial basis functions and thus, theoretically optimal convergence rates.}.

Studying the shape functions for the 12-node bi-quadratic transition element, we observe that the shape functions associated with nodes \#4, 11, and 12 are identical with the standard shape functions of an 8-node bi-quadratic Serendipity finite element. The other shape functions have been adjusted in order to ensure a C\textsuperscript{0}-continuous approximation of the displacement field between adjacent elements. Using this element type, a 1-irregular mesh refinement based on quadratic elements can be executed (see Fig.~\ref{fig:QuadtreeMesh}).
\subsection{Lagrange-Legendre: Mixed nodal/modal shape functions}
\label{app:LaLe}
The second case that is addressed in this section uses the proposed transition element to couple nodal and modal elements, i.e., we aim at a conformal/compatible coupling of an 8-node Serendipity finite element with two bi-quadratic \emph{p}-elements. The explanations regarding the derivation of the shape functions for this particular transition element will be shorter as the same steps are followed which have been discussed in detail in the previous section. 

One important aspect that needs to be taken into consideration when using high order hierarchic shape functions is to ensure a C\textsuperscript{0}-continuous approximation of the solution field. To this end, the local orientation has to be taken into account and different sign-factors have to be introduced. A comprehensive description of how to achieve a conformal coupling is provided in Ref.~\cite{PhDDuczek2014}. In the following, we assume that the orientation of adjacent edges is identical and therefore, no correction has to be executed.

\begin{figure}[b!]
	\begin{minipage}[t]{0.49\textwidth}
		\centering
		\includegraphics[clip,width=1\textwidth]{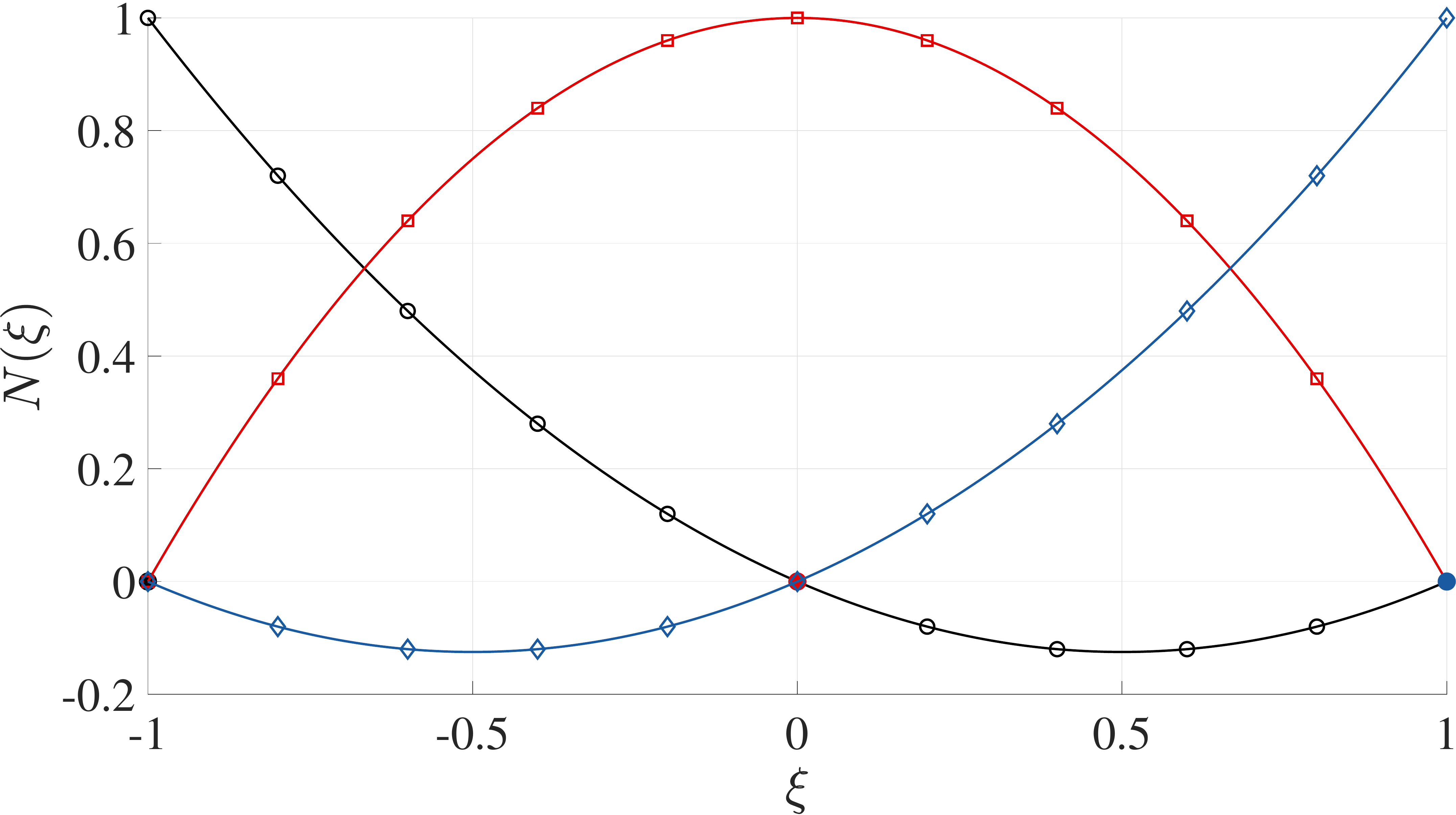}
		\caption{One-dimensional quadratic shape functions (3-node Lagrange element). Legend: \textcolor{Matlab1}{\rule[0.55ex]{5ex}{0.2ex}} $N_1(\xi)$ (\NewCirc), \textcolor{Matlab2}{\rule[0.55ex]{5ex}{0.2ex}} $N_2(\xi)$ ($\square$), \textcolor{Matlab3}{\rule[0.55ex]{5ex}{0.2ex}} $N_3(\xi)$ ($\Diamond$).}
		\label{fig:ShapeFuncQuad1d_LaLe}
	\end{minipage}
	\hfill
	\begin{minipage}[t]{0.49\textwidth}
		\centering
		\includegraphics[clip,width=1\textwidth]{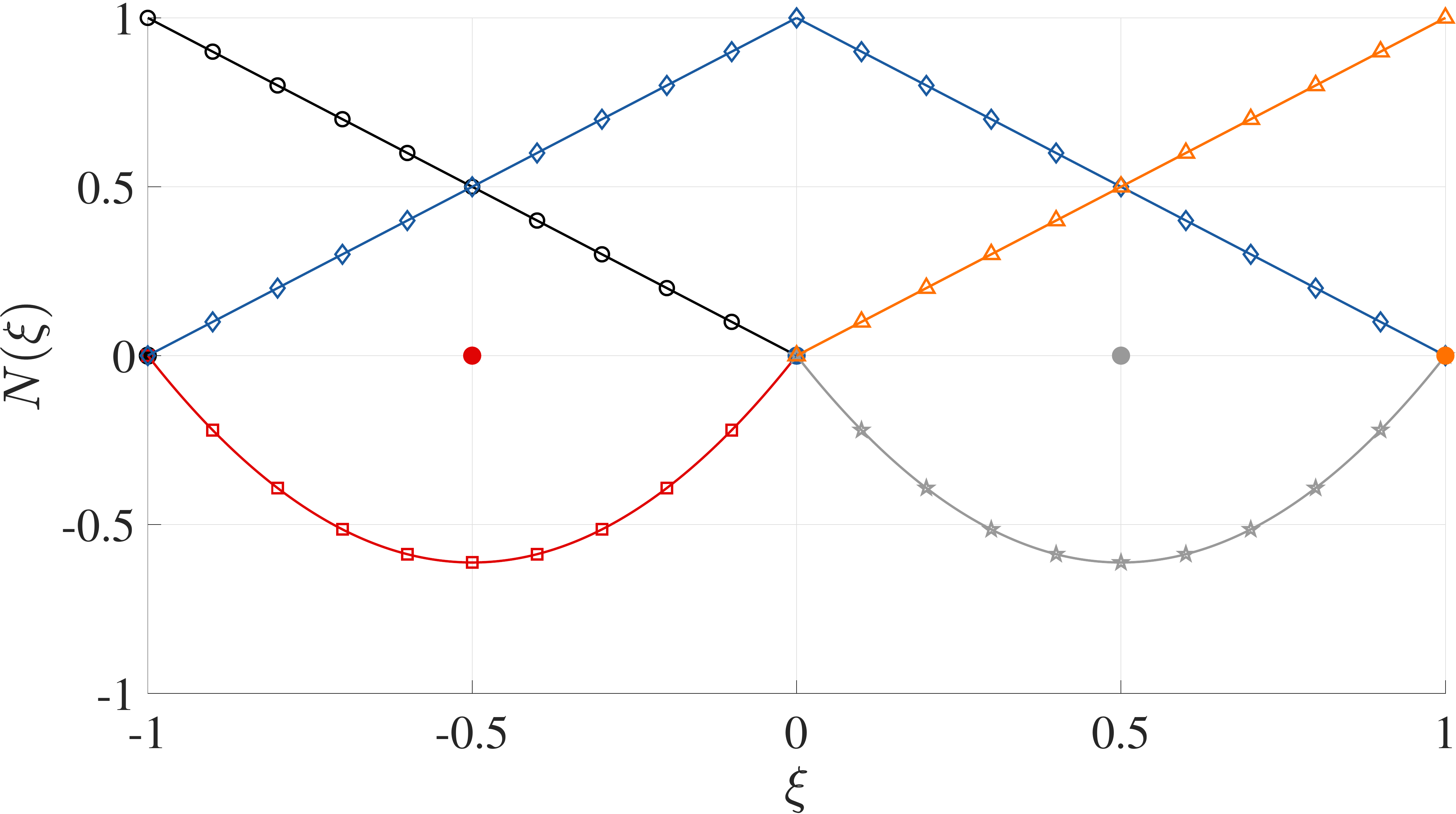}
		\caption{One-dimensional piecewise quadratic shape functions (based on 3-mode hierarchic elements). Legend: \textcolor{Matlab1}{\rule[0.55ex]{5ex}{0.2ex}} $N_1(\xi)$ (\NewCirc), \textcolor{Matlab2}{\rule[0.55ex]{5ex}{0.2ex}} $N_2(\xi)$ ($\square$), \textcolor{Matlab3}{\rule[0.55ex]{5ex}{0.2ex}} $N_3(\xi)$ ($\Diamond$), \textcolor{Matlab4}{\rule[0.55ex]{5ex}{0.2ex}} $N_4(\xi)$ (\NewStar), \textcolor{Matlab5}{\rule[0.55ex]{5ex}{0.2ex}} $N_5(\xi)$ ($\bigtriangleup$).}
		\label{fig:ShapeFuncPiecewiseQuad1d_LaLe}
	\end{minipage}
\end{figure}

The construction of the shape functions of this 12-node\footnote{Note that in this context it is conceptually wrong to refer to the transition element as a 12-node element. Considering, hierarchic \emph{p}-version shape functions we do not have physical nodes as all interior shape functions are only related to high order unknowns instead of displacements. However, for a lack of better notation we will retain this name in the remainder of the article.} transition element follows the procedure outlined in the previous section. The basis functions for the derivation are given by Eqs.~\eqref{eq:1DQuadShapeFunc_La}, \eqref{eq:1DQuadShapeFunc_Le} and \eqref{eq:QuadShapeFunc_Serendipity}. For the numbering of the nodes please see Fig.~\ref{fig:8NodeElem}. The interpolations along the four edges of the transition element are given as
\allowdisplaybreaks
\begin{align}
\text{E}_1:\; \bar{N}_i(\xi) &=
\begin{cases}
^\mathrm{Le}N^\mathrm{2}_{1}(\check{\xi}_1) & \text{for } i = 1, \text{ node \#1} \,,\\
^\mathrm{Le}N^\mathrm{2}_{2}(\check{\xi}_1) & \text{for } i = 2, \text{ node \#5} \,,\\
\begin{cases}
^\mathrm{Le}N^\mathrm{2}_{3}(\check{\xi}_1) \\
^\mathrm{Le}N^\mathrm{2}_{1}(\check{\xi}_2) \\
\end{cases}
& \text{for } i = 3, \text{ node \#6} \,, \\
^\mathrm{Le}N^\mathrm{2}_{2}(\check{\xi}_2) & \text{for } i = 4, \text{ node \#7} \,, \\
^\mathrm{Le}N^\mathrm{2}_{3}(\check{\xi}_2) & \text{for } i = 5, \text{ node \#2} \,,
\end{cases}
\\
\text{E}_2:\; \hat{N}_i(\eta) &=
\begin{cases}
^\mathrm{Le}N^\mathrm{2}_{1}(\check{\eta}_1) & \text{for } i = 1, \text{ node \#2} \,,\\
^\mathrm{Le}N^\mathrm{2}_{2}(\check{\eta}_1) & \text{for } i = 2, \text{ node \#8} \,,\\
\begin{cases}
^\mathrm{Le}N^\mathrm{2}_{3}(\check{\eta}_1) \\
^\mathrm{Le}N^\mathrm{2}_{1}(\check{\eta}_2) \\
\end{cases}
& \text{for } i = 3, \text{ node \#9} \,,\\
^\mathrm{Le}N^\mathrm{2}_{2}(\check{\eta}_2) & \text{for } i = 4, \text{ node \#10} \,, \\
^\mathrm{Le}N^\mathrm{2}_{3}(\check{\eta}_2) & \text{for } i = 5, \text{ node \#3} \,,
\end{cases}
\\
\text{E}_3:\; \tilde{N}_i(\xi) &=
\begin{cases}
^\mathrm{La}N^\mathrm{2}_{1}(\xi) & \text{for } i = 1, \text{ node \#4} \,,\\
^\mathrm{La}N^\mathrm{2}_{2}(\xi) & \text{for } i = 2, \text{ node \#11} \,,\\
^\mathrm{La}N^\mathrm{2}_{3}(\xi) & \text{for } i = 3, \text{ node \#3} \,,
\end{cases}
\\
\text{E}_4:\; \breve{N}_i(\eta) &=
\begin{cases}
^\mathrm{La}N^\mathrm{2}_{1}(\eta) & \text{for } i = 1, \text{ node \#1} \,,\\
^\mathrm{La}N^\mathrm{2}_{2}(\eta) & \text{for } i = 2, \text{ node \#12} \,,\\
^\mathrm{La}N^\mathrm{2}_{3}(\eta) & \text{for } i = 3, \text{ node \#4} \,.
\end{cases}
\end{align}
Considering the expressions to compute the shape functions on the divided edges $E_1$ and $E_2$, we have to introduce new coordinates $\check{\xi}_1$, $\check{\xi}_2$, $\check{\eta}_1$, and $\check{\eta}_2$. These coordinates are defined by a simple linear transformation from $[-1,0]$ or $[0,1]$ to the reference domain $[-1,1]$. The shape functions for the edges are exemplarily sketched in Figs.~\ref{fig:ShapeFuncQuad1d_LaLe} and \ref{fig:ShapeFuncPiecewiseQuad1d_LaLe}.

The functions living on the edges of the transition element are in the next step projected/blended into the interior of the element by using the projection operators $\mathcal{P}_{\xi}[\square]$, $\mathcal{P}_{\eta}[\square]$, and $\mathcal{P}_{\xi}[\mathcal{P}_{\eta}[\square]]$. First, we will derive the shape functions that are associated to the corner vertices:
\begin{alignat}{3}
N_1(\xi,\eta) &= N_{1}^{1}(\xi) \,^\mathrm{La}N^\mathrm{2}_{1}(\eta) &&+ N_{1}^{1}(\eta) \,^\mathrm{Le}N^\mathrm{2}_{1}(\check{\xi}_1) &&- N_{1}^{1}(\xi)N_{1}^{1}(\eta)\,, \\
N_2(\xi,\eta) &= N_{2}^{1}(\xi) \,^\mathrm{Le}N^\mathrm{2}_{1}(\check{\eta}_1) &&+ N_{1}^{1}(\eta) \,^\mathrm{Le}N^\mathrm{2}_{3}(\check{\xi}_2) &&- N_{2}^{1}(\xi)N_{1}^{1}(\eta)\,, \\
N_3(\xi,\eta) &= N_{2}^{1}(\xi) \,^\mathrm{Le}N^\mathrm{2}_{3}(\check{\eta}_2) &&+ N_{2}^{1}(\eta) \,^\mathrm{La}N^\mathrm{2}_{3}(\xi) &&- N_{2}^{1}(\xi)N_{2}^{1}(\eta)\,, \\
N_4(\xi,\eta) &= N_{1}^{1}(\xi) \,^\mathrm{La}N^\mathrm{2}_{3}(\eta) &&+ N_{2}^{1}(\eta) \,^\mathrm{La}N^\mathrm{2}_{1}(\xi) &&- N_{1}^{1}(\xi)N_{2}^{1}(\eta) = \!^\mathrm{s}N^{*}_{4}(\xi,\eta)  \,.
\end{alignat}
Second, the edge shape functions are listed:
\begin{alignat}{2}
&N_5(\xi,\eta)    &&= N_{1}^{1}(\eta) \,^\mathrm{Le}N^\mathrm{2}_{2}(\check{\xi}_1) = \!^\mathrm{s}N^{*}_{6}(\check{\xi}_1,\eta) \,, \\
&N_6(\xi,\eta)    &&= N_{1}^{1}(\eta) 
\begin{cases}
^\mathrm{Le}N^\mathrm{2}_{3}(\check{\xi}_1)\,, \\
^\mathrm{Le}N^\mathrm{2}_{1}(\check{\xi}_2)\,, \\
\end{cases} \\
&N_7(\xi,\eta)    &&= N_{1}^{1}(\eta) \,^\mathrm{Le}N^\mathrm{2}_{2}(\check{\xi}_2) = \!^\mathrm{s}N^{*}_{6}(\check{\xi}_2,\eta) \,, \\
&N_8(\xi,\eta)    &&= N_{2}^{1}(\xi) \,^\mathrm{Le}N^\mathrm{2}_{2}(\check{\eta}_1) = \!^\mathrm{s}N^{*}_{9}(\xi,\check{\eta}_1) \,, \\
&N_9(\xi,\eta)    &&= N_{2}^{1}(\xi)
\begin{cases}
^\mathrm{Le}N^\mathrm{2}_{3}(\check{\eta}_1)\,, \\
^\mathrm{Le}N^\mathrm{2}_{1}(\check{\eta}_2)\,, \\
\end{cases} \\
&N_{10}(\xi,\eta) &&= N_{2}^{1}(\xi) \,^\mathrm{Le}N^\mathrm{2}_{2}(\check{\eta}_2) = \!^\mathrm{s}N^{*}_{9}(\xi,\check{\eta}_2) \,, \\
&N_{11}(\xi,\eta) &&= N_{2}^{1}(\eta) \,^\mathrm{La}N^\mathrm{2}_{2}(\xi) = \!^\mathrm{s}N^{*}_{11}(\xi,\eta) \,, \\
&N_{12}(\xi,\eta) &&= N_{1}^{1}(\xi) \,^\mathrm{La}N^\mathrm{2}_{2}(\eta) = \!^\mathrm{s}N^{*}_{12}(\xi,\eta) \,.
\end{alignat}
\begin{figure}[b!]
	\centering
	\subfloat[Corner vertex \#1\label{fig:12NodeElem_SF1_LaLe}]{\includegraphics[clip,width=0.24\textwidth]{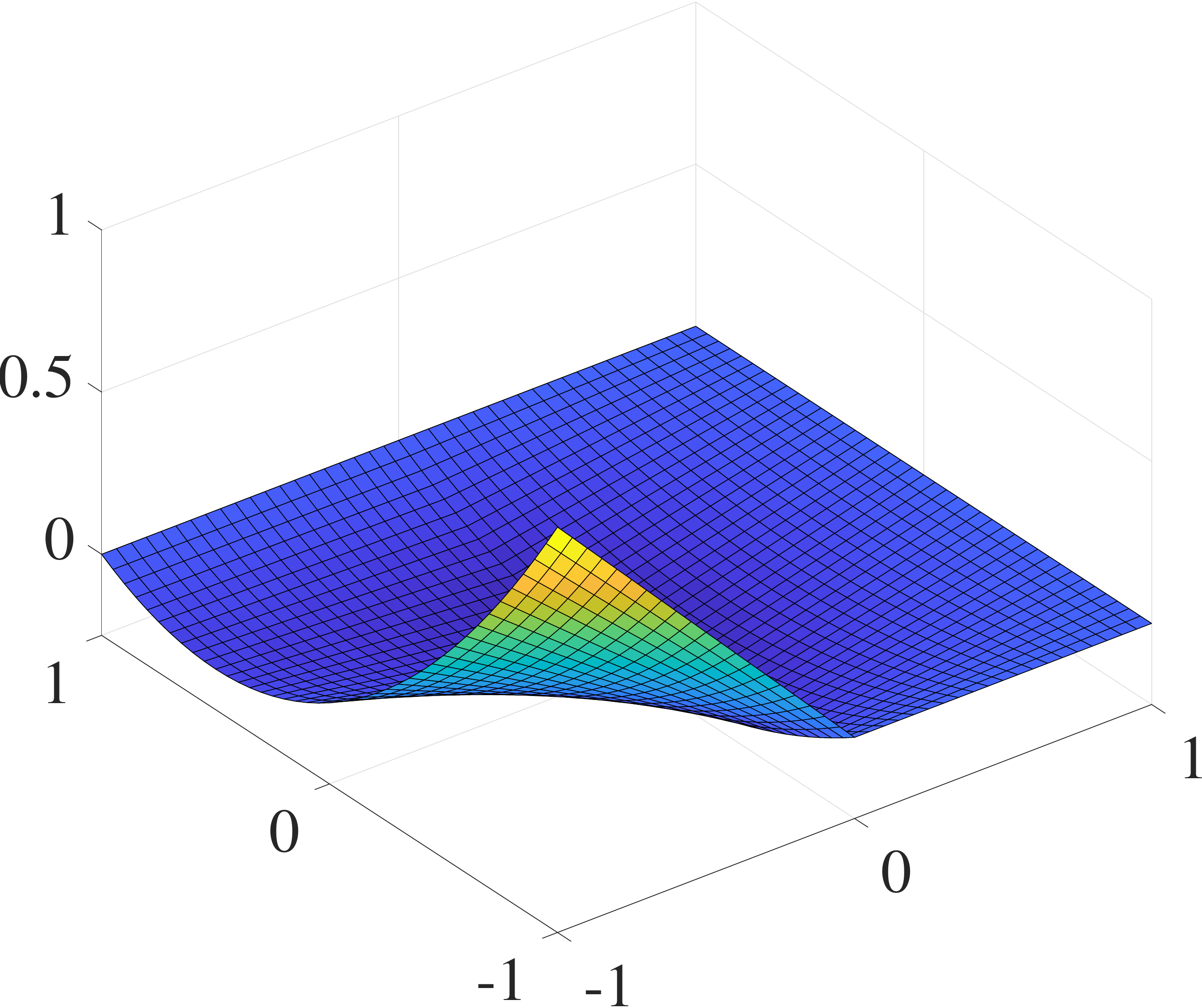}}
	\hfill
	\subfloat[Corner vertex \#2\label{fig:12NodeElem_SF2_LaLe}]{\includegraphics[clip,width=0.24\textwidth]{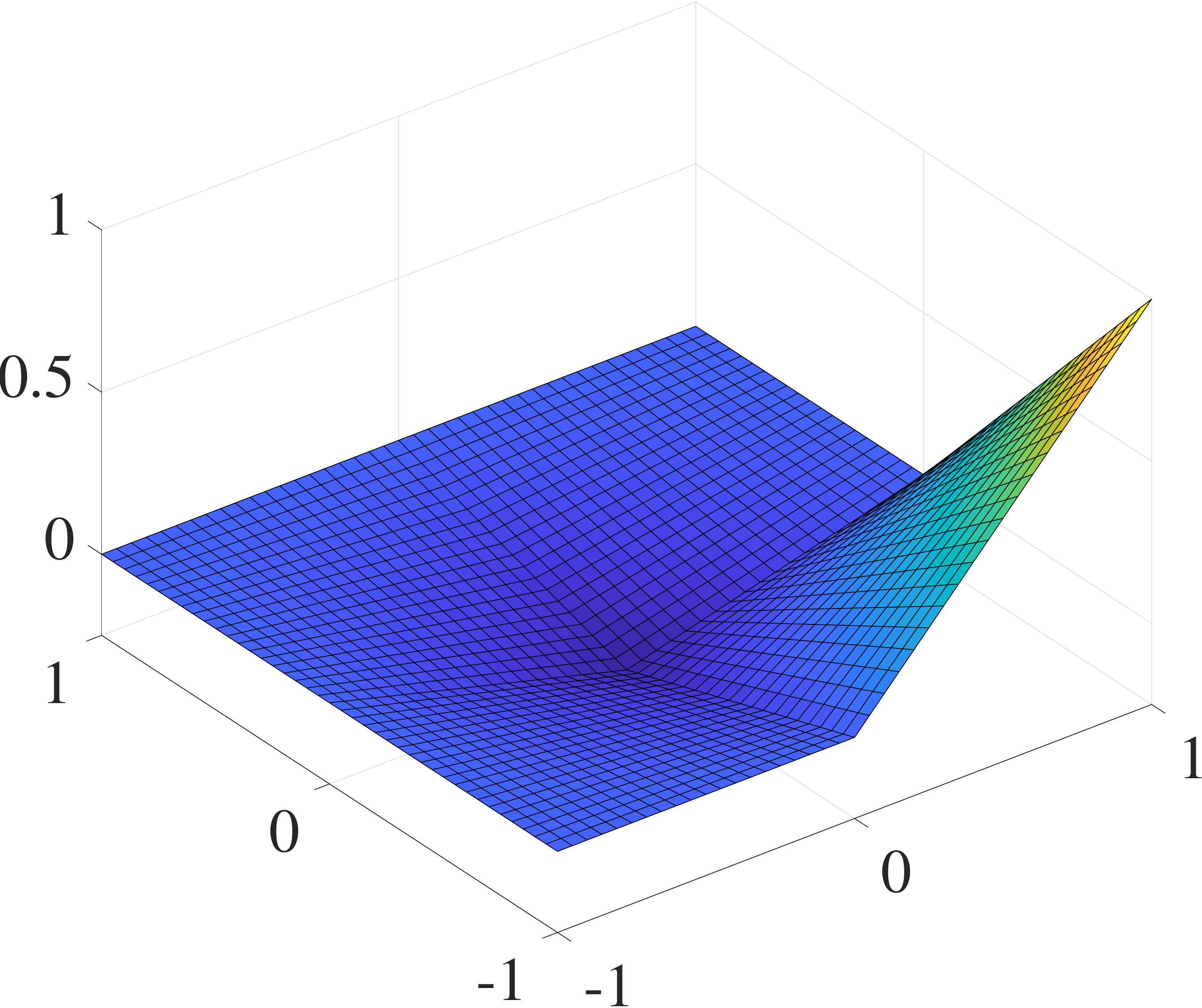}}
	\hfill
	\subfloat[Corner vertex \#3\label{fig:12NodeElem_SF3_LaLe}]{\includegraphics[clip,width=0.24\textwidth]{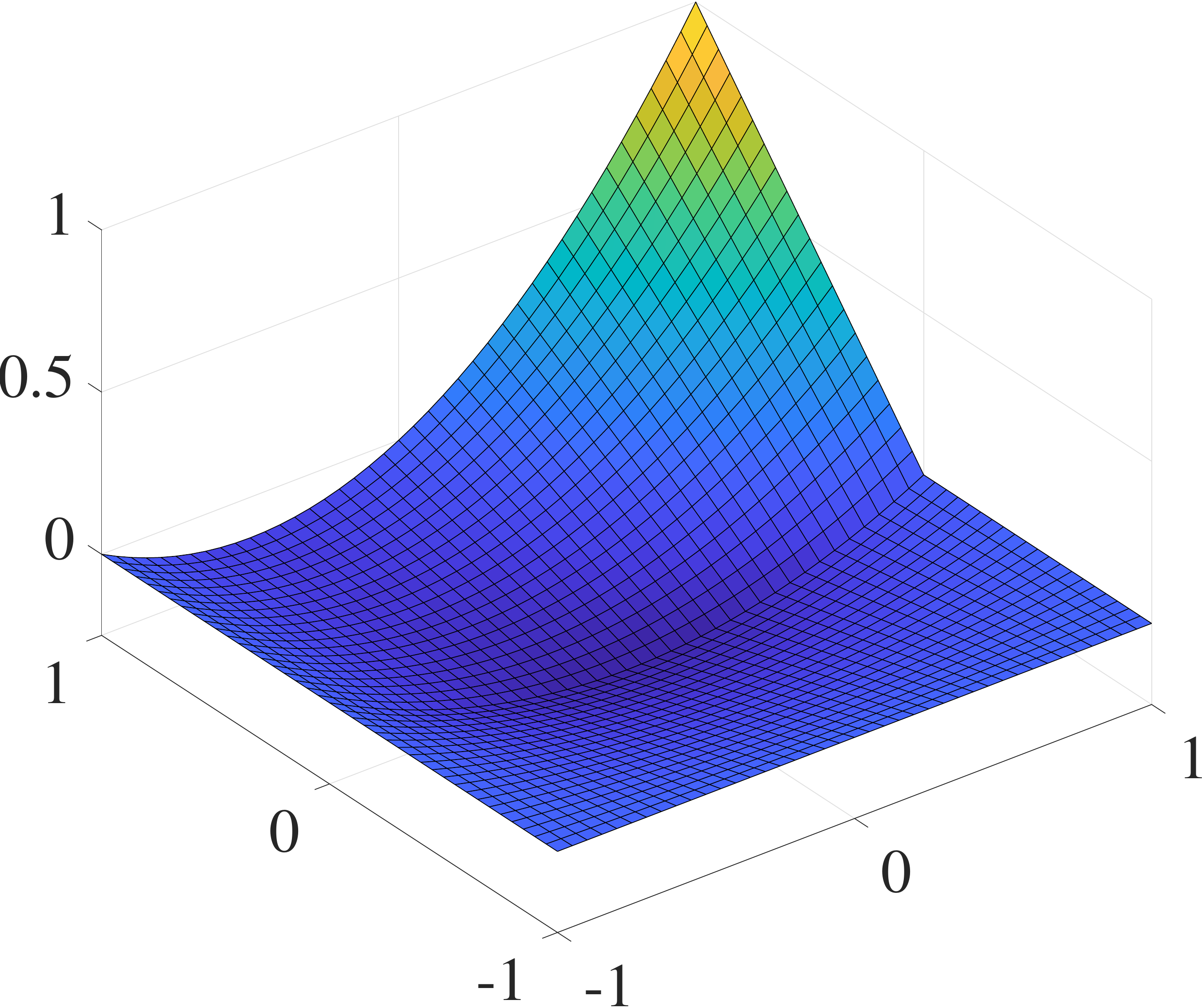}}
	\hfill
	\subfloat[Corner vertex \#4\label{fig:12NodeElem_SF4_LaLe}]{\includegraphics[clip,width=0.24\textwidth]{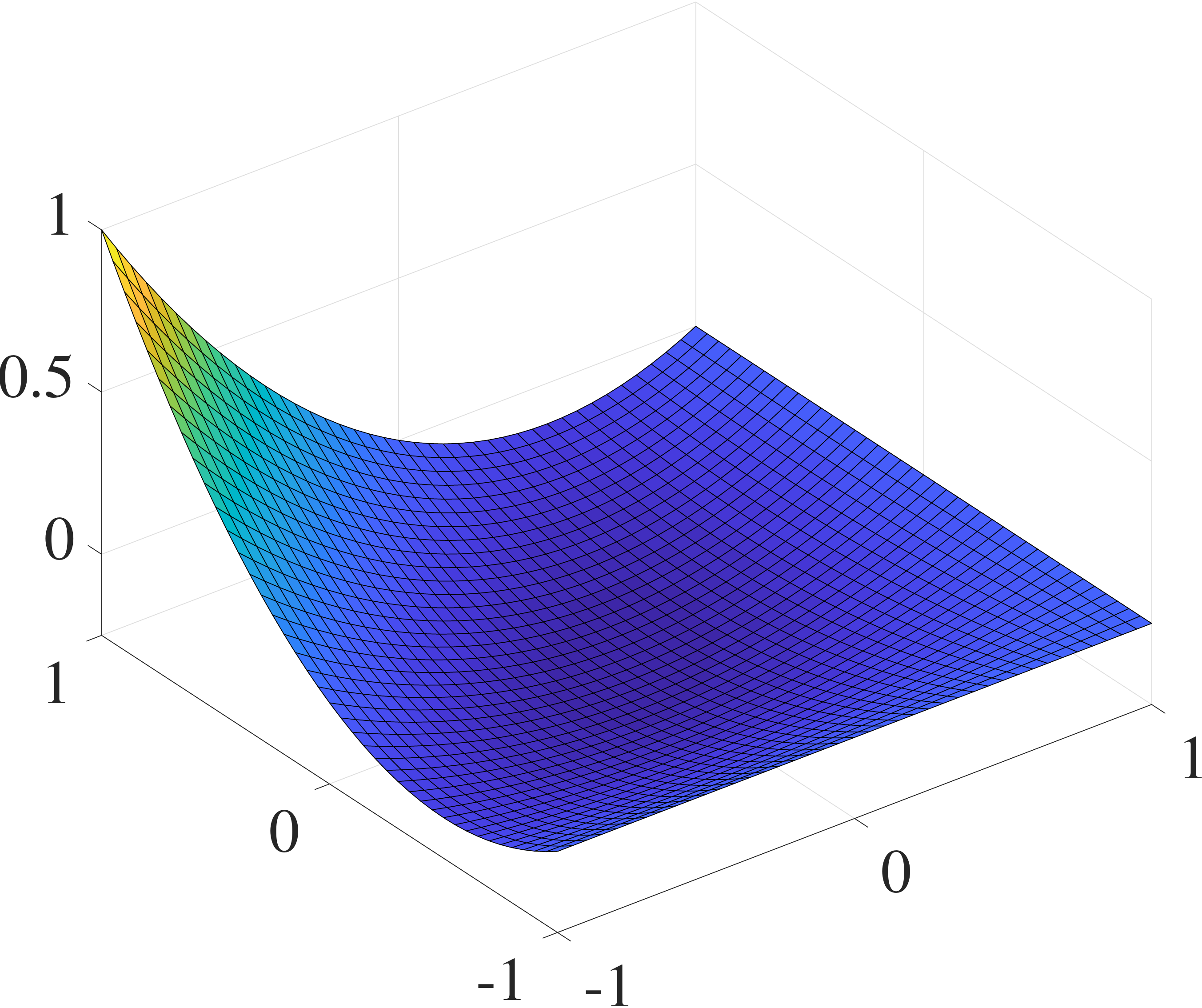}}\\
	\subfloat[Edge vertex \#5\label{fig:12NodeElem_SF5_LaLe}]{\includegraphics[clip,width=0.24\textwidth]{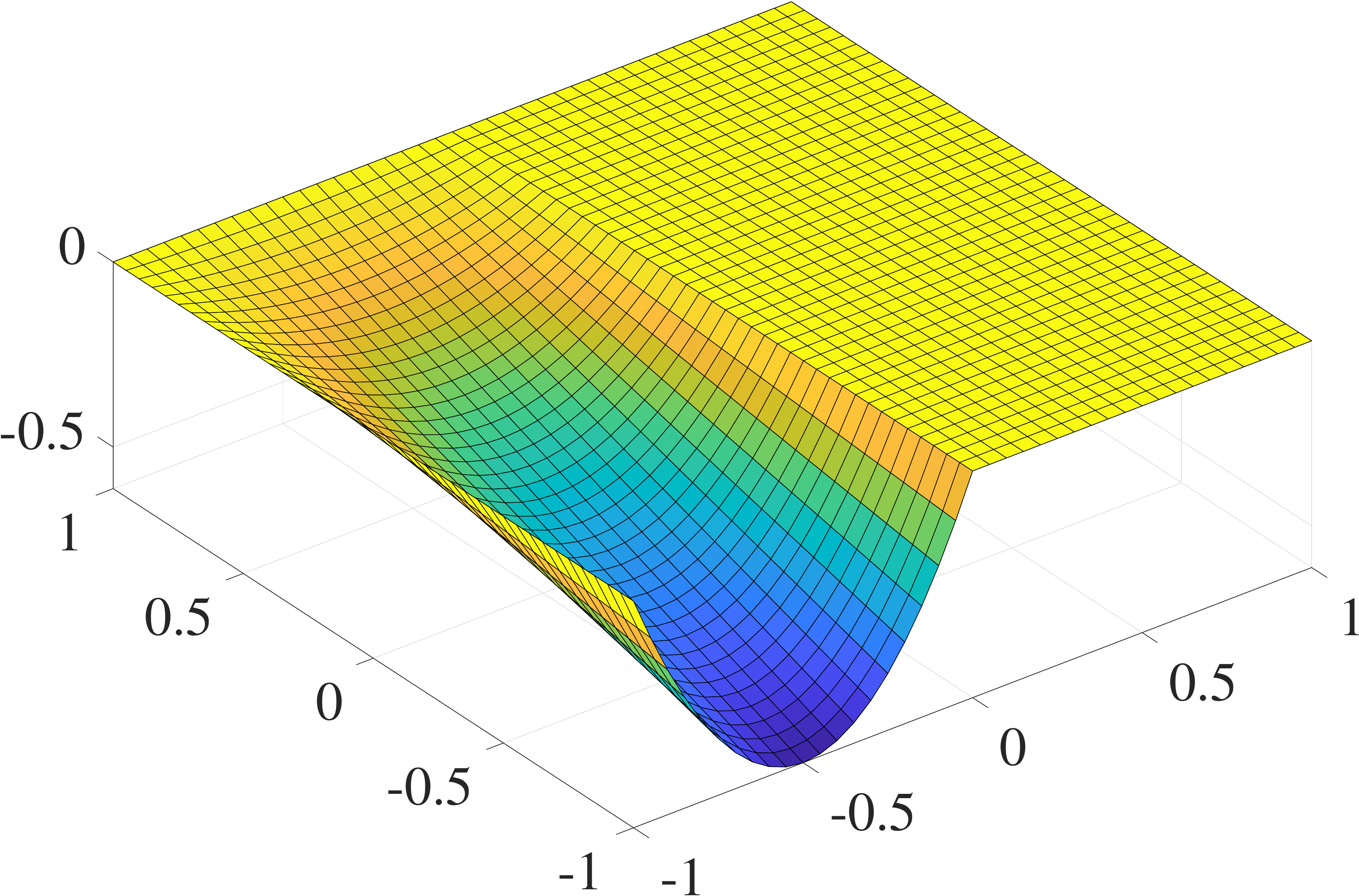}}
	\hfill
	\subfloat[Edge vertex \#6\label{fig:12NodeElem_SF6_LaLe}]{\includegraphics[clip,width=0.24\textwidth]{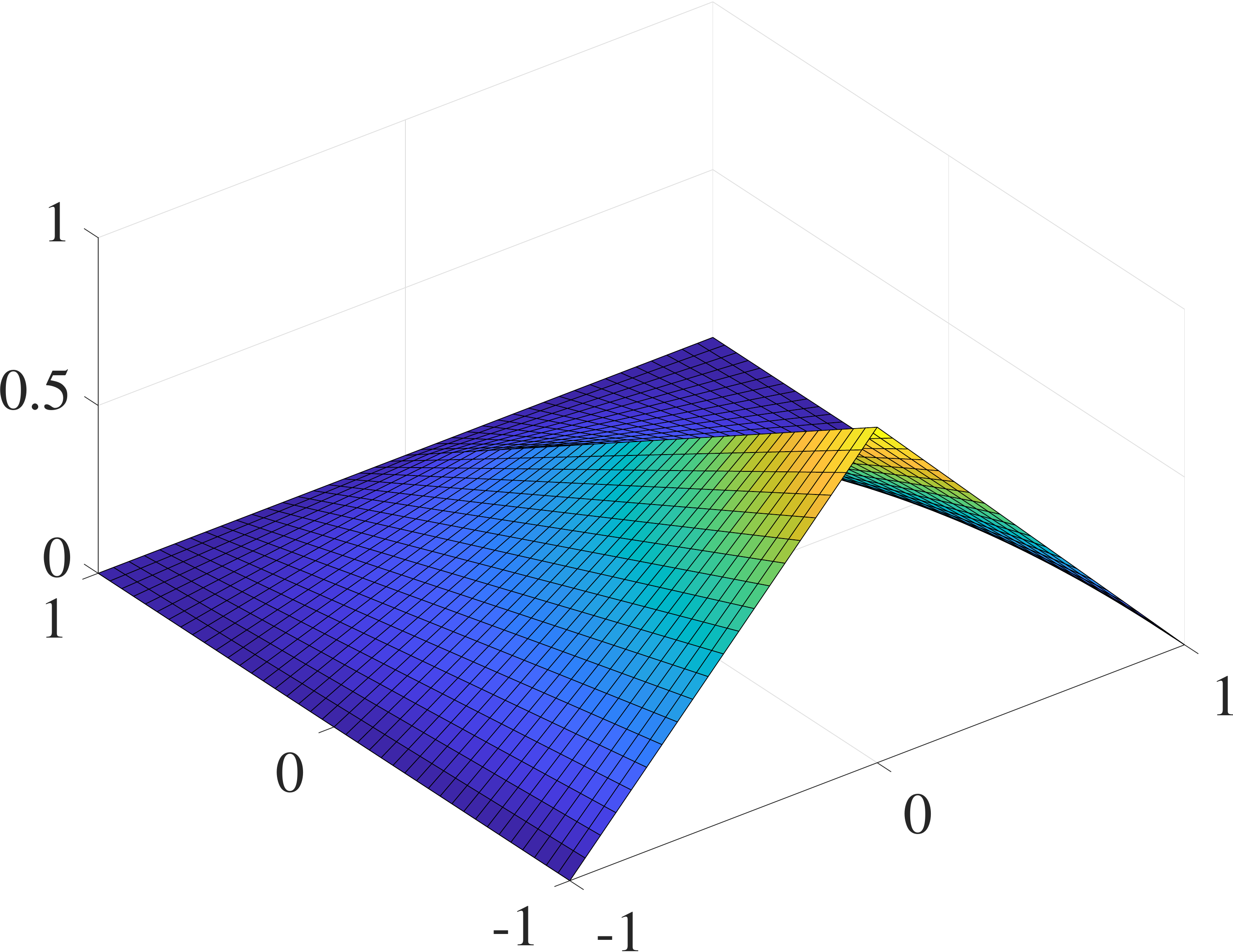}}
	\hfill
	\subfloat[Edge vertex \#7\label{fig:12NodeElem_SF7_LaLe}]{\includegraphics[clip,width=0.24\textwidth]{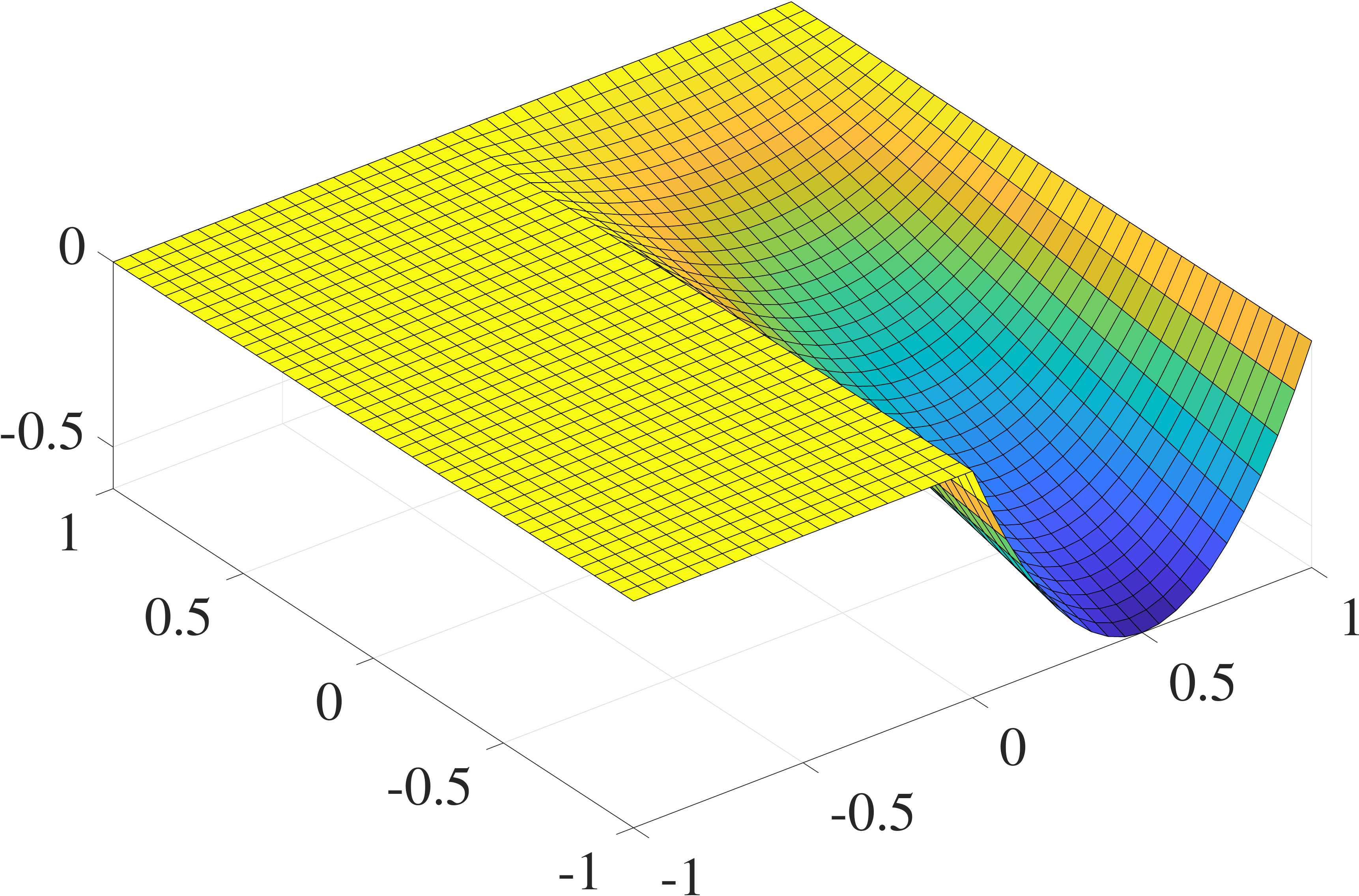}}
	\hfill
	\subfloat[Edge vertex \#8\label{fig:12NodeElem_SF8_LaLe}]{\includegraphics[clip,width=0.24\textwidth]{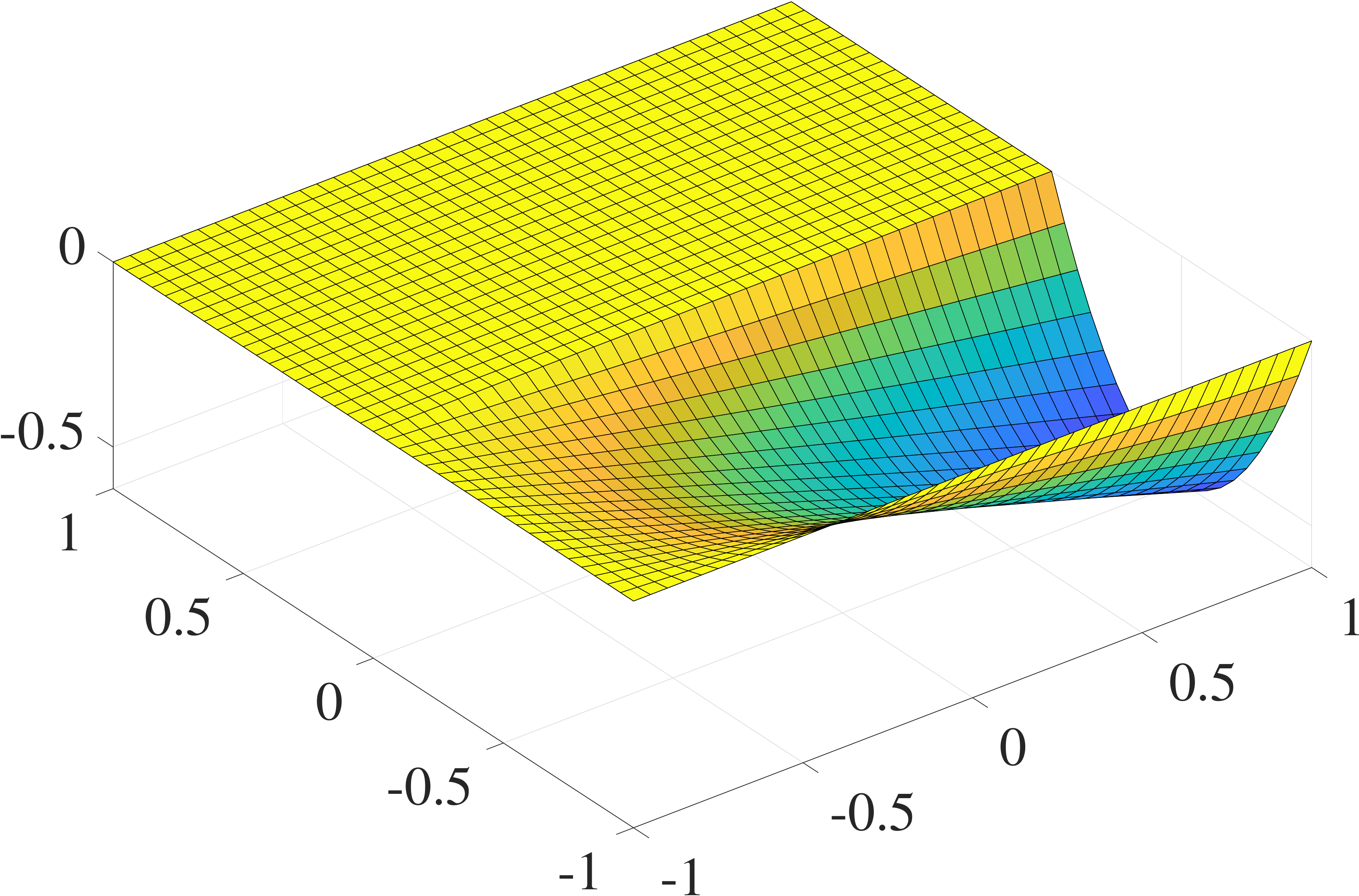}}\\
	\subfloat[Edge vertex \#9\label{fig:12NodeElem_SF9_LaLe}]{\includegraphics[clip,width=0.24\textwidth]{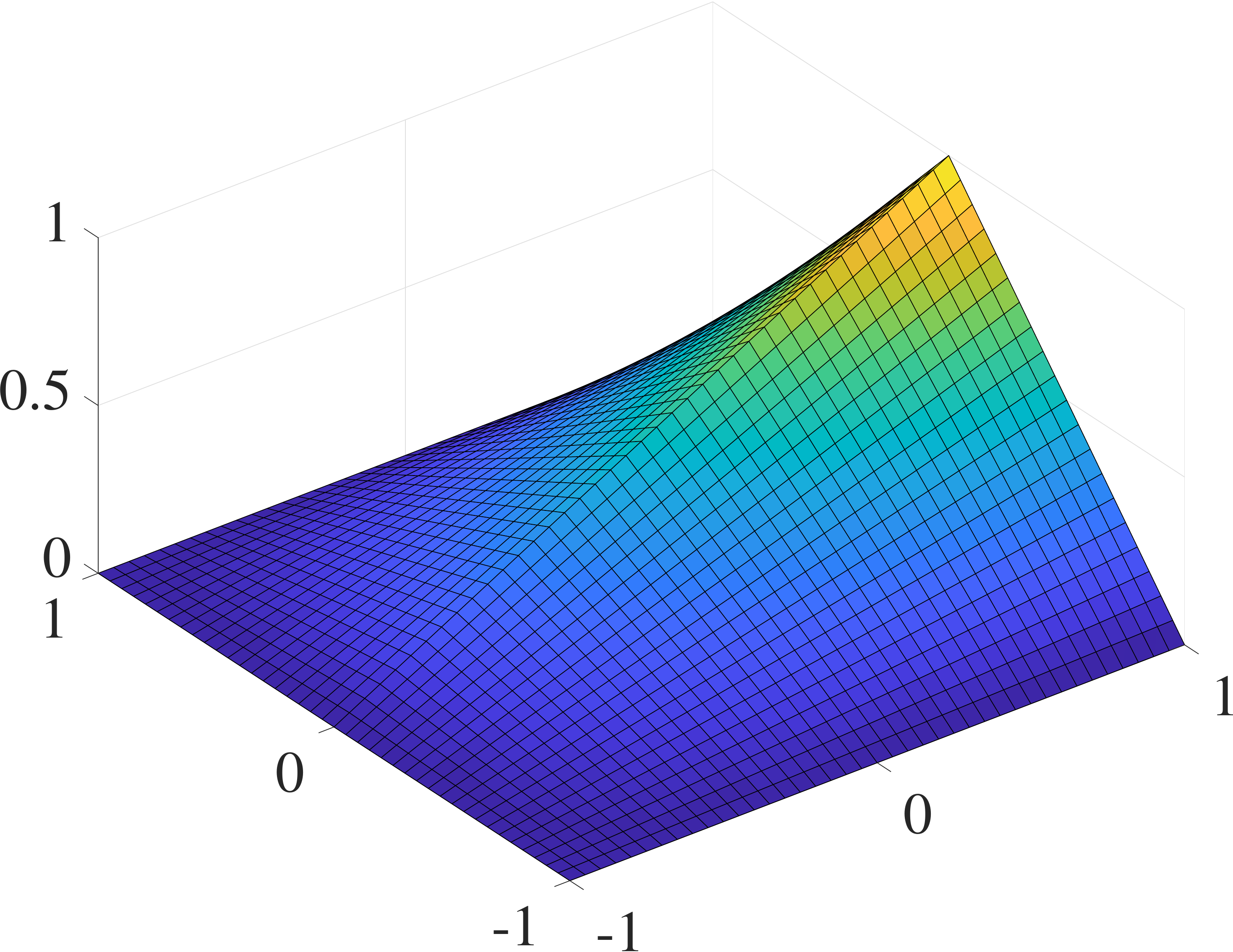}}
	\hfill
	\subfloat[Edge vertex \#10\label{fig:12NodeElem_SF10_LaLe}]{\includegraphics[clip,width=0.24\textwidth]{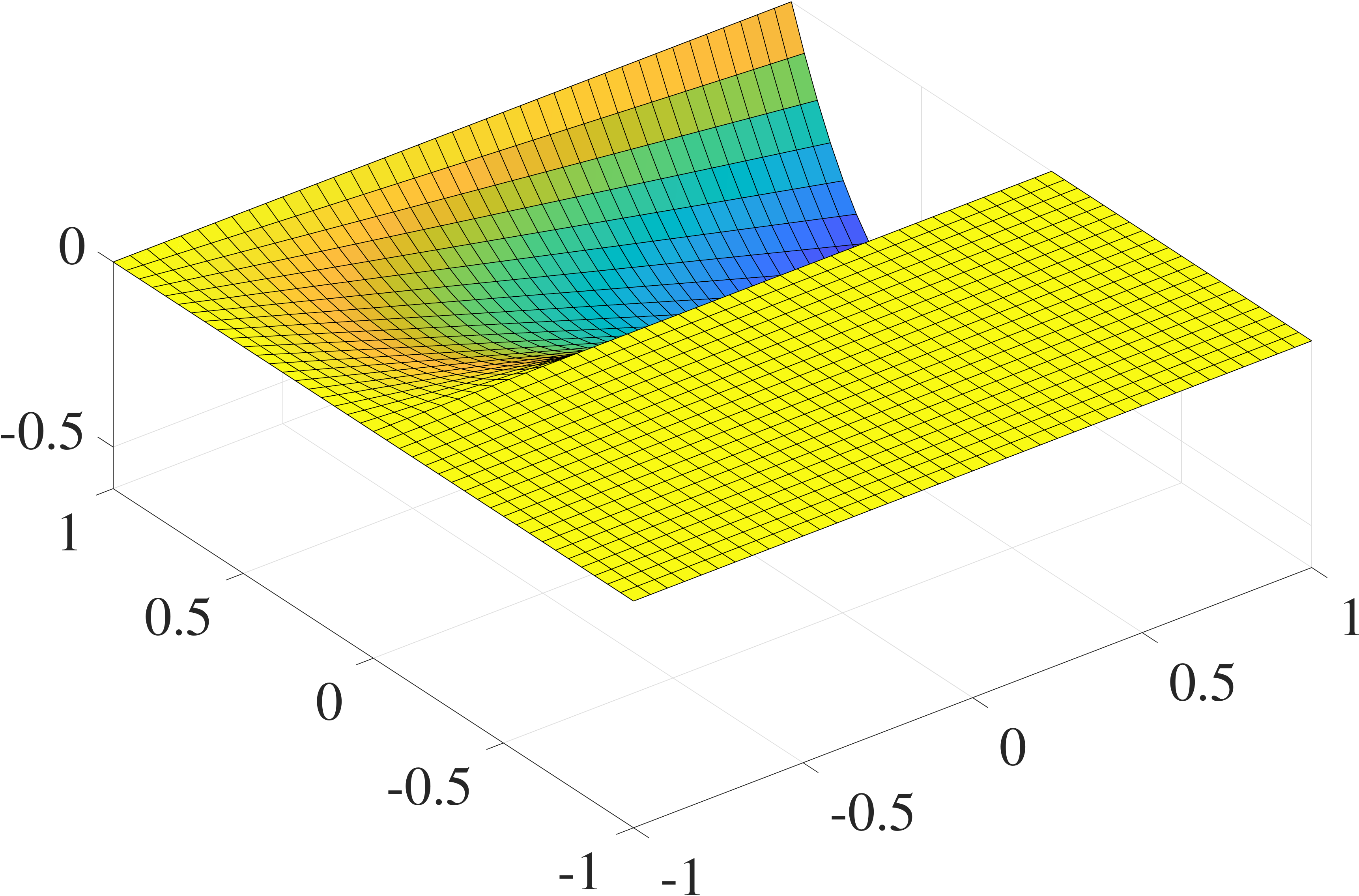}}
	\hfill
	\subfloat[Edge vertex \#11\label{fig:12NodeElem_SF11_LaLe}]{\includegraphics[clip,width=0.24\textwidth]{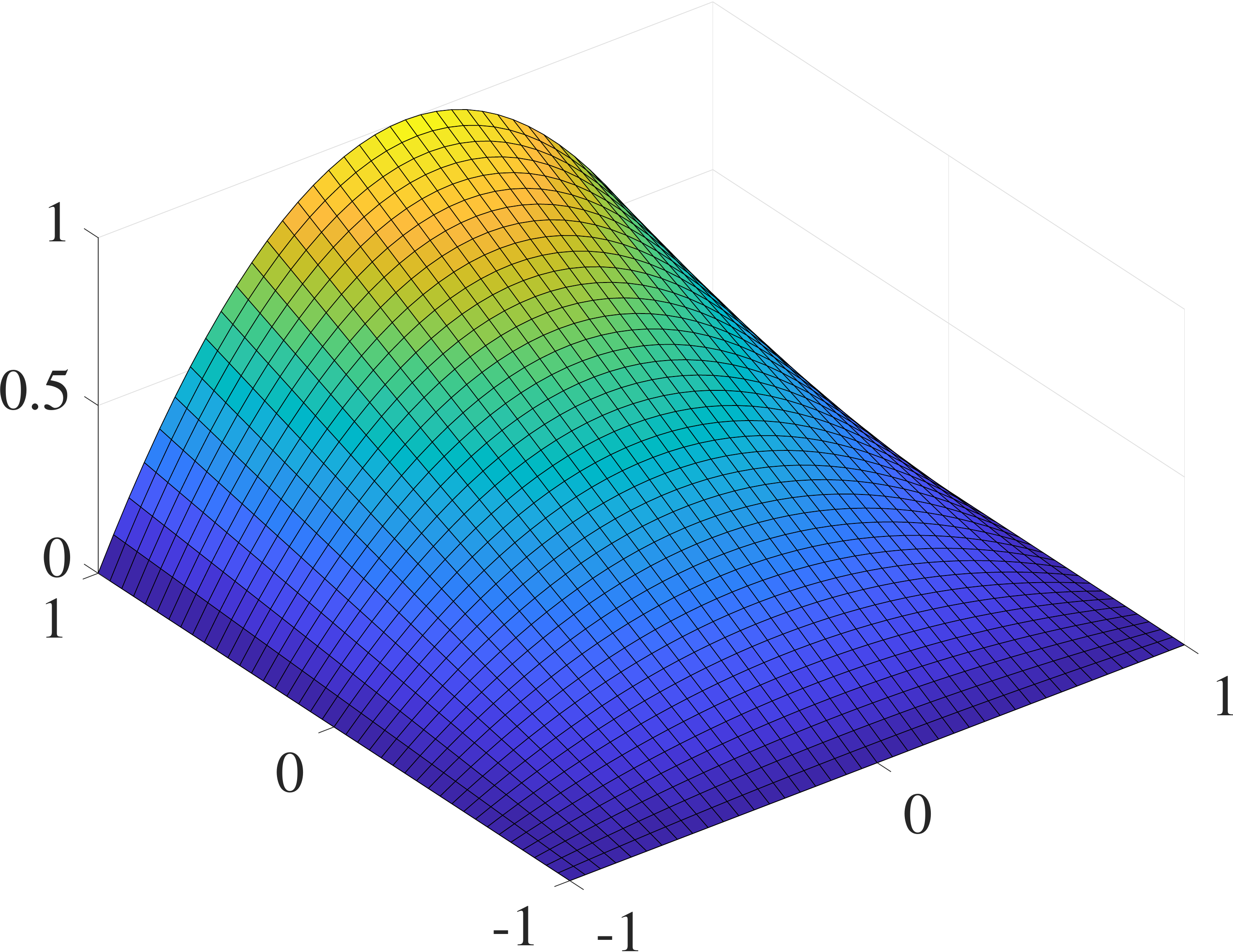}}
	\hfill
	\subfloat[Edge vertex \#12\label{fig:12NodeElem_SF12_LaLe}]{\includegraphics[clip,width=0.24\textwidth]{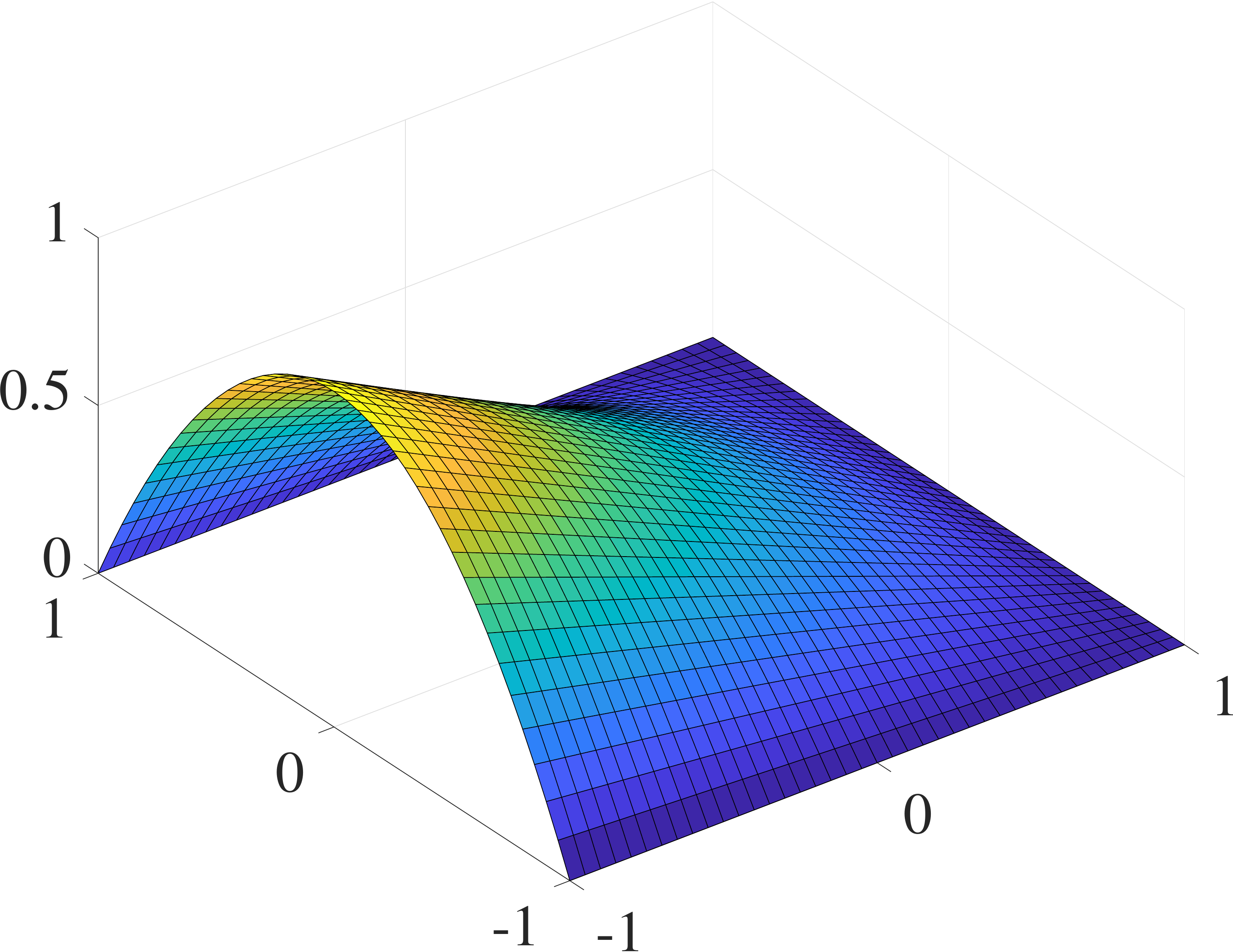}}\\
	\caption{Shape functions of the 12-node bi-quadratic transition element -- Lagrange- to Legendre-based shape functions.}
	\label{fig:12NodeElem_SF_LaLe}
\end{figure}%
All the shape functions of the derived transition element are illustrated in Fig.~\ref{fig:12NodeElem_SF_LaLe}. Here, we clearly observe the different sets of shape functions that have been employed for edges $E_1$ and $E_2$ or $E_3$ and $E_4$, respectively. Assuming that the local coordinate systems of the \emph{p}-edges are aligned this transition element is a means to combine different element types in one model.
\subsection{Legendre-Legendre: Modal shape functions}
The purpose of this section is to illustrate the coupling of hierarchic elements for local mesh refinement. The basis functions for the derivation are given by Eqs.~\eqref{eq:1DQuadShapeFunc_Le} and \eqref{eq:QuadShapeFunc_Trunk}. For the numbering of the nodes please see Fig.~\ref{fig:8NodeElem}. The interpolations along the four edges of the transition element are given as
\allowdisplaybreaks
\begin{align}
\text{E}_1:\; \bar{N}_i(\xi) &=
\begin{cases}
^\mathrm{Le}N^\mathrm{2}_{1}(\check{\xi}_1) & \text{for } i = 1, \text{ node \#1} \,,\\
^\mathrm{Le}N^\mathrm{2}_{2}(\check{\xi}_1) & \text{for } i = 2, \text{ node \#5} \,,\\
\begin{cases}
^\mathrm{Le}N^\mathrm{2}_{3}(\check{\xi}_1) \\
^\mathrm{Le}N^\mathrm{2}_{1}(\check{\xi}_2) \\
\end{cases}
& \text{for } i = 3, \text{ node \#6} \,, \\
^\mathrm{Le}N^\mathrm{2}_{2}(\check{\xi}_2) & \text{for } i = 4, \text{ node \#7} \,, \\
^\mathrm{Le}N^\mathrm{2}_{3}(\check{\xi}_2) & \text{for } i = 5, \text{ node \#2} \,,
\end{cases}
\\
\text{E}_2:\; \hat{N}_i(\eta) &=
\begin{cases}
^\mathrm{Le}N^\mathrm{2}_{1}(\check{\eta}_1) & \text{for } i = 1, \text{ node \#2} \,,\\
^\mathrm{Le}N^\mathrm{2}_{2}(\check{\eta}_1) & \text{for } i = 2, \text{ node \#8} \,,\\
\begin{cases}
^\mathrm{Le}N^\mathrm{2}_{3}(\check{\eta}_1) \\
^\mathrm{Le}N^\mathrm{2}_{1}(\check{\eta}_2) \\
\end{cases}
& \text{for } i = 3, \text{ node \#9} \,,\\
^\mathrm{Le}N^\mathrm{2}_{2}(\check{\eta}_2) & \text{for } i = 4, \text{ node \#10} \,, \\
^\mathrm{Le}N^\mathrm{2}_{3}(\check{\eta}_2) & \text{for } i = 5, \text{ node \#3} \,,
\end{cases}
\\
\text{E}_3:\; \tilde{N}_i(\xi) &=
\begin{cases}
^\mathrm{Le}N^\mathrm{2}_{1}(\xi) & \text{for } i = 1, \text{ node \#4} \,,\\
^\mathrm{Le}N^\mathrm{2}_{2}(\xi) & \text{for } i = 2, \text{ node \#11} \,,\\
^\mathrm{Le}N^\mathrm{2}_{3}(\xi) & \text{for } i = 3, \text{ node \#3} \,,
\end{cases}
\\
\text{E}_4:\; \breve{N}_i(\eta) &=
\begin{cases}
^\mathrm{Le}N^\mathrm{2}_{1}(\eta) & \text{for } i = 1, \text{ node \#1} \,,\\
^\mathrm{Le}N^\mathrm{2}_{2}(\eta) & \text{for } i = 2, \text{ node \#12} \,,\\
^\mathrm{Le}N^\mathrm{2}_{3}(\eta) & \text{for } i = 3, \text{ node \#4} \,.
\end{cases}
\end{align}
\begin{figure}[b!]
	\begin{minipage}[t]{0.49\textwidth}
		\centering
		\includegraphics[clip,width=1\textwidth]{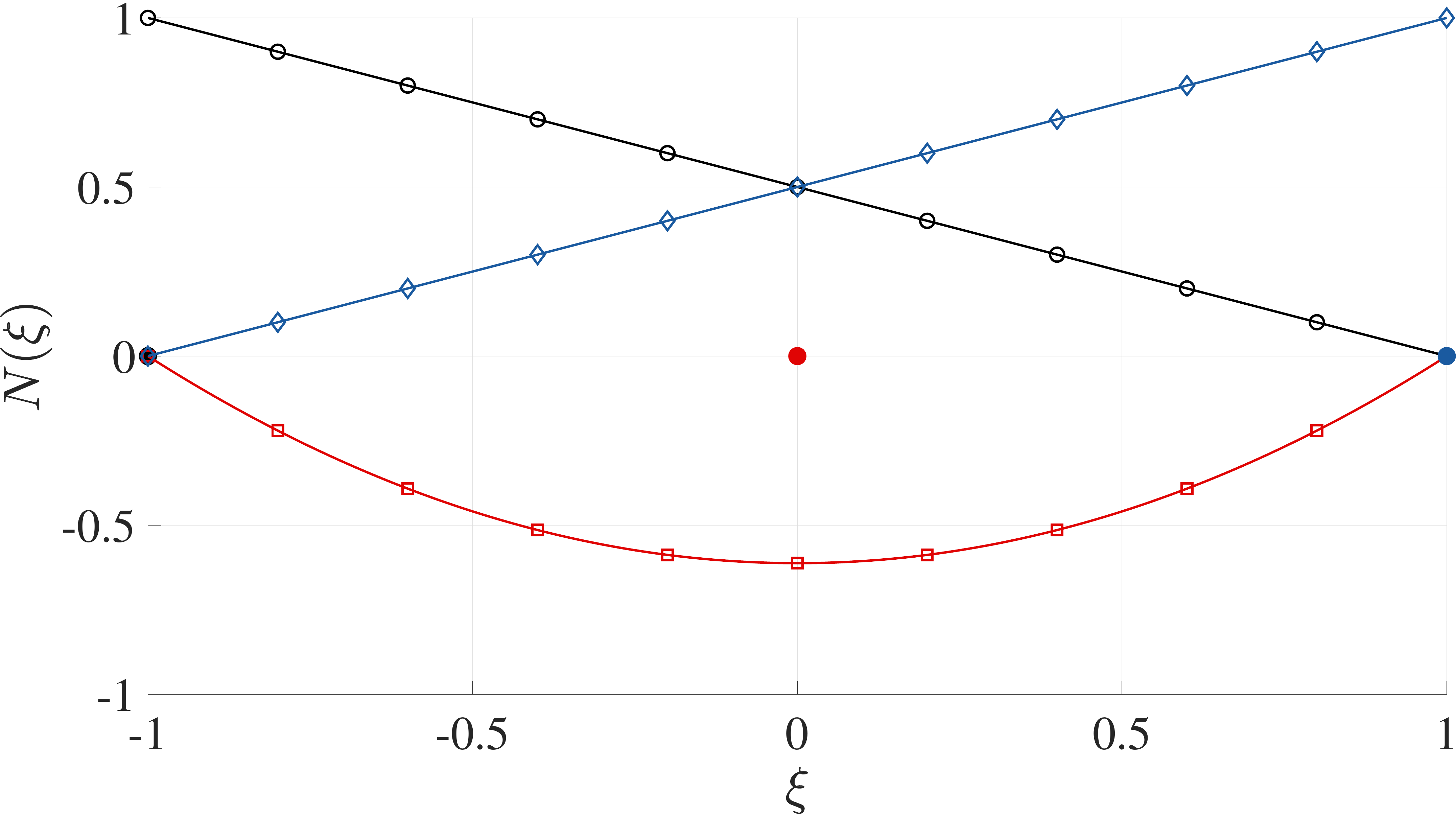}
		\caption{One-dimensional quadratic shape functions (3-mode hierarchic element). Legend: \textcolor{Matlab1}{\rule[0.55ex]{5ex}{0.2ex}} $N_1(\xi)$ (\NewCirc), \textcolor{Matlab2}{\rule[0.55ex]{5ex}{0.2ex}} $N_2(\xi)$ ($\square$), \textcolor{Matlab3}{\rule[0.55ex]{5ex}{0.2ex}} $N_3(\xi)$ ($\Diamond$).}
		\label{fig:ShapeFuncQuad1d_LeLe}
	\end{minipage}
	\hfill
	\begin{minipage}[t]{0.49\textwidth}
		\centering
		\includegraphics[clip,width=1\textwidth]{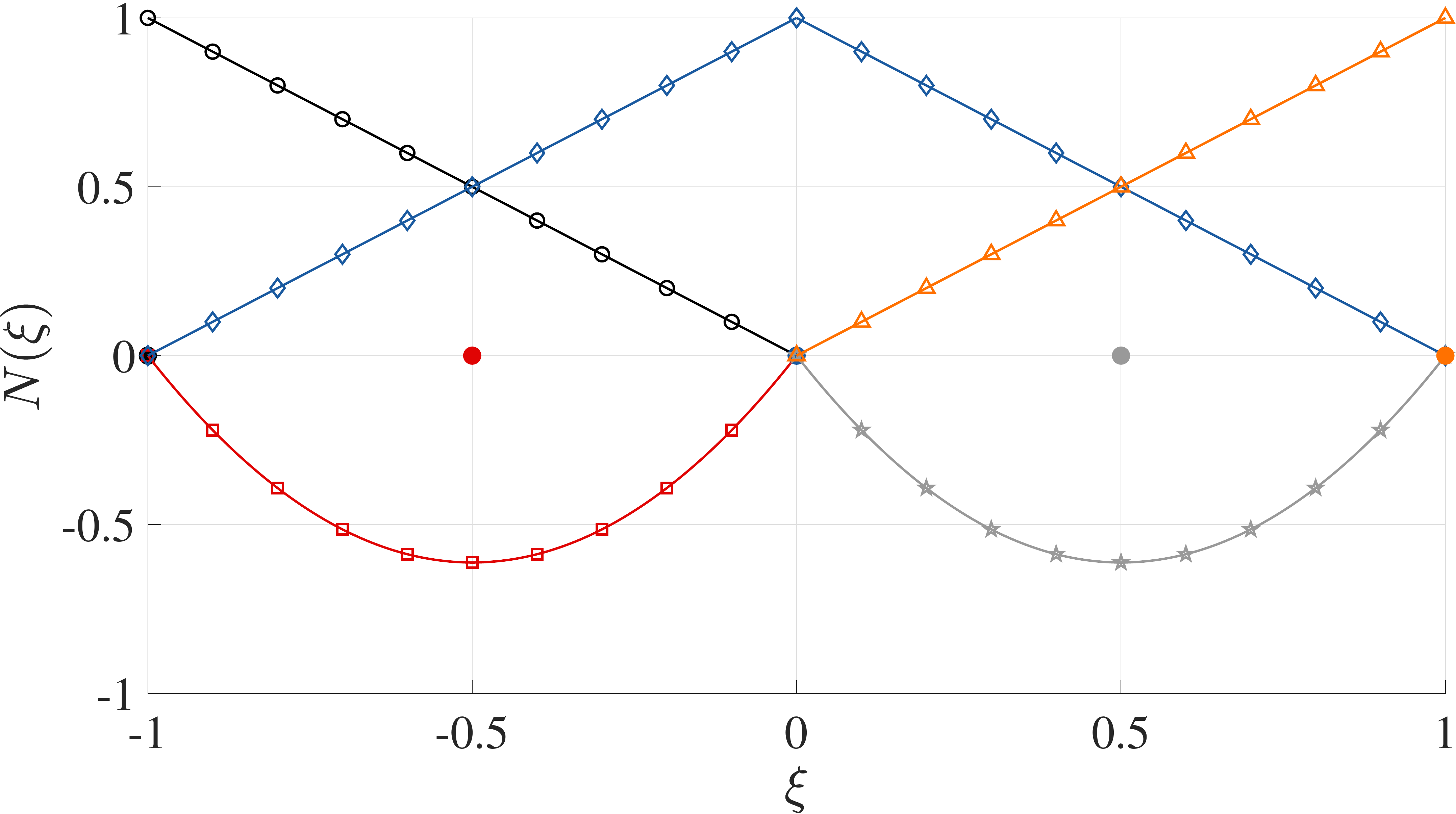}
		\caption{One-dimensional piecewise quadratic shape functions (based on 3-mode hierarchic elements). Legend: \textcolor{Matlab1}{\rule[0.55ex]{5ex}{0.2ex}} $N_1(\xi)$ (\NewCirc), \textcolor{Matlab2}{\rule[0.55ex]{5ex}{0.2ex}} $N_2(\xi)$ ($\square$), \textcolor{Matlab3}{\rule[0.55ex]{5ex}{0.2ex}} $N_3(\xi)$ ($\Diamond$), \textcolor{Matlab4}{\rule[0.55ex]{5ex}{0.2ex}} $N_4(\xi)$ (\NewStar), \textcolor{Matlab5}{\rule[0.55ex]{5ex}{0.2ex}} $N_5(\xi)$ ($\bigtriangleup$).}
		\label{fig:ShapeFuncPiecewiseQuad1d_LeLe}
	\end{minipage}
\end{figure}%
The shape functions for the edges are exemplarily sketched in Figs.~\ref{fig:ShapeFuncQuad1d_LeLe} and \ref{fig:ShapeFuncPiecewiseQuad1d_LeLe}. The functions living on the edges of the transition element are in the next step projected/blended into the interior of the element by using the projection operators $\mathcal{P}_{\xi}[\square]$, $\mathcal{P}_{\eta}[\square]$, and $\mathcal{P}_{\xi}[\mathcal{P}_{\eta}[\square]]$. First, we will derive the shape functions that are associated to the corner vertices:
\begin{alignat}{3}
N_1(\xi,\eta) &= N_{1}^{1}(\xi) \,^\mathrm{Le}N^\mathrm{2}_{1}(\eta) &&+ N_{1}^{1}(\eta) \,^\mathrm{Le}N^\mathrm{2}_{1}(\check{\xi}_1) &&- N_{1}^{1}(\xi)N_{1}^{1}(\eta)\,, \\
N_2(\xi,\eta) &= N_{2}^{1}(\xi) \,^\mathrm{Le}N^\mathrm{2}_{1}(\check{\eta}_1) &&+ N_{1}^{1}(\eta) \,^\mathrm{Le}N^\mathrm{2}_{3}(\check{\xi}_2) &&- N_{2}^{1}(\xi)N_{1}^{1}(\eta)\,, \\
N_3(\xi,\eta) &= N_{2}^{1}(\xi) \,^\mathrm{Le}N^\mathrm{2}_{3}(\check{\eta}_2) &&+ N_{2}^{1}(\eta) \,^\mathrm{Le}N^\mathrm{2}_{3}(\xi) &&- N_{2}^{1}(\xi)N_{2}^{1}(\eta)\,, \\
N_4(\xi,\eta) &= N_{1}^{1}(\xi) \,^\mathrm{Le}N^\mathrm{2}_{3}(\eta) &&+ N_{2}^{1}(\eta) \,^\mathrm{Le}N^\mathrm{2}_{1}(\xi) &&- N_{1}^{1}(\xi)N_{2}^{1}(\eta) = \!^\mathrm{t}N^{*}_{4}(\xi,\eta)  \,.
\end{alignat}
Second, the edge shape functions are listed:
\begin{alignat}{2}
&N_5(\xi,\eta)    &&= N_{1}^{1}(\eta) \,^\mathrm{Le}N^\mathrm{2}_{2}(\check{\xi}_1) = \!^\mathrm{t}N^{*}_{6}(\check{\xi}_1,\eta) \,, \\
&N_6(\xi,\eta)    &&= N_{1}^{1}(\eta) 
\begin{cases}
^\mathrm{Le}N^\mathrm{2}_{3}(\check{\xi}_1)\,, \\
^\mathrm{Le}N^\mathrm{2}_{1}(\check{\xi}_2)\,, \\
\end{cases} \\
&N_7(\xi,\eta)    &&= N_{1}^{1}(\eta) \,^\mathrm{Le}N^\mathrm{2}_{2}(\check{\xi}_2) = \!^\mathrm{t}N^{*}_{6}(\check{\xi}_2,\eta) \,, \\
&N_8(\xi,\eta)    &&= N_{2}^{1}(\xi) \,^\mathrm{Le}N^\mathrm{2}_{2}(\check{\eta}_1) = \!^\mathrm{t}N^{*}_{9}(\xi,\check{\eta}_1) \,, \\
&N_9(\xi,\eta)    &&= N_{2}^{1}(\xi)
\begin{cases}
^\mathrm{Le}N^\mathrm{2}_{3}(\check{\eta}_1)\,, \\
^\mathrm{Le}N^\mathrm{2}_{1}(\check{\eta}_2)\,, \\
\end{cases} \\
&N_{10}(\xi,\eta) &&= N_{2}^{1}(\xi) \,^\mathrm{Le}N^\mathrm{2}_{2}(\check{\eta}_2) = \!^\mathrm{t}N^{*}_{9}(\xi,\check{\eta}_2) \,, \\
&N_{11}(\xi,\eta) &&= N_{2}^{1}(\eta) \,^\mathrm{Le}N^\mathrm{2}_{2}(\xi) = \!^\mathrm{t}N^{*}_{11}(\xi,\eta) \,, \\
&N_{12}(\xi,\eta) &&= N_{1}^{1}(\xi) \,^\mathrm{Le}N^\mathrm{2}_{2}(\eta) = \!^\mathrm{t}N^{*}_{12}(\xi,\eta) \,.
\end{alignat}
\begin{figure}[htp!]
	\centering
	\subfloat[Corner vertex \#1\label{fig:12NodeElem_SF1_LeLe}]{\includegraphics[clip,width=0.24\textwidth]{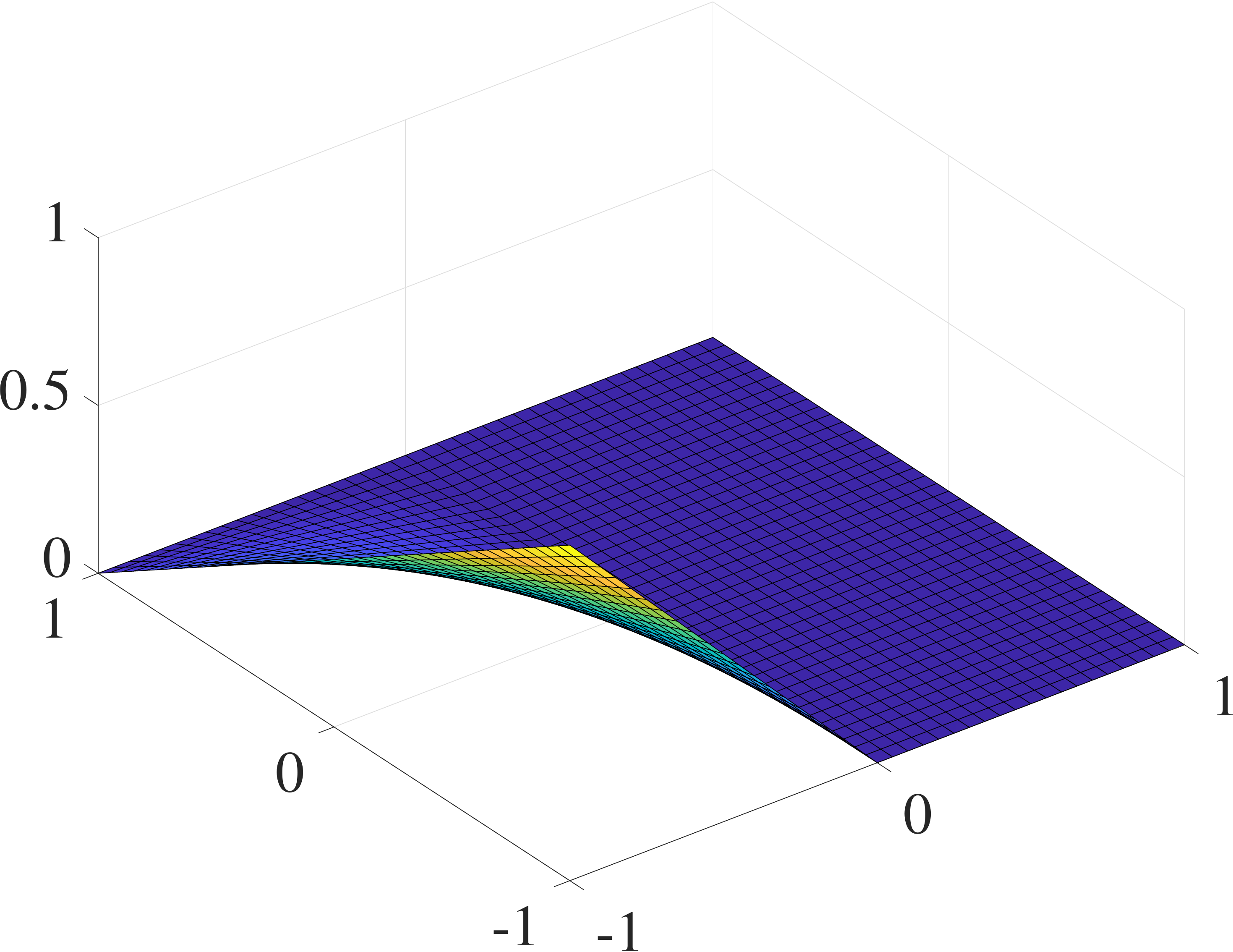}}
	\hfill
	\subfloat[Corner vertex \#2\label{fig:12NodeElem_SF2_LeLe}]{\includegraphics[clip,width=0.24\textwidth]{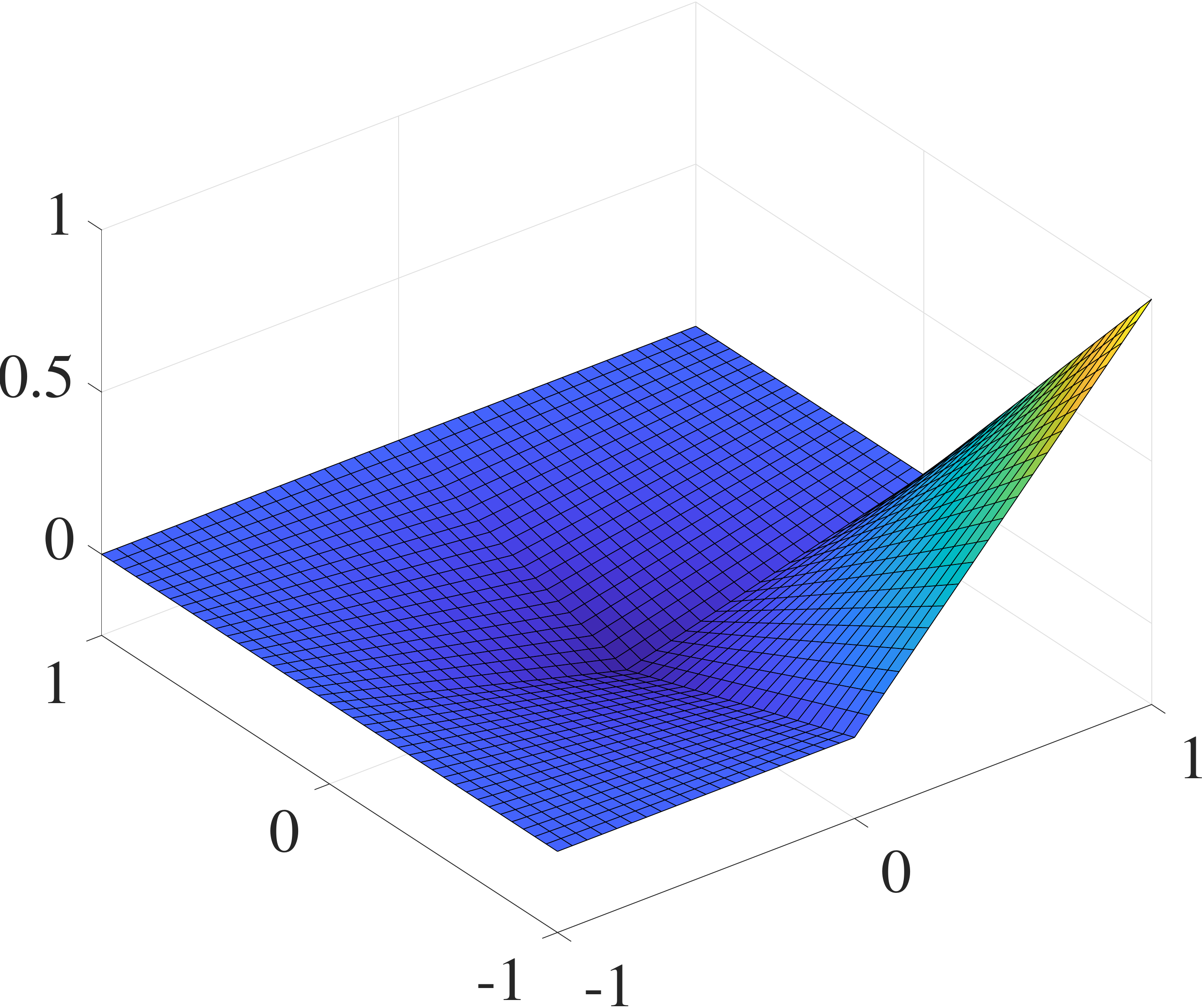}}
	\hfill
	\subfloat[Corner vertex \#3\label{fig:12NodeElem_SF3_LeLe}]{\includegraphics[clip,width=0.24\textwidth]{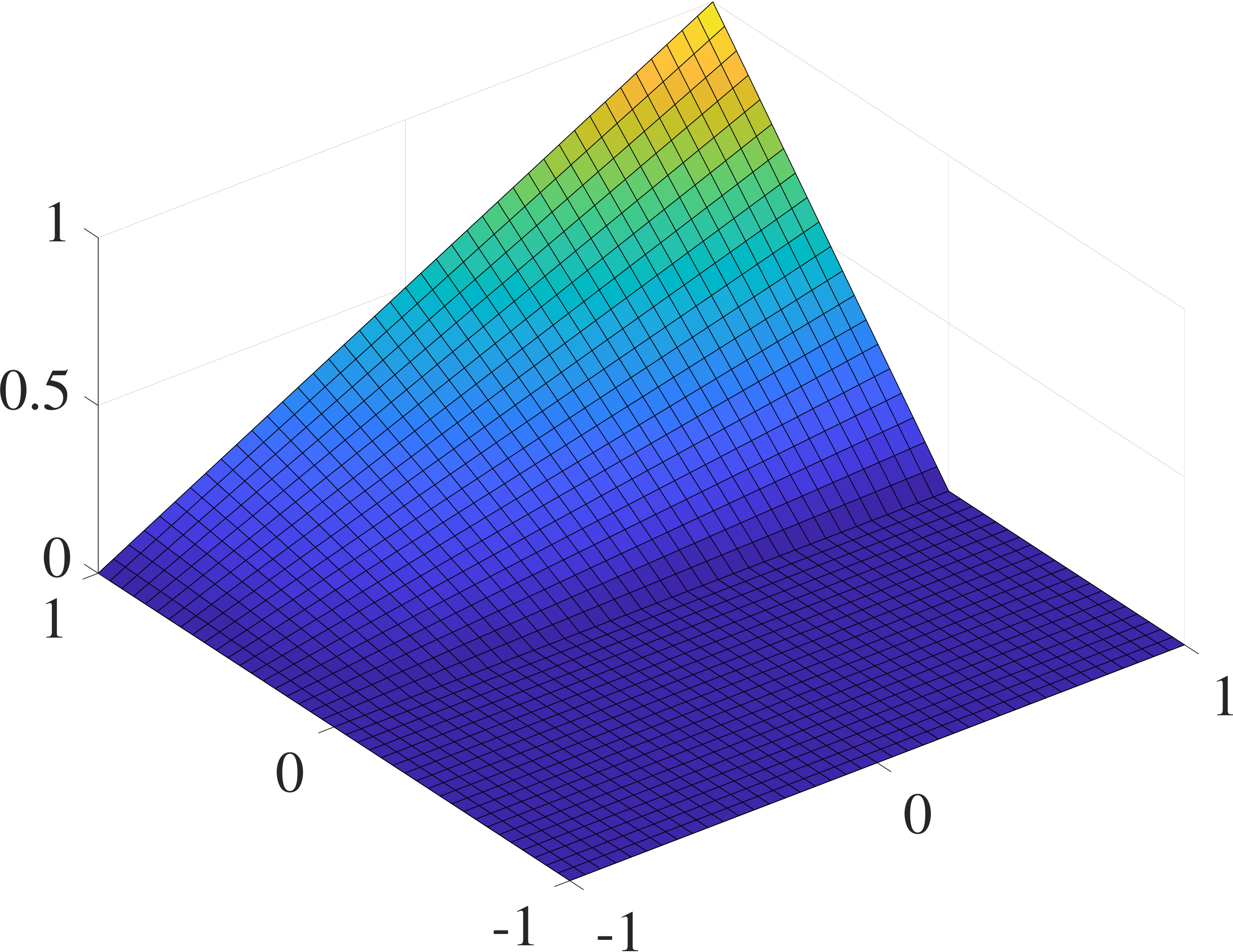}}
	\hfill
	\subfloat[Corner vertex \#4\label{fig:12NodeElem_SF4_LeLe}]{\includegraphics[clip,width=0.24\textwidth]{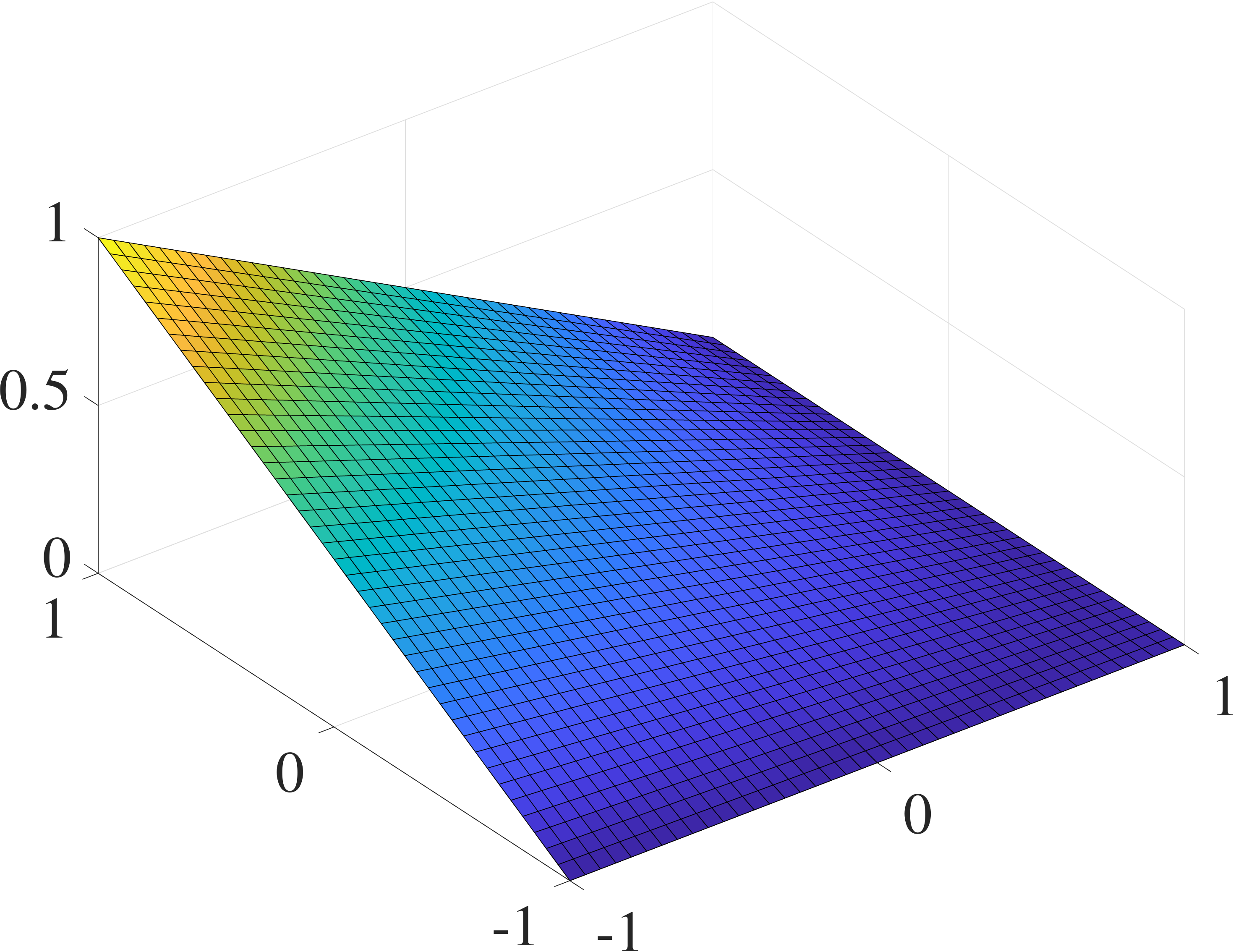}}\\
	\subfloat[Edge vertex \#5\label{fig:12NodeElem_SF5_LeLe}]{\includegraphics[clip,width=0.24\textwidth]{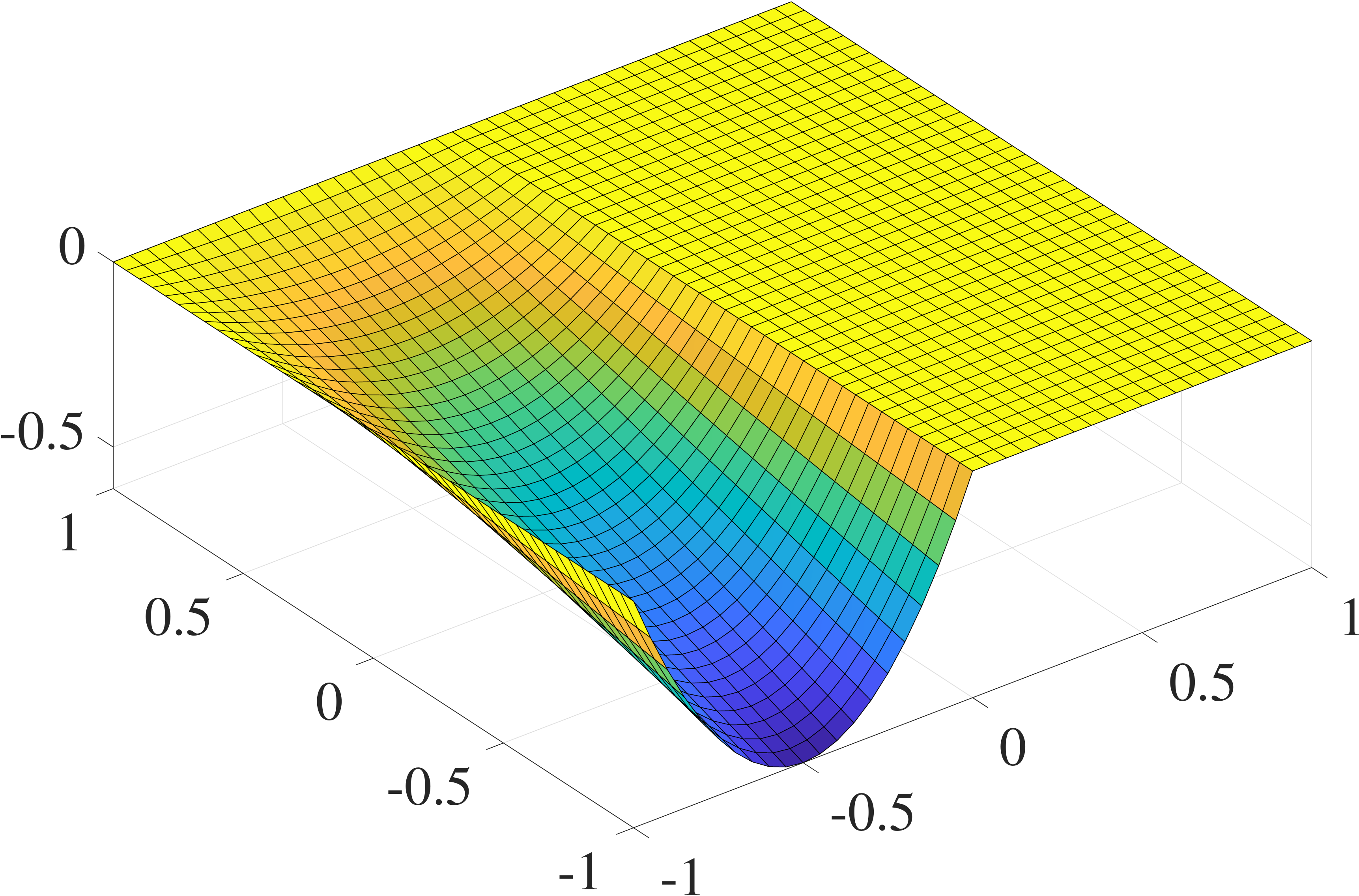}}
	\hfill
	\subfloat[Edge vertex \#6\label{fig:12NodeElem_SF6_LeLe}]{\includegraphics[clip,width=0.24\textwidth]{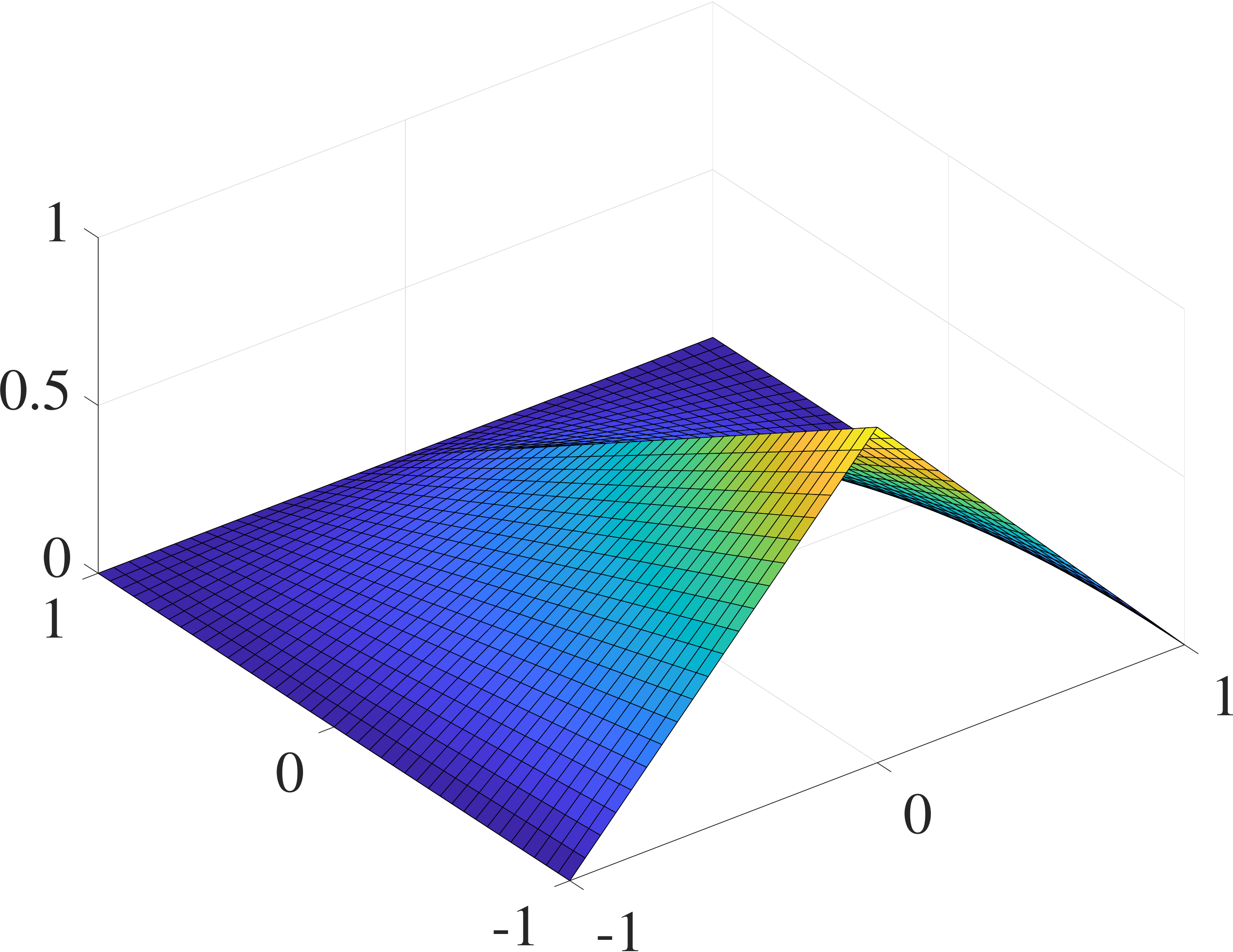}}
	\hfill
	\subfloat[Edge vertex \#7\label{fig:12NodeElem_SF7_LeLe}]{\includegraphics[clip,width=0.24\textwidth]{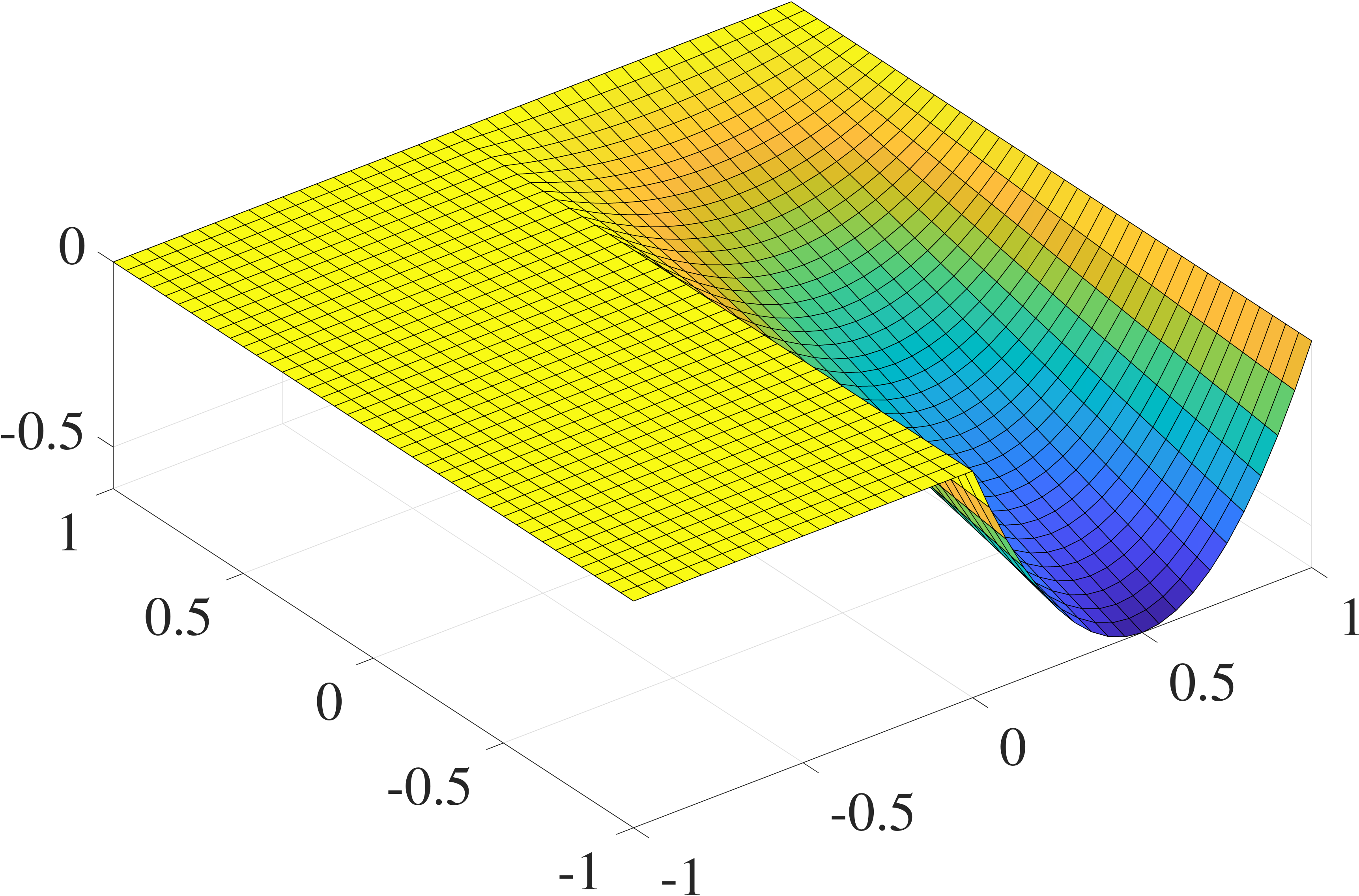}}
	\hfill
	\subfloat[Edge vertex \#8\label{fig:12NodeElem_SF8_LeLe}]{\includegraphics[clip,width=0.24\textwidth]{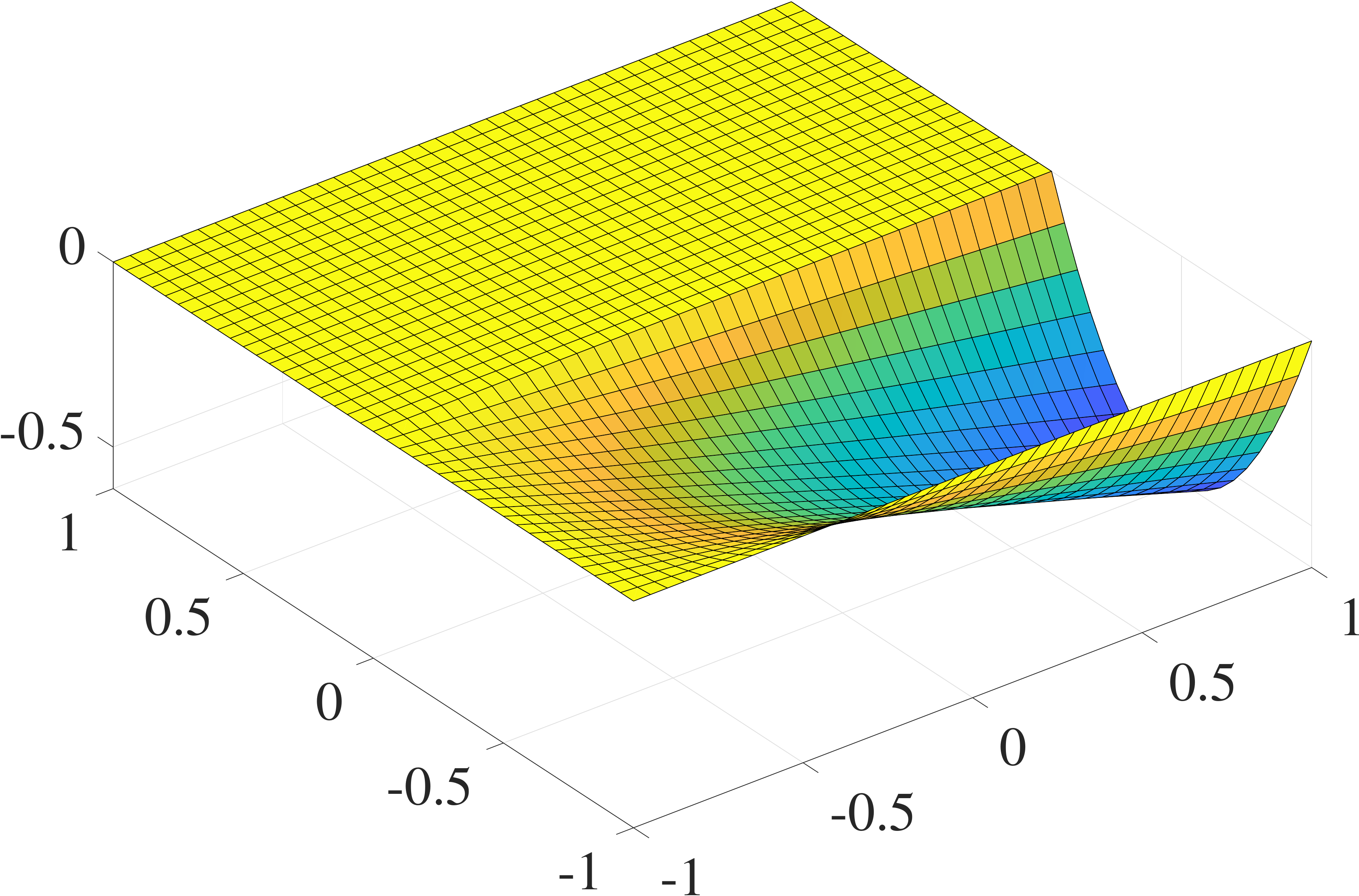}}\\
	\subfloat[Edge vertex \#9\label{fig:12NodeElem_SF9_LeLe}]{\includegraphics[clip,width=0.24\textwidth]{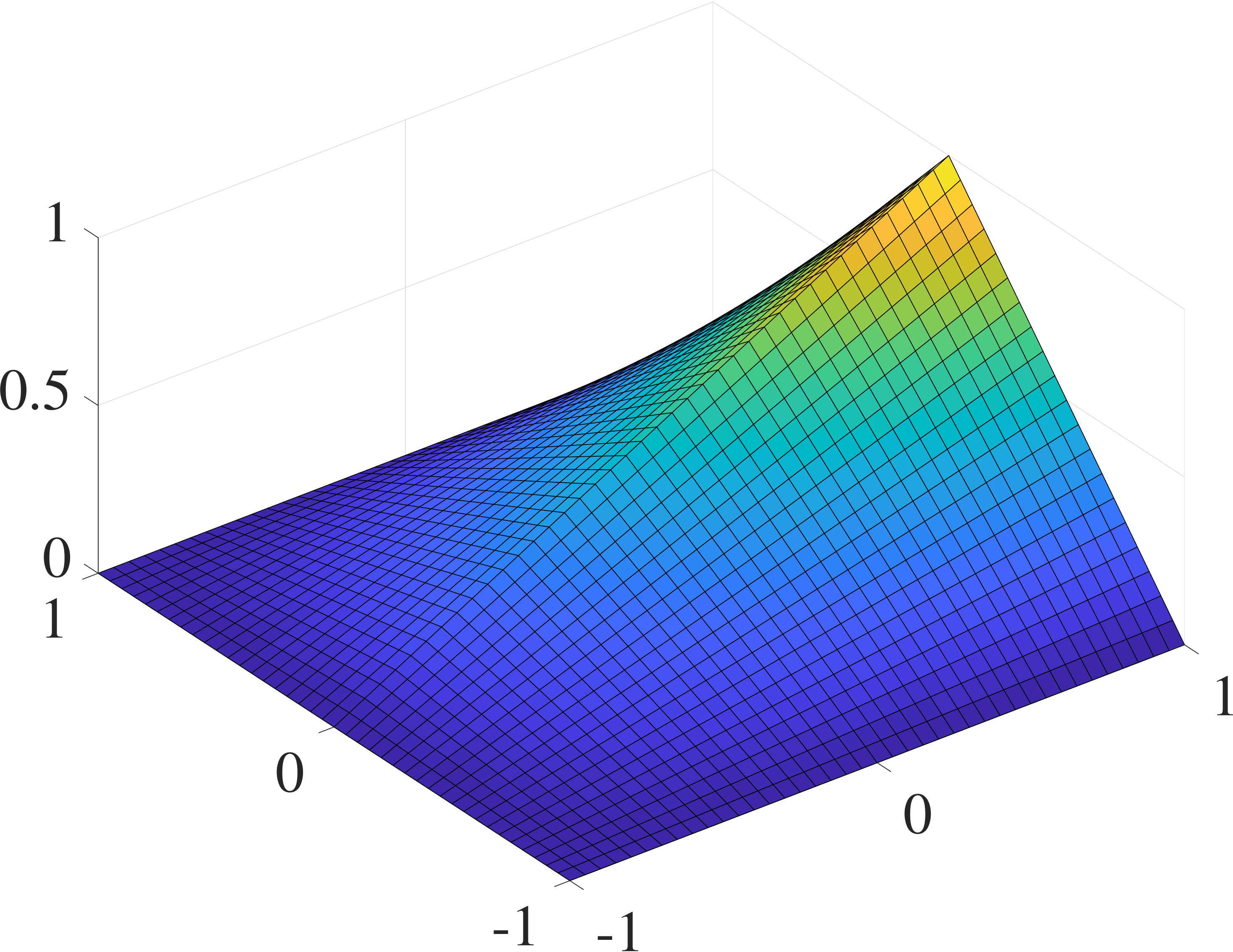}}
	\hfill
	\subfloat[Edge vertex \#10\label{fig:12NodeElem_SF10_LeLe}]{\includegraphics[clip,width=0.24\textwidth]{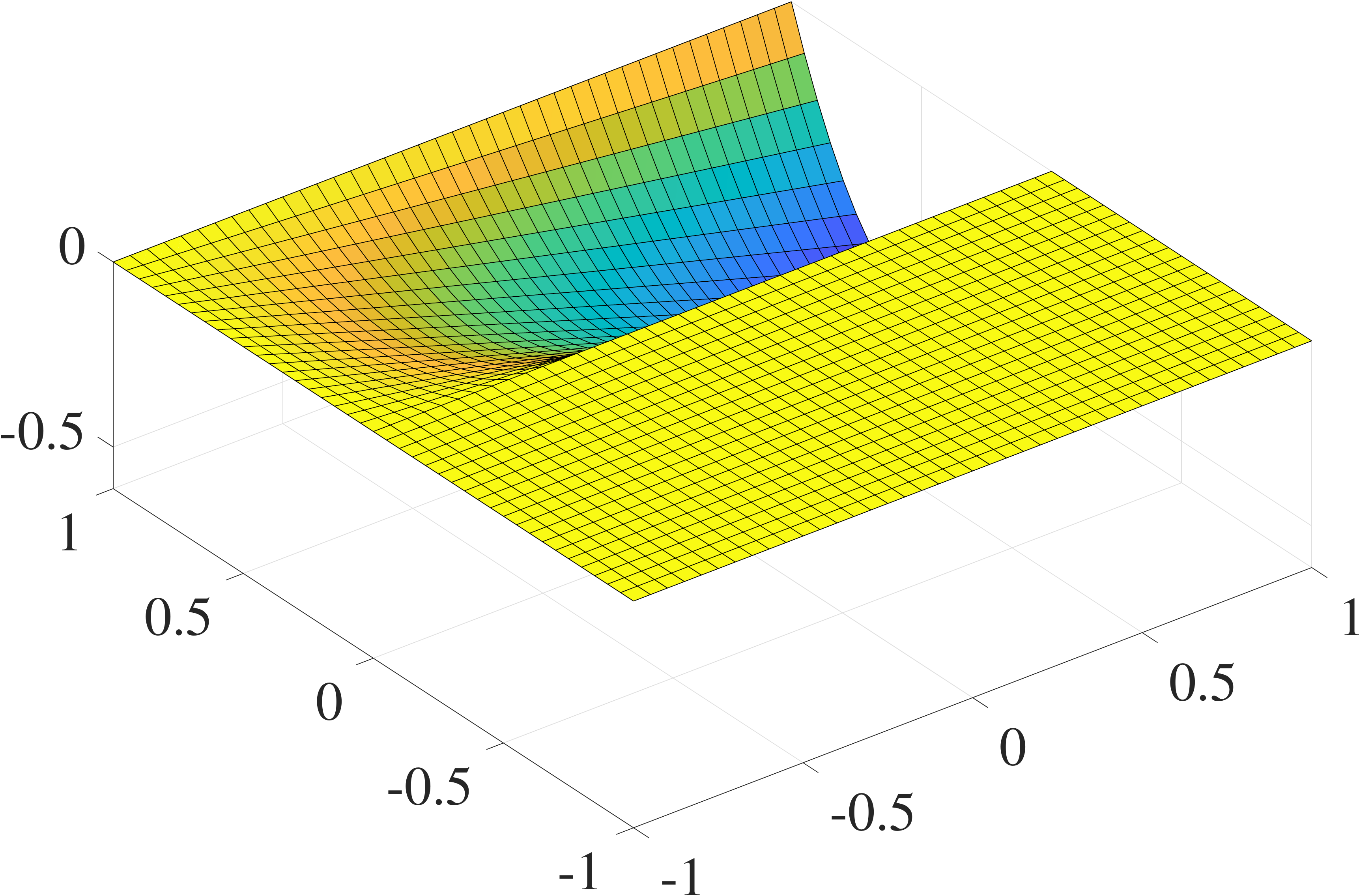}}
	\hfill
	\subfloat[Edge vertex \#11\label{fig:12NodeElem_SF11_LeLe}]{\includegraphics[clip,width=0.24\textwidth]{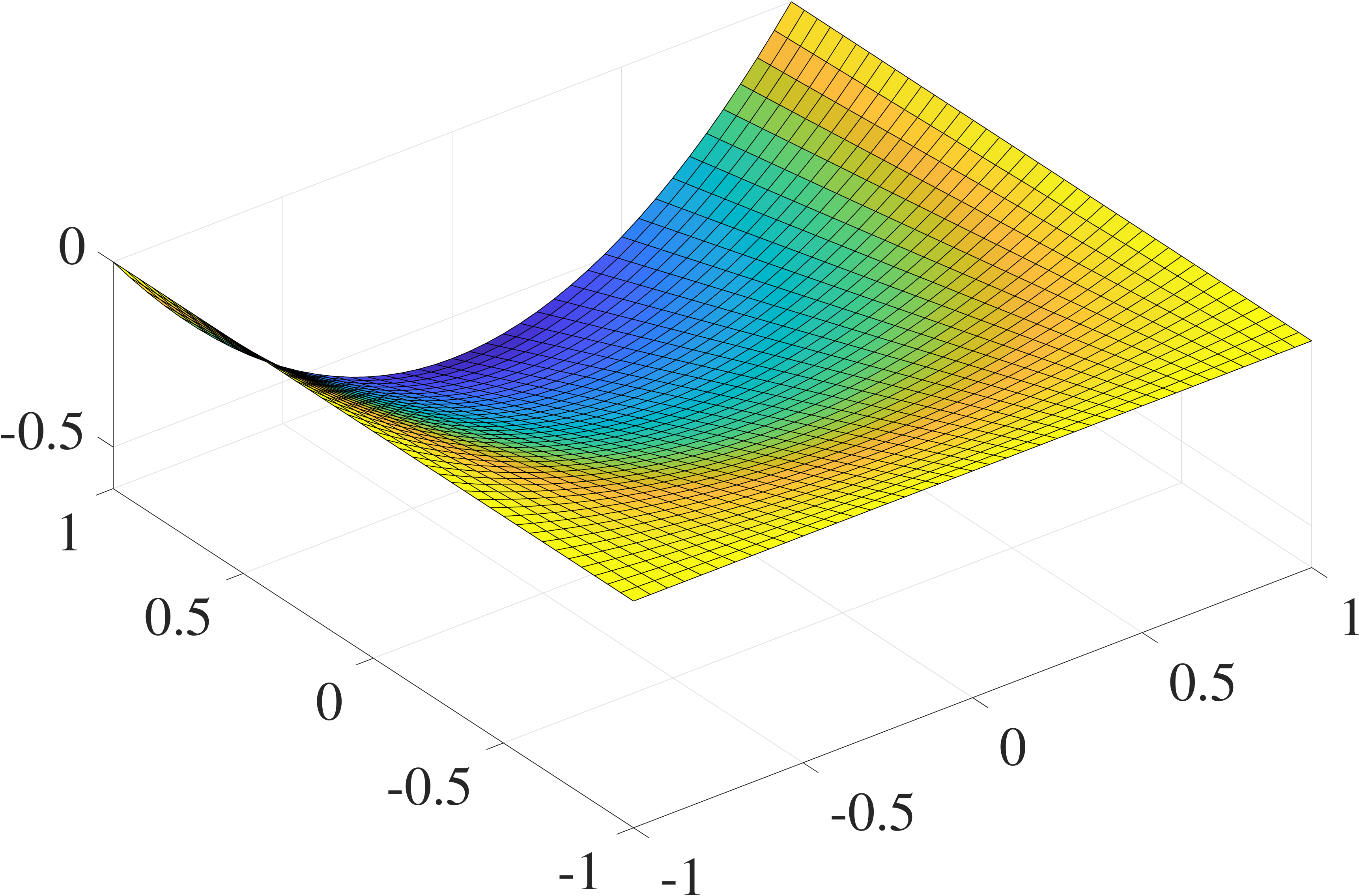}}
	\hfill
	\subfloat[Edge vertex \#12\label{fig:12NodeElem_SF12_LeLe}]{\includegraphics[clip,width=0.24\textwidth]{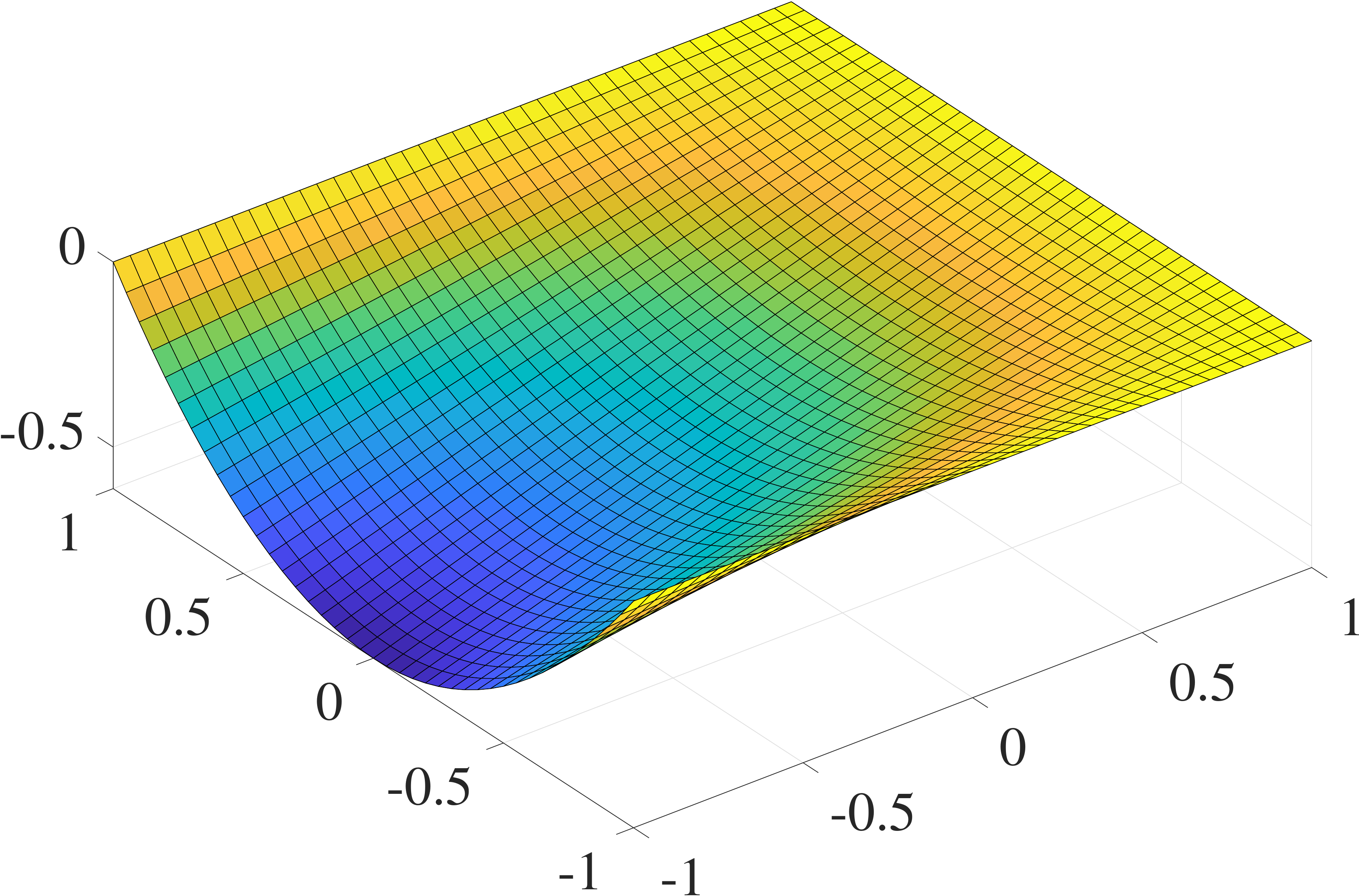}}\\
	\caption{Shape functions of the 12-node bi-quadratic transition element -- Legendre- to Legendre-based shape functions.}
	\label{fig:12NodeElem_SF_LeLe}
\end{figure}%
All the shape functions of the derived transition element are illustrated in Fig.~\ref{fig:12NodeElem_SF_LeLe}.
\section{Displacement field for the higher order patch test}
\label{sec:PatchTestHigh_Disp}
\begin{figure}[b!]
	\centering
	\includegraphics[clip,scale=0.75]{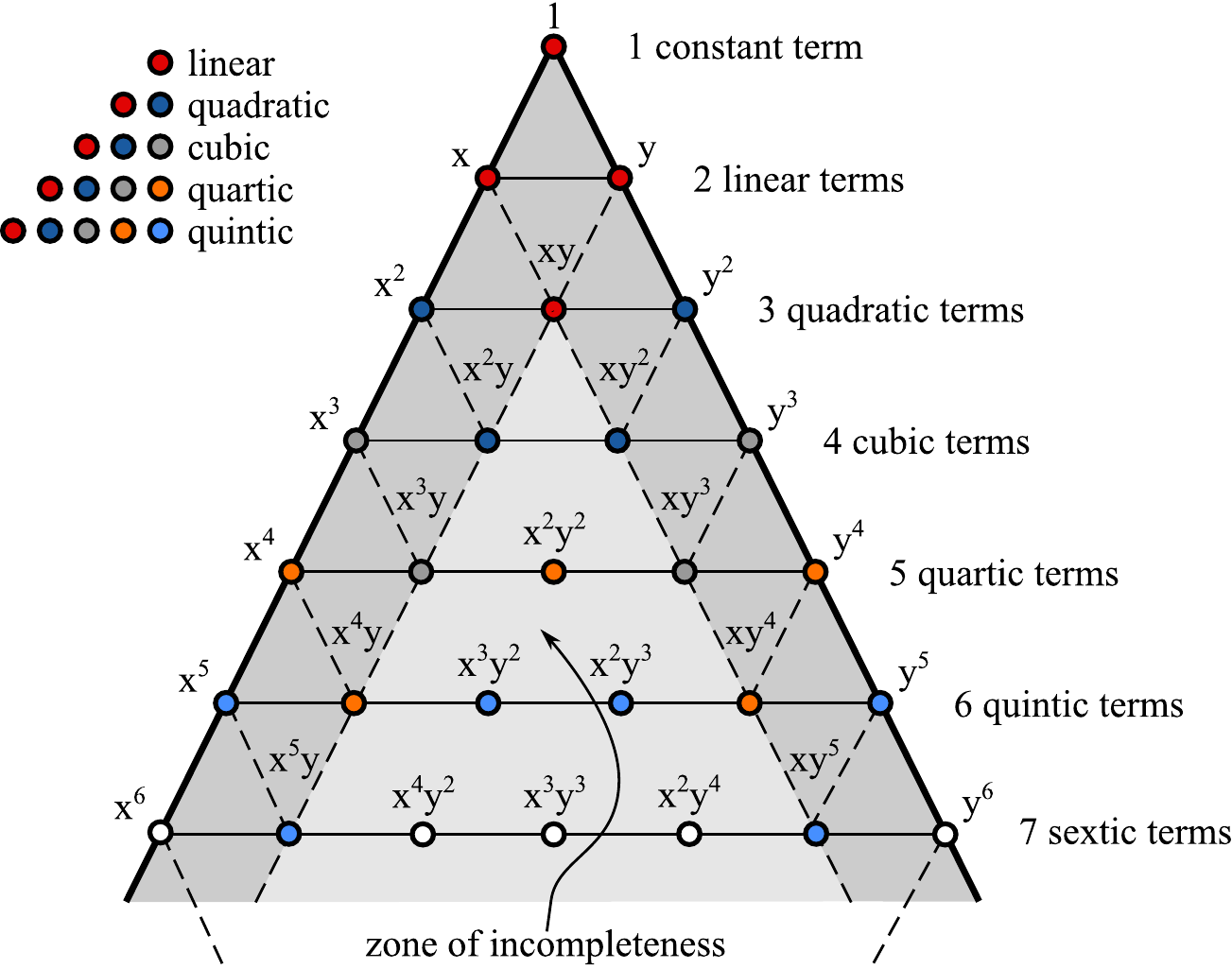}
	\caption{Pascal's triangle}
	\label{fig:PascalTri_app}
\end{figure}

In Sect.~\ref{sec:PatchTestHigh}, the higher order patch test was used to assess the convergence properties of the proposed transition elements. In this section, we provide a detailed explanation of how an admissible higher order displacement field is determined. The number of monomials that are required for a complete polynomial can be read from Pascal's triangle depicted in Fig.~\ref{fig:PascalTri_app}. For the higher order patch test, we prescribe a displacement field of the form
\begin{equation}
\mathbf{u} = \sum\limits_{i}^{} \mathbf{a}_iP^{p}_i(x,y)\,, \quad \forall\; i \in\{1,2,\ldots,n_\mathrm{p}\}
\label{eq:HighOrderDisplacementField_app}
\end{equation}
at the boundary of the computational domain, where $\mathbf{a}_i$ denotes unknown parameters of the polynomial functions and $n_\mathrm{p}$ represents the number of monomials in the extension. Note that the high order displacement field cannot be chosen arbitrarily but needs to fulfill the equilibrium conditions in order to be admissible for the intended purpose, i.e., the coefficients $\mathbf{a}_i$ in Eq.~\eqref{eq:HighOrderDisplacementField_app} have to be determined from
\begin{equation}
\mathrm{div}(\boldsymbol{\sigma}) = \mathbf{0}\,,
\end{equation}
with $\boldsymbol{\sigma}$ denoting the stress tensor. In our derivation, we assume that no body/volume loads are applied. In a two-dimensional framework, this equation simplifies to
\begin{alignat}{2}
\cfrac{\partial \sigma_\mathrm{x}}{\partial x} &+ \cfrac{\partial \tau_\mathrm{xy}}{\partial y} &&= 0\,, \label{eq:EC1_app}\\
\cfrac{\partial \tau_\mathrm{xy}}{\partial x} &+ \cfrac{\partial \sigma_\mathrm{y}}{\partial y} &&= 0\,. \label{eq:EC2_app}
\end{alignat}
Provided that plane stress conditions govern the behavior of the structure, the normal stresses $\sigma_\mathrm{x}$ and $\sigma_\mathrm{y}$ are defined as
\begin{align}
\sigma_\mathrm{x} &= \cfrac{E}{1-\nu^2} \left( \varepsilon_\mathrm{x} + \nu \varepsilon_\mathrm{y} \right)\,,\\
\sigma_\mathrm{y} &= \cfrac{E}{1-\nu^2} \left( \varepsilon_\mathrm{y} + \nu \varepsilon_\mathrm{x} \right)\,,
\end{align}
while the shear stress $\tau_\mathrm{xy}$ is given by
\begin{equation}
\tau_\mathrm{xy} = \cfrac{E}{2(1+\nu)}\,\gamma_\mathrm{xy}\,.
\end{equation}
In the previous equations, an isotropic, linear-elastic material behavior with Young's modulus $E$ and Poisson's ratio $\nu$ is implied. For the sake of completeness, we also provide the definitions of the normal and shear strains
\begin{alignat}{1}
\varepsilon_\mathrm{x} &= \cfrac{\partial u_\mathrm{x}}{\partial x}\,,\\
\varepsilon_\mathrm{y} &= \cfrac{\partial u_\mathrm{y}}{\partial y}\,,\\
\gamma_\mathrm{xy} &= \cfrac{\partial u_\mathrm{x}}{\partial y} + \cfrac{\partial u_\mathrm{y}}{\partial x}\,.
\end{alignat}
To derive an admissible high order displacement field, we have to implement the methodology discussed in the following:
\begin{enumerate}
	\item Set up the displacement field according to Eq.~\eqref{eq:HighOrderDisplacementField_app}:
	\begin{align}
	u_\mathrm{x}(x,y) &= \sum\limits_{i=1}^{n_\mathrm{p}} a_i P^p_i(x,y)\,,\\
	u_\mathrm{y}(x,y) &= \sum\limits_{i=1}^{n_\mathrm{p}} b_i P^p_i(x,y)\,.
	\end{align}
	\item Compute the second derivatives of the displacement field:
	\begin{align}
	\cfrac{\partial^2 u_\mathrm{x}(x,y)}{\partial x_i \partial x_j} &= \sum\limits_{k=4}^{n_\mathrm{p}} a_k \cfrac{\partial^2 P^p_k(x,y)}{\partial x_i \partial x_j}\,,\qquad i,j=1,2\,,\\
	\cfrac{\partial^2 u_\mathrm{y}(x,y)}{\partial x_i \partial x_j} &= \sum\limits_{k=4}^{n_\mathrm{p}} b_k \cfrac{\partial^2 P^p_k(x,y)}{\partial x_i \partial x_j}\,,\qquad i,j=1,2\,.
	\end{align}
%	\begin{align}
%	\cfrac{\partial^2 u_\mathrm{x}(x,y)}{\partial x^2} &= \sum\limits_{i=4}^{n_\mathrm{p}} a_i \cfrac{\partial^2 P^p_i(x,y)}{\partial x^2}\,,\\
%	\cfrac{\partial^2 u_\mathrm{x}(x,y)}{\partial y^2} &= \sum\limits_{i=4}^{n_\mathrm{p}} a_i \cfrac{\partial^2 P^p_i(x,y)}{\partial y^2}\,,\\
%	\cfrac{\partial^2 u_\mathrm{x}(x,y)}{\partial x \, \partial y} &= \sum\limits_{i=4}^{n_\mathrm{p}} a_i \cfrac{\partial^2 P^p_i(x,y)}{\partial x \, \partial y}\,,\\
%	\cfrac{\partial^2 u_\mathrm{y}(x,y)}{\partial x^2} &= \sum\limits_{i=4}^{n_\mathrm{p}} b_i \cfrac{\partial^2 P^p_i(x,y)}{\partial x^2}\,,\\
%	\cfrac{\partial^2 u_\mathrm{y}(x,y)}{\partial y^2} &= \sum\limits_{i=4}^{n_\mathrm{p}} b_i \cfrac{\partial^2 P^p_i(x,y)}{\partial y^2}\,,\\
%	\cfrac{\partial^2 u_\mathrm{y}(x,y)}{\partial x \, \partial y} &= \sum\limits_{i=4}^{n_\mathrm{p}} b_i \cfrac{\partial^2 P^p_i(x,y)}{\partial x \, \partial y}\,.
%	\end{align}
	%
	\item Substitute the derivatives of the displacement field into the equilibrium conditions Eqs.~\eqref{eq:EC1_app} and \eqref{eq:EC2_app}:
	\begin{align}
	\begin{split}
	\cfrac{E}{1-\nu^2} \left( \sum\limits_{i=4}^{n_\mathrm{p}} a_i \cfrac{\partial^2 P^p_i(x,y)}{\partial x^2} + \nu \sum\limits_{i=4}^{n_\mathrm{p}} b_i \cfrac{\partial^2 P^p_i(x,y)}{\partial x \, \partial y}\nu  \right) + \\ \cfrac{E}{2(1+\nu)} \left( \sum\limits_{i=4}^{n_\mathrm{p}} a_i \cfrac{\partial^2 P^p_i(x,y)}{\partial y^2} + \sum\limits_{i=4}^{n_\mathrm{p}} b_i \cfrac{\partial^2 P^p_i(x,y)}{\partial x \, \partial y} \right) &= 0\,,
	\end{split} \label{eq:EC11_app}
	\\
	\begin{split}
	\cfrac{E}{2(1+\nu)} \left( \sum\limits_{i=4}^{n_\mathrm{p}} a_i \cfrac{\partial^2 P^p_i(x,y)}{\partial x \, \partial y} +  \sum\limits_{i=4}^{n_\mathrm{p}} b_i \cfrac{\partial^2 P^p_i(x,y)}{\partial x^2} \right) + \\  
	\cfrac{E}{1-\nu^2} \left( \sum\limits_{i=4}^{n_\mathrm{p}} b_i \cfrac{\partial^2 P^p_i(x,y)}{\partial y^2} + \nu \sum\limits_{i=4}^{n_\mathrm{p}} a_i \cfrac{\partial^2 P^p_i(x,y)}{\partial x \, \partial y} \right) &= 0\,,
	\end{split} \label{eq:EC21_app}
	\end{align}
	It is important to notice that the results are independent of Young's modulus $E$. %since it is contained in all terms and thus, can simply be factored out and the whole equation is then divided by $E$.
	\item Sort the coefficients in Eqs.~\eqref{eq:EC11_app} and \eqref{eq:EC21_app} according to the powers of $x^iy^j$ (the ordering of the coefficients is chosen according to Pascal's triangle, see  Fig.~\ref{fig:PascalTri_app}); $i,j \in \{0,1,2,\ldots,p\}$:
	\begin{alignat}{7}
	c_0 &+ c_1x &&+ c_2y &&+ c_3x^2 &&+ c_4xy &&+ \ldots &&+ c_{n_\mathrm{p}-3}y^p &&= 0\,,\\
	d_0 &+ d_1x &&+ d_2y &&+ d_3x^2 &&+ d_4xy &&+ \ldots &&+ d_{n_\mathrm{p}-3}y^p &&= 0\,.
	\end{alignat}
	\item Solve for the unknown coefficients $a_i$ and $b_i$ of the polynomial displacement field:
	\begin{enumerate}
		\item Each new coefficient ($c_j$, $d_j$) is a linear combination (function) of the polynomial coefficients ($a_i$, $b_i$) and needs to be zero in order to generate an admissible function.
		\begin{align}
		c_j = f_1(a_i,b_i) &= 0 \\
		d_j = f_2(a_i,b_i) &= 0 
		\end{align}
		\item Due to the required second order derivatives, less than $2n_\mathrm{p}$ equations are available to determine all $2n_\mathrm{p}$ coefficients (at least the constant and linear terms vanish).
		\item We observe that in Eqs.~\eqref{eq:EC11_app} and \eqref{eq:EC21_app} all coefficients associated with the mixed second order derivatives ($\partial^2\square/\partial x\partial y$) are encountered twice, while the second order derivatives with respect to one coordinate ($\partial^2\square/\partial x^2$ or $\partial^2\square/\partial y^2$) are only contained once. Hence, we can prescribe all coefficients that only contain $x^i$ with $i > 1$ (or $y^i$ with $i > 1$) as well as the coefficients for the linear and constant terms (which are not included in the equations due to the second order derivatives); see also Pascal's triangle (Fig.~\ref{fig:PascalTri_app}).
		\begin{itemize}
			\item[-] In other words, compute all coefficients that are related to mixed polynomial terms $x^iy^j$ with $i,j > 0$ (interior terms in Pascal's triangle).
		\end{itemize}
	\end{enumerate}
\end{enumerate}
\begin{figure}[p]
	\centering
	\subfloat[$u_\mathrm{x}$, $p\,{=}\,1$]{\includegraphics[clip,width=0.315\textwidth]{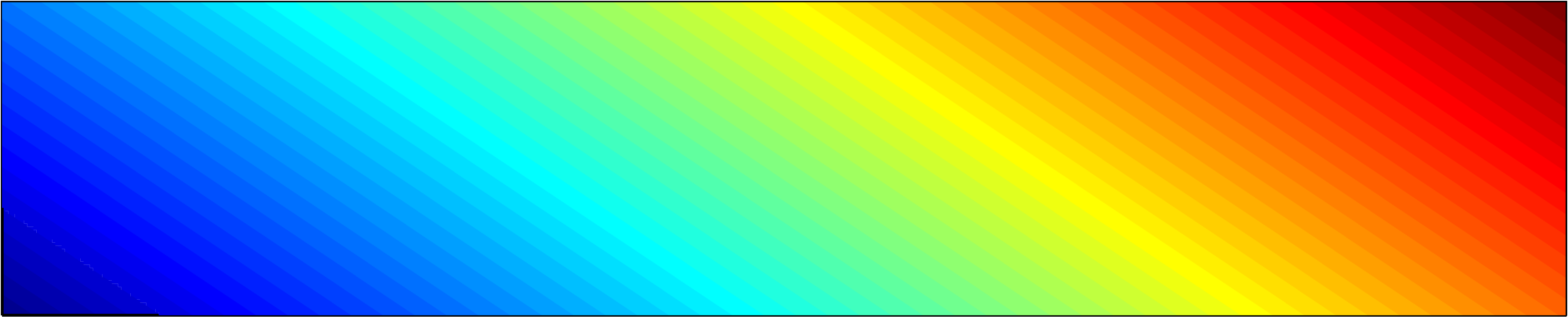}}
	\hfill
	\subfloat[$u_\mathrm{y}$, $p\,{=}\,1$]{\includegraphics[clip,width=0.315\textwidth]{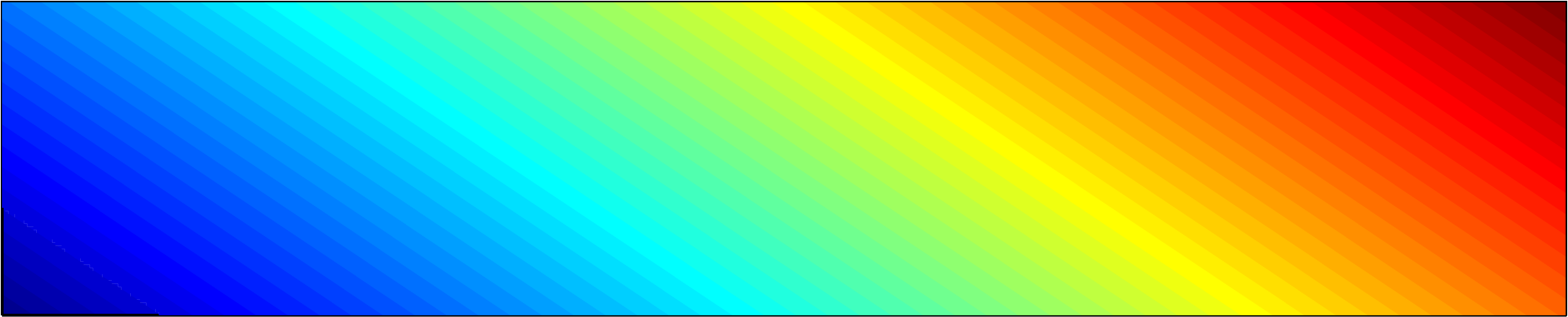}}
	\hfill
	\subfloat[$u_\mathrm{mag}$, $p\,{=}\,1$]{\includegraphics[clip,width=0.315\textwidth]{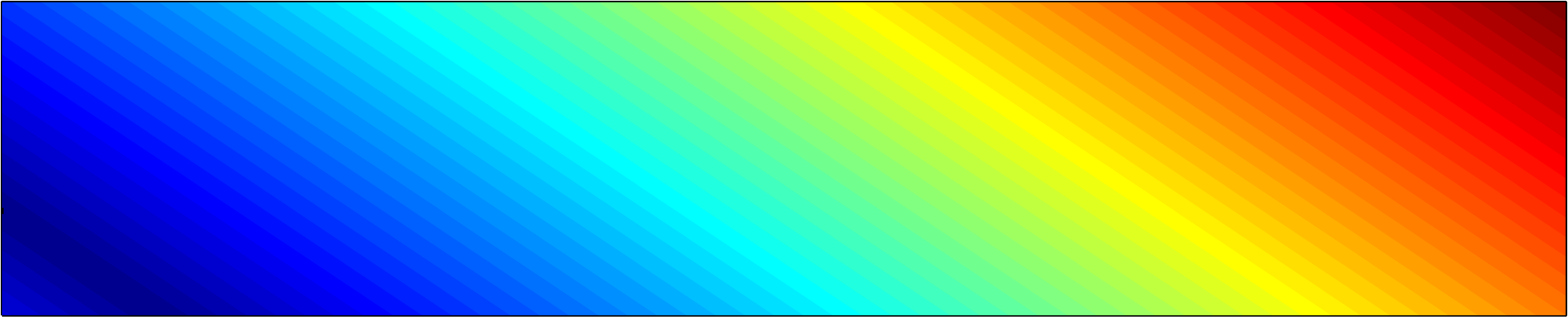}}
	\\
	\subfloat[$u_\mathrm{x}$, $p\,{=}\,2$]{\includegraphics[clip,width=0.315\textwidth]{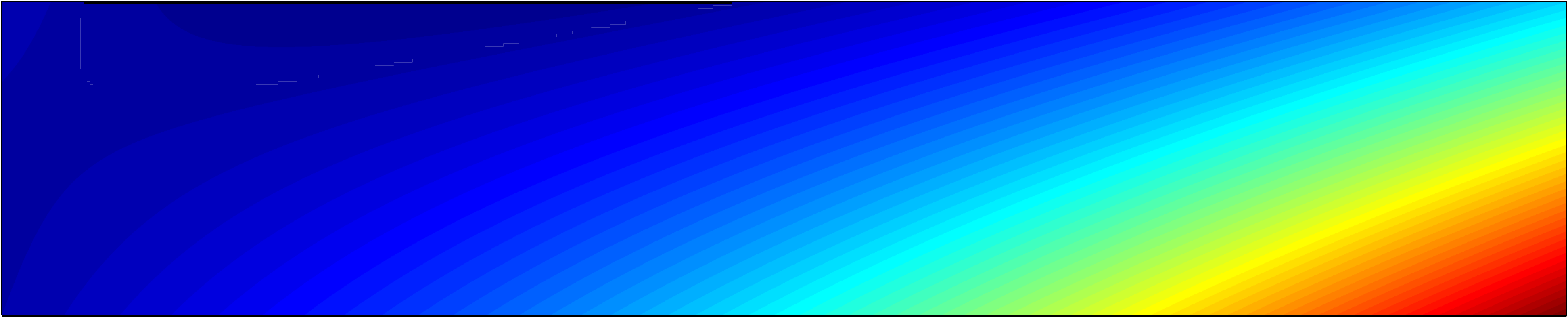}}
	\hfill
	\subfloat[$u_\mathrm{y}$, $p\,{=}\,2$]{\includegraphics[clip,width=0.315\textwidth]{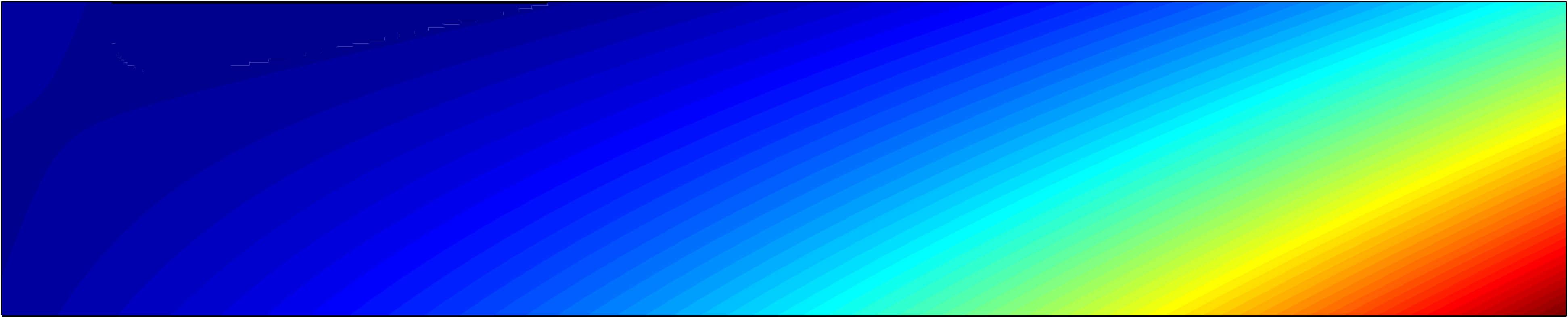}}
	\hfill
	\subfloat[$u_\mathrm{mag}$, $p\,{=}\,2$]{\includegraphics[clip,width=0.315\textwidth]{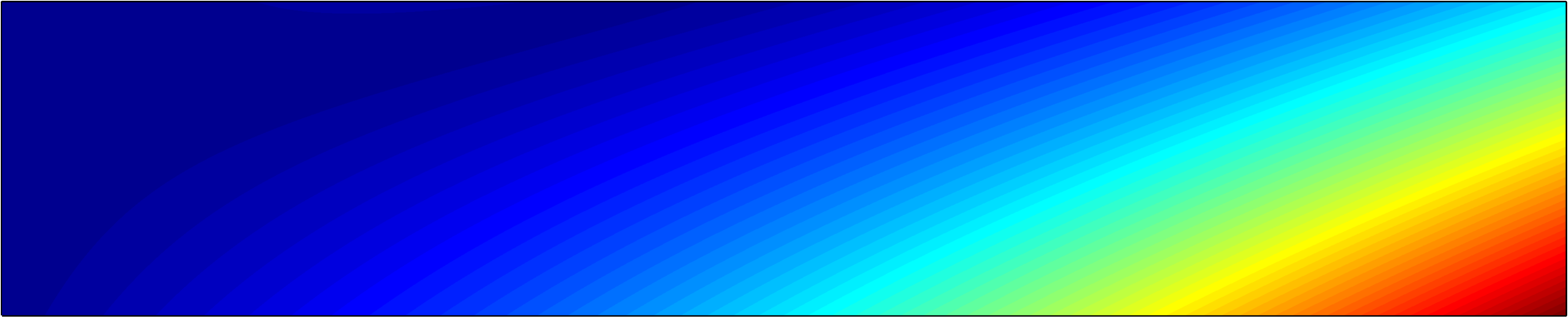}}
	\\
	\subfloat[$u_\mathrm{x}$, $p\,{=}\,3$]{\includegraphics[clip,width=0.315\textwidth]{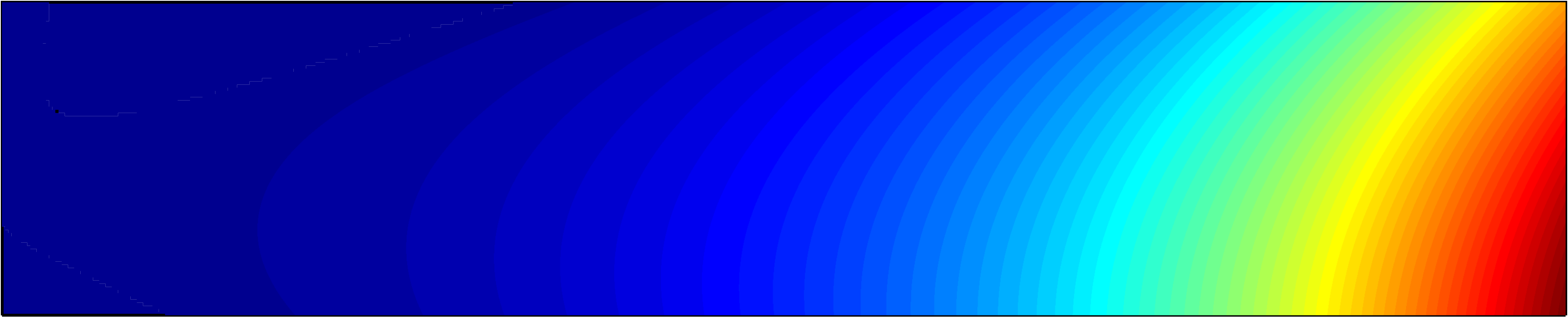}}
	\hfill
	\subfloat[$u_\mathrm{y}$, $p\,{=}\,3$]{\includegraphics[clip,width=0.315\textwidth]{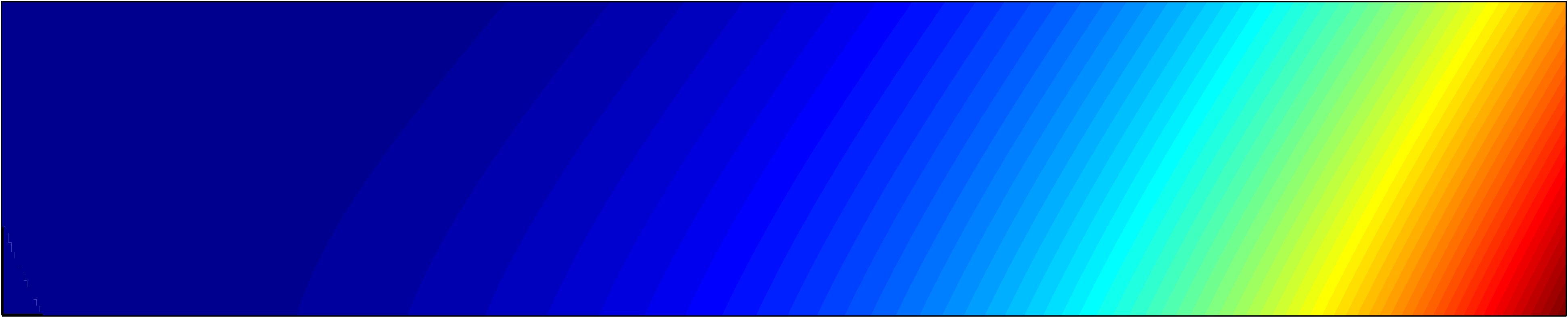}}
	\hfill
	\subfloat[$u_\mathrm{mag}$, $p\,{=}\,3$]{\includegraphics[clip,width=0.315\textwidth]{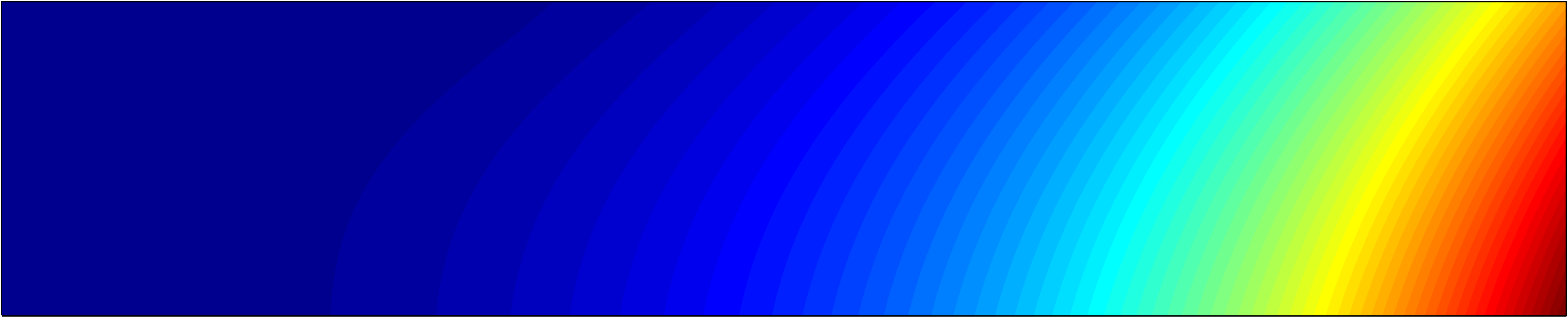}}
	\\
	\subfloat[$u_\mathrm{x}$, $p\,{=}\,4$]{\includegraphics[clip,width=0.315\textwidth]{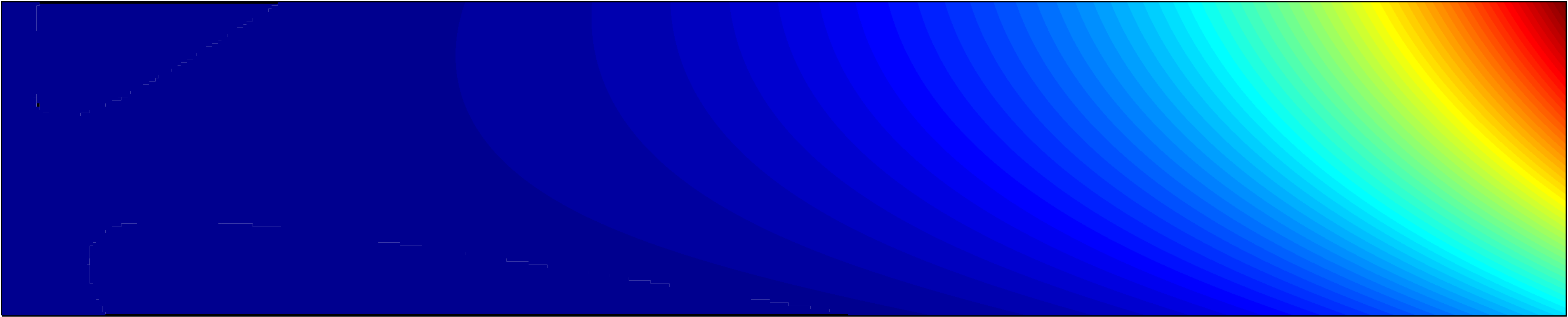}}
	\hfill
	\subfloat[$u_\mathrm{y}$, $p\,{=}\,4$]{\includegraphics[clip,width=0.315\textwidth]{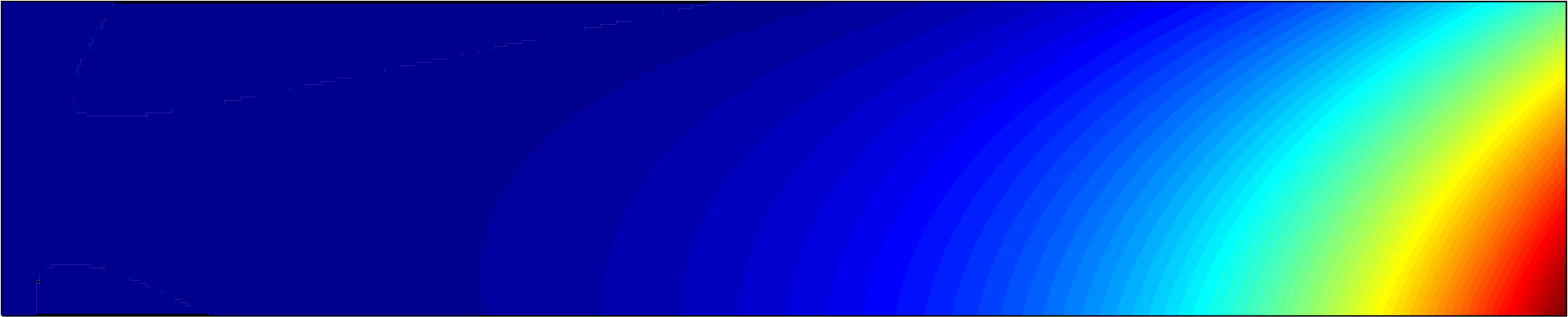}}
	\hfill
	\subfloat[$u_\mathrm{mag}$, $p\,{=}\,4$]{\includegraphics[clip,width=0.315\textwidth]{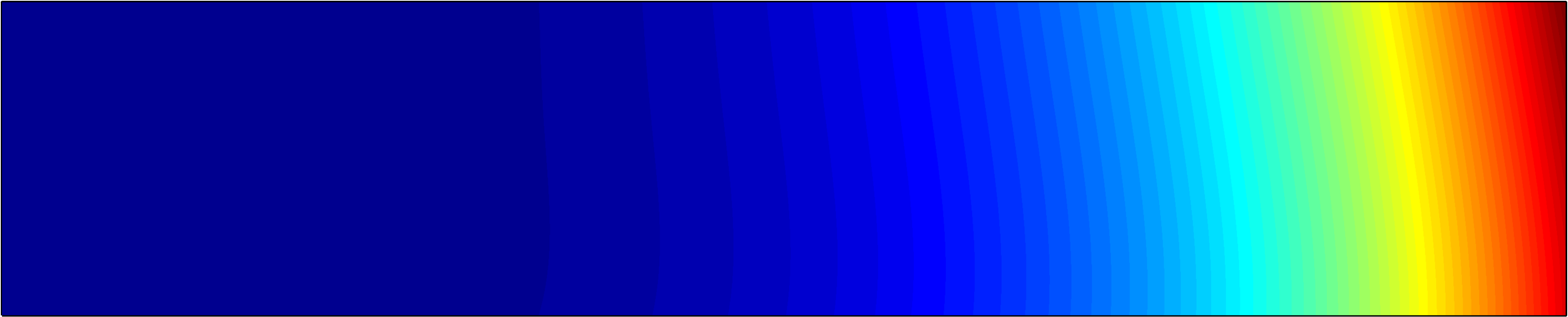}}
	\\
	\subfloat[$u_\mathrm{x}$, $p\,{=}\,5$]{\includegraphics[clip,width=0.315\textwidth]{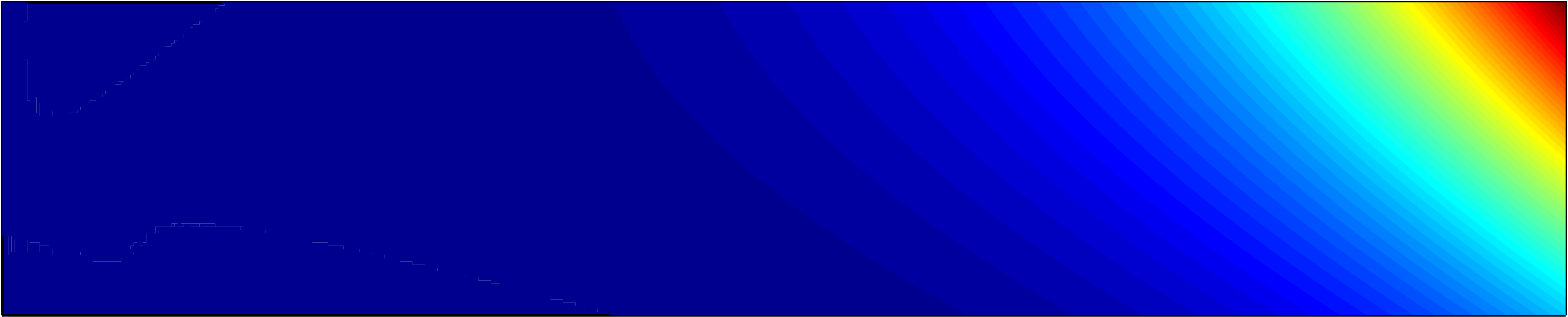}}
	\hfill
	\subfloat[$u_\mathrm{y}$, $p\,{=}\,5$]{\includegraphics[clip,width=0.315\textwidth]{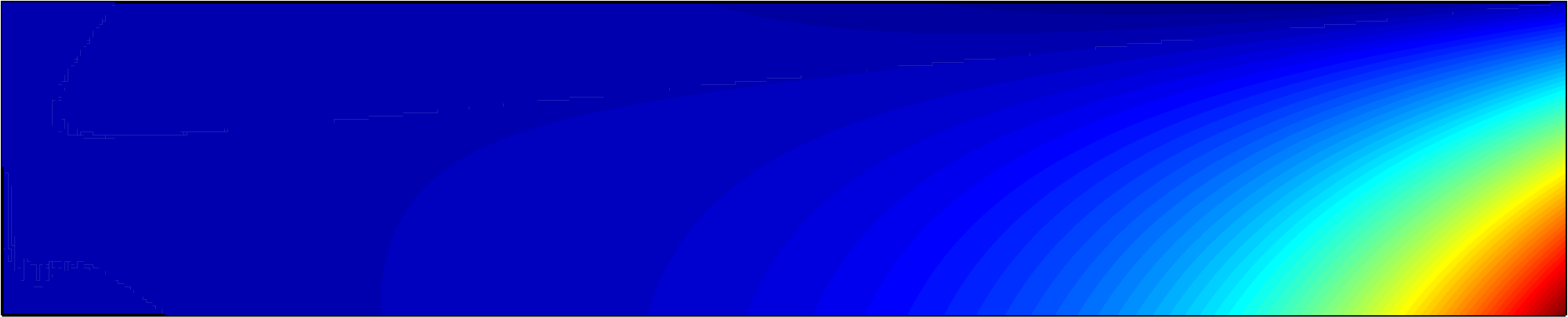}}
	\hfill
	\subfloat[$u_\mathrm{mag}$, $p\,{=}\,5$]{\includegraphics[clip,width=0.315\textwidth]{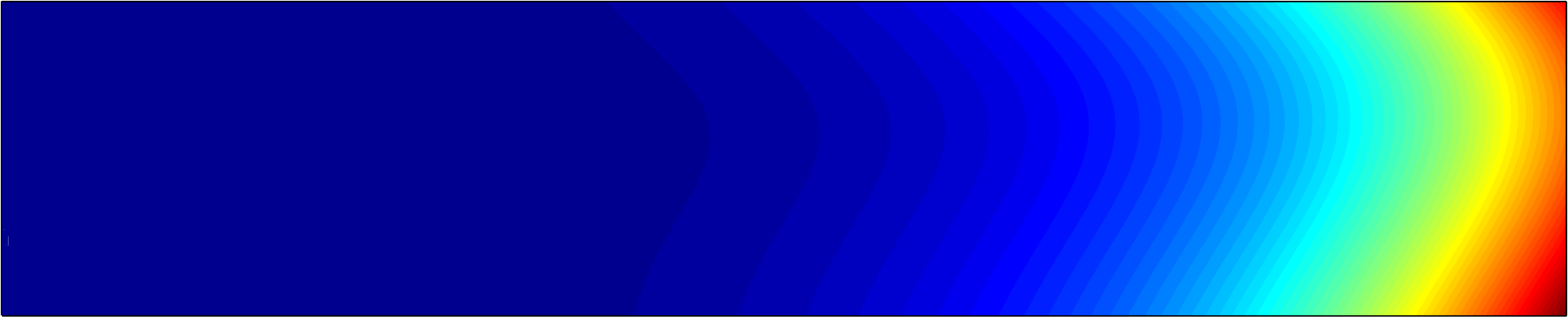}}
	\\
	\subfloat[$u_\mathrm{x}$, $p\,{=}\,6$]{\includegraphics[clip,width=0.315\textwidth]{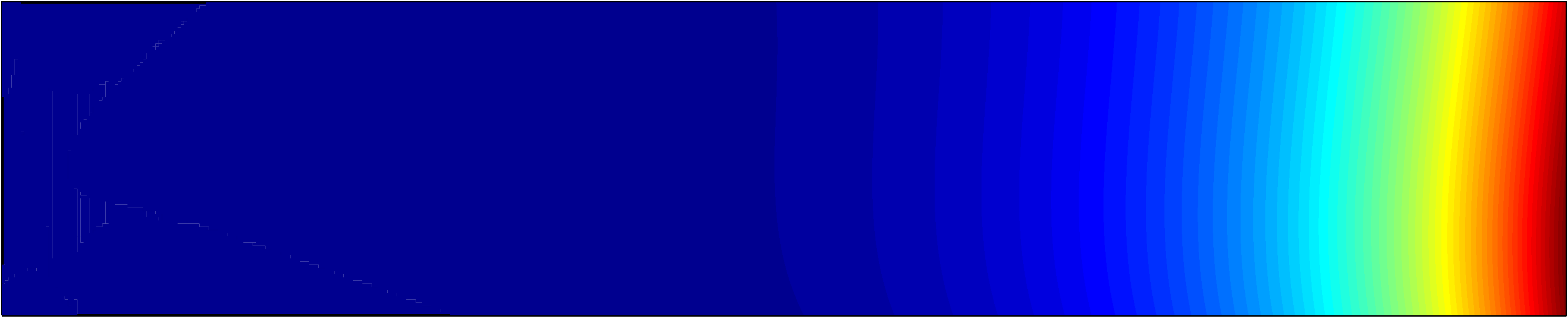}}
	\hfill
	\subfloat[$u_\mathrm{y}$, $p\,{=}\,6$]{\includegraphics[clip,width=0.315\textwidth]{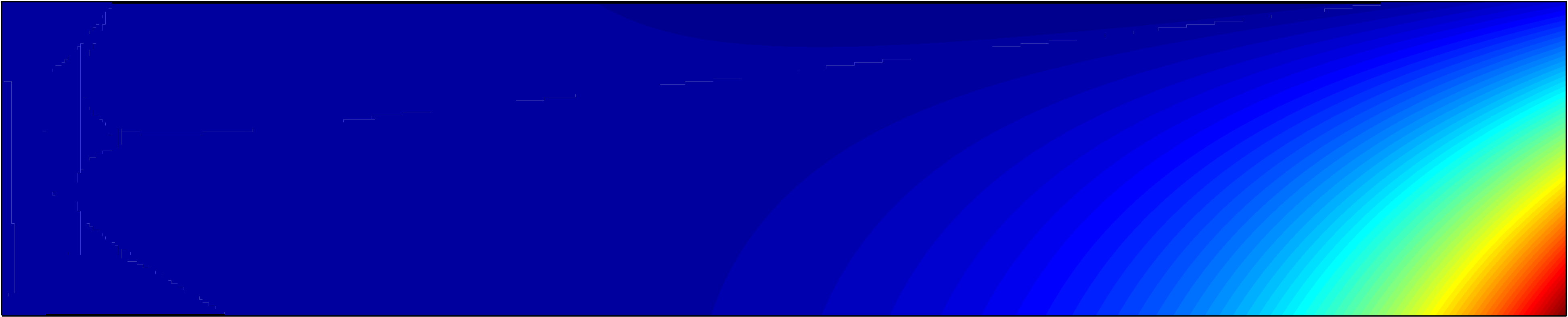}}
	\hfill
	\subfloat[$u_\mathrm{mag}$, $p\,{=}\,6$]{\includegraphics[clip,width=0.315\textwidth]{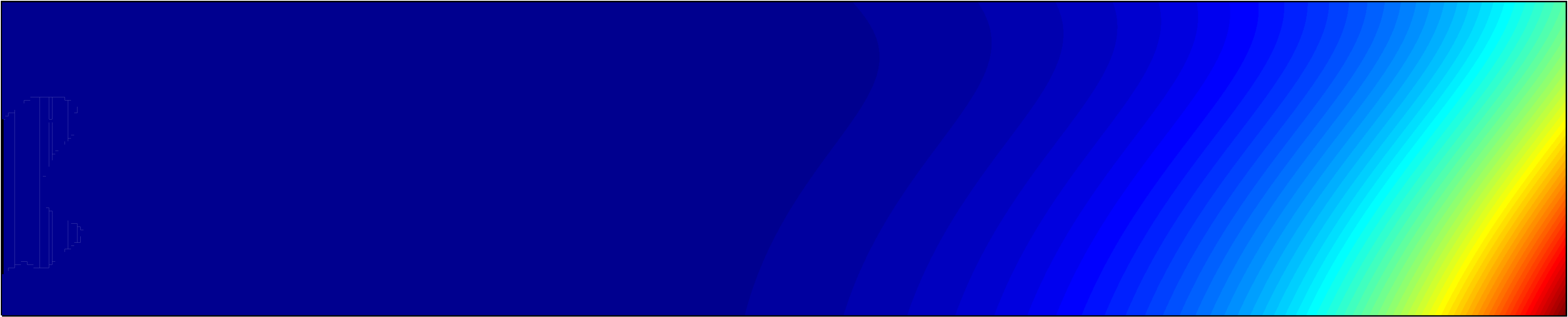}}
	\\
	\subfloat[$u_\mathrm{x}$, $p\,{=}\,7$]{\includegraphics[clip,width=0.315\textwidth]{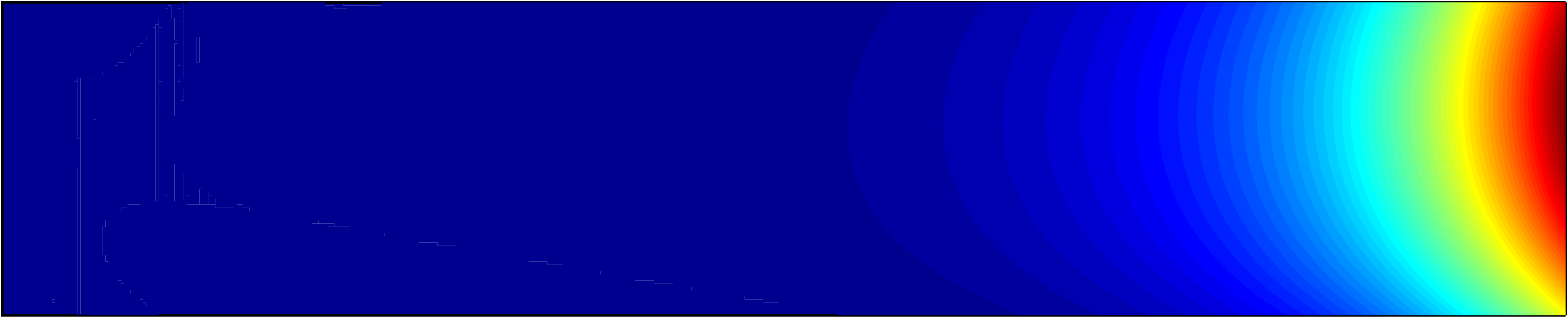}}
	\hfill
	\subfloat[$u_\mathrm{y}$, $p\,{=}\,7$]{\includegraphics[clip,width=0.315\textwidth]{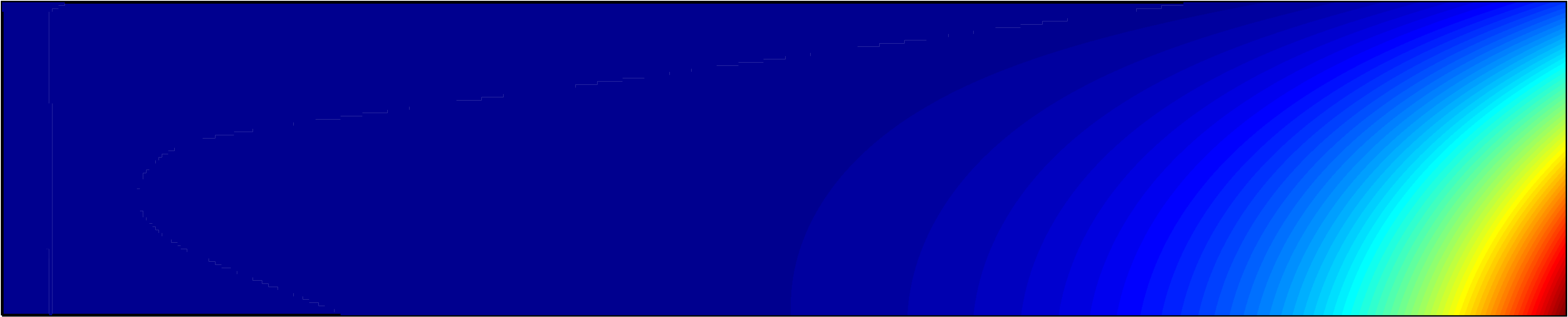}}
	\hfill
	\subfloat[$u_\mathrm{mag}$, $p\,{=}\,7$]{\includegraphics[clip,width=0.315\textwidth]{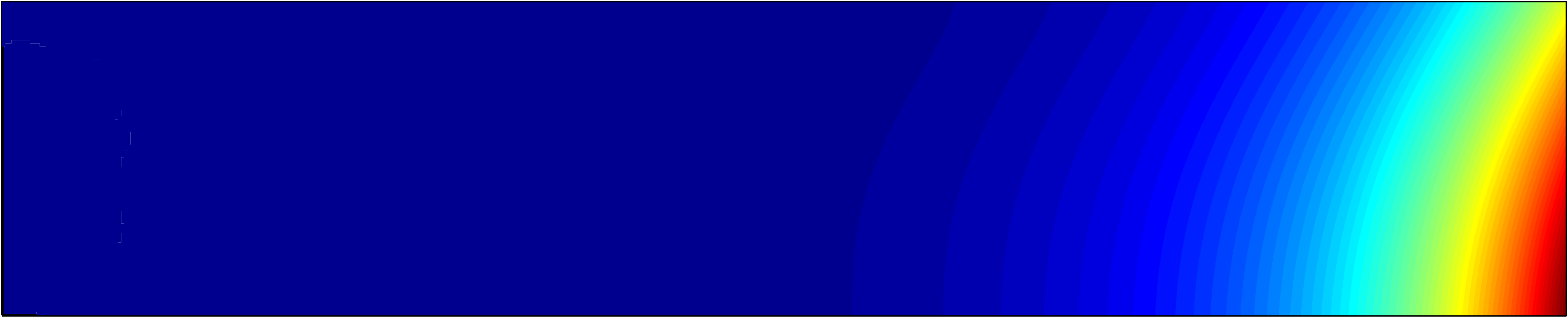}}
	\\
	\subfloat[$u_\mathrm{x}$, $p\,{=}\,8$]{\includegraphics[clip,width=0.315\textwidth]{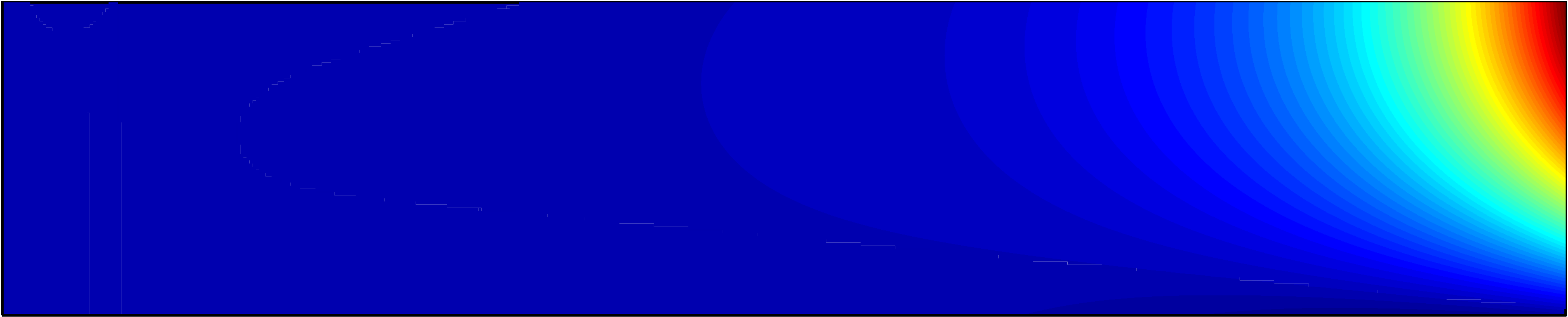}}
	\hfill
	\subfloat[$u_\mathrm{y}$, $p\,{=}\,8$]{\includegraphics[clip,width=0.315\textwidth]{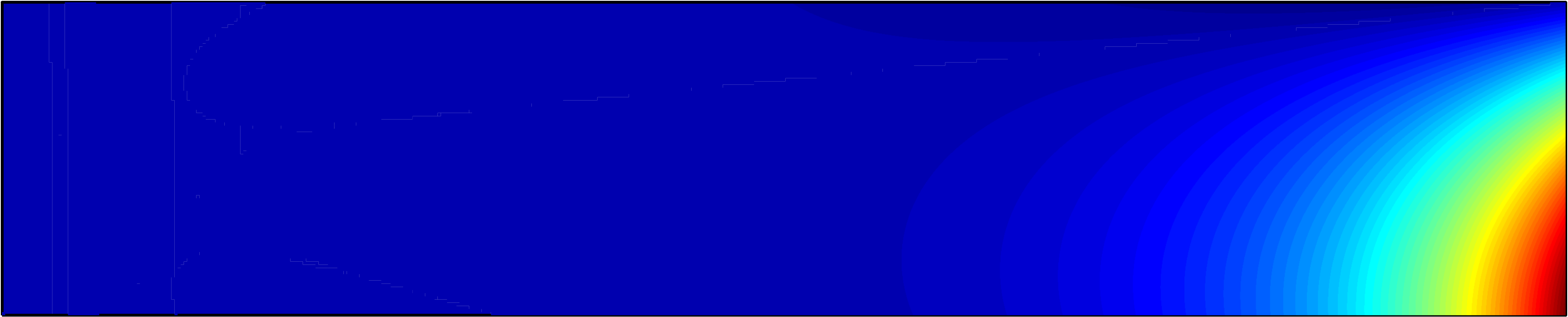}}
	\hfill
	\subfloat[$u_\mathrm{mag}$, $p\,{=}\,8$]{\includegraphics[clip,width=0.315\textwidth]{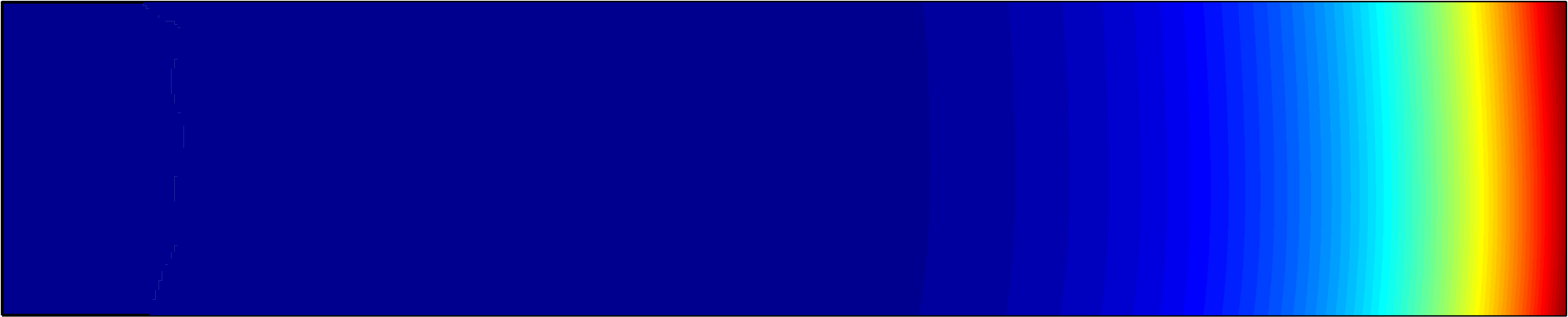}}
	\caption{Admissible displacement fields for the higher order patch test ($p\,{=}\,1,2,\ldots,8$), derived from the data provided in Table~\ref{tab:8OrderDisp_app}.}
	\label{fig:HigherOrderDispField_app}
\end{figure}%
\begin{table}[p]
	\centering
	\caption{Eighth order (admissible) displacement field for patch test version B. Note: all coefficients that are not included in the table below are prescribed as: $a_i\,{=}\,b_i\,{=}\,i\,.$ \label{tab:8OrderDisp_app}}
	\begin{tabular}{l|l||l|l}
		\toprule
		$u_x$ & & $u_y$ & \\[3pt]\hline
		$a_{5}\,{=}\,8\cfrac{\nu-4}{\nu+1}$			&	$a_{27}\,{=}\,-4\cfrac{17\nu+67}{\nu+1}$	&	$b_{5}\,{=}\,4\cfrac{3\nu-7}{\nu+1}$		&	$b_{27}\,{=}\,8\cfrac{21\nu-4}{\nu+1}$		\\[9pt]
		$a_{8}\,{=}\,3\cfrac{7\nu-13}{\nu+3}$		&	$a_{30}\,{=}\,21\cfrac{29\nu+5}{3\nu+5}$	&	$b_{8}\,{=}\,-3\cfrac{17\nu-3}{2\nu}$		&   $b_{30}\,{=}\,-21\cfrac{41\nu+17}{4\nu+2}$	\\[9pt]
		$a_{9}\,{=}\,-3\cfrac{17\nu+3}{2\nu}$		&	$a_{31}\,{=}\,-21\cfrac{123\nu+65}{4\nu+2}$	&	$b_{9}\,{=}\,6\cfrac{5\nu-2}{\nu+3}$		&   $b_{31}\,{=}\,-42\cfrac{11\nu+40}{3\nu+5}$	\\[9pt]
		$a_{12}\,{=}\,4\cfrac{11\nu+15}{\nu+1}$		&	$a_{32}\,{=}\,-70\cfrac{11\nu-25}{3\nu+5}$	&	$b_{12}\,{=}\,-4\cfrac{13\nu+9}{\nu+1}$		&   $b_{32}\,{=}\,35\cfrac{65\nu+29}{2\nu+1}$	\\[9pt]
		$a_{13}\,{=}\,-78$ 							&	$a_{33}\,{=}\,35\cfrac{65\nu+36}{2\nu+1}$	&	$b_{13}\,{=}\,-78$							&   $b_{33}\,{=}\,-35\cfrac{43\nu-15}{3\nu+5}$	\\[9pt]
		$a_{14}\,{=}\,-4\cfrac{13\nu+17}{\nu+1}$	&	$a_{34}\,{=}\,-21\cfrac{43\nu+115}{3\nu+5}$	&	$b_{14}\,{=}\,4\cfrac{15\nu+11}{\nu+1}$		&   $b_{34}\,{=}\,-21\cfrac{137\nu+65}{4\nu+2}$	\\[9pt]
		$a_{17}\,{=}\,5\cfrac{16\nu+37}{\nu+2}$		&	$a_{35}\,{=}\,-7\cfrac{137\nu+79}{4\nu+2}$	&	$b_{17}\,{=}\,-5\cfrac{11\nu+53}{3\nu+1}$	&   $b_{35}\,{=}\,14\cfrac{54\nu+25}{3\nu+5}$	\\[9pt]
		$a_{18}\,{=}\,-10\cfrac{11\nu-21}{3\nu+1}$	&	$a_{38}\,{=}\,8\cfrac{37\nu+41}{\nu+1}$		&	$b_{18}\,{=}\,-5\cfrac{37\nu+69}{\nu+2}$	&   $b_{38}\,{=}\,-8\cfrac{39\nu+35}{\nu+1}$	\\[9pt]
		$a_{19}\,{=}\,-5\cfrac{37\nu+79}{\nu+2}$	&	$a_{39}\,{=}\,-1092$						&	$b_{19}\,{=}\,-20\cfrac{13\nu-8}{3\nu+1}$	&   $b_{39}\,{=}\,-1092$						\\[9pt]
		$a_{20}\,{=}\,-10\cfrac{13\nu+29}{3\nu+1}$	&	$a_{40}\,{=}\,-56\cfrac{39\nu+43}{\nu+1}$	&	$b_{20}\,{=}\,5\cfrac{21\nu+37}{\nu+2}$		&   $b_{40}\,{=}\,56\cfrac{41\nu+37}{\nu+1}$	\\[9pt]
		$a_{23}\,{=}\,4\cfrac{33\nu-17}{\nu+1}$		&	$a_{41}\,{=}\,2870$ 						&	$b_{23}\,{=}\,-8\cfrac{4\nu+29}{\nu+1}$		&   $b_{41}\,{=}\,2870$							\\[9pt]
		$a_{24}\,{=}\,-80$							&	$a_{42}\,{=}\,56\cfrac{14\nu+45}{\nu+1}$ 	&	$b_{24}\,{=}\,-80$							&   $b_{42}\,{=}\,-56\cfrac{43\nu+39}{\nu+1}$	\\[9pt]
		$a_{25}\,{=}\,-80\cfrac{4\nu_21}{3\nu+3}$	&	$a_{43}\,{=}\,-1204$ 						&	$b_{25}\,{=}\,-40\cfrac{17\nu-33}{3\nu+3}$	&   $b_{43}\,{=}\,-1204$						\\[9pt]
		$a_{26}\,{=}\,-170$							&	$a_{44}\,{=}\,-8\cfrac{43\nu+47}{\nu+1}$	&	$b_{26}\,{=}\,-170$							&   $b_{44}\,{=}\,8\cfrac{45\nu+41}{\nu+1}$		\\
		\bottomrule
	\end{tabular}
\end{table}%
In the following, we apply the outlined procedure to determine admissible (polynomial) displacement fields of orders $p\,{=}\,2,3,\ldots,8$. Due to the second-order derivatives that are involved in the equilibrium conditions, any linear displacement field is admissible and therefore, all three coefficients can assume arbitrary values. Note that for higher orders, the situation changes a little bit. Taking, for example, a complete fifth-order polynomial, we have 21 unknown coefficients and therefore, we need to determine the values of 42 parameters (two-dimensional case). Fortunately, as discussed above, we only need to solve for the mixed terms and can prescribe all other coefficients. Therefore, we need to find 20 conditions (10 mixed terms) to uniquely determine our function. For the sake of simplicity, we prescribe $a_i = i$ and $b_i = i$ as coefficients related to the terms that do not contain a mixed product\footnote{Remark: The coefficients are numbered according to their position in Pascal's triangle starting from the top and continuing downwards and from left to right. Hence, a polynomial function has the following form: $P_\mathrm{p} = a_1 + a_2x + a_3y + a_4x^2 + \ldots$}. The results for a complete eighth order polynomial are listed in Table~\ref{tab:8OrderDisp_app}. The higher order (admissible) displacement field is hierarchic in the sense that the coefficients of a complete polynomial of order $p$ also contain those of order $p\,{-}\,1$. Thus, all lower order displacement fields ($p\,{\le}\,8$) can be recovered from the data compiled in Table~\ref{tab:8OrderDisp_app}. For the computational model shown in Fig.~\ref{fig:PatchTestHigherOrder}, the corresponding admissible displacement fields for complete polynomials of orders $p\,{=}\,1,2,\ldots,8$ are depicted in Fig.~\ref{fig:HigherOrderDispField_app}.